\pdfoutput=1 
\documentclass[reqno, 12pt]{amsart}
\usepackage{amsmath}%
\usepackage{amsfonts}%
\usepackage{amssymb}%
\usepackage{graphicx}
\usepackage{epsfig}

\usepackage{url}

\usepackage{ifthen}
\newboolean{BKMRK}
\setboolean{BKMRK}{true}    

\ifthenelse {\boolean{BKMRK}}
  { \usepackage{hyperref} }
  {  }

\ifthenelse {\boolean{BKMRK}}
  { \hypersetup{colorlinks=true, linkcolor=red, citecolor=red, pdfstartview=FitH} }
  {  }

\def\picDblWidth{3.20in}    

\pdfoutput=1    

\newtheorem{theorem}{Theorem}
\theoremstyle{plain}

\numberwithin{equation}{section}

\begin{document}
\title{Experiments With Zeta Zeros and Perron's Formula}
\author{Robert Baillie}
\address{State College, PA 16803}


\date{\today}


\begin{abstract}



Of what use are the zeros of the Riemann zeta function?  We can use sums involving zeta zeros to count the primes up to $x$.  Perron's formula leads to sums over zeta zeros that can count the squarefree integers up to $x$, or tally Euler's $\phi$ function and other arithmetical functions.  This is largely a presentation of experimental results.

\end{abstract}

\maketitle

\setcounter{tocdepth}{1}  
\smaller[2]    
\tableofcontents
\larger[2]

\section{Introduction}

Most mathematicians know that the zeros of the Riemann zeta function have something to do with the distribution of primes. What is less well known that we can use the values of the zeta zeros to calculate, quite accurately in some cases, the values of interesting arithmetical functions.  For example, consider Euler's totient function, $\phi(n)$, which counts the number of integers from 1 through $n$ that are relatively prime to $n$.  The summatory function of $\phi$, that is, the partial sum of $\phi(n)$ for $1 \leq n \leq x$,
\[
\Phi (x)=\sum _{n \leq x} \phi(n)
\]
is an irregular step function that increases roughly as $x^2$.  In fact, a standard result in analytic number theory is that \cite[Theorem 330, p. 268]{HardyAndWright}
\[
  \Phi(x) = \frac{3 x^2}{\pi^2} + O( x \log x ).
\]
But such an estimate mainly tells us how fast $\Phi(x)$ increases.  $\frac{3 x^2}{\pi^2}$ is a smoothly increasing function and does not replicate the details of the step function $\Phi(x)$ at all.
But by using a sum whose terms include zeta zeros, we can get a function that appears to closely match the step function $\Phi(x)$.  Our approximation rises rapidly at the integers, where $\Phi(x)$ has jump discontinuities, and is fairly level between integers, where $\Phi(x)$ is constant.

As another example, consider $\sigma(n)$, the sum of the (positive) divisors of $n$.  The summatory function
\[
T (x)=\sum _{n \leq x} \sigma(n^2).
\]
increases roughly as $x^3$, but zeta zeros give us a good approximation to $T(x)$.  Both of these sums come to us courtesy of Perron's formula, a remarkable theorem that takes integrals of zeta in the complex plane and yields approximations to the sums of integer arithmetical functions.  We then use Cauchy's theorem to estimate these integrals by summing residues; these residues are expressions that involve zeta zeros.

Asymptotic estimates of the type
\begin{equation}\label{E:Tx1}
T(x) = f(x) + O(g(x))
\end{equation}
are known for some of the summatory functions we consider here.  See, for example, section ~\ref{S:Powerful}.  It often turns out that the method used here (summing residues) easily produces
\begin{equation}\label{E:Tx2}
T(x) \simeq f(x) + \text{sum over zeta zeros},
\end{equation}
and even gives us the exact coefficients in $f(x)$, although it does not yield the error estimate $g(x)$.  Assuming that equation \eqref{E:Tx2} holds, equation \eqref{E:Tx1} is essentially estimating the size of a sum over zeta zeros.

Using zeta zeros, we can even do a pretty good job of computing the value of $\pi(x)$, the number of primes $\leq x$, although Perron's formula is not used in this case.  But all these examples suggest that zeta zeros have the power to predict the behavior of many arithmetical functions.

This paper is long on examples and short on proofs.  The approach here is experimental.  We will use standard results of analytic number theory, with appropriate references to those works, instead of repeating the proofs here.  In many cases, we will give approximations to summatory functions that are not proven.  Proofs that these expressions do, in fact, approach the summatory functions could make interesting papers.  Further, this paper is still incomplete, and should be considered a work in progress.

Many thanks are due to professors Harold Diamond and Robert Vaughan, who have patiently answered numerous questions from the author.

\section{The Riemann Zeta Function and its Zeros}
Here is a brief summary of some basic facts about the Riemann zeta function and what is known about its zeros.  As is customary in analytic number theory, we let $s = \sigma + i t$ denote a complex variable.  $\Re(s)$ and $\Im(s)$ are the real and imaginary parts of $s$.

The Riemann zeta function is defined for $\Re(s) > 1$ by
\[
\zeta (s)=\sum _{n=1}^{\infty } \frac{1}{n^s}.
\]

Zeta can be analytically continued to the entire complex plane via the functional equation \cite[p. 13]{Titchmarsh},

\[
\zeta (s)=2^s \pi ^{s-1} \zeta (1-s) \sin \left(\frac{\pi s}{2}\right) \Gamma (1-s).
\]

Zeta has a simple pole (that is, a pole of order 1) at $s = 1$, where the residue is 1.

The real zeros of zeta are $s = -2, -4, -6,$ \ldots .  They are often referred to as the ``trivial'' zeros.  There are also complex zeros (the so-called ``nontrivial'' zeros).  These zeros occur in conjugate pairs, so if $a + b i$ is a zero, then so is $a - b i$.

When we count the complex zeros, we number them starting at the real axis, going up.  So, the ``first'' three zeros are approximately $1/2 + 14.135 i$, $1/2 + 20.022 i$, and $1/2 + 25.011 i$.  The imaginary parts listed here are approximate, but the real parts are \textit{exactly} $1/2$.  Gourdon \cite{Gourdon} has shown that the first $10^{13}$ complex zeros all have real part equal to $1/2$.  The Riemann Hypothesis is the claim that all complex zeros lie on the so-called ``critical'' line $Re(s) = 1/2$.  This is still unproven.  However, we do know \cite{HardyAndLittlewood} that the zeta function has infinitely many zeros on the critical line.  We even know \cite{Conrey} that over 40\% of all of the complex zeros lie on this line.

\section{Zeta Zeros and the Distribution of Primes}

Let $c$ be the least upper bound of the real part of the complex zeta zeros.  We know that $1/2 \leq c \leq 1$.  Then it follows that \cite[p. 430]{MontgomeryAndVaughan}

\[
\pi(x)= \text{li}(x) + O\left(x^c \log(x)\right)
\]

\noindent where
\[
\text{li}(x)=\int_2^x \frac{1}{\log (t)} \, dt.
\]

This shows a connection between the real parts of the zeta zeros and the size of the error term in the estimate for $\pi(x)$. However, this sort of estimate says nothing about how $\pi(x)$ is connected to the \emph{imaginary} parts of the zeta zeros.  We will see that the imaginary parts of the zeros can be used to calculate the values of $\pi(x)$ and many other number-theoretic functions.

This topic falls under the category of ``explicit formulas''.  Two examples of explicit formulas are commonly found in number theory books: a formula for $\pi(x)$, and a formula for Chebyshev's psi function.  We will examine these, then consider many more examples.



\ifthenelse {\boolean{BKMRK}}
  { \section{The Riemann-von Mangoldt Formula for \texorpdfstring{$\pi(x)$}{pi(x)}}\label{S:rVonM} }
  { \section{The Riemann-von Mangoldt Formula for $\pi(x)$}\label{S:rVonM} }

von Mangoldt expanded on the work of Riemann and derived a remarkable explicit formula for $\pi(x)$ that is based on sums over the complex zeta zeros.

First, note that, if $x$ is prime, then $\pi(x)$ jumps up by 1, causing $\pi(x)$ to be discontinuous.  $\pi(x)$ is continuous for all other $x$.

For reasons that will become clear in a moment, we will define a function $\pi_0(x)$, to be a slight variation of $\pi(x)$.  These two functions agree except where $\pi(x)$ is discontinuous:
\[
\pi_0(x)=
  \begin{cases}
    \pi(x)-1/2  &\text{if $x$ is prime,}\\
    \pi(x) &\text{otherwise.}
  \end{cases}
\]
When $\pi(x)$ jumps up by 1, $\pi_0(x)$ equals the average of the values of $\pi(x)$ before and after the jump.  We can also write $\pi_0(x)$ using limits:
\[
\pi _0(x)=\frac{1}{2} \left(\lim_{t\to x^-} \, \pi (t)+\lim_{t\to x^+} \, \pi (t)\right).
\]

The Riemann-von Mangoldt formula \cite[pp. 44-55]{RieselBook} and \cite[eq. 5]{RieselPaper} states that

\begin{equation}\label{E:r-vonm}
\pi _0(x)=\sum _{n=1}^{\infty } \frac{\mu (n)}{n}f\left(x^{1/n}\right)
\end{equation}

where

\begin{equation}\label{E:r-vonm2}
f(x)=\sum _{n=1}^{\infty } \frac{\pi_0\left(x^{1/n}\right)}{n}
=\text{li}(x)
-\sum _{k=1}^{\infty} \text{li}\left(x^{\rho_k}\right)
+\int_x^{\infty }\frac{1}{\left(t^3-t\right) \log (t)} \, dt-\log (2).
\end{equation}

Here, $\mu(n)$ is the M\"{o}bius mu function.  If the prime factorization of $n$ is
\[
n=p_1^{a_1} \cdot p_2^{a_2} \cdots p_k^{a_k},
\]
then $\mu(n)$ is defined as follows:

\begin{equation}\label{E:muDefinition}
\mu(n)=
  \begin{cases}
    1 &\text{if $n = 1$}\\
    0  &\text{if any $a_1, a_2, \cdots, a_k > 1$,}\\
    (-1)^{k} &\text{otherwise.}
  \end{cases}
\end{equation}
That is, $\mu(n) = 0$ if $n$ is divisible by the square of a prime; otherwise $\mu(n)$ is either $-1$ or $+1$, according to whether $n$ has an odd or even number of distinct prime factors.

Several comments about the Riemann-von Mangoldt formula are in order.

First, notice that the sum in ~\eqref{E:r-vonm} converges to $\pi_0(x)$, not to $\pi(x)$.  That's why we needed to define $\pi_0(x)$.

Second, the sum in ~\eqref{E:r-vonm} is not an infinite sum: it terminates as soon as $x^{1/n} < 2$, because as soon as $n$ is large enough to make $x^{1/n} < 2$, we will have $f(x^{1/n}) = 0$.

Third, the summation over the complex zeta zeros is not absolutely convergent, so the value of the sum may depend on the order in which we add the terms. So, by convention, the sum is computed as follows.  For a given $N$, we add the terms for the first $N$ conjugate \emph{pairs} of zeta zeros, where the zeros are taken in increasing order of their imaginary parts.  Then, we let $N$ approach $\infty$:

\[
\lim_{N\to \infty } \, \left(\sum _{k=1}^N
   \text{li}(x^{\rho_k})
 + \text{li}(x^{\rho_k})^{*}
  \right).
\]
where the asterisk denotes the complex conjugate.  Here is an equivalent way to write this sum:
\[
\lim_{T\to \infty } \, \left(\sum _{|t| \leq T} \text{li}\left(x^{\rho}\right)\right)
\]
where $t$ is the imaginary part of the complex zero $\rho$.

Fourth, as a practical matter, when we compute the sum over zeta zeros, we need to use only those zeta zeros that have positive imaginary part.  Suppose that $a + b i$ is a zeta zero with $b > 0$, and that the corresponding term in the series is $A + B i$.  When evaluated at the conjugate zeta zero, $a - b i$, the term will have the value $A - B i$.  These terms are conjugates of each other, so when they are added, the imaginary parts cancel, and the sum is just $2A$.  Therefore, we can efficiently compute the sum as
\[
2 \Re( \sum _{k=1}^{\infty} \text{li}\left(x^{\rho_k}\right) ),
\]
where now the sum is over only those zeta zeros that have positive imaginary part.

\begin{figure*}[ht]
  \centerline{
    \mbox{\includegraphics[width=\picDblWidth]{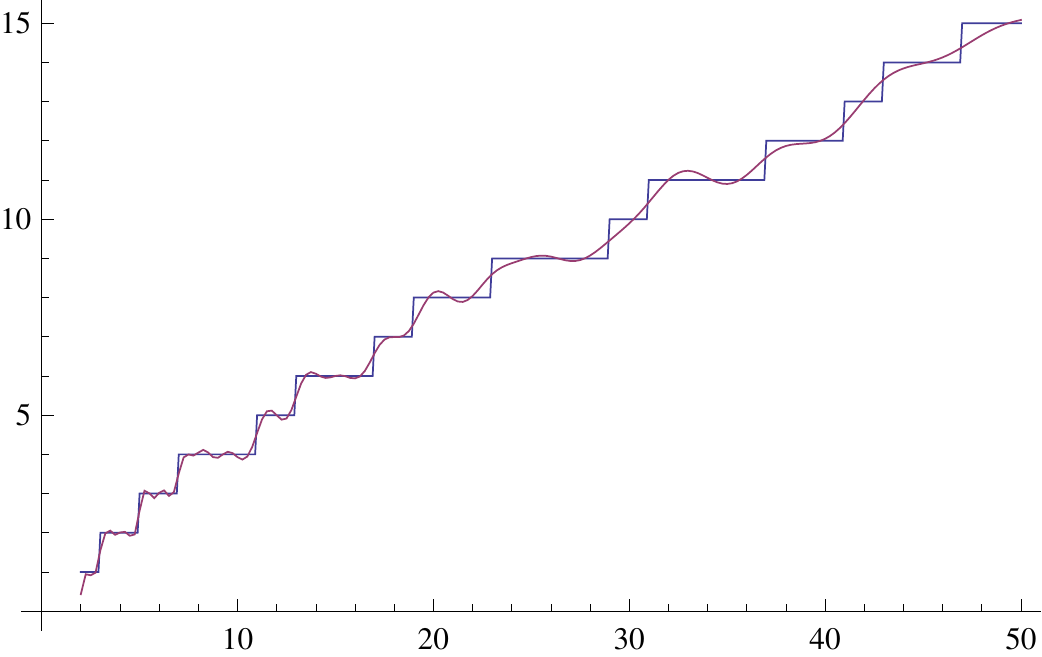}}
    \mbox{\includegraphics[width=\picDblWidth]{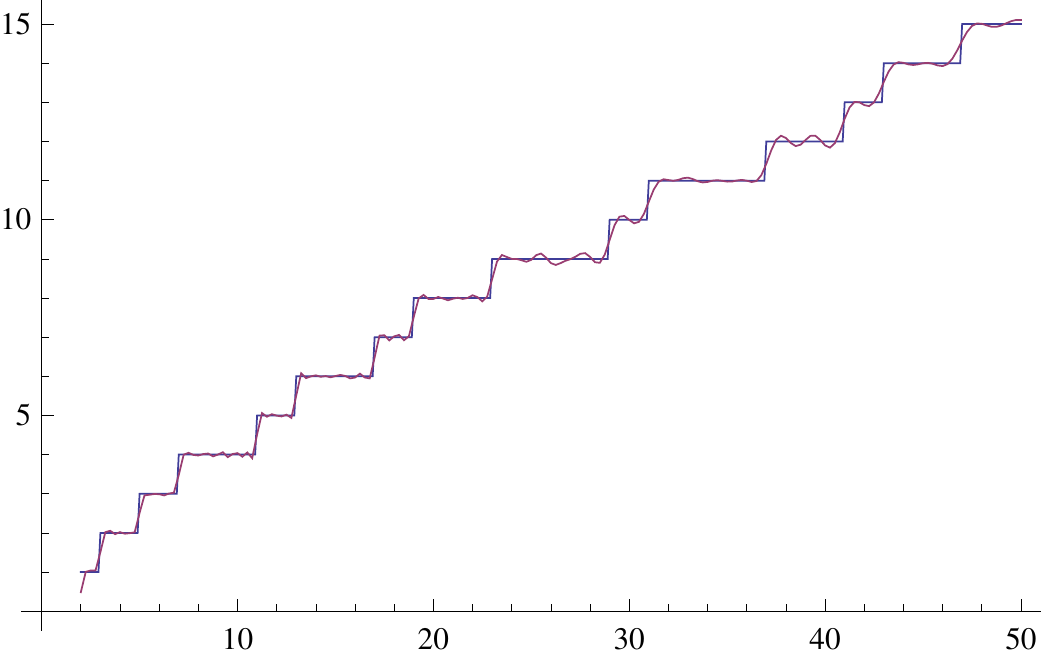}}
  }
  \caption{approximating $\pi(x)$ with $10$ and $50$ pairs of zeta zeros}
  \label{fig:twopix}
\end{figure*}

Can ~\eqref{E:r-vonm} and ~\eqref{E:r-vonm2} really be used to approximate the graph of $\pi(x)$?  Yes! Figure ~\ref{fig:twopix} shows the step function $\pi(x)$ and some approximations to $\pi(x)$.  To make the left-hand graph, we used $10$ pairs of zeta zeros in the sum in equation ~\eqref{E:r-vonm2}.  The step function $\pi(x)$ is displayed as a series of vertical and horizontal line segments.  Of course, the vertical segments aren't really part of the graph.  However, they are useful because they help show that, when $x$ is prime, ~\eqref{E:r-vonm} is close to the midpoints of the jumps of the step function.

In the right-hand graph, we used the first 50 pairs of zeta zeros.  For larger values of $x$, this gives a better approximation.

\begin{figure*}[ht]
  \centerline{
    \mbox{\includegraphics[width=\picDblWidth]{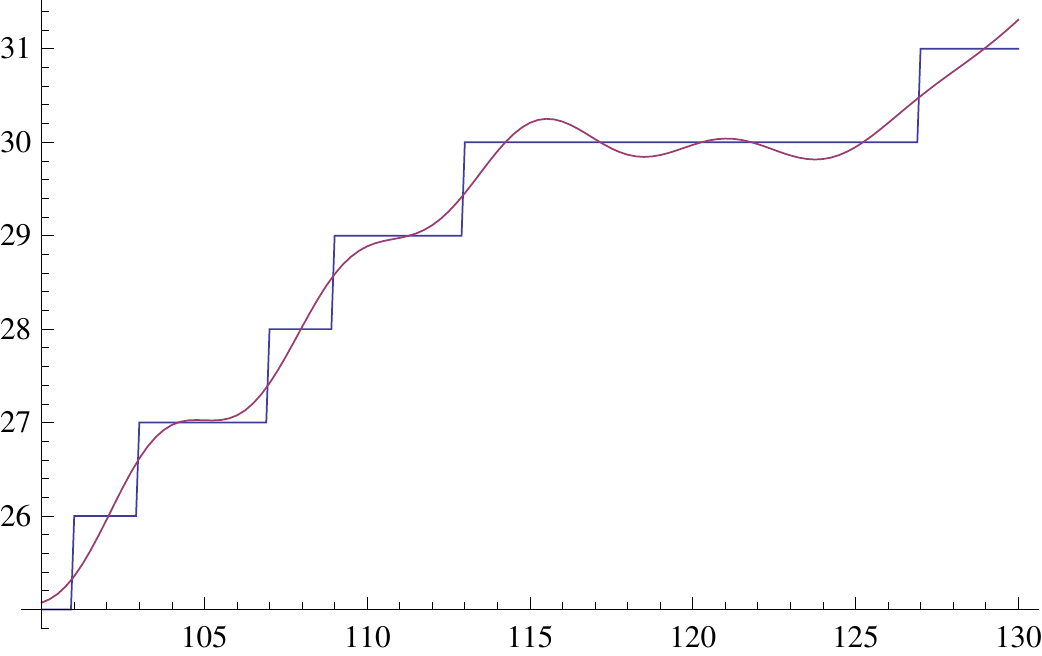}}
    \mbox{\includegraphics[width=\picDblWidth]{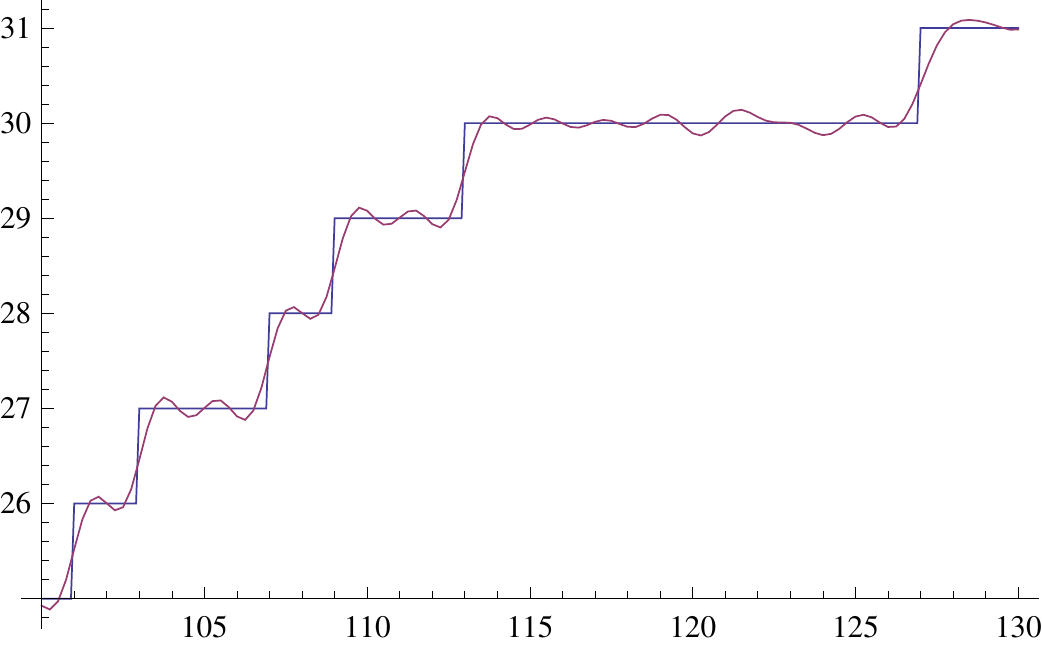}}
  }
  \caption{using 50 and 200 pairs of zeros to approximate $\pi(x)$ for $100 \leq x \leq 130$}
  \label{fig:riemann-vonm100130}
\end{figure*}

In Figure ~\ref{fig:riemann-vonm100130}, we have plotted ~\eqref{E:r-vonm} for $100 \leq x \leq 130$, using $50$ and $200$ pairs of zeta zeros in the sum in ~\eqref{E:r-vonm2}.  This is an interesting interval because, after the prime at $x = 113$, there is a relatively large gap of 14, with the next prime at $x = 127$.  Notice how the graphs tend to level off after $x = 113$ and remain fairly constant until the next prime at $x = 127$.  It is quite amazing that the sums ~\eqref{E:r-vonm} and ~\eqref{E:r-vonm2} involving zeta zeros somehow 'know' that there are no primes between 113 and 127!  This detailed tracking of the step function $\pi(x)$ (or $\pi_0(x)$) is due solely to the sums over zeta zeros in ~\eqref{E:r-vonm2}.

Stan Wagon's book \cite[pp. 512-519]{Wagon} also has an interesting discussion of this topic.

In section \ref{S:CountingPrimesInAP}, we see how to use zeros of Dirichlet $L$-functions to count primes in arithmetic progressions.

\begin{figure*}[ht]
  \centerline{
    \mbox{\includegraphics[width=\picDblWidth]{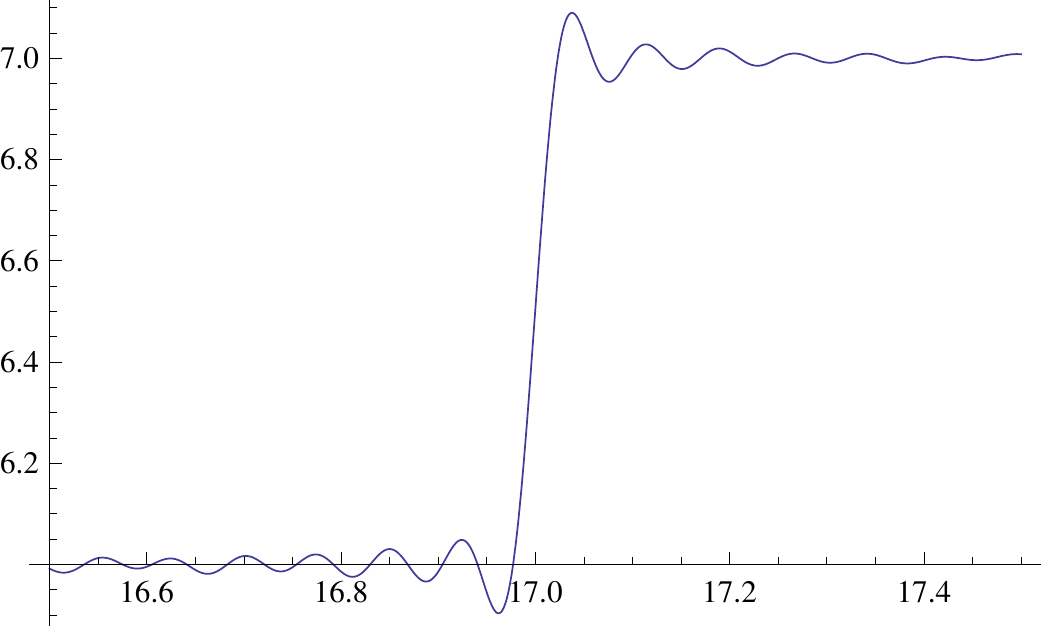}}
    \mbox{\includegraphics[width=\picDblWidth]{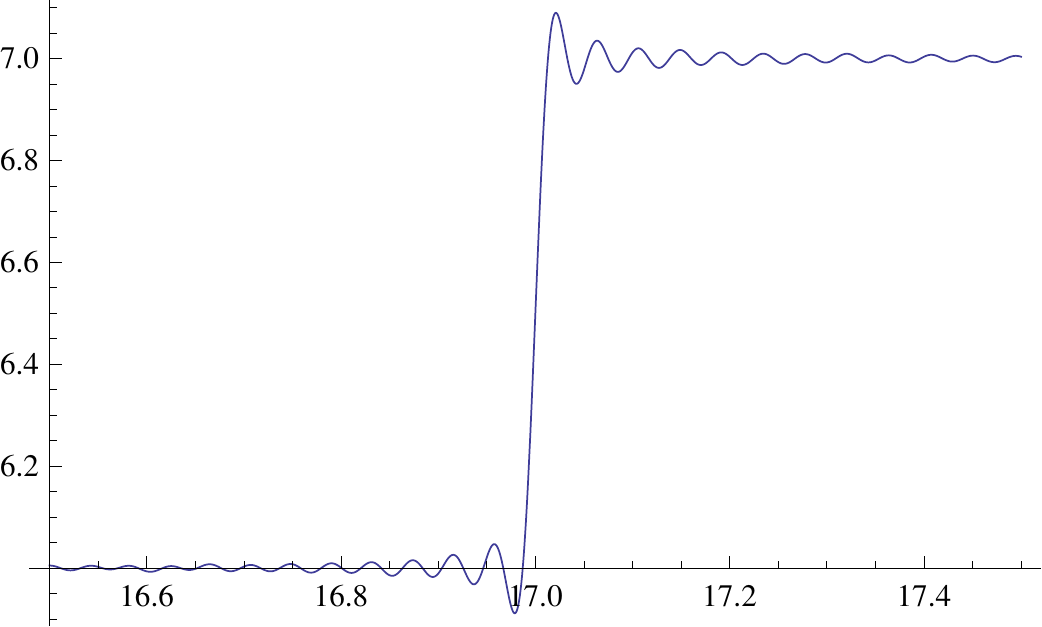}}
  }
  \caption{using 1000 and 2000 pairs of zeros to approximate $\pi(x)$ near the jump at $x = 17$}
  \label{fig:riemann-vonm17-1k}
\end{figure*}

\section{The Gibbs Phenomenon}\label{S:Gibbs}


\begin{figure*}[ht]
  \centerline
  {
    \mbox{\includegraphics[width=\picDblWidth]{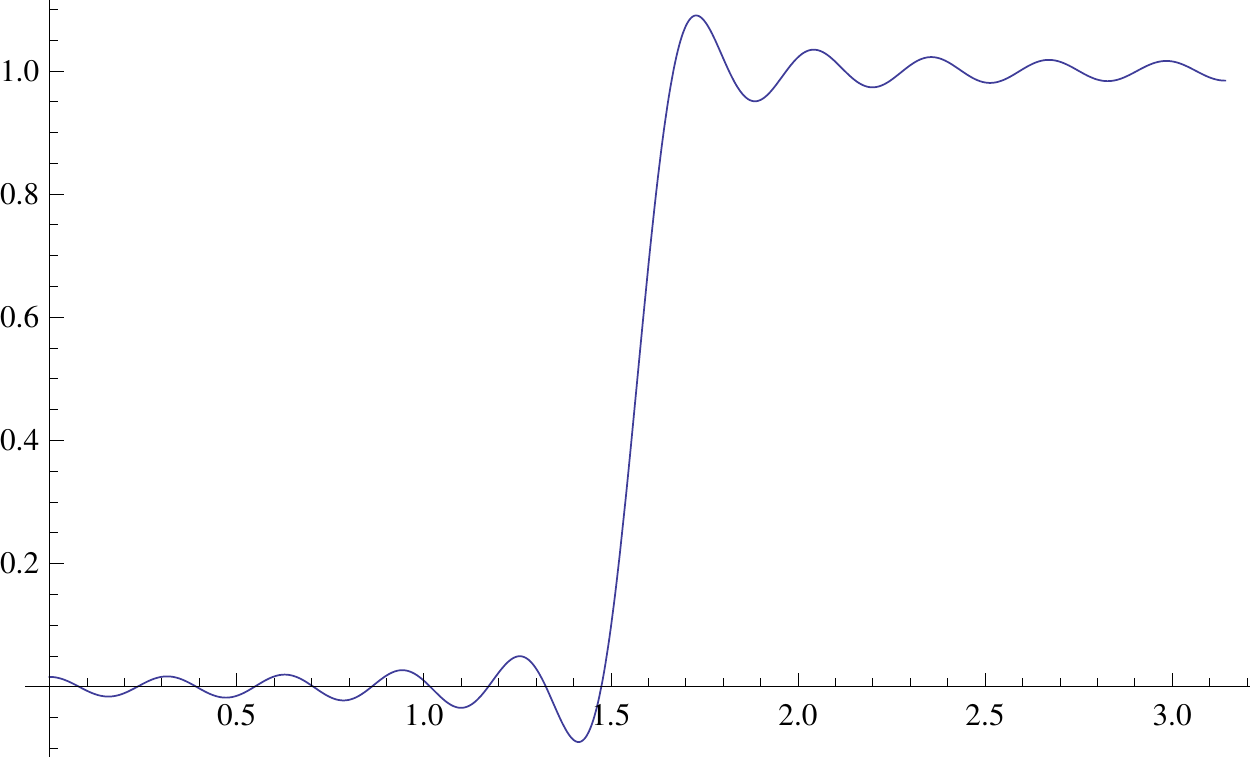}}
    \mbox{\includegraphics[width=\picDblWidth]{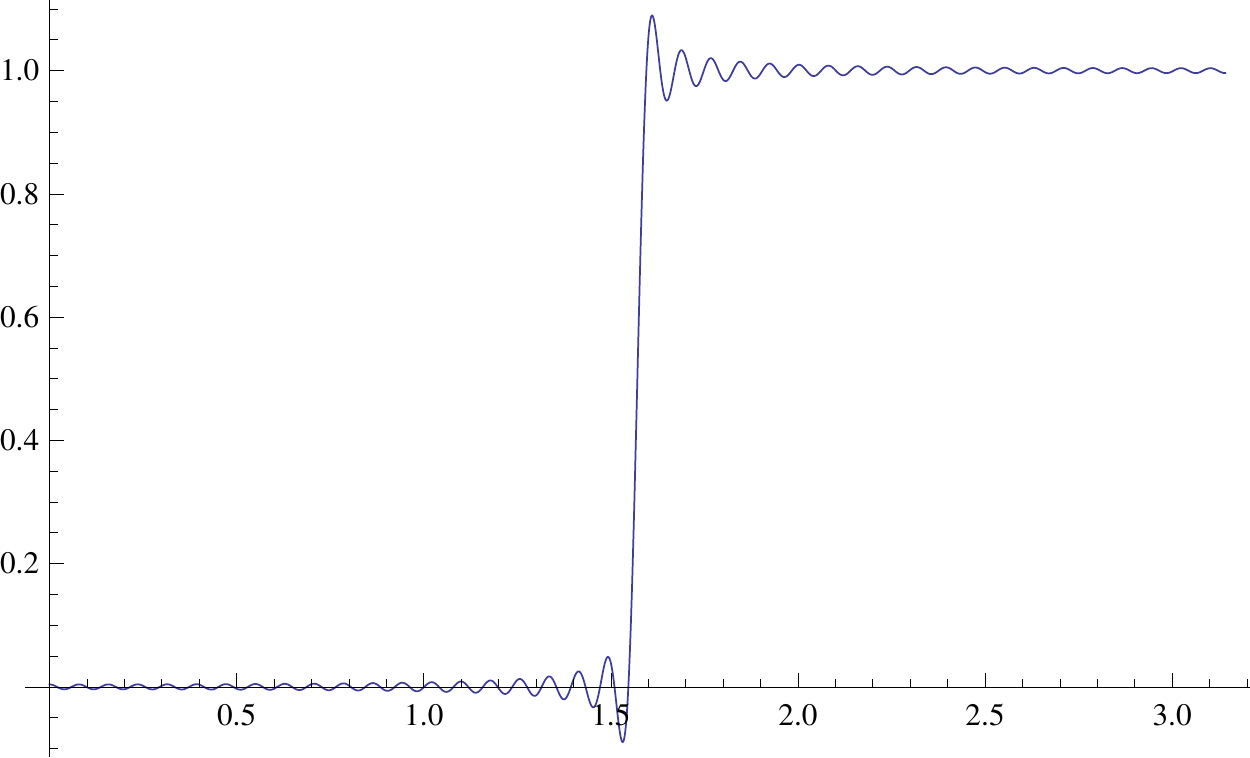}}
  }
  \caption{Partial sums of 10 and 40 terms of a Fourier series with unit jump}
  \label{fig:gibbs1}
\end{figure*}

When a Fourier series represents a step function, the $n^{th}$ partial sum of the Fourier series will undershoot and overshoot before and after the jumps of the step function.  If the step function jumps up or down by 1 unit, there will be an undershoot and an overshoot in the \emph{infinite} sum, each in the amount of
\begin{equation}\label{E:gibbsConstant}
g = -\frac{1}{2} + \frac{1}{\pi} \int_{0}^{\pi} \frac{\sin(t)}{t} \, dt \simeq 0.0894898722
\end{equation}
before after the jump.  This is called the Gibbs phenomenon.  Figure ~\ref{fig:gibbs1} shows the overshoot in the partial sum of the first $n$ terms of a Fourier series that has unit jumps at $x = \pi/2$, $x = 3\pi/2$, $x = 5\pi/2$, etc.  In this figure, the first peak after the jump occurs at about $x = 1.72788$, $y = 1.08991$ for $n = 10$, and at about $x = 1.61007$, $y = 1.08952$ for $n = 40$.  (For any given $n$, the first peak occurs at $x = \frac{\pi}{2} + \frac{\pi}{2 n}$).  If the step function has a jump of $c$, then the undershoot and overshoot are both $c$ times the constant in equation \eqref{E:gibbsConstant}.

Figure ~\ref{fig:riemann-vonm17-1k} shows a closeup of the approximation ~\eqref{E:r-vonm} to $\pi(x)$ near $x = 17$.  (There's nothing special about 17).  This looks very much like the partial sum of an Fourier series for a step function, including a Gibbs-like phenomenon before and after the jump.  In the right-hand side of Figure ~\ref{fig:riemann-vonm17-1k}, which is a sum over 2000 pairs of zeta zeros, the minimum value just before the jump is about $6 - 0.08885$; the maximum value just after the jump is about $7 + .08983$.  This overshoot is visible in the some of the graphs in the sections below.

\section{Rules for Residues}\label{S:RulesForResidues}
Later, we will need to compute the residues of functions that have poles of orders 1 through 4.  Formulas for poles of orders 1 through 3 can be found in advanced calculus books (see, for example, Kaplan \cite[p. 564]{Kaplan} for the formulas for orders 1 and 2 and Kaplan \cite[p. 575, exercise 8]{Kaplan} for the formula for order 3).  For convenience, we state these formulas here.  Their complexity increases rapidly as the order of the pole increases.  However, Kaplan \cite[p. 564]{Kaplan} also gives an algorithm for computing the residue at a pole of any order.  Instead of writing out the formula for order 4, we will give a general \textit{Mathematica} procedure that implements that algorithm.

If $f(z) = \frac{A(z)}{B(z)}$ where $A(z)$ and $B(z)$ are analytic in a neighborhood of $s$, and where $A(s) \neq 0$ and $B(z)$ has a zero of multiplicity 1 at $z = s$, then the residue of $f(z)$ at $z = s$ is
\begin{equation}\label{E:stdResidueFormula1}
res( f(z), s ) = \frac{A(s)}{B'(s)}.
\end{equation}

If the conditions are the same as above, but $B(z)$ has a zero of multiplicity 2, then the residue is
\begin{equation}\label{E:stdResidueFormula2}
res( f(z), s ) = \frac{6 A'(s) B''(s)-2 A(s) B^{(3)}(s)}{3 B''(s)^2}.
\end{equation}

If the conditions are the same as above, but $B(z)$ has a zero of multiplicity 3, then the residue is
\begin{equation}\label{E:stdResidueFormula3}
 \frac{120 B^{(3)}(s)^2 A''(s)-60 B^{(4)}(s) B^{(3)}(s) A'(s)-12 A(s)
   B^{(5)}(s) B^{(3)}(s)+15 A(s) B^{(4)}(s)^2}{40 B^{(3)}(s)^3}.
\end{equation}

Here is \textit{Mathematica} code that computes the formula for any order $m$, $m > 0$.
\smaller
\begin{verbatim}
residueFormula[m_Integer?Positive, z_Symbol] :=
Module[
  (* compute the residue of A[z]/B[z] where
     B has a zero of order m and A != 0 *)
  { aa, a, bb, b, k, ser, res },
  aa[z] = Sum[a[k] z^k, {k, 0, m - 1}];   (* need z^m through z^(m-1) *)
  bb[z] = Sum[b[k] z^k, {k, m, 2*m - 1}]; (* need z^m through z^(2*m-1) *)
  ser = Series[aa[z]/bb[z], {z, 0, -1}];  (* need z^m through z^(-1) *)
  res = Coefficient[ser, z, -1];  (* residue = coefficient of z^(-1) *)
  (* replace a[ ] and b[ ] with derivatives *)
  res = res /. {a[k_] -> D[A[z], {z, k}]/k!, b[k_] -> D[B[z], {z, k}]/k!};
  Simplify[res]
]
\end{verbatim}
\larger
For example, the following command will generate the expression in residue rule \eqref{E:stdResidueFormula2}:
\begin{verbatim}
  residueFormula[2, s].
\end{verbatim}
In various sections below, we will first define $A(s)$ and $B(s)$, then use \verb+residueFormula+ to expand the result.  We do this, for example, in section \ref{S:lambdaTimesTau}.

Finally, later (sections \ref{S:TallyingTheSquarefreeDivisors} and \ref{S:lambdaTimesTau}), we will need a formula for the residue at $z = 1$ of expressions like
\[
\frac{\zeta(z)^2}{\zeta(2z)}.
\]
where the residue (here, of order 2) arises not from a zero in the denominator, but from a pole in the numerator at $z = 1$.  If $f(z)$ has a pole of order $n$ at $z = s$, then the residue at $z = s$ is \cite[p. 564]{Kaplan}
\begin{equation}\label{E:stdResidueFormulaN}
\frac{1}{(n-1)!}  \lim_{z\to s} \, \frac{d ^{n-1}\left(f(z)(z - s)^n\right)}{d z^{n-1}}.
\end{equation}

\section{Perron's Formula}

Given a sequence of numbers $a_n$ for $n = 1, 2, \dots$, a Dirichlet series is a series of the form
\[
\sum _{n=1}^{\infty } \frac{a_n}{n^s}
\]
where $s$ is a complex variable. If a Dirichlet series converges, then it converges in some half-plane $Re(s) > c_0$.

Now, suppose we have a Dirichlet series
\begin{equation}\label{E:perron1}
\sum _{n=1}^{\infty } \frac{a_n}{n^s}=F(s)
\end{equation}
that holds for $Re(s) > c_0$.

Perron's formula inverts \eqref{E:perron1} to give a formula for the summatory function of the $a_n$, as a function of $x$:
\begin{equation}\label{E:perron2}
\sideset{}{'} \sum _{n=1}^x a_n
 = \lim_{T\to \infty } \, \frac{1}{2 \pi  i}
   \int_{c-i T}^{c+i T} F(s) \frac{x^s}{s} \, ds.
\end{equation}

This holds for $x > 0$ and for $c > c_0$.
Whenever we see an $x$ as the upper limit on a summation, we mean that the sum extends up to $\lfloor x\rfloor$.

The prime on the summation means that, when $x$ is an integer, instead of adding $a_x$ to the sum, we add $\frac{1}{2}a_x$.  Let $G(x)$ be the sum on the left-hand side of \eqref{E:perron2}.  Suppose $x$ is an integer and that $a_x > 0$.  Then $G(x-.1)$ does not include the term $a_x$ at all, $G(x)$ has $\frac{1}{2}a_x$ as its final term, and $G(x+.1)$ has $a_x$ as its final term.  Thus, at $x$, $G$ has a jump discontinuity, and the value of $G(x)$ is the average of, for example, $G(x-.1)$ and $G(x+.1)$.  Of course, this argument also holds with .1 replaced by any $\epsilon$ between 0 and 1.

To summarize, the left side of \eqref{E:perron2} is a step function that may jump when $x$ is an integer.  If this step function jumps from $a$ to a new value $b$, the right side converges to the midpoint, $(a+b)/2$, as $T$ approaches infinity.  The important point is that the integral on the right side of \eqref{E:perron2} equals the sum on the left, \textit{even where the sum is discontinuous}.

We will estimate the integral in \eqref{E:perron2} by integrating around a rectangle in the complex plane. Recall that, under appropriate conditions, the integral around a closed path is equal to $2 \pi i$ times the sum of the residues at the poles that lie inside that path.  We will use the sum of the residues to approximate the line integral.

In all of our examples, the $a_n$ will be integers.  It is quite remarkable that the line integral \eqref{E:perron2} in the complex plane can be used to sum the values on the left side of \eqref{E:perron2}.  It is not at all obvious that a complex line integral has anything whatever to do with sums of the integers $a_n$.  Nevertheless, Perron's formula is proved in \cite[pp. 138-140]{MontgomeryAndVaughan}, among other places.  For the reader who is still skeptical, we will give numerical examples illustrating that the left and right sides of \eqref{E:perron2} can, indeed, be quite close when a modest number of poles is taken into account.

Perron's formula was first published by Oskar Perron in 1908 \cite{Perron}.

The plan for most of the rest of this paper is as follows.  In each section, we will apply Perron's formula to an arithmetic function whose Dirichlet series yields a sum over zeta zeros.  This means that the Dirichlet series must have the form
\[
  \sum_{n=1}^{\infty } \frac{a_n}{n^s}=\frac{F(s)}{G(s)}
\]
where $G(s)$ has factors of the zeta function.  The point in having $\zeta(s)$ (or $\zeta(2s)$, etc.) in the denominator is that, each zeta zero will give rise to a pole in the right-hand side.  When we sum the residues at these poles, we get one term in the sum for each zeta zero.

There are two sources we will frequently cite that conveniently gather together many such arithmetic functions: \cite{Gould} and \cite{McCarthy}.  The former lists the Dirichlet series in a convenient form.  The latter derives them, along with many asymptotic formulas.  Below, we will work with some of the most common arithmetic functions.  These two references also discuss many more arithmetic functions which we do not consider, such as generalizations of von Mangoldt's function and generalizations of divisor functions.

In some cases, Perron's formula leads to a sum of terms involving zeta zeros, which, like Perron's integral, is provably close to the summatory function.  In many cases, the sum appears to be close to the summatory function, but no proof has been written down.  Finally, in other cases, the sum over zeta zeros appears to diverge.


\ifthenelse {\boolean{BKMRK}}
  { \section{Using Perron's Formula to Compute \texorpdfstring{$\psi(x)$}{psi(x)}} }
  { \section{Using Perron's Formula to Compute $\psi(x)$} }

We will now present one of the best-known applications of Perron's formula.  This example is discussed in detail by Davenport \cite[pp. 104-110]{Davenport}, Ivi\'{c} \cite[p. 300]{Ivic}, and Montgomery and Vaughan \cite[pp. 400-401]{MontgomeryAndVaughan}.

Chebychev's psi ($\psi$) function is important in number theory.  Many theorems about $\psi(x)$ have a corresponding theorem about $\pi(x)$, but $\psi$ is easier to work with.  $\psi$ is most easily defined in terms of the von Mangoldt Lambda function, which is defined as follows:


\[
\Lambda(n)=
  \begin{cases}
    \log{p}  &\text{if $n$ is the $k^{th}$ power of a prime (where $k > 0$),}\\
    0 &\text{otherwise.}
  \end{cases}
\]

$\psi(x)$ is the summatory function of $\Lambda$:

\[
\psi (x)=\sum _{n=1}^x \Lambda (n).
\]

It is to be understood that the upper limit on the summation is really the greatest integer not exceeding $x$, that is,
 $\lfloor x\rfloor$.

One can easily check that $\psi(x)$ is the log of the least common multiple of the integers from $1$ through $x$.  For example, for $n = 1, 2, \ldots, 10$, the values of $\Lambda(n)$ are 0, $\log(2)$, $\log(3)$, $\log(2)$, $\log(5)$, 0, $\log(7)$, $\log(2)$, $\log(3)$, and 0, respectively.  The sum of these values is
\[
\psi(10) = 3 \cdot \log(2) + 2 \cdot \log(3) + \log(5) + \log(7) = \log(2^3 \cdot 3^2 \cdot 5 \cdot 7) = \log(2520)
\]
where 2520 is the least common multiple of the integers from 1 through 10.

\begin{figure*}[ht]
  \mbox{\includegraphics[width=\picDblWidth]{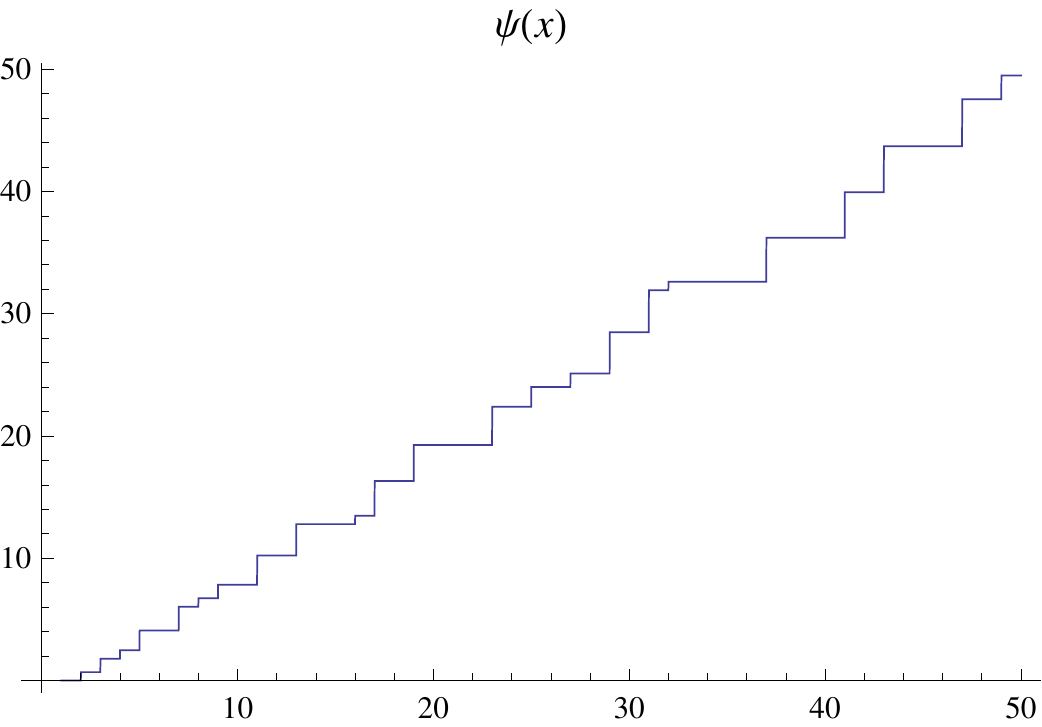}}
  \caption{The graph of $\psi(x)$}
  \label{fig:psiFunction}
\end{figure*}

Figure ~\ref{fig:psiFunction} shows a graph of $\psi(x)$.

The Dirichlet series for $\Lambda$ is \cite[Theorem 294, p. 253]{HardyAndWright}, \cite[p. 105]{Davenport}

\begin{equation}\label{E:lambdaDirichletSeries}
\sum _{n=1}^{\infty } \frac{\Lambda (n)}{n^s}=-\frac{\zeta'(s)}{\zeta (s)}.
\end{equation}

This series converges for $\Re(s) > 1$.

When we apply Perron's formula to this Dirichlet series, we will get a formula that approximates the corresponding summatory function, $\psi(x)$.  But recall that, in Perron's formula \eqref{E:perron2}, the summation includes half of the last term if $x$ is an integer that causes the step function to jump.  So, Perron's formula actually gives a formula for $\psi_0$, a slight variant of $\psi$:

\[
\psi_0(x)=
  \begin{cases}
    \psi(x) - \frac{1}{2}\Lambda(x)  &\text{if $x$ is the $k^{th}$ power of a prime (where $k > 0$),}\\
    \psi(x) &\text{otherwise.}
  \end{cases}
\]
That is, where $\psi(x)$ jumps up by the amount $\Lambda(x)$, $\psi_0(x)$ equals the average of the two values $\psi(x)$ and $\psi(x-1)$.  Recall that we made a similar adjustment to $\pi(x)$ and called it $\pi_0(x)$.

Given the Dirichlet series \eqref{E:lambdaDirichletSeries}, Perron's formula states that
\begin{equation}\label{E:psiIntegral}
\psi_0(x)= \lim_{T\to \infty } \, \frac{1}{2 \pi  i}
   \int_{c-i T}^{c+i T} -\frac{\zeta'(s)}{\zeta(s)} \frac{x^s}{s} \, ds.
\end{equation}

The computation of $\psi_0(x)$ (or $\psi(x)$) now reduces to computing a line integral on the segment from $c-i T$ to $c + i T$ in the complex plane. Here is our plan for computing that integral.  First, we will choose $c$.  Then, we will form a rectangle whose right side is the line segment from $c - i T$ to $c + i T$.  Next, we will adjust the rectangle to make it enclose some of the poles of our integrand.  (We must make sure that the rectangle does not pass through any zeta zeros.)   The left side of the rectangle will extend from, say, $a + i T$ to $a - i T$.  We will integrate around this rectangle in the counterclockwise direction.

The integral around the rectangle is
\begin{equation}\label{E:psiIntegralRectangle}
  \frac{1}{2 \pi i} \int_{c - i T}^{c + i T} I(s) \, ds +
  \frac{1}{2 \pi i} \int_{c + i T}^{a + i T} I(s) \, ds +
  \frac{1}{2 \pi i} \int_{a + i T}^{a - i T} I(s) \, ds +
  \frac{1}{2 \pi i} \int_{a - i T}^{c - i T} I(s) \, ds.
\end{equation}
where the integrand $I(s)$ is
\begin{equation}\label{E:psiIntegrand}
- \frac{\zeta '(s)}{\zeta (s)} \frac{x^s}{s}.
\end{equation}

The exact value of this integral can be computed easily: it is the sum of the residues at the poles that lie inside the rectangle.

It will turn out that the last three integrals in \eqref{E:psiIntegralRectangle} are small, so we have, approximately,
\[
\psi_0(x) \simeq \frac{1}{2 \pi  i} \int_{c-i T}^{c+i T} I(s) \, ds \simeq \text{sum of residues}.
\]

Our task now is to locate the poles of the integrand, determine the residues at those poles, then adjust the rectangle so it encloses those poles.

First, the integrand has a pole at $s = 0$.

Second, it is known \cite[Theorems 281 and 283, p. 247]{HardyAndWright} that
\[
\zeta (s)=\frac{1}{s-1} + O(1)
\]
and that
\[
\zeta '(s)=\frac{-1}{(s-1)^2} + O(1).
\]
So,
\[
\frac{\zeta '(s)}{\zeta (s)}=\frac{1}{s-1} + O(1).
\]
Therefore, the integrand has a pole at $s = 1$.

Finally, the denominator is 0 at every value of $s$ such that $\zeta (s) = 0$, that is, at each zero of the zeta function.  We will see that every zeta zero, both real and complex, gives rise to a term in a sum for $\psi(x)$.

Next, we list the residues at these poles.  At $s = 0$, the residue is $-\frac{\zeta'(0)}{\zeta(0)} = -\log(2 \pi) \approx -1.837877$.

At $s = 1$, the residue is $x$.

To compute the residue of this integrand at a typical zeta zero, we will use residue formula \eqref{E:stdResidueFormula1}.  But to use this formula, we must assume that if $\rho$ is a zero of zeta, then it is a simple zero, that is, a zero of multiplicity 1.  All of the trillions of complex zeta zeros that have been computed are, indeed, simple.  So are the real zeros.  In any case, since our goal is to get a formula that involves a sum over only a modest number of zeros, this assumption is certainly valid for the zeros we will be considering.

To get the residue of \eqref{E:psiIntegrand} at $s = \rho$, we let $B(s) = \zeta(s)$ and let $A(s)$ equal everything else in the integrand \eqref{E:psiIntegrand}.  Thus, $A(s) = -\zeta'(s) x^s/s$.   With these choices, $A(s)/B(s)$ is equal to the integrand.  The residue at the zeta zero $s = \rho$ will be
\[
\frac{A(\rho)}{B'(\rho)} = -\frac{\zeta '(\rho ) x^{\rho }}{\rho } \frac{1}{\zeta '(\rho )} = -\frac{x^{\rho}}{\rho}.
\]

Notice that the above argument holds equally well for both complex and real zeta zeros.

We will now adjust our rectangle to make it enclose some of the above poles.

First, we'll make sure that the rectangle around which we integrate encloses the poles at $s = 0$ and $s = 1$, where the residues are $-\frac{\zeta '(0)}{\zeta (0)}$ and $x$, respectively.  We can do this, for example, by moving the right boundary of the rectangle to the line $\Re(s)=2$.

Next, we'll adjust the left-hand boundary of the rectangle so it encloses the first $M$ real zeros $-2k = -2, -4, \cdots, -2M$.  The sum of the residues at the corresponding poles is
\[
- \sum _{k=1}^M \frac{x^{-2 k}}{-2 k}.
\]

We also adjust the top and bottom boundaries of the rectangle to make it enclose the first $N$ \emph{pairs} of complex zeros.  Note that, for a given $N$, there will be a range of values of $T$ such that there are $N$ zeros whose imaginary parts, in absolute value, are less than $T$.  We can choose $N$, or we can choose $T$, whichever is more convenient.

The $k^{th}$ pair of conjugate zeros gives rise to the two terms
\[
-\frac{x^{\rho_k}}{\rho_k}
\]
(where the bar denotes the complex conjugate), and
\[
-\frac{x^{\overline{\rho_k}}}{\overline{\rho_k}}.
\]

It is easy to verify that these terms are conjugates of each other.  When we add these terms, their imaginary parts will cancel and their real parts will add.  Therefore, the sum of these two conjugate terms will be
\[
-2 \Re( \frac{x^{\rho_k}}{\rho_k})
\]

The sum of the residues at the complex zeros of zeta will be
\[
-2 \Re(\sum _{k=1}^N \frac{x^{\rho_k}}{\rho_k}).
\]

When we take into account all of these poles, the sum of the above residues is
\begin{equation}\label{E:psiZetaZeroSumApprox}
\psi_0(x)\simeq x-\frac{\zeta '(0)}{\zeta (0)}-2
   \Re\left(\sum _{k=1}^N \frac{x^{\rho _k}}{\rho _k}\right)
   -\sum _{k=1}^{M } \frac{x^{-2 k}}{-2 k}.
\end{equation}
We can extend the left-hand side of the rectangle as far to the left as we wish, enclosing an arbitrarily large number of real zeta zeros.  The second sum converges quickly to $\frac{1}{2}\log(1-\frac{1}{x^2})$ as $M$ approaches $\infty$.  It can be proven that it is also valid to let $N$ approach $\infty$ in \eqref{E:psiZetaZeroSumApprox}.  We then get the following exact formula for $\psi_0(x)$:
\begin{equation}\label{E:psiZetaZeroSumProved}
\psi_0(x) = x - \frac{\zeta '(0)}{\zeta (0)}
 - 2 \Re\left(\sum _{k=1}^{\infty} \frac{x^{\rho _k}}{\rho _k}\right)
 - \sum _{k=1}^{\infty} \frac{x^{-2 k}}{-2 k}.
\end{equation}
See, for example, Montgomery and Vaughan \cite[p. 401]{MontgomeryAndVaughan}.  Further, Davenport \cite[Equation 11, p. 110]{Davenport} gives the following bound for the error as a function of $x$ and $T$:
\begin{equation}\label{E:psiZetaZeroSumWithError}
\psi_0(x) = x-\frac{\zeta '(0)}{\zeta (0)}
   - \sum _{|t_k| \leq T} \frac{x^{\rho_k}}{\rho_k}
   - \sum _{k=1}^{\infty } \frac{x^{-2 k}}{-2 k}
   + O( \frac{x \log^2(x T)}{T} ),
\end{equation}
where $x$ is an integer and $t_k$ is the imaginary part of $\rho_k$, the $k^{th}$ complex zeta zero.  Note that, for a fixed $x$, the error term approaches 0 as $T$ approaches $\infty$, consistent with \eqref{E:psiZetaZeroSumProved}.

The second sum in \eqref{E:psiZetaZeroSumProved} is too small to have any visible effect on the graphs that we will display, so we will use $M = 0$ when we graph this formula.

The dominant term, $x$, shows roughly how fast $\psi(x)$ increases.  In fact, it is known that $\psi(x)$ is asymptotic to $x$ (written $\psi(x) \thicksim x$), which means that $\lim_{x\to \infty } \, \frac{\psi (x)}{x}=1$.

The first two residues, $x$ and $-\frac{\zeta '(0)}{\zeta (0)}$, arising from poles at $s = 0$ and $s = 1$, give a good linear approximation to $\psi(x)$.  That approximation, which comes from equation \eqref{E:psiZetaZeroSumProved} with $N = 0$ and $M = 0$, is shown at the top left of Figure ~\ref{fig:psiApproxGraphs}.  This smooth approximation, while good on a large scale, misses the details, namely, the jumps at the primes and the prime powers.  To replicate those details, the approximation needs some of the terms with zeta zeros in equation \eqref{E:psiZetaZeroSumProved}.

Ingham \cite[pp. 77-80]{Ingham} also proves equation \eqref{E:psiZetaZeroSumProved}.  On page 80, he makes the intriguing remark,
\begin{quote}
The 'explicit formula' [\eqref{E:psiZetaZeroSumProved}] suggests that there are connections between the numbers $p^m$ (the discontinuities of $\psi_0(x)$) and the numbers $\rho$.  But no relationship essentially more explicit than [\eqref{E:psiZetaZeroSumProved}] has ever been established between these two sets of numbers.
\end{quote}
This is still true today.

It is also worth pointing out that in \eqref{E:psiZetaZeroSumProved}, the sum over complex zeta zeros does not require the calculation of the zeta function itself.  In the approximations in later sections of this paper, the sums \emph{do} involve zeta or its derivatives, which makes it harder to prove the existence of an exact formula like \eqref{E:psiZetaZeroSumProved}.

\subsection{Numerical Experiments}

We will now work through some numerical examples.  Given a rectangle, we will compute line integrals on each of its sides. We will see that the integral on the first segment, from $c - i T$ to $c + i T$, is significant, and that the integrals on the other three sides are small.  This means that the first integral alone is approximately equal to the integral around the entire rectangle, which, in turn, equals the sum of the residues inside the rectangle. (Note that Perron's formula, equation \eqref{E:psiIntegral}, includes a factor of $1/(2 \pi i)$, so the right side of \eqref{E:psiIntegral} is the sum of the residues, not $2 \pi i$ times the sum of the residues.)  As we integrate around the rectangle in a counterclockwise direction, we'll let $I_1$, $I_2$, $I_3$, and $I_4$ be the line integrals on the right, top, left, and bottom sides, respectively.

The Dirichlet series \eqref{E:lambdaDirichletSeries} converges for $\Re(s) > 1$, so in the integral in \eqref{E:psiIntegral}, we must have $c > 1$.  Let's take $c = 2$.  Let's extend the left-hand side of the rectangle to $\Re(s) = -1$, the rectangle encloses no real zeros.  Then we will have $M = 0$ in \eqref{E:psiZetaZeroSumProved}.

Equation \eqref{E:psiIntegral} suggests that, for a given $x$, if $T$ is large, then the sum of the residues may be close to $\psi_0(x)$. We will experiment with different values of $T$ in \eqref{E:psiIntegral}, that is, different values of $N$ in equation \eqref{E:psiZetaZeroSumProved}.

Let's start with $x = 10$.  Then $\psi_0(x) = \psi(x) = \log(2520) \approx 7.832$.  If $T = 14$, then the lower-right corner of the rectangle is at $2 - 14 i$, and its top left corner is at $-1 + 14 i$.  The first complex zeta zero is at $s \approx 1/2 + 14.135 i$, so this rectangle encloses no complex zeta zeros.  Therefore, we will have $N = 0$ in equation \eqref{E:psiZetaZeroSumProved}.  Numerical integration tells us that $I_1 \approx 8.304$, $I_2 \approx -.069 - 0.305 i$, $I_3 \approx -0.004$, and $I_4 \approx -.069 + 0.305 i$.  The sum of these four integrals is
\[
I_1 + I_2 + I_3 + I_4 \approx 8.162.
\]

As a check on our numerical integrations, we can substitute $x = 10$, $N = 0$, and $M = 0$ into equation \eqref{E:psiZetaZeroSumProved}.  This, in effect, computes the sum of the four integrals by using residues.  We get 8.162, as expected.  Cauchy was right: we \emph{can} compute a contour integral by adding residues!

Now let's extend the rectangle in the vertical direction to enclose one pair of complex zeta zeros.  We will now have $N = 1$ in equation \eqref{E:psiZetaZeroSumProved}.  We can, for example, take $T = 15$, because the first two zeta zeros are $1/2 + 14.135 i$ and $1/2 + 20.022 i$. A computation of the integral over the sides of this rectangle gives $I_1 \approx 7.602$, $I_2 \approx .078 - 0.213 i$, $I_3 \approx -0.006$, and $I_4 \approx .078 + 0.213 i$.  The total is $\approx 7.751$, a little closer to our goal, $\psi_0(x) \approx 7.832$.

Because the second zeta zero has imaginary part roughly $20.022$, we could have enclosed one pair of zeros by extending our rectangle up to $T = 20$ instead of to $T = 15$.  Had we done this, we would have found that $I_1 \approx 7.518$, $I_2 \approx .119 + 0.087 i$, $I_3 \approx -0.005$, and $I_4 \approx .119 - 0.087 i$.  However, the \emph{sum} of these four integrals is unchanged ($\approx 7.751$), because the rectangle with $T = 20$ encloses exactly the same poles as the rectangle with $T = 15$.

Let's set $T = 100$ and extend the rectangle again, so its lower-right and upper-left corners are $2 - 100 i$ and $-1 + 100 i$.  There are now $N = 29$ zeta zeros with positive imaginary part less than $T$.  The integrals over the sides of this rectangle are $I_1 \approx 7.815$, $I_2 \approx -.014 - 0.018 i$, $I_3 \approx -0.005$, and $I_4 \approx -.014 + 0.018 i$.  The sum is $\approx 7.782$, not far from the goal of $7.832$.

In Figure ~\ref{fig:psiApproxGraphs} below, we show the step function $\psi(x)$, along with approximations to $\psi(x)$ using formula \eqref{E:psiZetaZeroSumProved} with $M = 0$.  To create the four graphs, we set $N = 0$, $N = 1$, $N = 29$, and $N = 100$. We computed \eqref{E:psiZetaZeroSumProved} for $x = 2$ through $x = 50$.
Notice that, if $x$ is a value that causes the step function to jump, the smooth curve is usually quite close to the midpoint of the $y$ values, that is, $\psi_0(x)$.

\begin{figure*}[ht]
  \centerline{
    \mbox{\includegraphics[width=\picDblWidth]{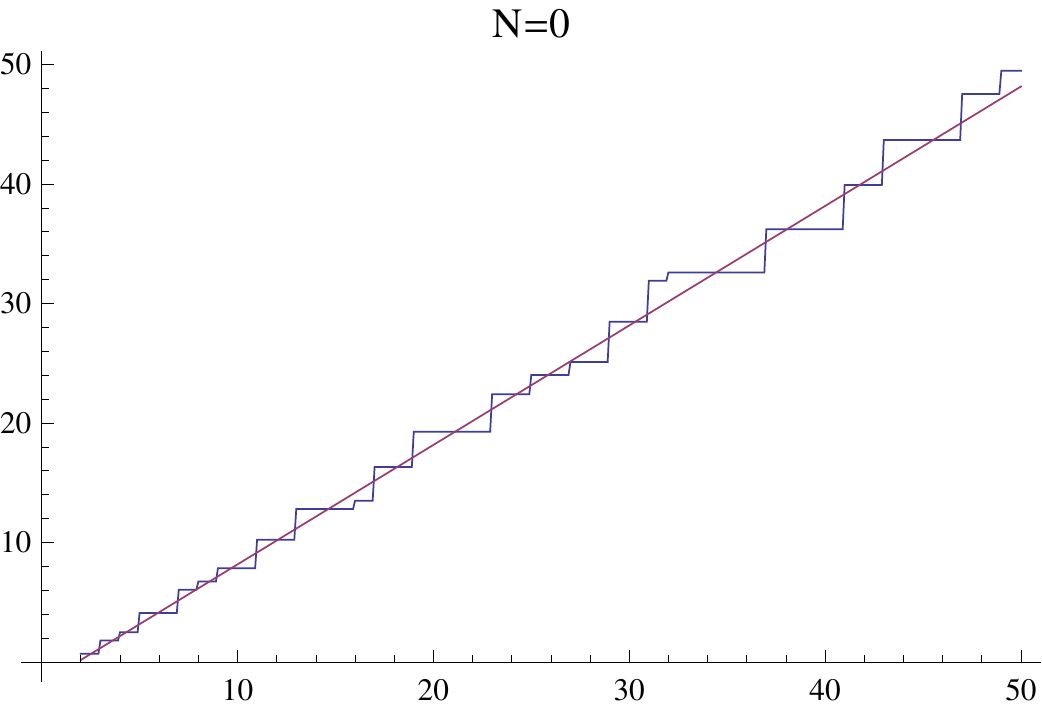}}
    \mbox{\includegraphics[width=\picDblWidth]{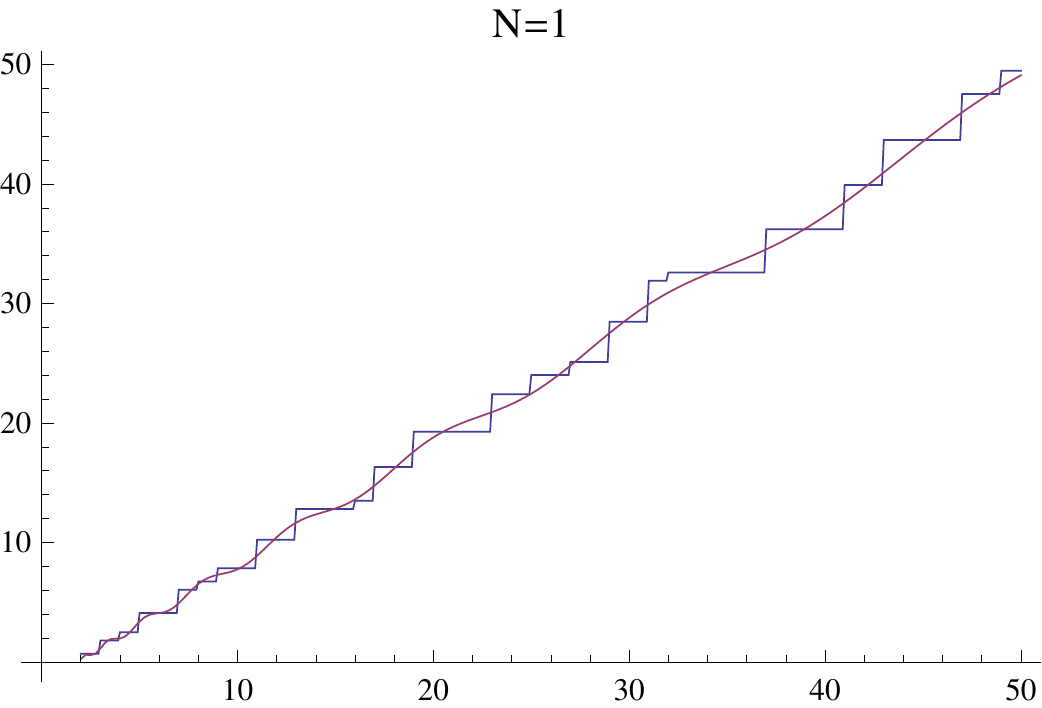}}
  }
  \centerline{
    \mbox{\includegraphics[width=\picDblWidth]{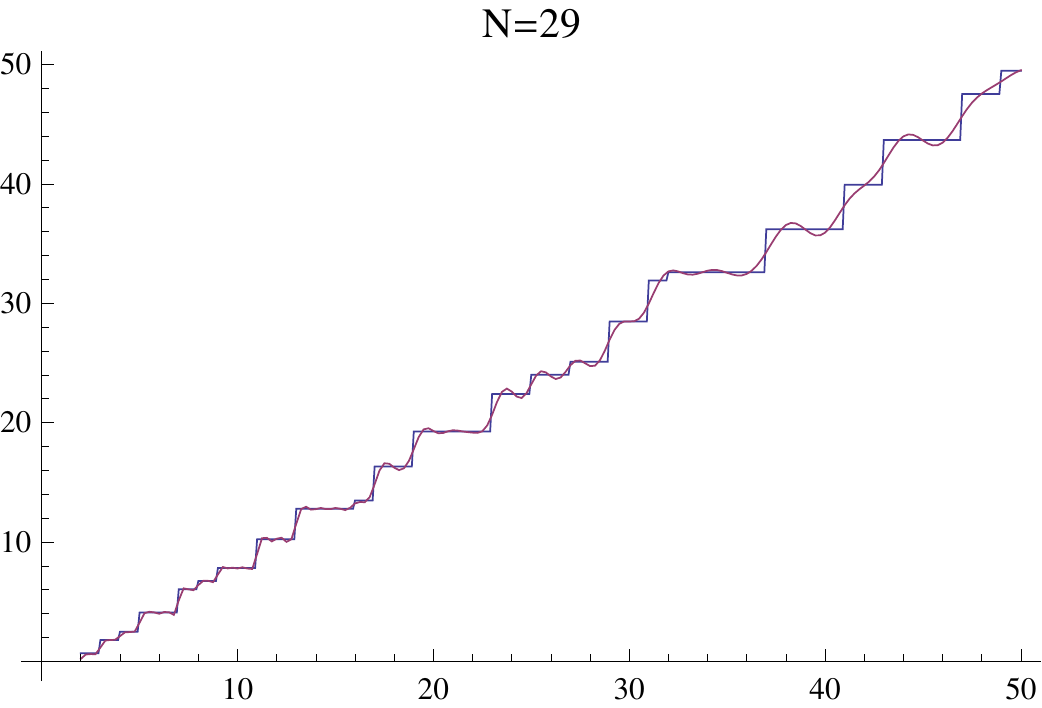}}
    \mbox{\includegraphics[width=\picDblWidth]{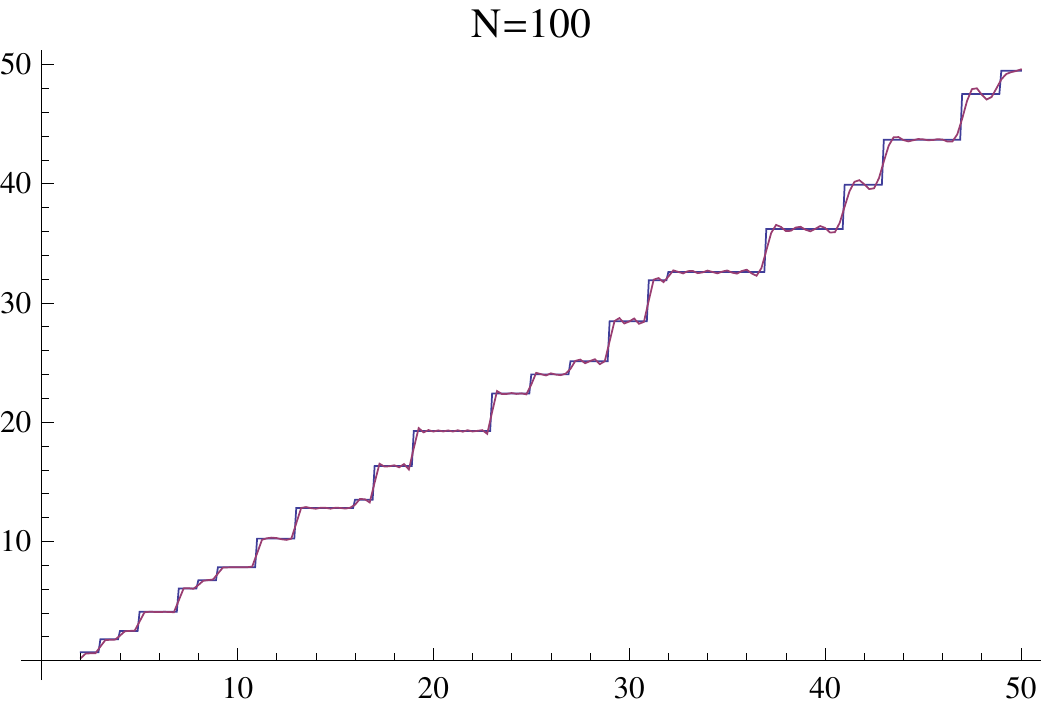}}
  }
  \caption{approximating $\psi(x)$ with $N = 0$, $N = 1$, $N = 29$, and $N = 100$ pairs of zeta zeros}
  \label{fig:psiApproxGraphs}
\end{figure*}


In Figure ~\ref{fig:psiApproxGraphs}, notice that for a fixed $N$ (or $T$), the agreement gets worse as $x$ increases.  Second, for a given $x$, as $N$ (or $T$) increases, we will usually see better agreement between the plotted points and the step function.

\begin{figure*}[ht]
    \mbox{\includegraphics[width=\picDblWidth]{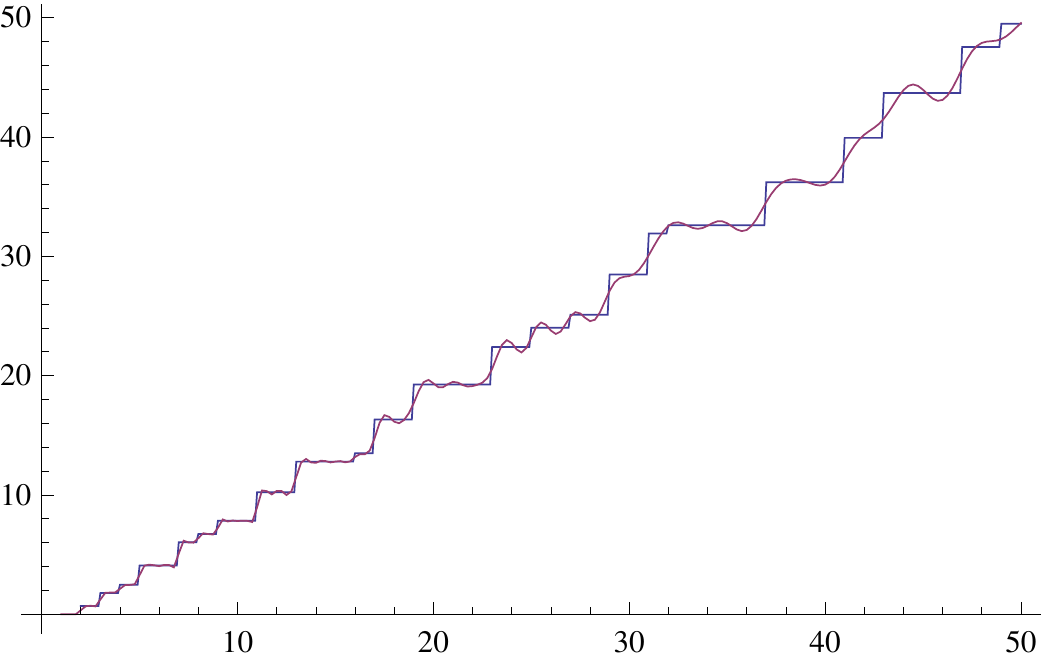}}
  \caption{Approximating $\psi(x)$ with an integral: $c = 2$, $T = 100$}
  \label{fig:psiIntegralApproximation}
\end{figure*}

What would happen if we used this integral approximation,
\begin{equation}\label{E:psiIntegralApprox}
\psi_0(x) \simeq \frac{1}{2 \pi  i}
   \int_{c-i T}^{c+i T} -\frac{\zeta'(s)}{\zeta(s)} \frac{x^s}{s} \, ds,
\end{equation}
based on \eqref{E:psiIntegral}, instead of the summing residues?  In Figure ~\ref{fig:psiIntegralApproximation}, we set $c = 2$ and $T = 100$ in \eqref{E:psiIntegralApprox} and perform the integration for many $x \le 50$.  Recall that a rectangle with $T = 100$ encloses $N = 29$ pairs of complex zeta zeros, so in this sense, it is fair to compare this approximation with the third graph in Figure ~\ref{fig:psiApproxGraphs}.
We can see that, at least in this example and for these $x$, the approximations are about equally good.


Incidentally, there's nothing special about $c = 2$ here.  If we had computed the integral using $c = 1.1$ or $c = 1.5$ instead, the graphs at this scale would look about the same.  If we had used $c = 2.5$, the graph would be a little bumpier.

\begin{figure*}[ht]
  \mbox{\includegraphics[width=\picDblWidth]{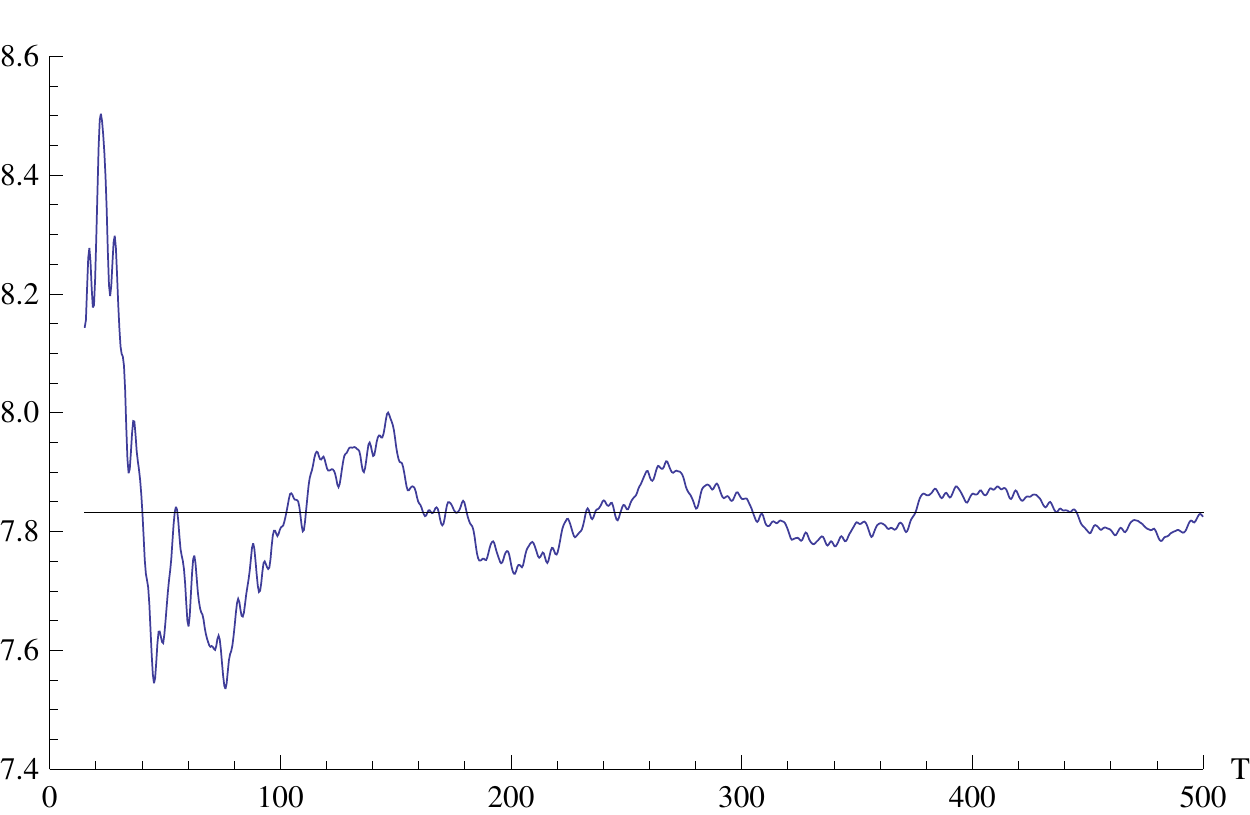}}
  \caption{$x = 10.5$, $c = 1.1$}
  \label{fig:PsiIntVsT}
\end{figure*}

\begin{figure*}[ht]
  \mbox{\includegraphics[width=\picDblWidth]{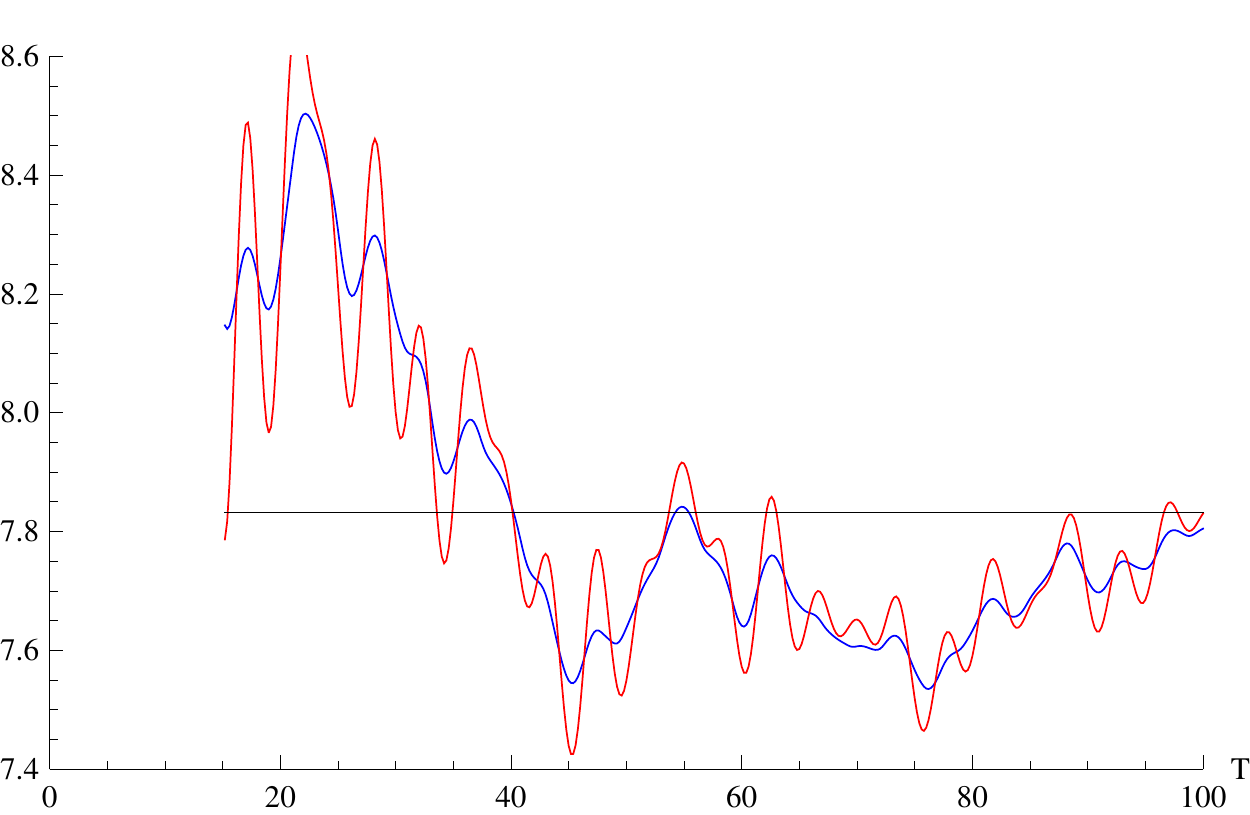}}
    \caption{$x = 10.5:  c = 1.1$ (blue) and $c = 2.0$ (red)}
  \label{fig:Psic1Vc2}
\end{figure*}

Here is another way to visualize how close the Perron integral \eqref{E:psiIntegralApprox} is to the summatory function.  Take $x = 10.5$.  $\psi_0(x) \simeq 7.832$.  Set $c = 1.1$ in equation \eqref{E:psiIntegralApprox}.  Figure ~\ref{fig:PsiIntVsT} shows how close the integral is to 7.832 as a function of $T$.

Again taking $x = 10.5$, figure ~\ref{fig:Psic1Vc2} shows the integrals with $c = 1.1$ and $c = 2.0$.

\begin{figure*}[ht]
  \centerline
  {
    \mbox{\includegraphics[width=\picDblWidth]{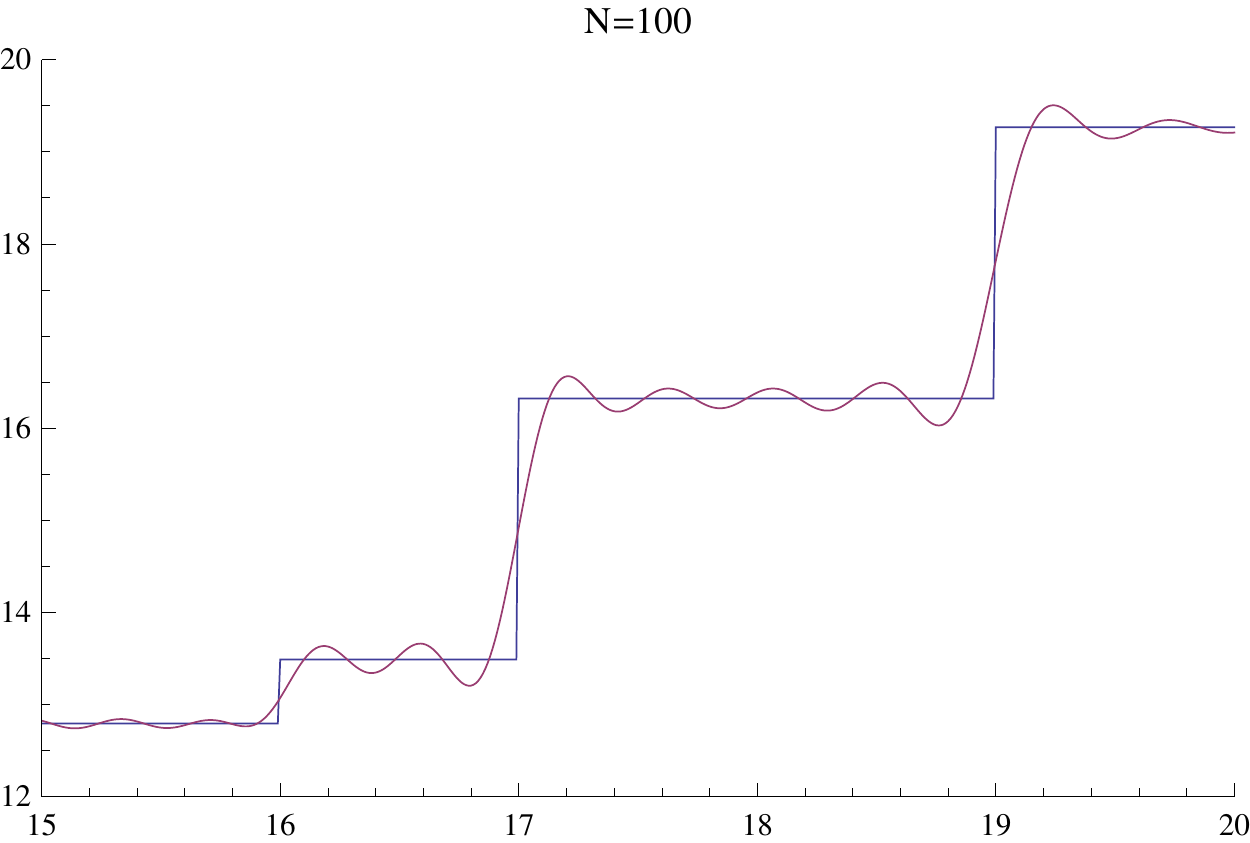}}
    \mbox{\includegraphics[width=\picDblWidth]{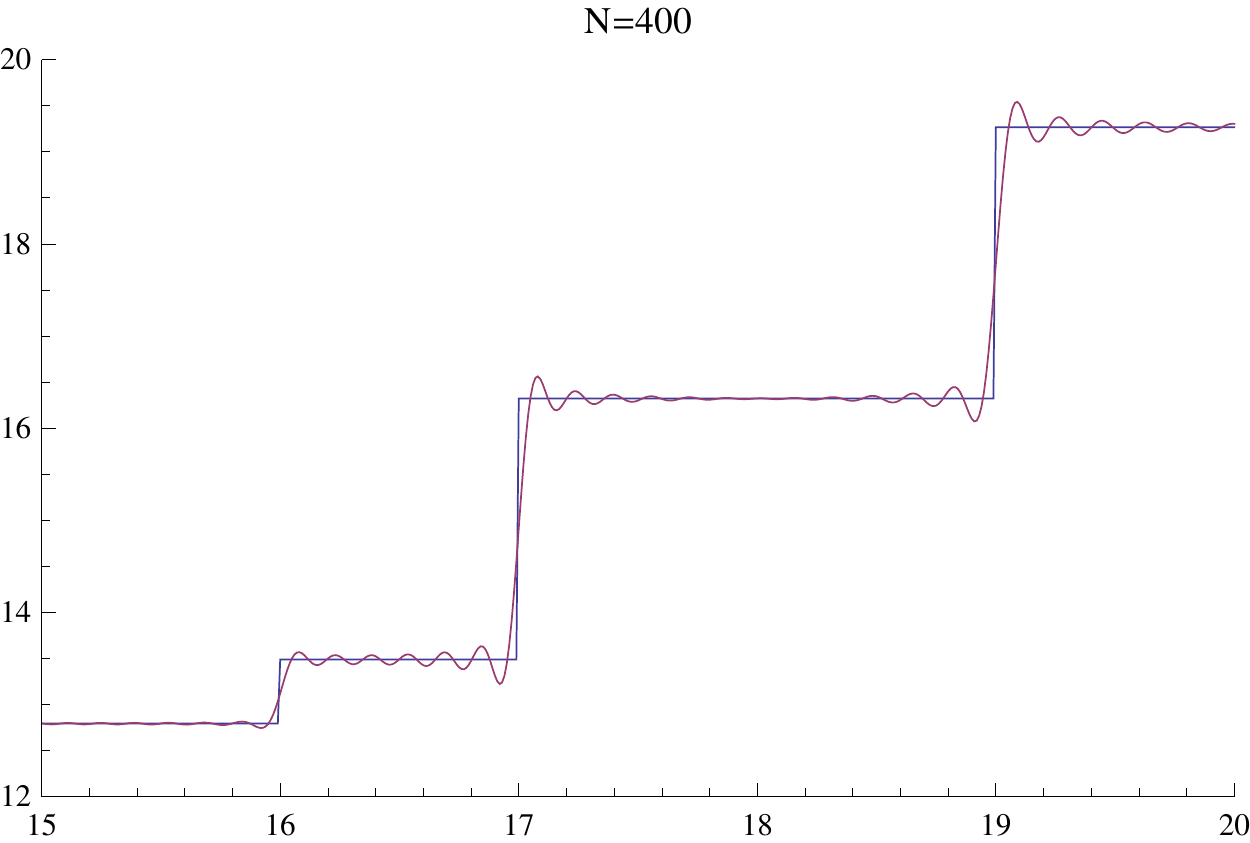}}
  }
  \caption{The Gibbs phenomenon in approximations to $\psi(x)$}
  \label{fig:psiGibbs}
\end{figure*}

Notice that the graphs in Figure ~\ref{fig:psiApproxGraphs} do not show any sign of Gibbs phenomenon.  However, this is due only to the scale of the graphs.  Figure ~\ref{fig:psiGibbs} shows graphs similar to those of Figure ~\ref{fig:psiApproxGraphs}, but with $N = 100$ and $N = 400$, and for $x$ between 15 and 20.  The Gibbs phenomenon at $x = 17$ is quite clear: the step function jumps by $\psi(17.1) - \psi(16.9) = \log(17) \approx 2.833$.  At $x = 16$, the Gibbs phenomenon is less discernible because the jump in the step function there is only $\psi(16.1) - \psi(15.9) = \log(2) \approx .693$.

As with Fourier series, the undershoot and overshoot seem to be related to the size of the jump.  For example, let's see what happens at $x = 16$, $x = 17$, and $x = 97$ and $x = 127$, where the step function $\psi(x)$ has jumps of size $ \log(2)$, $\log(17)$, $ \log(97)$ and $\log(127)$.  We'll use equation \eqref{E:psiZetaZeroSumProved} with $N = 2000$ pairs of zeta zeros.  For each of these values of $x$, let $y_{min}$ be the $y$-value at the minimum just before the jump, and let $y_{max}$ be the $y$-value at the peak just after the jump.  Taking into account the size of the jump, the number that corresponds to the Gibbs constant $g$ in equation \eqref{E:gibbsConstant} is:
\[
d = \frac{y_{max} - y_{min} - c}{2} \frac{1}{c}
\]
where $c$ is the size of the jump.  (The final division by $c$ normalizes for the size of the jump.)  For $x = 16$, $x = 17$, $x = 97$, and $x = 127$, $d$ is about
$0.08965$, $0.08951$, $0.08957$, and $0.08974$.  2000 terms are not enough to see the small jump of size $\log(2)$ at $x = 128$.

\bigskip




\ifthenelse {\boolean{BKMRK}}
  { \section{Computing the Mertens Function \texorpdfstring{$M(x)$}{M(x)}} }
  { \section{Computing the Mertens Function $M(x)$} }

Recall the definition \eqref{E:muDefinition} of the M\"{o}bius $\mu$ function.  The summatory function of $\mu(n)$ is the Mertens function, usually denoted by $M(x)$:
\[
M(x)=\sum _{n=1}^{x} \mu (n).
\]

The Mertens function has an interesting connection to the Riemann Hypothesis.  In 1897, Mertens conjectured that $|M(x)|<\sqrt{x}$ for all $x > 1$.  It is known \cite{KotnikandteRiele} that the Mertens conjecture is true for $x <  10^{14}$.  If this conjecture were true, then the Riemann Hypothesis would also have to be true.  But in 1985, Odlyzko and te Riele proved that the Mertens conjecture is false \cite{OdlyzkoAndteRiele}.  (This does not disprove the Riemann Hypothesis.)   Odlyzko and te Riele did not give a specific $x$ for which the Mertens conjecture is false.  In 2006, Kotnik and te Riele proved \cite{KotnikandteRiele} that there is some $x < 1.59 \times 10^{40}$ for which the Mertens conjecture is false.  The 1985 results used sums with 2000 pairs of complex zeta zeros, each computed to 100 decimals.  The 2006 results used sums with 10000 pairs of zeros, with each zero computed to 250 decimals.

Here is the Dirichlet series involving $\mu(n)$:
\begin{equation}\label{E:muDirichlet}
  \sum _{n=1}^{\infty } \frac{\mu (n)}{n^s}=\frac{1}{\zeta (s)}.
\end{equation}
This holds for $\Re(s)>1$.

From this Dirichlet series, Perron's formula suggests an approximation to the summatory function.  As usual, the prime in the summation and the subscript 0 below means that, where the step function $M(x)$ has a jump discontinuity, $M_0(x)$ is the average of the function values just before and after the jump.
\begin{equation}\label{E:MertensIntegral}
\sideset{}{'} \sum _{n=1}^x \mu(n) =
M_0(x)= \lim_{T\to \infty } \, \frac{1}{2 \pi  i}
   \int_{c-i T}^{c+i T} \frac{1}{\zeta(s)} \frac{x^s}{s} \, ds.
\end{equation}

We will approximate this line integral with an integral over a rectangle that encloses poles of the integrand in \eqref{E:MertensIntegral}.  We will first locate the poles of the integrand, then compute the residues at those poles.  The poles occur at the places where the denominator is 0, that is, $s = 0$, and at the zeros of zeta.  The residue at $s = 0$ is $\frac{1}{\zeta (0)}=-2$.

To get the residue at $\rho_k$, the $k^{th}$ complex zeta zero, we will use \eqref{E:stdResidueFormula1} for the residue of the quotient $A(s)/B(s)$ where
\[
  \frac{A(s)}{B(s)} = \frac{1}{\zeta(s)} \frac{x^s}{s}.
\]

To use \eqref{E:stdResidueFormula1}, we will set $A(s) = \frac{x^s}{s}$ and $B(s) = \zeta(s)$.  Then $B'(s) = \zeta'(s)$, and the residue of $\frac{1}{\zeta(s)} \frac{x^s}{s}$ at $s = \rho_k$ is
\[
\frac{x^{\rho _k}}{\rho _k \zeta '\left(\rho _k\right)}.
\]

At the real zeros $s = -2$, $-4$, $-6 \dots$, we get a similar expression, but with $\rho_k$ replaced by $-2k$.

The alleged approximation that uses the sum of the residues over the first $N$ complex pairs of zeros and the first $M$ real zeros, is
\begin{equation}\label{E:MertensApproximation}
M_0(x) \simeq -2
 + \sum _{k=1}^N x^{\rho _k} \frac{1}{\rho _k \zeta '\left(\rho _k\right)}
 + \sum _{k=1}^M x^{-2 k} \frac{1}{(-2 k) \zeta '(-2 k)}.
\end{equation}

The second sum is small for large $x$.  As we will see, only about 3 terms of this series are large enough to affect the graph, and then, only for small $x$.  In Section \ref{S:TitchmarshTheorem} will discuss to what extent this approximation can be proved.

Note that we have not proved that the expression on the right of \eqref{E:MertensApproximation} is close to $M_0(x)$.

We will form a rectangle that encloses the pole at $s = 0$, along with the poles at the first $N$ complex zeros, and the poles at the first $M$ real zeros.  We will then integrate around the rectangle.  We will see that the line integrals on the top, left, and bottom sides are small, so that the line integral in \eqref{E:MertensIntegral} can be approximated by the sum of integrals on all four sides, which, in turn, equals the sum of the residues of the poles that lie inside the rectangle.

\subsection{Numerical Experiments}

\begin{figure*}[ht]
  \centerline{
    \mbox{\includegraphics[width=\picDblWidth]{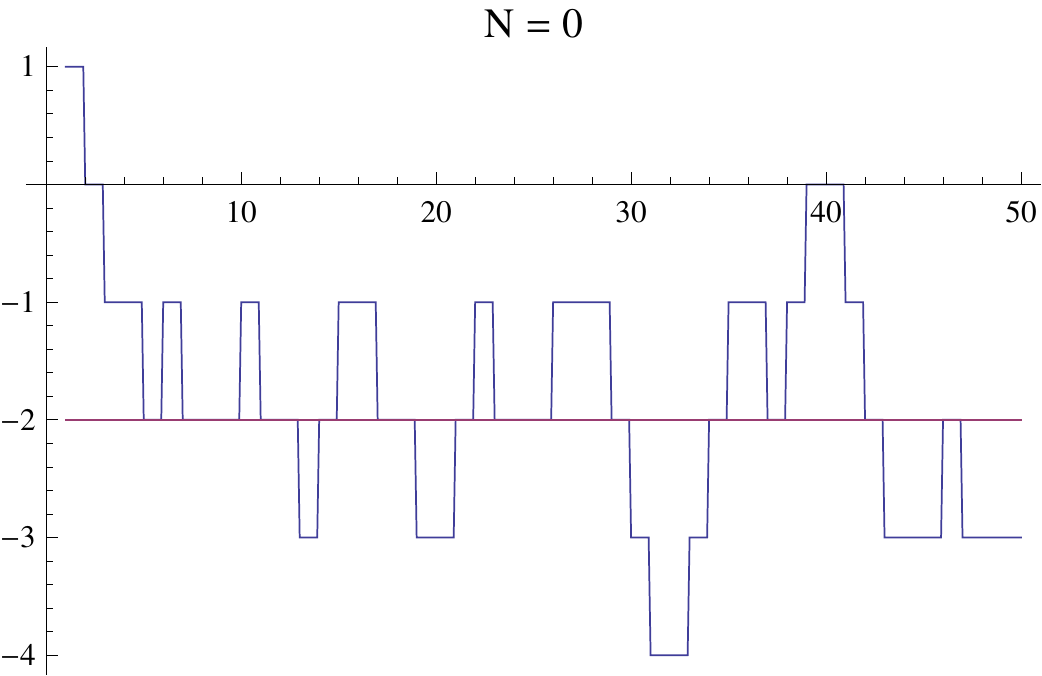}}
    \mbox{\includegraphics[width=\picDblWidth]{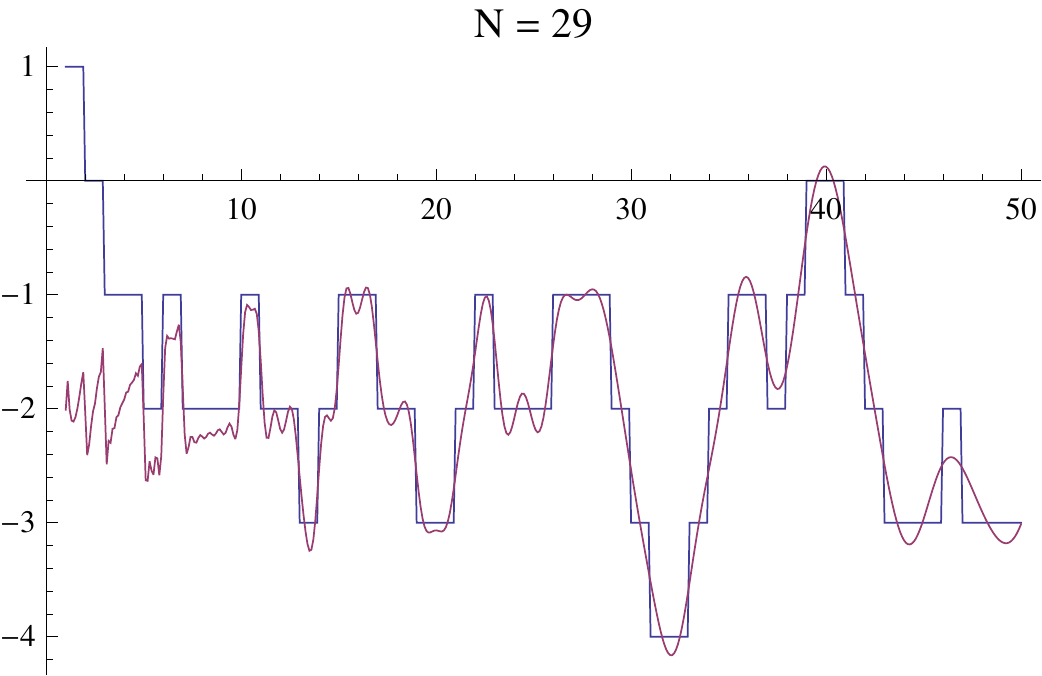}}
  }
  \centerline{
    \mbox{\includegraphics[width=\picDblWidth]{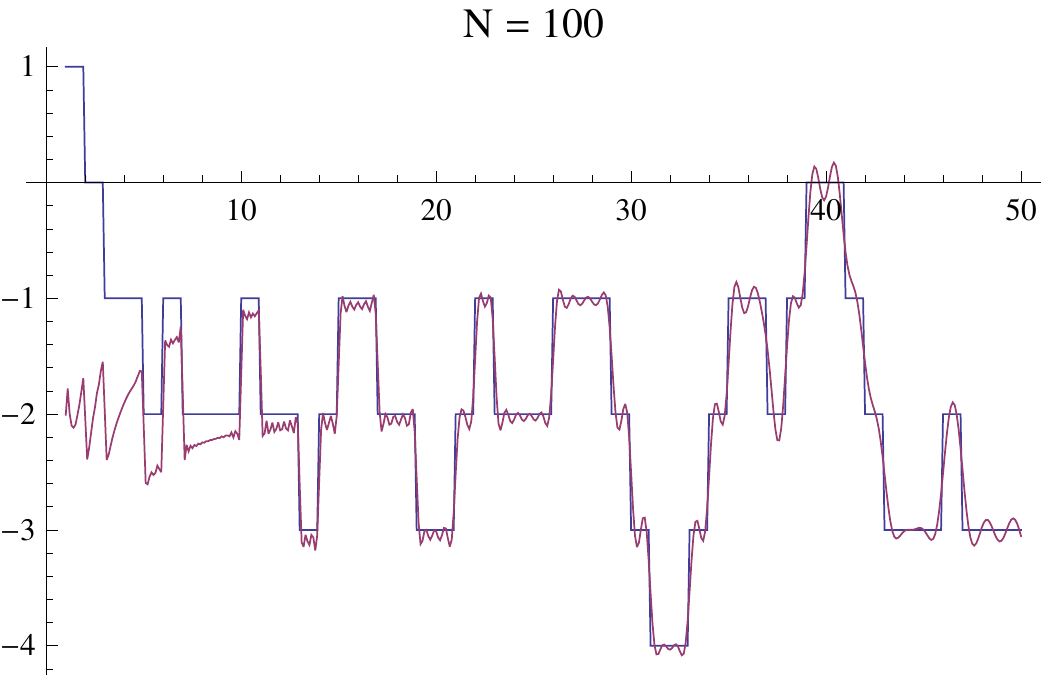}}
    \mbox{\includegraphics[width=\picDblWidth]{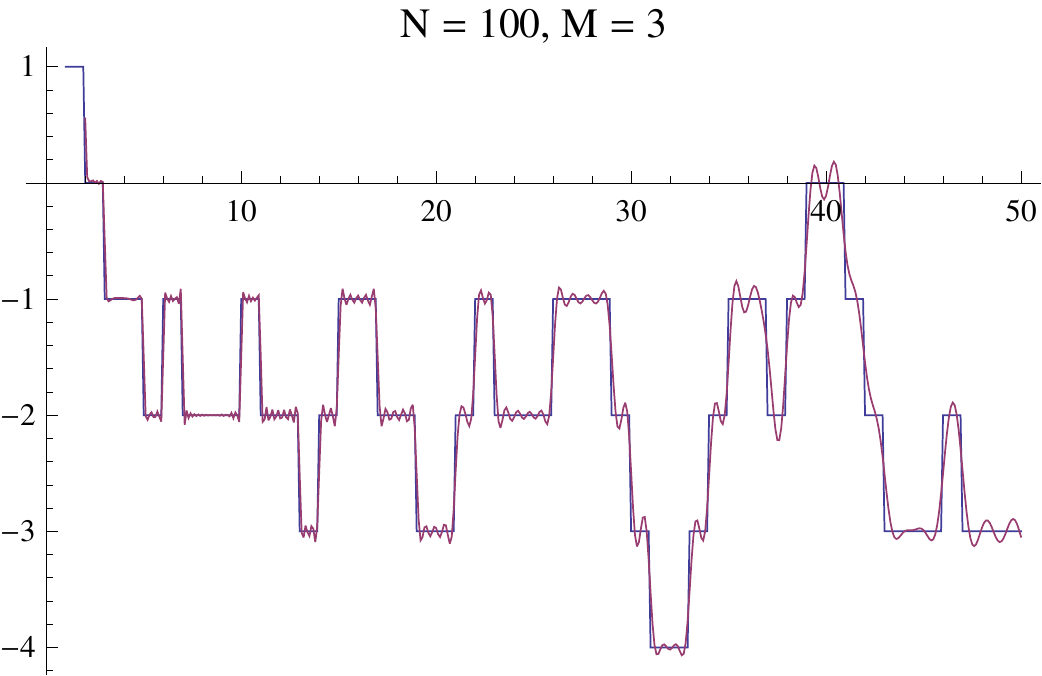}}
  }
  \caption{approximating $M_0(x)$ with $N = 0$, $N = 29$, and $N = 100$ pairs of zeta zeros}
  \label{fig:mxApproxGraphs}
\end{figure*}

In the first three graphs in Figure ~\ref{fig:mxApproxGraphs}, we set $M = 0$, and used $N = 0$, $N = 29$, and $N = 100$ pairs of complex zeta zeros in \eqref{E:MertensApproximation}.  Notice that, even with $N = 100$, the approximation is not very good for small x.  So, in the fourth graph, we also include the sum of the residues of the first $M = 3$ pairs of real zeta zeros.  This makes the agreement quite good, even for small $x$.

Now let's look at numerical examples.  Let's set $c = 2$, and use \eqref{E:MertensApproximation} with $N = 29$ and $M = 3$.  A rectangle that encloses the pole at $s = 0$, the poles at the three real zeros $s = -2$, $s = -4$, and $s = -6$, and the 29 pairs of complex zeros could extend, for example, from $2 - 100 i$ to $-7 + 100 i$.  If $x = 32$, then one can see from the graph that $M(x)=-4$, and, since this $x$ is divisible by a prime power, $\mu(x) = 0$, so that $M(x)$ is continuous at this $x$.  The integrals along the right, top, left, and bottom sides are, respectively, $I_1 \approx -4.506$, $I_2 \approx .180 - 0.325 i$, $I_3 \approx 0$, and $I_4 \approx .180 + 0.325 i$.  The total is $\approx -4.145$.

Notice that, for a fixed $N$ and $M$, \eqref{E:MertensApproximation} diverges more and more from the step function as $x$ increases.  We can remedy this by using a larger value of $N$ in \eqref{E:MertensApproximation}.


\begin{figure*}[ht]
  \centerline{
    \mbox{\includegraphics[width=\picDblWidth]{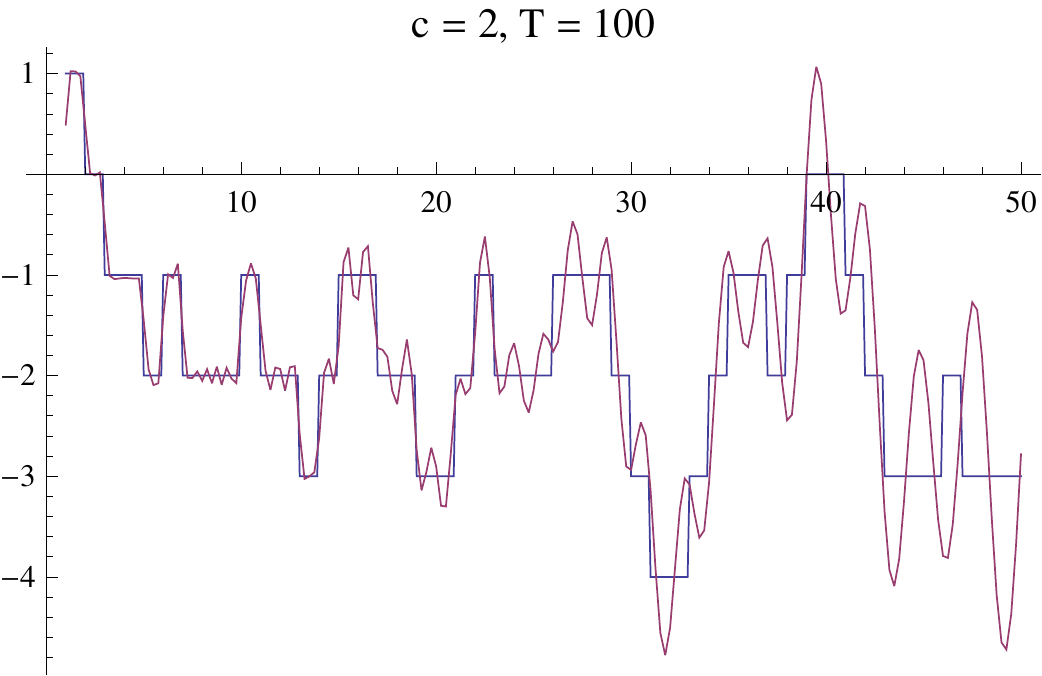}}
    \mbox{\includegraphics[width=\picDblWidth]{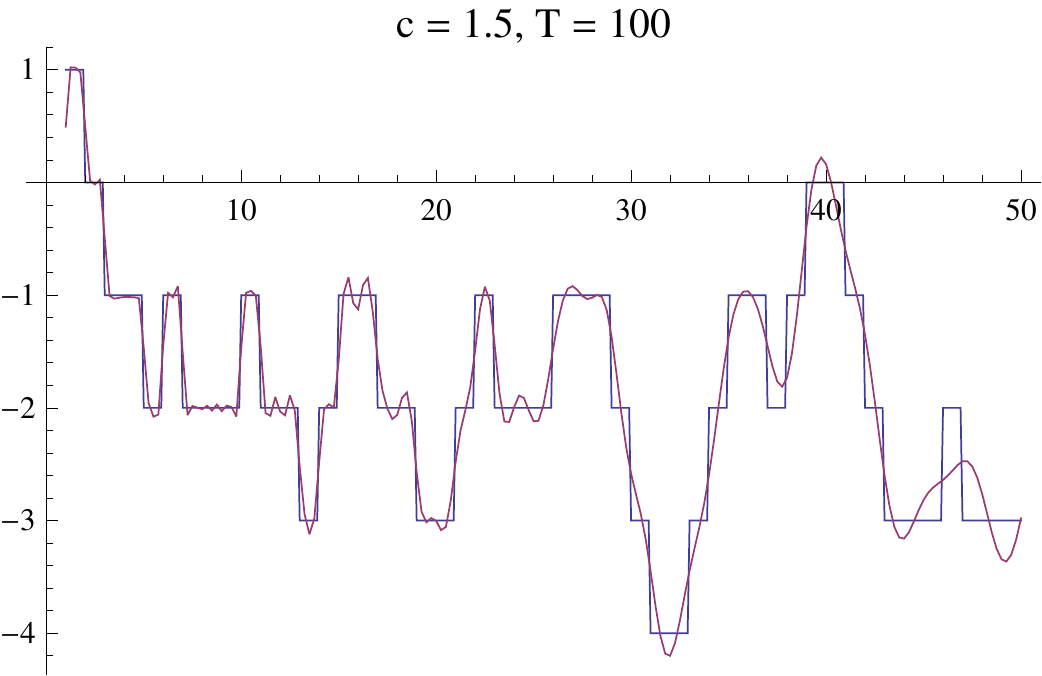}}
  }
  \caption{approximating $M_0(x)$ with integrals: $T = 100$, $c = 2$ and $c = 1.5$}
  \label{fig:mxIntegralApproximationA}
\end{figure*}

In Figure ~\ref{fig:mxApproxGraphs}, we drew approximations based on the sum of the residues given by equation \eqref{E:MertensApproximation}.  What would happen if we chose $c$ and $T$, then, for each $x$, we plotted the \textit{integral}
\begin{equation}\label{E:MxIntegral}
  \frac{1}{2 \pi  i}
   \int_{c - i T}^{c + i T} \frac{1}{\zeta(s)} \frac{x^s}{s} \, ds
\end{equation}
instead of the \textit{sum} of the residues?

\begin{figure*}[ht]
  \centerline{
    \mbox{\includegraphics[width=\picDblWidth]{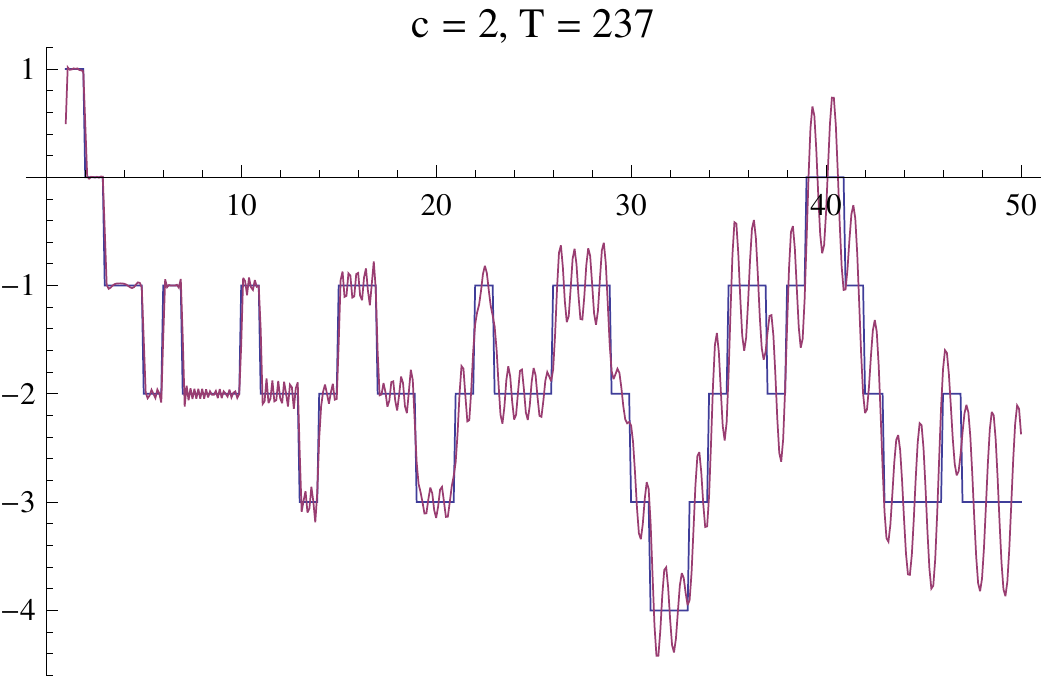}}
    \mbox{\includegraphics[width=\picDblWidth]{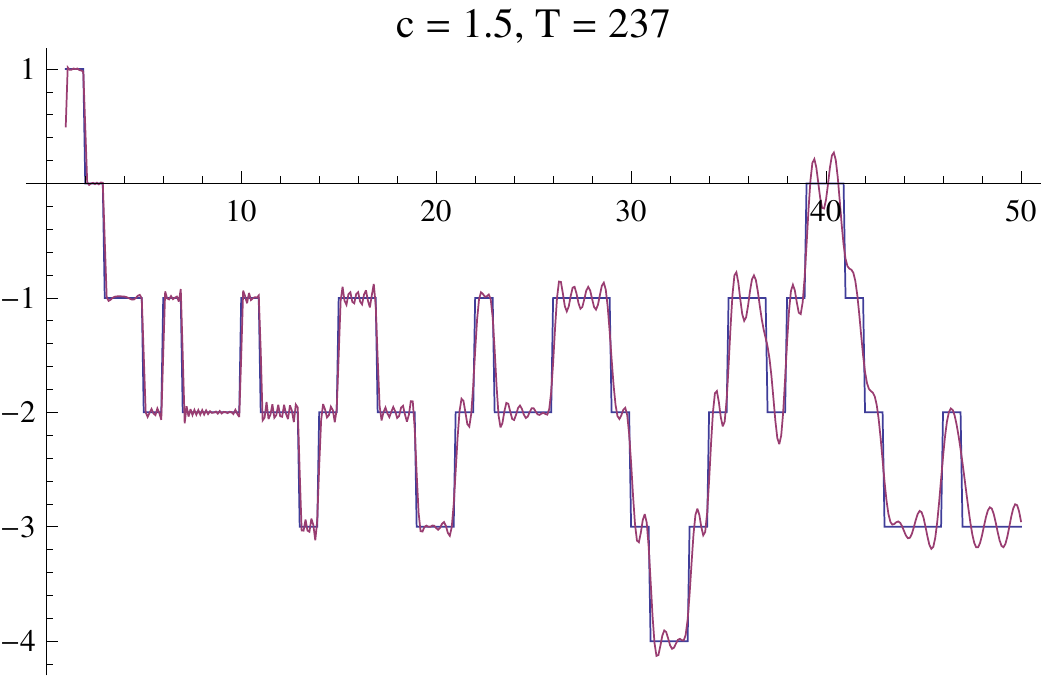}}
  }
  \caption{approximating $M_0(x)$ with integrals: $T = 237$, $c = 2$ and $c = 1.5$}
  \label{fig:mxIntegralApproximationB}
\end{figure*}

On the left side of Figure ~\ref{fig:mxIntegralApproximationA}, we show the values we get by substituting $c = 2$ and $T = 100$ in \eqref{E:MxIntegral}.  The integral \eqref{E:MxIntegral} is clearly following the step function, but the agreement with the step function is not very good.  On the right side of Figure ~\ref{fig:mxIntegralApproximationA}, we plot the integral using $c = 1.5$.  With $c = 1.5$, the agreement is much better.  In fact, for the $x$ in range of this graph, this integral approximation is about as good as the sum of residues when $N = 29$ pairs of zeros are included.

Since the integral approximation is so bad with $c = 2$, $T = 100$, one might ask whether we can get a better approximation by taking a larger $T$.  So, let's try a $T$ that corresponds to, say, 100 pairs of zeros; to do this, we can set $T = 237$.  The graph is in figure ~\ref{fig:mxIntegralApproximationB}.  The overall level of agreement isn't that different, but using a larger $T$ introduces higher-frequency terms into the graphs.  What does improve the approximation is to use smaller a value of $c$ (1.5 instead of 2).

\begin{figure*}[ht]
  \mbox{\includegraphics[width=\picDblWidth]{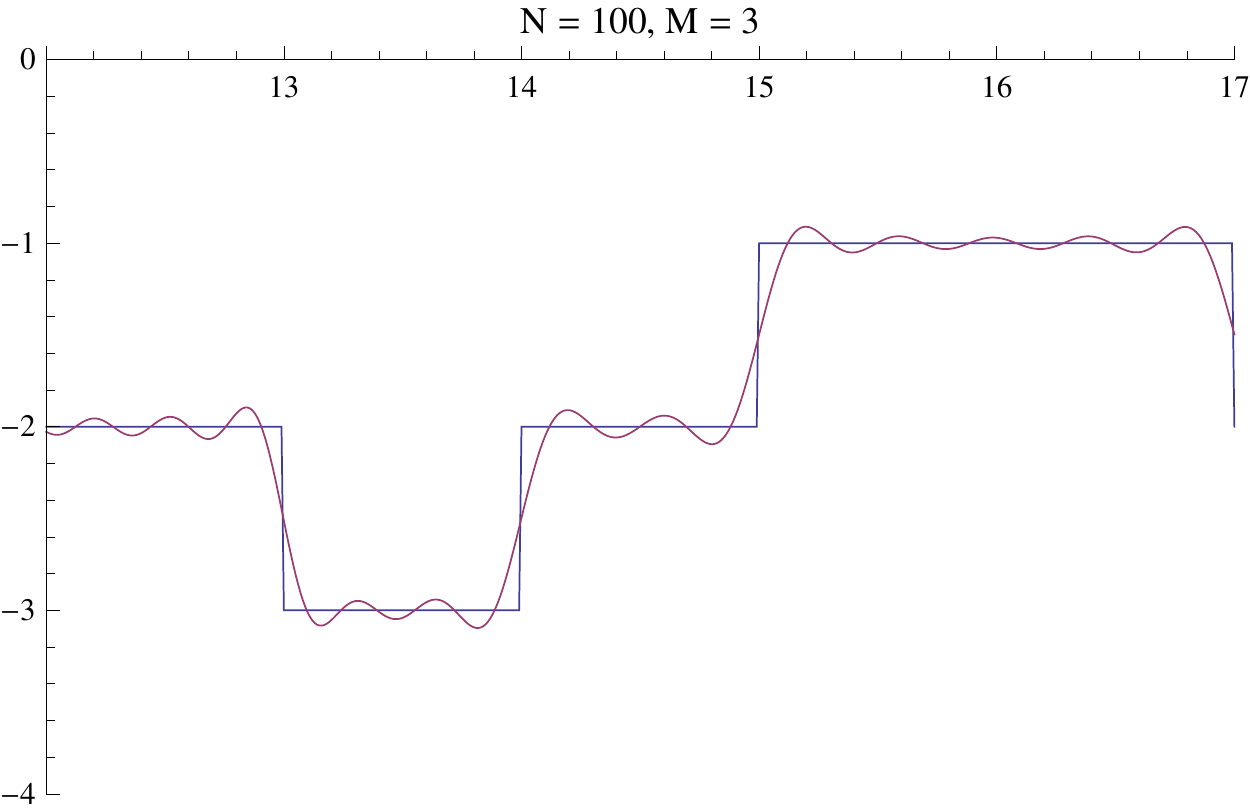}}
  \caption{Gibbs phenomenon in $M(x)$: $N = 100$, $M = 3$}
  \label{fig:mxGibbs}
\end{figure*}

Figure ~\ref{fig:mxGibbs} shows the Gibbs phenomenon in the graph of $M(x)$.

\subsection{Theorems About $M_0(x)$}\label{S:TitchmarshTheorem}

It certainly appears from the above graphs that the sums in \eqref{E:MertensApproximation} provide a good approximation to the step function $M_0(x)$.  Here we summarize what has actually been proven.



Assuming the Riemann Hypothesis, and assuming that all zeta zeros are simple (i.e., are of multiplicity 1), Titchmarsh \cite[p. 318]{TitchmarshFirstEd} proved the following Theorem:

\begin{theorem}
There is a sequence $T_{\nu}$, $\nu \leq T_{\nu} \leq \nu + 1$, such that
\begin{equation}\label{E:provedSumM0}
M_0(x) =
 -2 +
 \lim_{\nu \to \infty } \,
   \sum_{|\gamma| < T_{\nu}} \frac{x^{\rho}}{\rho \zeta'(\rho)}
 + \sum_{k=1}^{\infty} \frac{ (-1)^{k-1} (2 \pi / x)^{2 k} }{ (2 k)! k \zeta(2k + 1) }.
\end{equation}
\end{theorem}
Note: because
\[
\zeta'(-2 k) = (-1)^k \frac{(2 k)!}{2 (2 \pi)^{2 k}} \zeta(2 k + 1),
\]
the $k^{th}$ term in the second sum is the same as the $k^{th}$ term in equation \eqref{E:MertensApproximation}, namely,
\[
x^{-2 k} \frac{1}{(-2 k) \zeta '(-2 k)}.
\]

In Titchmarsh's Theorem, $\gamma$ is the imaginary part of the zeta zero $\rho$.  The proof in \cite{TitchmarshFirstEd} sums the residues inside a rectangle, the sides of which must not pass through any zeta zeros.  This means that the $T_{\nu}$ in this Theorem cannot be the imaginary part of a zeta zero.  In addition, the proof makes use of a theorem which says that, given $\epsilon > 0$, each interval $[\nu, \nu + 1]$ contains a $T = T_{\nu}$ such that
\[
\frac{1}{\zeta(\sigma + i T)} = O(T^{\epsilon})
\]
for $\frac{1}{2} \leq \sigma \leq 2$.

Bartz \cite{BartzMx2} proves a result similar to Titchmarsh's Theorem, but without assuming the Riemann Hypothesis.  Her result is also more general, in that it works even if the zeta zeros are not simple.  In the easiest case, where the zeros \emph{are} simple, she proves that there is an absolute positive constant $T_0$ and a sequence $T_{\nu}$, $2^{\nu-1} T_0 \leq T_{\nu} \leq 2^{\nu} T_0$, (for $\nu \geq 1$), such that \eqref{E:provedSumM0} holds.

One would \emph{like} to have this result
\[
M_0(x) = -2
 + \sum _{k=1}^{\infty} x^{\rho _k} \frac{1}{\rho _k \zeta '\left(\rho _k\right)}
 + \sum _{k=1}^{\infty} x^{-2 k} \frac{1}{(-2 k) \zeta '(-2 k)}.
\]
but note that neither Titchmarsh's nor Bartz's result is this strong.


\section{Counting the Squarefree Integers}\label{S:CountingTheSquarefreeIntegers}
A squarefree number is a positive integer that is not divisible by the square of a prime. The first 10 squarefree numbers are 1, 2, 3, 5, 6, 7, 10, 11, 13, and 14.

The number of squarefree integers from 1 through $x$ is usually denoted by $Q(x)$. For example, $Q(10) = 7$, because 4, 8, and 9 are divisible by squares of primes, so they don't contributes to the count. The graph of $Q(x)$ is an irregular step function.

Define $q(n) = 1$ if $n$ is squarefree and 0 otherwise ($q(1) = 1$).  Then $Q(x)$ is the summatory function of $q(n)$.

The following Dirichlet series holds for $\Re(s) > 1$ \cite[p. 255]{HardyAndWright}:
\[
  \sum _{n=1}^{\infty } \frac{q(n)}{n^s}=\frac{\zeta (s)}{\zeta (2 s)}.
\]

Perron's formula tells us that, if $c > 1$,
\begin{equation}\label{E:QxIntegral}
\sideset{}{'} \sum _{n=1}^x q(n) = Q_0(x)
 = \lim_{T\to \infty } \, \frac{1}{2 \pi i}
   \int_{c-i T}^{c+i T} \frac{\zeta(s)}{\zeta(2s)} \frac{x^s}{s} \, ds
\end{equation}
where $Q_0(x)$ equals $Q(x)$ except where $Q(x)$ has a jump discontinuity (that is, at the squarefree integers), and, at those jumps, $Q_0(x)$ is the average of the values before and after the jump.

The integrand has a pole at $s = 0$, where the residue is 1.  Recall that zeta has a pole of order 1 at $s = 1$.  Therefore, the integrand has a pole at $s = 1$.  The residue there is $6 x / \pi^2$.

Finally, the integrand has poles wherever $2s$ is a zero of zeta.  If $\rho$ is any zeta zero, then the pole occurs at $s = \rho / 2$.  To get the residue at $s = \rho / 2$, we will use equation \eqref{E:stdResidueFormula1} for the residue of the quotient $A(s)/B(s)$
\[
  \frac{A(s)}{B(s)} = \frac{\zeta(s)}{\zeta(2s)} \frac{x^s}{s},
\]
where $A(s) = \zeta(s) \frac{x^s}{s}$ and $B(s) = \zeta(2s)$.  Then $B'(s) = 2 \zeta(2s)$, and the residue at the zeta zero $s = \rho / 2$ is
\[
\frac{A(s)}{B'(s)} = 
  \frac{\zeta(\frac{\rho}{2})}{2 \zeta'(\rho)} \frac{x^\frac{\rho}{2}}{\frac{\rho}{2}} =
  \frac{x^\frac{\rho}{2}}{\rho} \frac{\zeta(\frac{\rho}{2})}{\zeta'(\rho)}.
\]

This holds whether $\rho$ is a real or a complex zeta zero.  Each zeta zero gives rise to one term of this form.  Therefore, the sum over the residues at the first $N$ complex zeros and the first $M$ real zeros, plus the residues at $s = 0$ and $s = 1$, gives this alleged approximation

\begin{equation}\label{E:QxSum}
Q_0(x)\simeq
 1 + \frac{6 x}{\pi ^2}
 +2 \Re\left(\sum_{k=1}^N x^{\frac{\rho _k}{2}} \frac{
   \zeta \left(\frac{\rho_k}{2}\right)}{\rho _k \zeta '\left(\rho_k\right)}\right)
 + \sum _{k=1}^M x^{-\frac{2 k}{2}} \frac{ \zeta\left(\frac{-2 k}{2}\right)}{(-2 k) \zeta '(-2 k)}.
\end{equation}

In the second sum, only the terms with $M \leq 3$ are large enough to affect the graph.  Moreover, if $k$ is even, then the term in the second sum
\[
\frac{x^{-\frac{2 k}{2}} \zeta\left(\frac{-2 k}{2}\right)}{(-2 k) \zeta '(-2 k)} =
 \frac{x^{-k} \zeta(-k)}{(-2 k) \zeta '(-2 k)}
\]
is zero, because $\zeta(-k) = 0$ since $-k$ is a zeta zero if $k$ is even.  Therefore, the second sum reduces to a sum over odd $k \leq M$.


If the zeta zeros are simple, Bartz \cite{BartzSquareFree} proves that there is an absolute positive constant $T_0$ and a sequence $T_{\nu}$, $2^{\nu-1} T_0 \leq T_{\nu} \leq 2^{\nu} T_0$, (for $\nu \geq 1$), such that
\begin{equation}\label{E:CountingSquareFree}
Q_0(x)=
 1 + \frac{6 x}{\pi ^2} +
 \lim_{\nu \to \infty } \,
   \sum_{|\gamma| < T_{\nu}} x^{\frac{\rho _k}{2}} \frac{
   \zeta \left(\frac{\rho_k}{2}\right)}{\rho _k \zeta '\left(\rho_k\right)}
 + \sum _{k=1}^{\infty} x^{-\frac{2 k}{2}} \frac{ \zeta\left(\frac{-2 k}{2}\right)}{(-2 k) \zeta '(-2 k)},
\end{equation}
where $\gamma$ is the imaginary part of the zeta zero $\rho$.  Bartz has a similar result \cite{BartzCubeFree} for cube-free integers.

\subsection{Numerical Experiments}

\begin{figure*}[ht]
  \centerline{
    \mbox{\includegraphics[width=\picDblWidth]{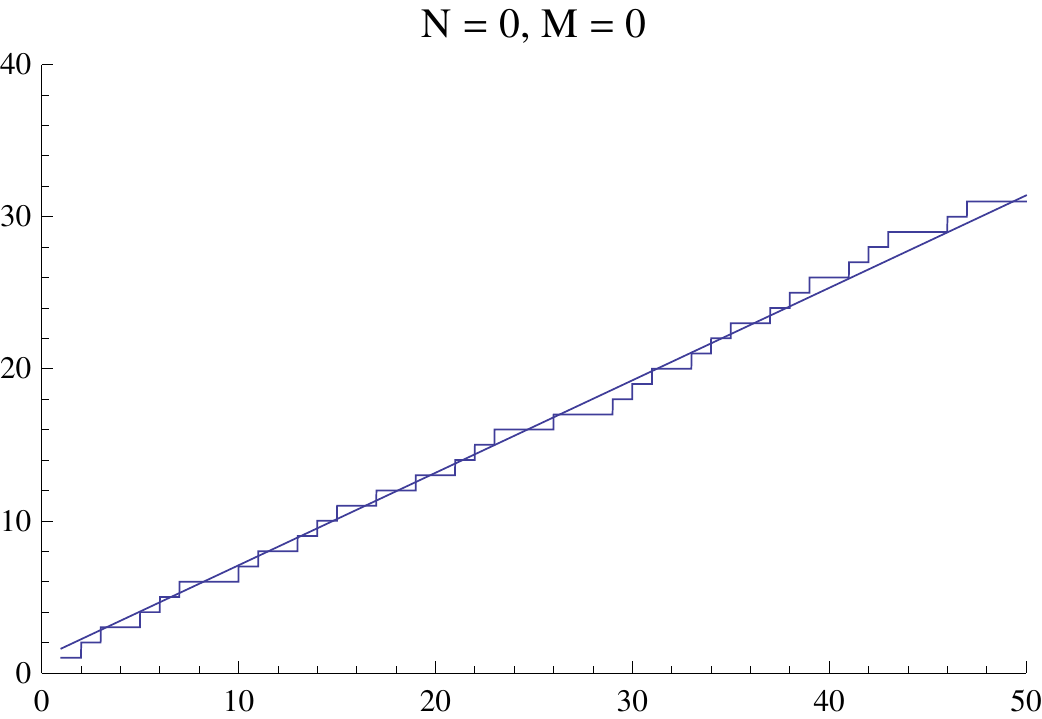}}
    \mbox{\includegraphics[width=\picDblWidth]{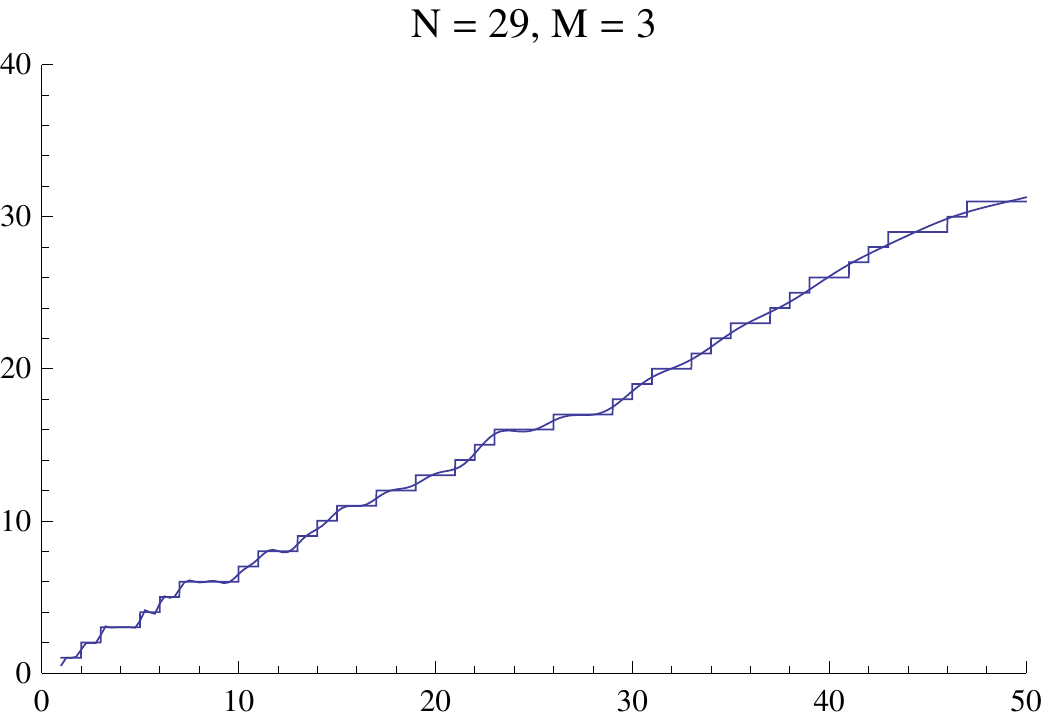}}
  }
  \caption{Approximating $Q_0(x)$ with sums of residues}
  \label{fig:qxSumApprox1}
\end{figure*}

\begin{figure*}[ht]
  \centerline{
    \mbox{\includegraphics[width=\picDblWidth]{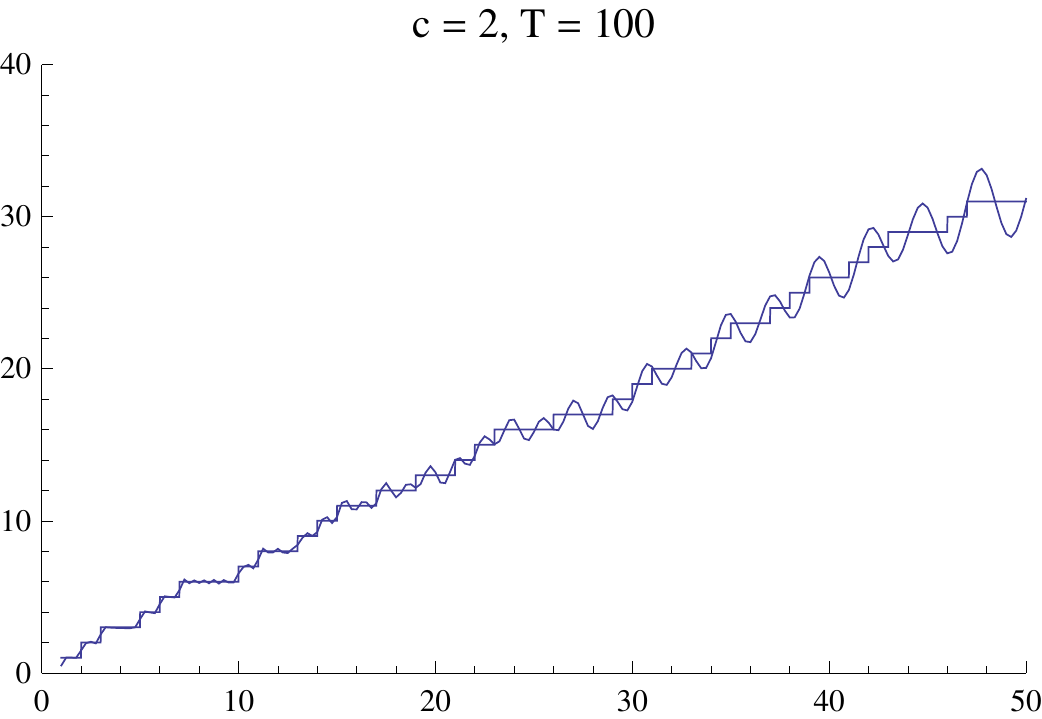}}
    \mbox{\includegraphics[width=\picDblWidth]{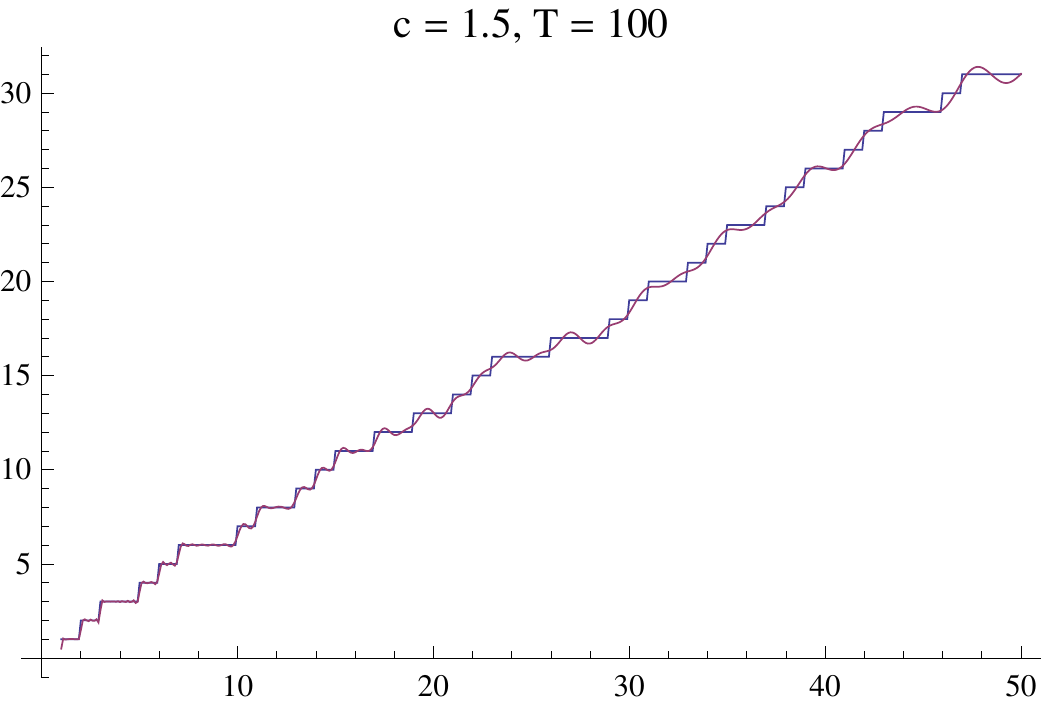}}
  }
  \caption{Approximating $Q_0(x)$ with integrals; $c = 2$ and $c = 1.5$}
  \label{fig:qxIntegralApprox29}
\end{figure*}

Figure ~\ref{fig:qxSumApprox1} shows graphs of two approximations to $Q_0(x)$ using \eqref{E:QxSum}.  The graph on the left shows the linear approximation that comes from the first two terms of \eqref{E:QxSum}.  The graph on the right shows \eqref{E:QxSum} with $N = 29$, $M = 3$.  The three terms in the second sum have a noticeable affect on the graph only for $x < 10$.

The left side of Figure ~\ref{fig:qxIntegralApprox29} shows a graph of the integral approximation
\[
Q_0(x) \simeq  \frac{1}{2 \pi i}
   \int_{c-i T}^{c+i T} \frac{\zeta(s)}{\zeta(2s)} \frac{x^s}{s} \, ds.
\]
with $c = 2$ and $T = 100$.  As noted above, a rectangle with $T = 100$ will enclose 29 pairs of complex zeros.  The right side shows the integral approximation, this time, with $c = 1.5$.

\ifthenelse {\boolean{BKMRK}}
  { \section{Tallying the Euler \texorpdfstring{$\phi$}{phi} Function} }
  { \section{Tallying the Euler $\phi$ Function} }

The Euler phi (or "totient") function is defined as follows.  If $n$ is a positive integer, then $\phi(n)$ is the number of integers from 1 through $n$ that are relatively prime to $n$.  In other words, $\phi(n)$ is the number of integers $k$ from 1 through $n$ such that the greatest common divisor of $k$ and $n$ is 1.  For example, $\phi(10) = 4$, and if $n$ is prime, then $\phi(n) = n - 1$.  We will let $\Phi(x)$ denote the summatory function of $\phi(n)$, so that
\[
  \Phi(x) = \sum _{n=1}^x \phi(n).
\]

Then
\[
  \Phi_0(x) = \sideset{}{'} \sum _{n=1}^x \phi(n)
\]
so that $\Phi_0(x)$ equals $\Phi(x)$ except where $\Phi(x)$ has a jump discontinuity (this happens at every positive integer since $\phi(n) > 0$).  Further, at those jumps, $\Phi_0(x)$ is the average of the values of $\Phi(x)$ before and after the jump.

The following Dirichlet series holds for $\Re(s) > 2$ \cite[Theorem 287, p. 250]{HardyAndWright}:
\begin{equation}\label{E:PhiDirichletSeries}
  \sum _{n=1}^{\infty } \frac{\phi(n)}{n^s}=\frac{\zeta (s-1)}{\zeta (s)}.
\end{equation}

From this Dirichlet series, Perron's formula tells us that, if $c > 2$, then
\begin{equation}\label{E:PhiIntegral}
 \Phi_0(x)
 = \lim_{T\to \infty } \, \frac{1}{2 \pi i}
   \int_{c-i T}^{c+i T} \frac{\zeta(s-1)}{\zeta(s)} \frac{x^s}{s} \, ds.
\end{equation}

The integrand has a pole at $s = 0$ where the residue is $1/6$.  $\zeta(s)$ has a pole of order 1 at $s = 1$, so, because of the $\zeta(s-1)$ in the numerator, the integrand has a pole at $s - 1 = 1$, that is, at $s = 2$.  The residue at $s = 2$ is $3 x^2 / \pi^2$.

The integrand also has poles at every real and complex zero of zeta.  We will apply residue formula \eqref{E:stdResidueFormula1} with $B(s) = \zeta(s)$ and $A(s)$ equal to everything else in the integrand, that is, $A(s) = \zeta(s-1) x^s/s$.  The residue at the zeta zero $s = \rho$ will be
\[
\frac{A(\rho)}{B'(\rho)} = \frac{\zeta (\rho - 1) x^{\rho }}{\rho } \frac{1}{\zeta '(\rho )}
 = \frac{x^{\rho } \zeta(\rho - 1) }{\rho \zeta '(\rho ) }.
\]

Therefore, if we create a rectangle that encloses the first $N$ pairs of complex zeta zeros, the first $M$ real zeros, and the poles at $s = 0$ and $s = 2$, then the integral around the rectangle is just the sum of the residues, which suggests that

\begin{equation}\label{E:PhiSum}
\Phi_0(x) \simeq
\frac{1}{6} + \frac{3 x^2}{\pi ^2}
  + 2 \Re\left(\sum _{k=1}^N x^{\rho _k} \frac{ \zeta
   \left(\rho _k-1\right)}{\rho _k \zeta '\left(\rho_k\right)}\right)
  + \sum _{k=1}^M x^{-2 k} \frac{ \zeta (-2 k-1)}{(-2 k) \zeta '(-2 k)}.
\end{equation}
A quick calculation will show that, for $x > 1.5$, the second sum is too small to visibly affect our graphs, so we will take $M = 0$ in the graphs below.

\begin{figure*}[ht]
  \mbox{\includegraphics[width=\picDblWidth]{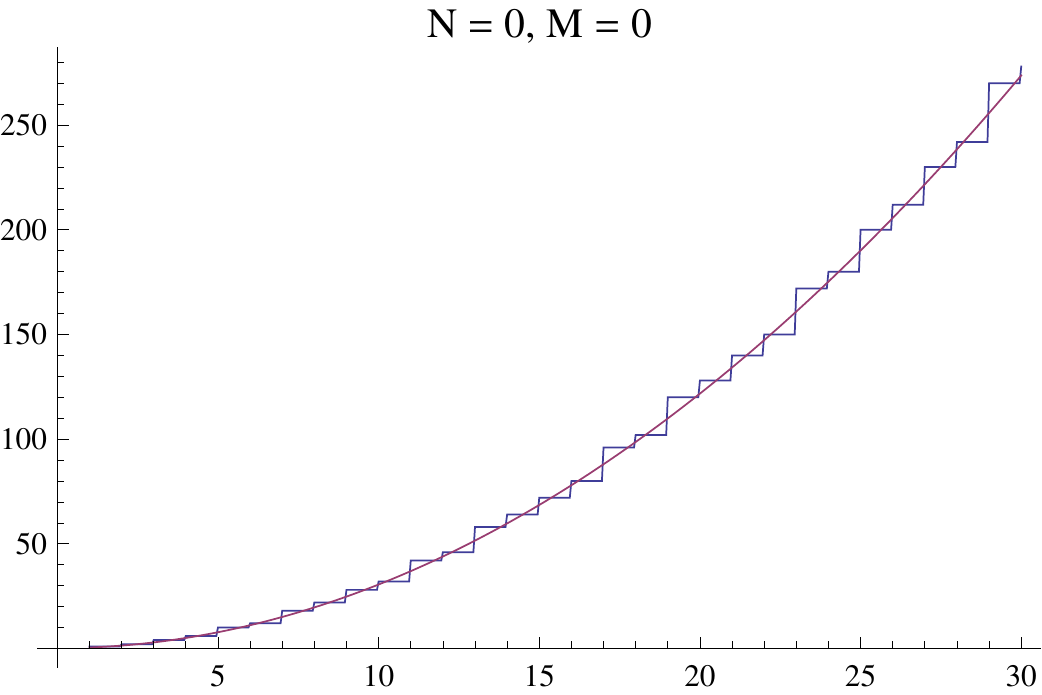}}
  \caption{Approximating $\Phi(x)$ with the quadratic $\frac{1}{6} + \frac{3 x^2}{\pi ^2}$}
  \label{fig:phiSumApprox0}
\end{figure*}

Notice that, if no zeta zeros are included in \eqref{E:PhiSum} (that is, if $N = M = 0$), then the estimate for $\Phi(x)$ is given by the quadratic
\[
  \Phi(x) \simeq \frac{1}{6} + \frac{3 x^2}{\pi ^2}.
\]

Figure ~\ref{fig:phiSumApprox0} shows this quadratic approximation to $\Phi(x)$.  This approximation looks quite good.  Not only that, in Figure ~\ref{fig:phiSumApprox0}, the value of this quadratic appears to be between $\Phi(x-1)$ and $\Phi(x)$ for \textit{every} integer $x$.  In fact, the first integer $x$ where
\[
\Phi(x-1) \leq \frac{1}{6} + \frac{3 x^2}{\pi ^2} \leq \Phi(x)
\]
is \textit{false}, is $x = 820$: $\Phi(819) = 204056$, $\Phi(820) = 204376$, but $\frac{1}{6} + \frac{3 \cdot 820^2}{\pi ^2} \simeq 204385.25831$.  There is only one such $x$ less than $10^3$, but there are 36 of them less than $10^4$, 354 of them less than $10^5$, 3733 of them less than $10^6$, and 36610 such $x$ less than $10^7$.

Compare this to the standard estimate with an error term in \cite[Theorem 330, p. 268]{HardyAndWright} that
\begin{equation}\label{E:PhiSumWithOEstimate}
  \Phi(x) = \sum _{n=1}^x \phi(n) = \frac{3 x^2}{\pi^2} + O( x \log x ).
\end{equation}
If \eqref{E:PhiSum} is a valid approximation to $\Phi_0(x)$ (or $\Phi(x)$), \eqref{E:PhiSumWithOEstimate} implies that
\[
 \Re\left(\sum _{k=1}^{\infty} x^{\rho _k} \frac{ \zeta
   \left(\rho _k-1\right)}{\rho _k \zeta '\left(\rho_k\right)}\right) = O( x \log x ).
\]

\begin{figure*}[ht]
  \centerline{
    \mbox{\includegraphics[width=\picDblWidth]{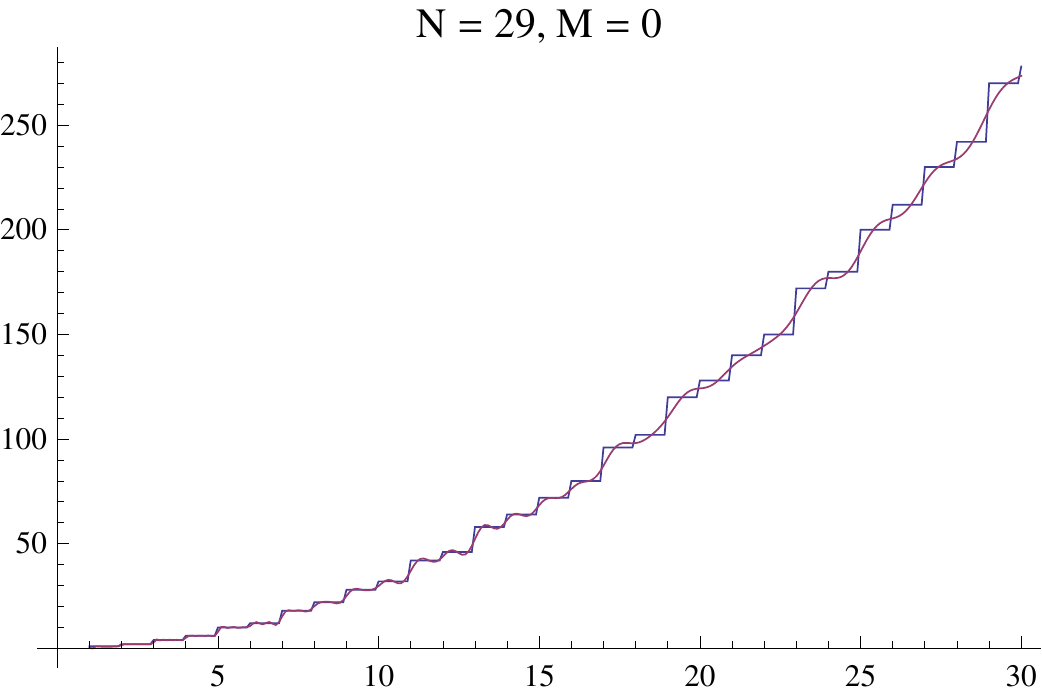}}
    \mbox{\includegraphics[width=\picDblWidth]{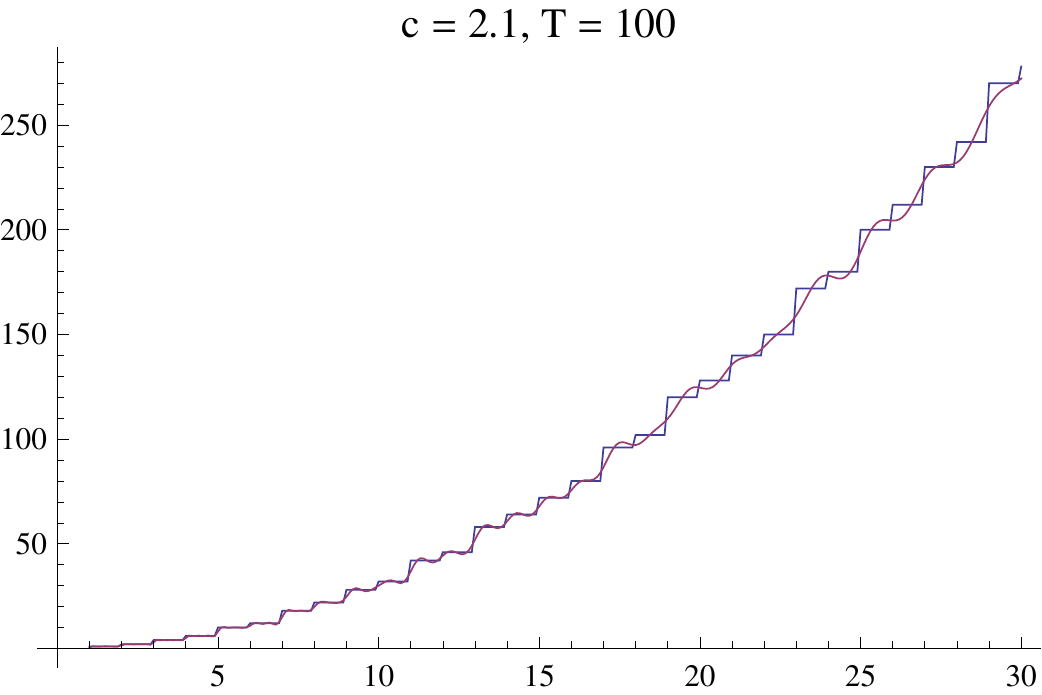}}
  }
  \caption{Approximating $\Phi_0(x)$ with a sum and an integral}
  \label{fig:phiSumAndIntegralApprox29}
\end{figure*}

Here is the integral approximation that we get from equation \eqref{E:PhiIntegral}:
\begin{equation}\label{E:PhiIntegral2}
  \Phi_0(x) \simeq
   \frac{1}{2 \pi i}
   \int_{c-i T}^{c+i T} \frac{\zeta(s-1)}{\zeta(s)} \frac{x^s}{s} \, ds.
\end{equation}

The graphs in Figure ~\ref{fig:phiSumAndIntegralApprox29} show that the sum \eqref{E:PhiSum} and the integral \eqref{E:PhiIntegral2} seem to give reasonable approximations to the step function $\Phi(x)$ (or $\Phi_0(x)$), at least for these $x$.  Quite a bit of analysis would be required to estimate how close \eqref{E:PhiSum} and \eqref{E:PhiIntegral2} are to $\Phi_0(x)$.

R\c{e}ko\'{s} \cite{Rekos} discusses this explicit formula related to Euler's $\phi$ function:
\[
  f(z) = \lim_{n \to \infty }
  \sum_{0 < \Im(\rho) < T_n} \frac{e^{\rho z} \zeta(\rho - 1)}{ \zeta'(\rho)}.
\]

\ifthenelse {\boolean{BKMRK}}
  { \section{Tallying the Liouville \texorpdfstring{$\lambda$}{lambda} Function} }
  { \section{Tallying the Liouville $\lambda$ Function} }

The Liouville lambda function is defined as follows.  For a positive integer $n$, factor $n$ into primes.  If $\Omega(n)$ is the number of such factors counting multiplicity, then $\lambda(n) = (-1)^{\Omega(n)}$.  For example, if $n = 12$, then $\lambda(n) = (-1)^3 = -1$, because 12 has three prime factors: 2, 2, and 3.  We will denote the summatory function of $\lambda$ by $L(x)$:
\[
  L(x) = \sum _{n=1}^x \lambda(n).
\]
It is not hard to see that $L(x)$ is the number of integers in the range $1 \leq n \leq x$ that have an even number of prime factors minus the number of integers in that range that have an odd number of prime factors.  This is because each $n$ with an even number of prime factors contributes $+1$ to the sum $L(x)$, while each $n$ with an odd number of prime factors contribute $-1$.

The following Dirichlet series holds for $\Re(s) > 1$ \cite[Theorem 300, p. 255]{HardyAndWright}, \cite[eq. D-53]{Gould}, \cite[p. 255]{McCarthy}:
\[
\sum _{n=1}^{\infty } \frac{\lambda (n)}{n^s}=\frac{\zeta (2 s)}{\zeta (s)}.
\]

The analysis is similar to those in the previous sections.

From this Dirichlet series, Perron's formula tells us that, if $c > 1$, then
\begin{equation}\label{E:LxIntegral}
 L_0(x)
 = \lim_{T\to \infty } \, \frac{1}{2 \pi i}
   \int_{c-i T}^{c+i T} \frac{\zeta(2s)}{\zeta(s)} \frac{x^s}{s} \, ds.
\end{equation}

The integrand has a pole at $s = 0$ where the residue is $1$.  $\zeta(s)$ has a pole of order 1 at $s = 1$, so, because of the $\zeta(2s)$ in the numerator, the integrand has a pole at $2s = 1$, that is, at $s = 1/2$.  The residue at $s = 1/2$ is
$\sqrt{x} / \zeta \left(\frac{1}{2}\right) \simeq -0.684765 \sqrt x$.

The integrand also has poles at every real and complex zero of zeta.  We will apply residue formula \eqref{E:stdResidueFormula1} with $B(s) = \zeta(s)$ and $A(s)$ equal to everything else in the integrand, that is, $A(s) = \zeta(2s) x^s/s$.  The residue at the zeta zero $s = \rho$ will be
\[
\frac{A(\rho)}{B'(\rho)} = \frac{\zeta (2 \rho) x^{\rho }}{\rho } \frac{1}{\zeta '(\rho )}
 = \frac{x^{\rho } \zeta(2 \rho) }{\rho \zeta '(\rho ) }.
\]

This expression holds for any zeta zero, whether real or complex.  Recall that the negative even integers are the real zeros of zeta.  Therefore, if $\rho = -2k$ is any of these real zeros of zeta, then the $\zeta(2 \rho)$ in the numerator is always 0, because $2 \rho$ is another real zero.  So, in the sum of the residues that correspond to the real zeros, every term will be 0.

Therefore, we can create a rectangle that encloses the first $N$ pairs of complex zeta zeros and the poles at $s = 0$ and $s = 1/2$.  When we integrate around the rectangle, the integral is just the sum of the residues, which suggests that we might have the approximation

\begin{equation}\label{E:LxSum}
L_0(x) \simeq
  1 + \frac{\sqrt{x}}{\zeta \left(\frac{1}{2}\right)}
  + 2 \Re\left(\sum _{k=1}^N x^{\rho _k} \frac{ \zeta \left(2
   \rho_k\right)}{\rho _k \zeta '\left(\rho_k\right)}\right).
\end{equation}

\begin{figure*}[ht]
  \centerline{
    \mbox{\includegraphics[width=\picDblWidth]{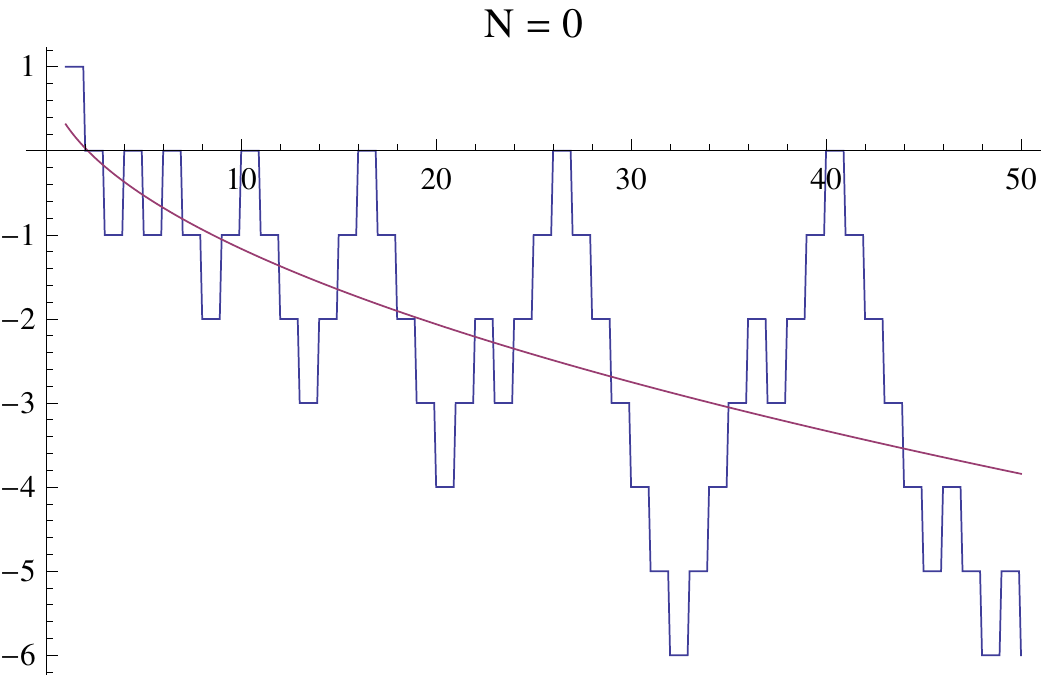}}
    \mbox{\includegraphics[width=\picDblWidth]{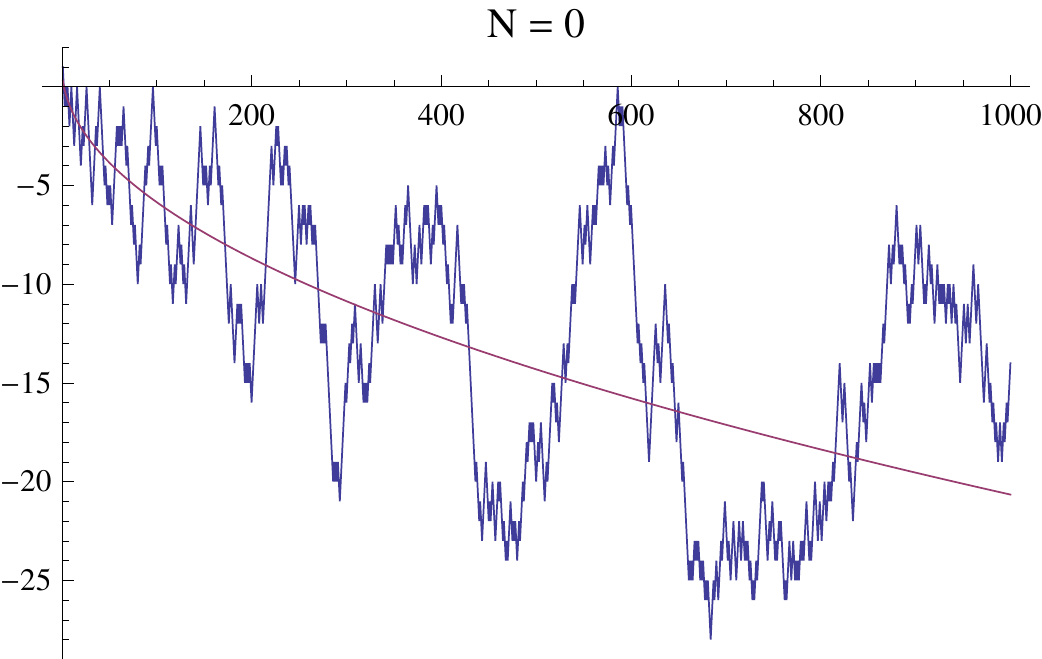}}
  }
  \caption{Approximating $L_0(x)$ with $1 + \frac{\sqrt x}{\zeta(1/2)}$}
  \label{fig:LxSumApprox0}
\end{figure*}

Figure ~\ref{fig:LxSumApprox0} shows the approximation to the step function based on the first two terms of equation \eqref{E:LxSum}.  As you can see, (unlike the corresponding approximation in the previous section), these first two terms completely obscure the behavior of the step function.

\begin{figure*}[ht]
  \centerline{
    \mbox{\includegraphics[width=\picDblWidth]{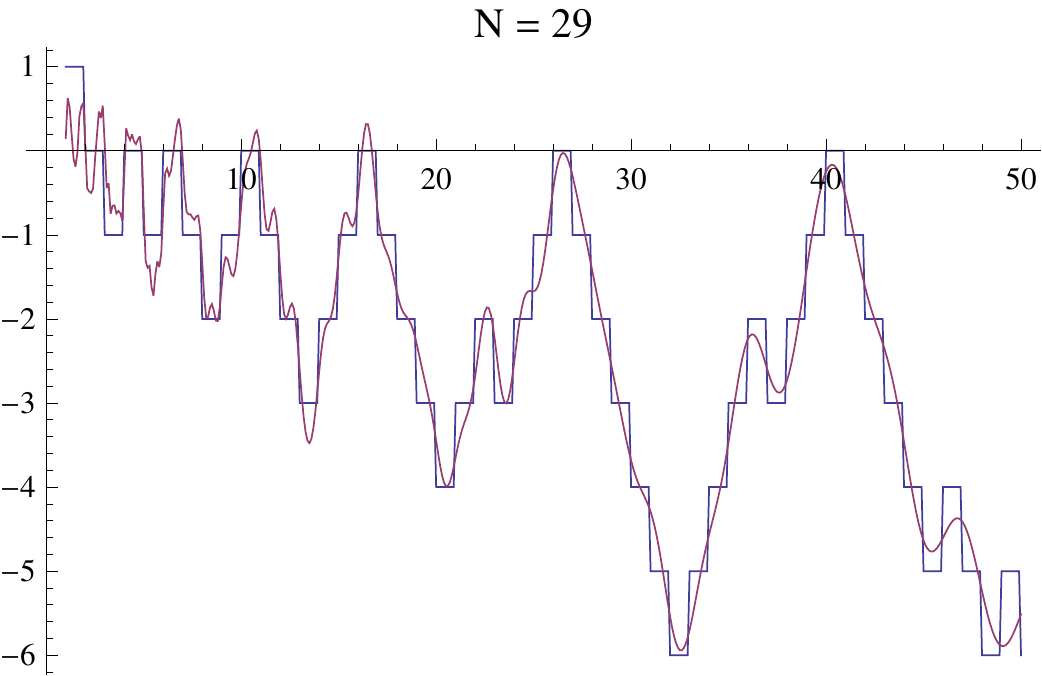}}
    \mbox{\includegraphics[width=\picDblWidth]{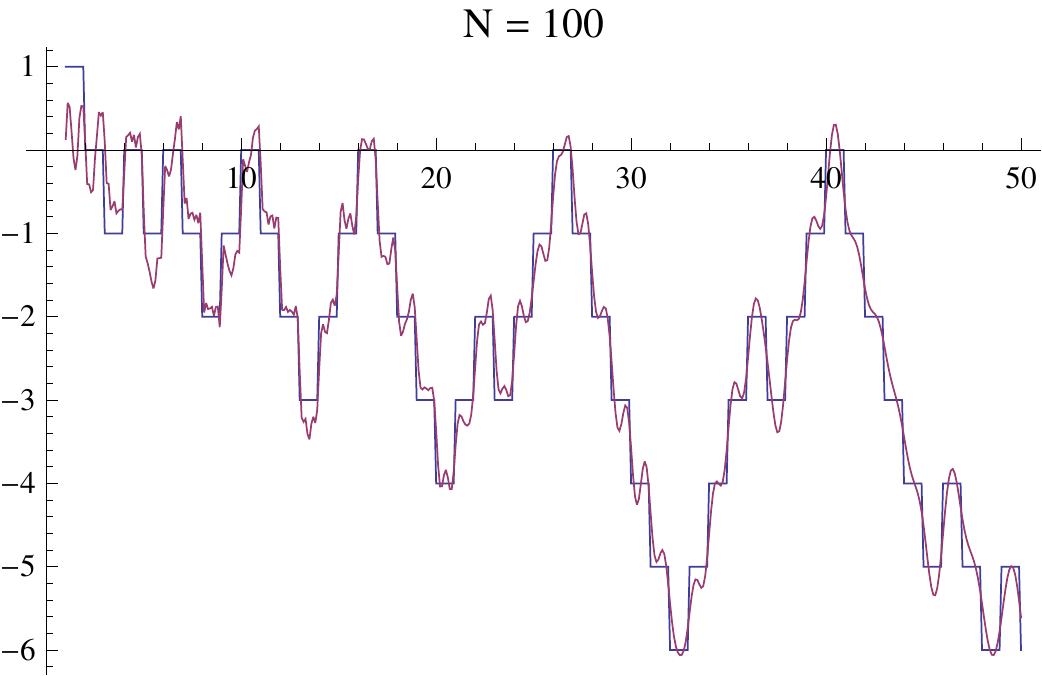}}
  }
  \caption{Approximating $L_0(x)$ with sums over $N = 29$ and $N = 100$ pairs of zeta zeros}
  \label{fig:LxSumApproximations}
\end{figure*}

\begin{figure*}[ht]
  \mbox{\includegraphics[width=\picDblWidth]{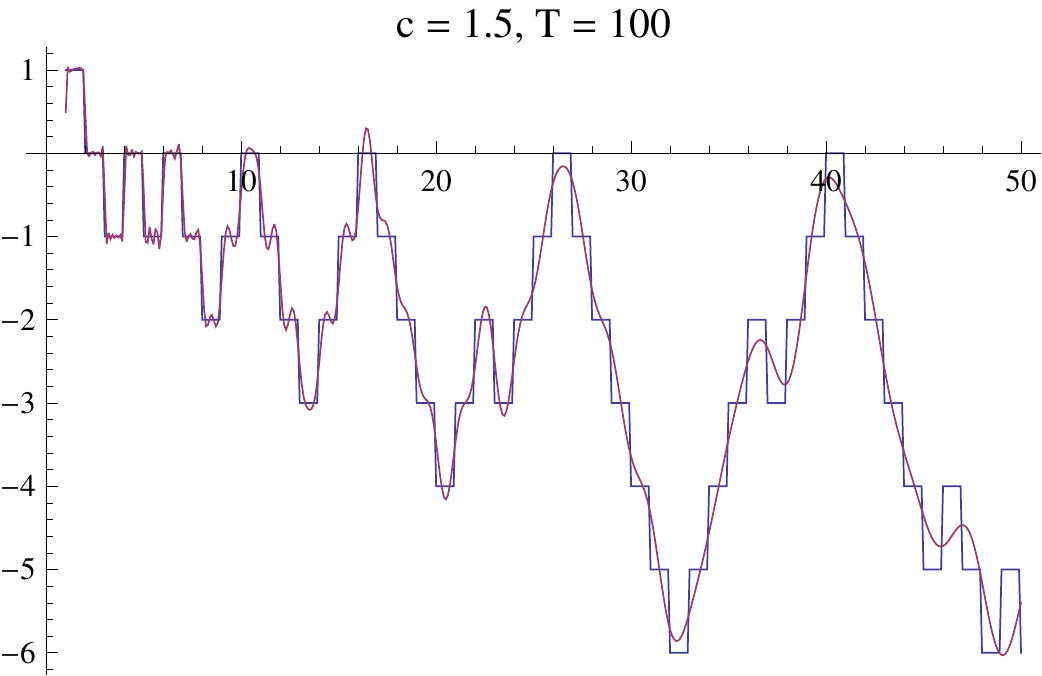}}
  \caption{Approximating $L_0(x)$ with an integral; $c = 1.5$, $T = 100$}
  \label{fig:LxIntegralApproximation}
\end{figure*}

\begin{figure*}[ht]
  \mbox{\includegraphics[width=\picDblWidth]{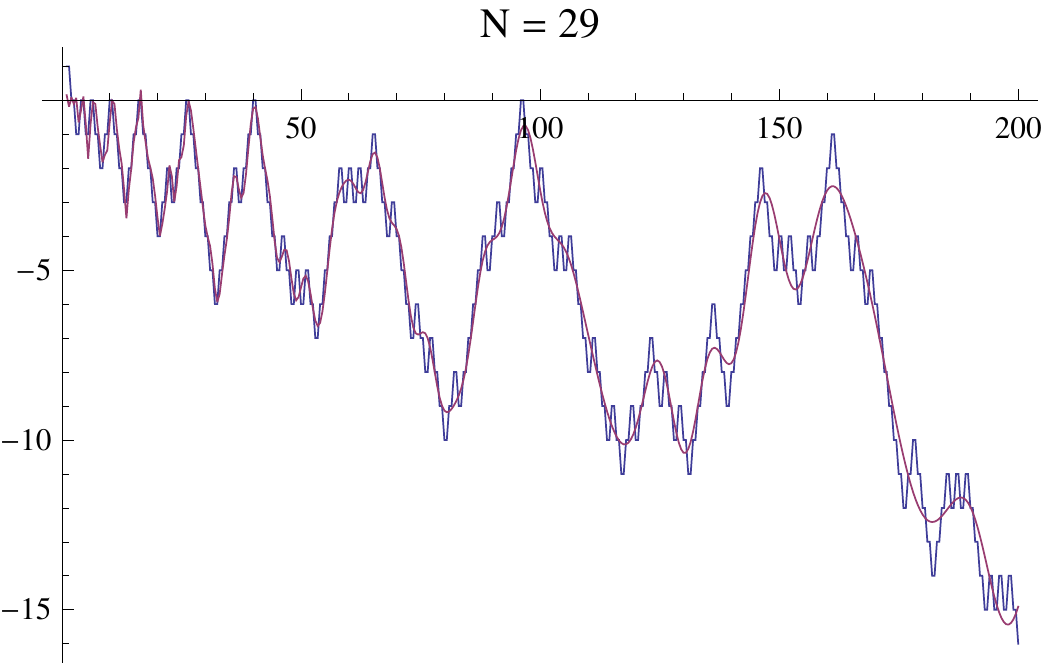}}
  \caption{Approximating $L_0(x)$ with a sum over $N = 29$ pairs of zeta zeros}
  \label{fig:LxSumApproximationx200}
\end{figure*}

You can see from the graph on the left of Figure ~\ref{fig:LxSumApprox0} that for $1 < x \leq 50$, the step function $L(x)$ is less than or equal to 0.  In the graph on the right, $L(96) = L(586) = 0$, so for $1 < x \leq 1000$, $L(x)$ is also less than or equal to 0.  In 1919, George P\'{o}lya conjectured that $L(x) \leq 0$ for all $x > 1$.  However, in 1958, this conjecture was proven to be false.  We now know \cite{Tanaka1} that the smallest counterexample $x > 1$ for which $L(x) > 0$ is $x = 906150257$.  The first ten values of $x$ for which $L(x) = 0$ are 2, 4, 6, 10, 16, 26, 40, 96, 586, and 906150256 \cite[Seq. A028488]{Sloane}.  It is not known \cite{Tanaka1} whether $L(x)$ changes sign infinitely often, but it is known (see \cite{Tanaka1} and \cite{Tanaka2}) that
\[
L(x) - \frac{\sqrt{x}}{\zeta \left(\frac{1}{2}\right)}
\]
does change sign infinitely often.


For small $x$, the approximations in Figure ~\ref{fig:LxSumApproximations}, are not terribly close to the step function $L_0(x)$, even when 100 pairs of zeros are used.  However, the integral approximation shown in Figure ~\ref{fig:LxIntegralApproximation} is a lot better.

To investigate this discrepancy, let's set $x = 3.5$.  The value of the step function at $x = 3.5$ is $L(3.5) = -1$.  Let's integrate around the rectangle that extends from $1.5 - 100 i$ to $-1 + 100 i$.  The integrals along the right, top, left, and bottom sides are: $I_1 \approx -1.01558$, $I_2 \approx 0.00623 - 0.01804 i$, $I_3 \approx .35419$, and $I_4 \approx 0.00623 + 0.01804 i$.  The sum of these four integrals is $ \approx -0.64893 $.  Note that $I_1$, the integral from $1.5 - 100 i$ to $1.5 + 100 i$, is fairly close to $L(3.5)$.  If the integral along the left side ($I_3$) was, say, 0, then the sum would be $ \approx -1.00312$, much closer to $L(3.5)$

Finally, Figure ~\ref{fig:LxSumApproximationx200} shows that, for $x \leq 200$, the sum over $N = 29$ pairs of zeros tracks the irregular step function $L(x)$ pretty closely.  As usual, as $x$ increases, for a fixed $N$, the fine details of the step function gradually get lost.

\section{Tallying the Squarefree Divisors}\label{S:TallyingTheSquarefreeDivisors}

For positive integers $n$, $\nu(n)$ is defined to be the number of distinct prime factors of $n$.  Therefore, if the prime factorization of $n$ is
\[
n=p_1^{a_1} \cdot p_2^{a_2} \cdots p_k^{a_k},
\]
then $\nu(n) = k$.  It is not hard to see that $2^{\nu(n)}$ is the number of \textit{squarefree} divisors of $n$.  This example shows why:  Let $n = 60 = 2^2 \cdot 3 \cdot 5$.  $\nu(60) = 3$, because 60 has three distinct prime factors.  We can list all the squarefree divisors of 60 by taking all  $2^3 = 8$ combinations of the presence or absence of the three primes, that is, $2^a \cdot 3^b \cdot 5^c$, where $a$, $b$, and $c$ independently take the values 0 and 1.

The following Dirichlet series holds for $\Re(s) > 1$ \cite[Theorem 301, p. 255]{HardyAndWright}, \cite[eq. D-9]{Gould}.
\[
\sum _{n=1}^{\infty } \frac{2^{\nu (n)}}{n^s}=\frac{\zeta(s)^2}{\zeta (2 s)}.
\]

Let's denote the summatory function of $2^{\nu(n)}$ by $T(x)$, so that
\[
T(x) = \sum _{n=1}^x 2^{\nu (n)}.
\]

Perron's formula says that, for $c > 1$,
\begin{equation}\label{E:TwoNuIntegral}
 T_0(x)
 = \lim_{T\to \infty } \, \frac{1}{2 \pi i}
   \int_{c-i T}^{c+i T} \frac{\zeta(s)^2}{\zeta(2 s)} \frac{x^s}{s} \, ds.
\end{equation}
where the subscript 0 has the usual meaning.  The integrand
\[
  \frac{\zeta(s)^2}{\zeta(2 s)} \frac{x^s}{s}
\]
has poles at $s = 0$ and at every place where $2s$ is a zeta zero.  It also has a pole of order 2 at $s = 1$, due to the presence of zeta squared in the numerator.  The residue at $s = 0$ is $-1/2$.  The residue at $s = 1$ is calculated as follows.

According to \eqref{E:stdResidueFormulaN}, the residue can be obtained by differentiating the product of $(s - 1)^2$ and the above integrand, then taking the limit as $s$ approaches 1.  When we differentiate the product
\[
(s - 1)^2 \frac{\zeta(s)^2}{\zeta(2 s)} \frac{x^s}{s},
\]
with respect to s, we get
\begin{align*}
 &-\frac{(s-1)^2 x^s \zeta (s)^2}{s^2 \zeta (2 s)}
 -\frac{2 (s-1)^2 x^s \zeta(s)^2 \zeta '(2 s)}{s \zeta (2 s)^2}
 +\frac{(s-1)^2 x^s \zeta (s)^2 \log (x)}{s\zeta (2 s)} \\
  & \quad  
 + \frac{2 (s-1)^2 x^s \zeta (s) \zeta '(s)}{s \zeta (2 s)}+\frac{2 (s-1) x^s \zeta(s)^2}{s \zeta (2 s)}.
\end{align*}

The limits of the first three terms as $s$ approaches 1 are easy because \cite[Theorem 281, p. 247]{HardyAndWright}
\[
\lim_{s\to 1} \, (s-1) \zeta (s) = 1.
\]
$\zeta(2) = \pi^2/6$, so the limits of the first three terms are $-6x/\pi^2$, $-72 x \zeta '(2)/\pi ^4$, and $6 x \log (x)/\pi ^2$.

For the fourth and fifth terms, we will need two facts.  First \cite[p. 166]{Finch},
\begin{equation}\label{ZetaAndStieltjesConstants}
\zeta(s)=\frac{1}{(s-1)} + \gamma + \sum_{n=1}^{\infty} \frac{(-1)^n}{n!} \gamma_n (s - 1)^n
\end{equation}
where $\gamma \simeq .57721566$ is Euler's constant, and $\gamma_n$ are the Stieltjes constants.  Second \cite[Theorem 283, p. 247]{HardyAndWright},
\[
\zeta '(s)=\frac{-1}{(s-1)^2} + O(1).
\]
The sum of the 4th and 5th terms above is
\[
\frac{2 x^s (s-1) \zeta (s)}{s \zeta(2 s)} (\zeta (s) + (s - 1) \zeta '(s)).
\]
The limit of
\[
\frac{2 x^s (s-1) \zeta (s)}{s \zeta(2 s)}
\]
is $12 x/\pi^2$, and we claim that the limit of
\[
\zeta (s) + (s - 1) \zeta '(s)
\]
is $\gamma$.  To see this, we write out the first few terms of the Laurent series for $\zeta$ and $\zeta'$:
\begin{align*}
&\zeta (s) + (s - 1) \zeta '(s) \\
& = \frac{1}{s - 1} + \gamma + (\text{powers of $(s-1)$}) + (s - 1)\left(\frac{-1}{{s - 1}^2} + O(1)\right) \\
& =  \gamma + (\text{powers of $(s-1)$}) + (s - 1) O(1).
\end{align*}
As $s$ approaches 1, this expression approaches $\gamma$, as claimed.  Combining these limits, the residue of the integrand at $s = 1$ is
\[
\frac{-6x}{\pi^2}  -\frac{72 x \zeta '(2)}{\pi ^4} + \frac{6 x \log (x)}{\pi ^2} + \frac{12 x \gamma}{\pi^2} \\
= \frac{6 x \left(-12 \zeta '(2)+2 \gamma  \pi ^2-\pi ^2\right)}{\pi ^4} + \frac{6 x \log (x)}{\pi ^2}.
\]

Next, if $\rho$ is any real or complex zeta zero, we can derive an expression for the residue at $s = \rho/2$ using the same method as in section \ref{S:CountingTheSquarefreeIntegers} where we counted the squarefree integers.  Again, we will use equation \eqref{E:stdResidueFormula1} for the residue of the quotient $A(s)/B(s)$
\[
  \frac{A(s)}{B(s)} = \frac{\zeta(s)^2}{\zeta(2s)} \frac{x^s}{s},
\]
where $A(s) = \zeta(s)^2 \frac{x^s}{s}$ and $B(s) = \zeta(2s)$.  Then $B'(s) = 2 \zeta(2s)$, and the residue at the zeta zero $s = \rho / 2$ is
\[
\frac{A(s)}{B'(s)}
  = \frac{\zeta(\frac{\rho}{2})^2}{2 \zeta'(\rho)} \frac{x^\frac{\rho}{2}}{\frac{\rho}{2}}
  = \frac{x^{ \frac{\rho}{2} } \zeta( \frac{\rho}{2} )^2 }{\rho \zeta'(\rho) }.
\]

Putting all these residues together, we get this alleged approximation for $T(x)$ that involves the first $N$ pairs of complex zeta zeros and the first $M$ real zeros:
\begin{equation}\label{E:TwoNuSumApprox}
T_0(x) \simeq
 a_1 + a_2 x+a_3 x \log (x)+
 2 \Re\left(\sum _{k=1}^N x^{\frac{\rho _k}{2}} \frac{
   \zeta \left(\frac{\rho _k}{2}\right){}^2}{\rho _k \zeta
   '\left(\rho _k\right)}\right)
   +\sum _{k=1}^{M } x^{-\frac{2 k}{2}}  \frac{
     \zeta \left(\frac{-2 k}{2}\right)^2}{(-2 k) \zeta '(-2 k)}
\end{equation}
where
\[
a_1 = -1/2,
\]
\[
a_2 = \frac{6 \left(-12 \zeta '(2)+2 \gamma  \pi ^2-\pi ^2\right)}{\pi ^4} \simeq 0.78687,
\]
and
\[
a_3 = 6/\pi^2 \simeq .60793.
\]

The second sum is too small to visibly affect the graphs, so we'll set $M = 0$ when we draw the graphs below.

\begin{figure*}[ht]
  \mbox{\includegraphics[width=\picDblWidth]{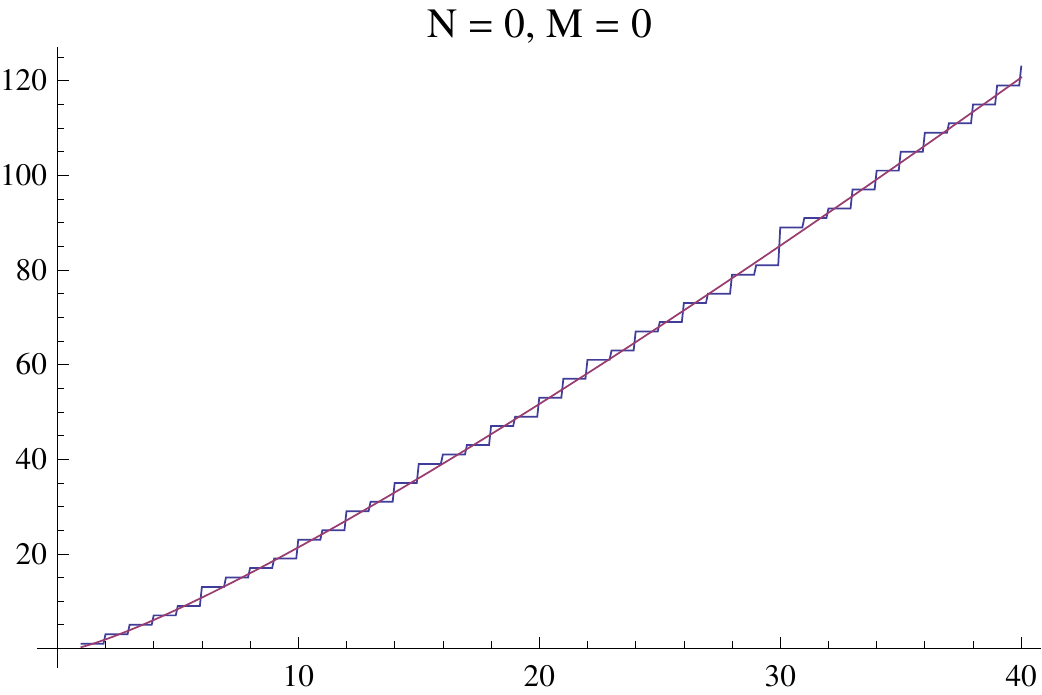}}
  \caption{Approximating the sum of $2^{\nu(n)}$ with the first three terms of \eqref{E:TwoNuSumApprox}}
  \label{fig:TwoNuSumApprox0}
\end{figure*}

\begin{figure*}[ht]
  \centerline{
    \mbox{\includegraphics[width=\picDblWidth]{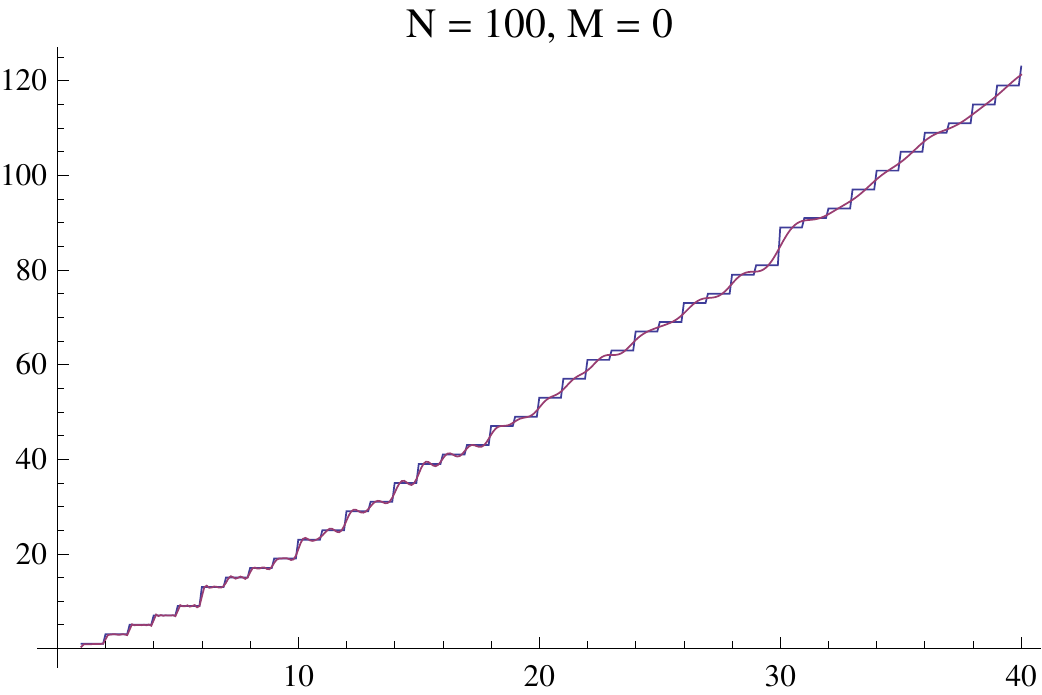}}
    \mbox{\includegraphics[width=\picDblWidth]{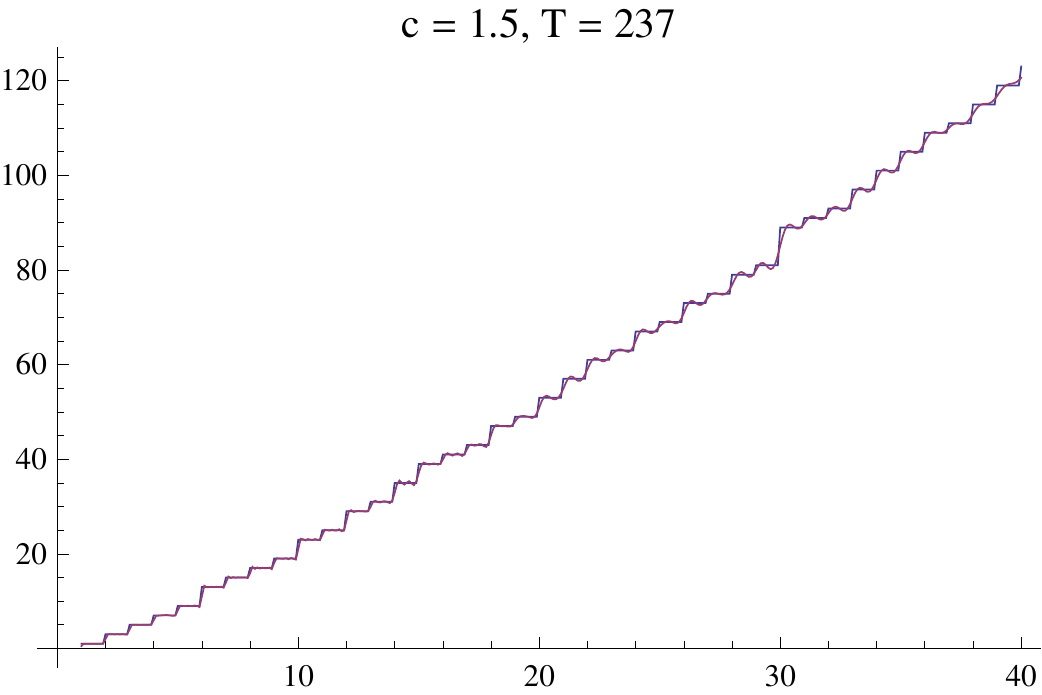}}
  }
  \caption{Approximating the sum of $2^{\nu(n)}$ with a sum and an integral}
  \label{fig:TwoNuSumApproximations}
\end{figure*}

\begin{figure*}[ht]
  \centerline{
    \mbox{\includegraphics[width=\picDblWidth]{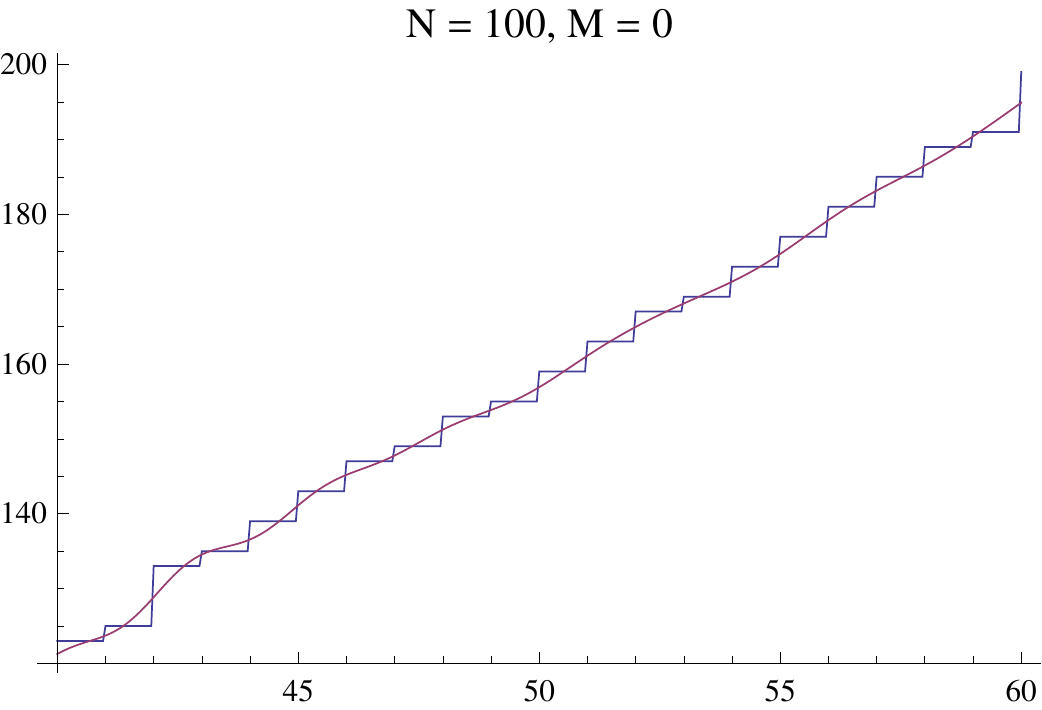}}
    \mbox{\includegraphics[width=\picDblWidth]{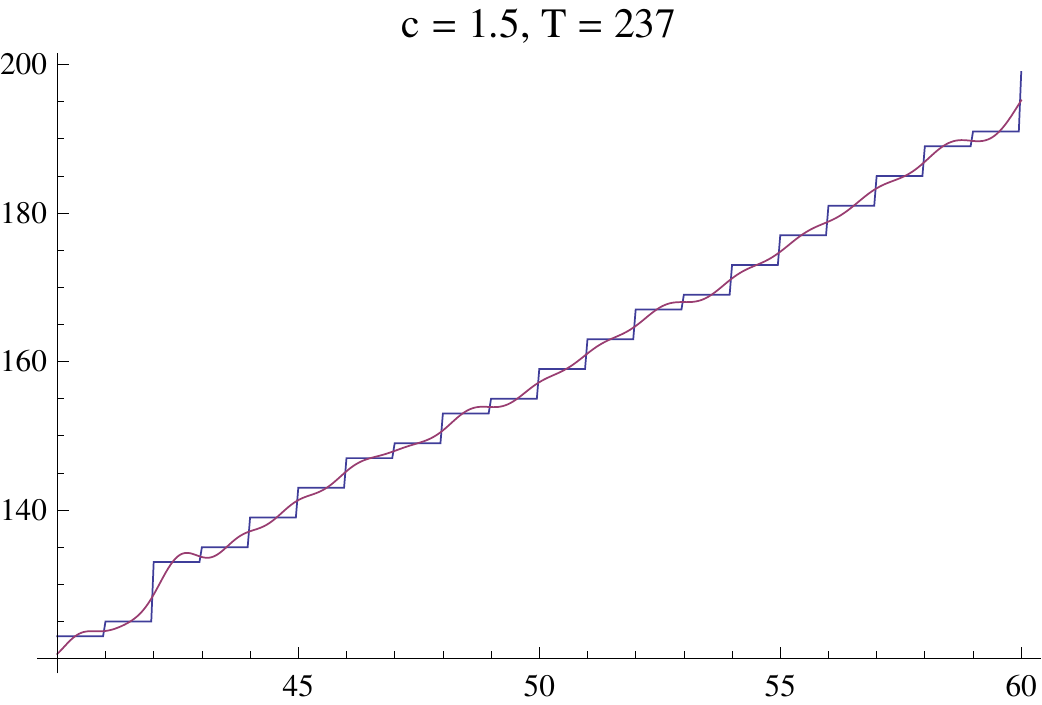}}
  }
  \caption{Approximating the sum of $2^{\nu(n)}$ with the same sum and integral for $40 \leq x \leq 60$}
  \label{fig:TwoNuSumApproximationsx40-60}
\end{figure*}

We will now use \eqref{E:TwoNuSumApprox} to draw a graph of $T_0(x)$.  Figure ~\ref{fig:TwoNuSumApprox0} shows the approximation to the step function based on the first three terms of equation \eqref{E:TwoNuSumApprox}.  Figure ~\ref{fig:TwoNuSumApproximations} shows the approximations using \eqref{E:TwoNuSumApprox} with $N = 100$, $M = 0$, and the integral
\[
 T_0(x) \simeq
   \frac{1}{2 \pi i}
   \int_{c-i T}^{c+i T} \frac{\zeta(s)^2}{\zeta(2 s)} \frac{x^s}{s} \, ds.
\]
with $c = 1.5$ and $T = 237$.  This value of $T$ is large enough for the corresponding rectangle to enclose 100 pairs of complex zeros.  The approximations are quite good in this range of $x$.  Figure ~\ref{fig:TwoNuSumApproximationsx40-60} shows the same approximations, but for the range $40 \leq x \leq 60$.

\textit{Mathematica} can compute the residue at the second-order pole at $s = 1$, as follows:
\begin{verbatim}
    res = Residue[Zeta[s]^2/Zeta[2 s] * x^s/s, {s, 1}]
    Collect[res, Log[x], Simplify].
\end{verbatim}
The second command separates out the term that contains $\log(x)$.  The result is
\[
\frac{6 x \log (x)}{\pi ^2} +
\frac{6 x \left((2 \gamma -1) \pi ^2-12 \zeta '(2)\right)}{\pi^4}.
\]

Assuming the zeta zeros are simple, Wiertelak \cite[Theorem 4]{Wiertelak} proves that, for $x > 1/4$,
\begin{equation}\label{E:WiertelakEq}
T_0(x) =
 a_1 + a_2 x+a_3 x \log (x)+
 \sum_{0 < \Im(\rho) < T_n}
   x^{\frac{\rho _k}{2}} \frac{ \zeta(\rho/2)^2 }{ \rho \zeta'(\rho) }
  + \frac{1}{\pi^2 x} \sum _{k=1}^{\infty } \frac{ \zeta(2k+2)^2 x^{-2 k} }{ (2k+1) \binom{4k+2}{2k+1} \zeta(4k+3) }
\end{equation}
for a certain sequence $T_n$, where $a_1$, $a_2$, and $a_3$ are the values given above.  Note: because
\[
\zeta'(-2 k) = (-1)^k \frac{(2 k)!}{2 (2 \pi)^{2 k}} \zeta(2 k + 1),
\]
the second sum in \eqref{E:WiertelakEq} has the same as form the second sum in equation \eqref{E:TwoNuSumApprox}.

If we have the Dirichlet series for an arithmetical function, the procedure for deriving an approximation should be clear by now.  So, from now on, we will give fewer details about the derivations of the residues except where new issues arise, such as poles of order 2, 3, or 4.  In each section below, we'll use the same symbol, $T(x)$, for the respective summatory function of that section, and $T_0(x)$ for the summatory function modified in the usual way.

\ifthenelse {\boolean{BKMRK}}
  { \section{Computing sigma Sums \texorpdfstring{$\sigma(n^2)$}{sigma(n squared)} and \texorpdfstring{$\sigma(n)^2$}{sigma(n) squared} } }
  { \section{Computing Sigma Sums $\sigma(n^2)$ and $\sigma(n)^2$} }

$\sigma(n)$ is the sum of the (positive) divisors of the positive integer $n$.  Here we consider two sums involving the $\sigma$ function that can be approximated using Perron's formula.

\subsection{The Sum of $\sigma(n^2)$}
This Dirichlet series holds for $s > 3$ \cite[set $k = 1$ in eq. D-51]{Gould}, \cite[eq. 5.39, p. 237]{McCarthy}:
\[
\sum _{n=1}^{\infty} \frac{\sigma (n^2)}{n^s} = \frac{\zeta(s) \zeta(s - 1) \zeta(s - 2)}{\zeta (2s - 2)}.
\]
Note that there is a typo in \cite[eq. D-51]{Gould}: the denominator on the right side should be $\zeta(2(s - k))$.  Let $T_0(x)$ denote the (modified) summatory function of $\sigma(n^2)$.  Then, based on Perron's formula, we have this integral approximation
\begin{equation}\label{E:SigmaSum1IntegralApprox}
T_0(x) \simeq  \frac{1}{2 \pi i}
   \int_{c-i T}^{c+i T} \frac{\zeta(s) \zeta(s - 1) \zeta(s - 2)}{\zeta(2s - 2)} \frac{x^s}{s} \, ds.
\end{equation}

We will not derive an estimate that shows how close this integral is to the sum for given values of $x$ and $T$.  However, we will see that, for modest values of $x$ and $T$, the integral follows the step function rather closely.

The integrand has poles at $s = 0$, $s = 1$, $s = 2$, and $s = 3$.  There are also poles at those $s$ such that $2s - 2 = \rho$, where $\rho$ is any zeta zero, that is, at each $s = \rho/2 + 1$.  The sum of the residues, which we hope will be close to $T_0(x)$, is 

\begin{equation}\label{E:SigmaSum1Approx}
T_0(x) \simeq
a_1 + a_2 x + a_3 x^2 + a_4 x^3
+ 2 \Re\left(\sum _{k=1}^N x^{\frac{\rho _k}{2}+1} \frac{
   \zeta \left(\frac{\rho _k}{2}+1\right) \zeta \left(\frac{\rho
   _k}{2}\right) \zeta \left(\frac{\rho _k}{2}-1\right)}{2
   \left(\frac{\rho _k}{2}+1\right) \zeta '\left(\rho
   _k\right)}\right)
\end{equation}
where
\[
a_1=\frac{1}{48}\simeq 0.0208333,
\]
\[
a_2=-\frac{1}{12}\simeq -0.083333,
\]
\[
a_3=-\frac{1}{4}=-0.25,
\]
and
\[
a_4=\frac{5 \zeta (3)}{\pi ^2}\simeq 0.608969.
\]

The terms $a_1$, $a_2 x$, $a_3 x^2$, and $a_4 x^3$ are the residues of the integrand at $s = 0$, $s = 1$, $s = 2$, and $s = 3$, respectively.  One would expect to have another sum similar to the one in \eqref{E:SigmaSum1Approx}, but over the real zeros instead of the complex zeros.  However, this sum is always 0: If $k \ge 1$, then if $\rho = -2k$ is the $k^{th}$ real zeta zero, then at least one of $\rho/2 + 1$, $\rho/2$, or $\rho/2 - 1$ will also be a real zeta zero, which causes the $k^{th}$ term to drop out.

The factor of 2 in the denominator of the terms in \eqref{E:SigmaSum1Approx} arises when we compute the residue of $A(s)/B(s)$, where
\[
A(s) = \zeta(s) \zeta(s-1) \zeta(s-2) \frac{x^s}{s}
\]
and
\[
B(s) = \zeta(2s - 2).
\]
When we use residue rule \eqref{E:stdResidueFormula1}, we compute the derivative of $B(s)$ with respect to $s$.  This produces the factor of 2 that appears in \eqref{E:SigmaSum1Approx}.

\begin{figure*}[ht]
  \mbox{\includegraphics[width=\picDblWidth]{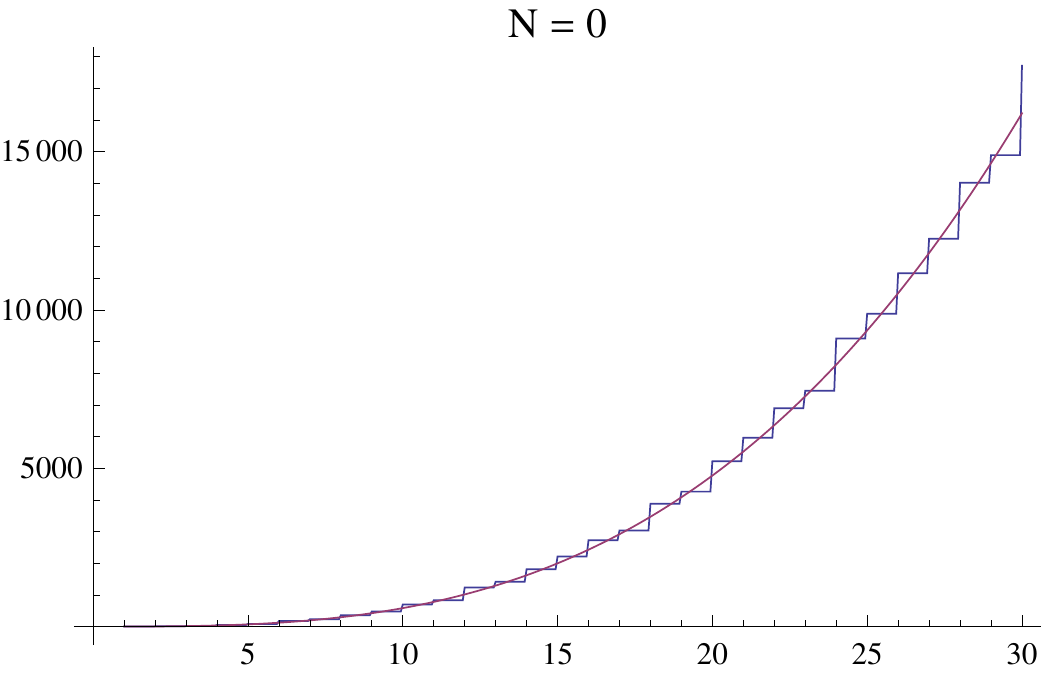}}
  \caption{Approximating the sum of $\sigma(n^2)$ with the first four terms of \eqref{E:SigmaSum1Approx}}
  \label{fig:SigmaSum1Approx0}
\end{figure*}

\begin{figure*}[ht]
  \centerline{
    \mbox{\includegraphics[width=\picDblWidth]{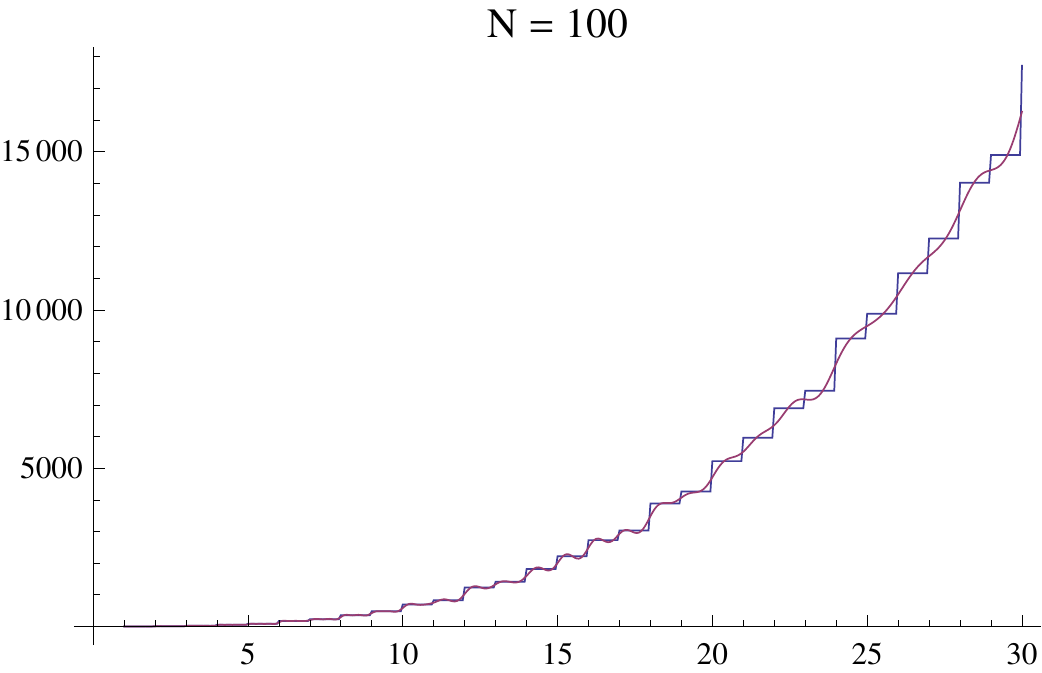}}
    \mbox{\includegraphics[width=\picDblWidth]{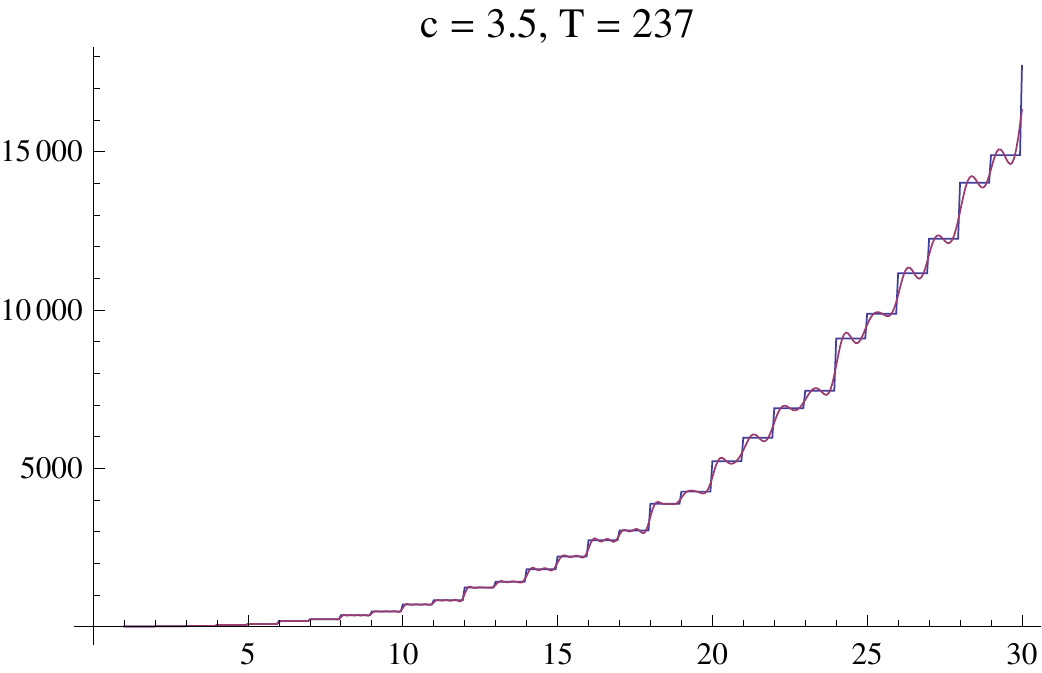}}
  }
  \caption{Approximating the sum of $\sigma(n^2)$ using a sum and an integral}
  \label{fig:SigmaSum1Approximations}
\end{figure*}

Figure ~\ref{fig:SigmaSum1Approx0} shows the polynomial approximation to the summatory function $T_0(x)$ that comes from the first four terms of \eqref{E:SigmaSum1Approx}.

Figure ~\ref{fig:SigmaSum1Approximations} shows approximations to $T_0(x)$ using a sum \eqref{E:SigmaSum1Approx} and an integral \eqref{E:SigmaSum1IntegralApprox}.  Note that the approximating sum increases roughly as $x^3$.  Moreover, the increases from one value of $x$ to another can be both quite large, and quite irregular.  From $x = 23$ to $x = 24$, the step function changes from 7449 to 9100, an increase of 1651.  From $x = 24$ to $x = 25$, the step function changes from 9100 to 9881, an increase of 781.  In spite of these rapid and irregular steps, the integral approximation tracks the step function quite well.

\subsection{The Sum of $\sigma(n)^2$}

This Dirichlet series holds for $s > 3$ \cite[set $a = b = 1$ in Theorem 305, p. 256]{HardyAndWright}:
\[
\sum _{n=1}^{\infty} \frac{\sigma (n)^2}{n^s} = \frac{\zeta(s) \zeta(s - 1)^2 \zeta(s - 2)}{\zeta (2s - 2)}.
\]
This Dirichlet series is similar to the one for the sum of $\sigma(n^2)$, except here, there is a pole of order 2 at $s = 2$.  The integral approximation to the summatory function $T_0(x)$ is
\begin{equation}\label{E:SigmaSum2IntegralApprox}
T_0(x) \simeq  \frac{1}{2 \pi i}
   \int_{c-i T}^{c+i T} \frac{\zeta(s) \zeta(s - 1)^2 \zeta(s - 2)}{\zeta(2s - 2)} \frac{x^s}{s} \, ds.
\end{equation}

From the residues at $s = 0$, $s = 1$, $s = 2$, $s = 3$, and at each $s = \rho/2 + 1$, we get an expression that may approximate $T_0(x)$:

\begin{equation}\label{E:SigmaSum2Approx}
T_0(x) \simeq
 a_1 + a_2 x + a_3 x^2 + a_4 x^2 \log (x) + a_5 x^3
 + 2 \Re\left(\sum _{k=1}^N x^{\frac{\rho _k}{2}+1} \frac{
   \zeta \left(\frac{\rho _k}{2}+1\right) \zeta \left(\frac{\rho
   _k}{2}\right){}^2 \zeta \left(\frac{\rho _k}{2}-1\right)}{2
   \left(\frac{\rho _k}{2}+1\right) \zeta '\left(\rho_k\right)}\right)
\end{equation}
where
\[
a_1=-\frac{1}{576}\simeq -0.001736,
\]
\[
a_2=\frac{1}{24}\simeq 0.041667,
\]
\[
a_3=\frac{1}{8} \left(\frac{12 \zeta '(2)}{\pi ^2}-4 \gamma +1-2
   \log (2 \pi )\right)\simeq -0.765567,
\]
\[
a_4=-\frac{1}{4}=-0.25,
\]
and
\[
a_5=\frac{5 \zeta (3)}{6}\simeq 1.00171.
\]
The terms $a_1$ and $a_2 x$ are the residues of the integrand at $s = 0$ and $s = 1$, respectively.  The residue at $s = 2$ gives rise to $a_3 x^2$ and $a_4 x^2 \log(x)$.  The term $a_5 x^3$ is the residue at $s = 3$.

\begin{figure*}[ht]
  \mbox{\includegraphics[width=\picDblWidth]{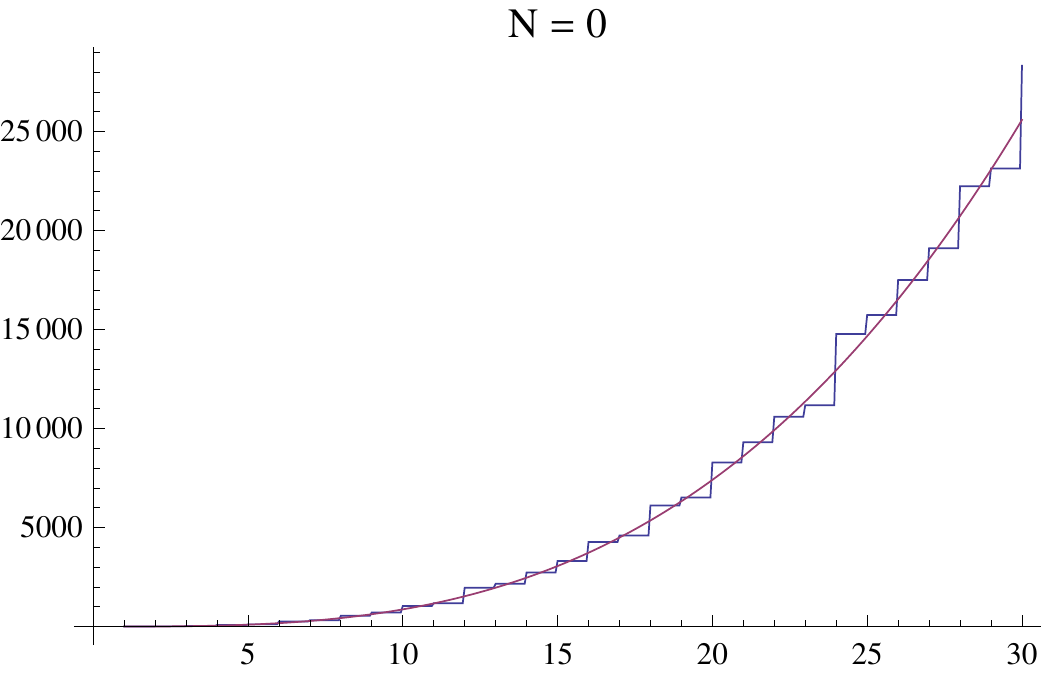}}
  \caption{Approximating the sum of $\sigma(n)^2$ with the first five terms of \eqref{E:SigmaSum2Approx}}
  \label{fig:SigmaSum2Approx0}
\end{figure*}

\begin{figure*}[ht]
  \centerline{
    \mbox{\includegraphics[width=\picDblWidth]{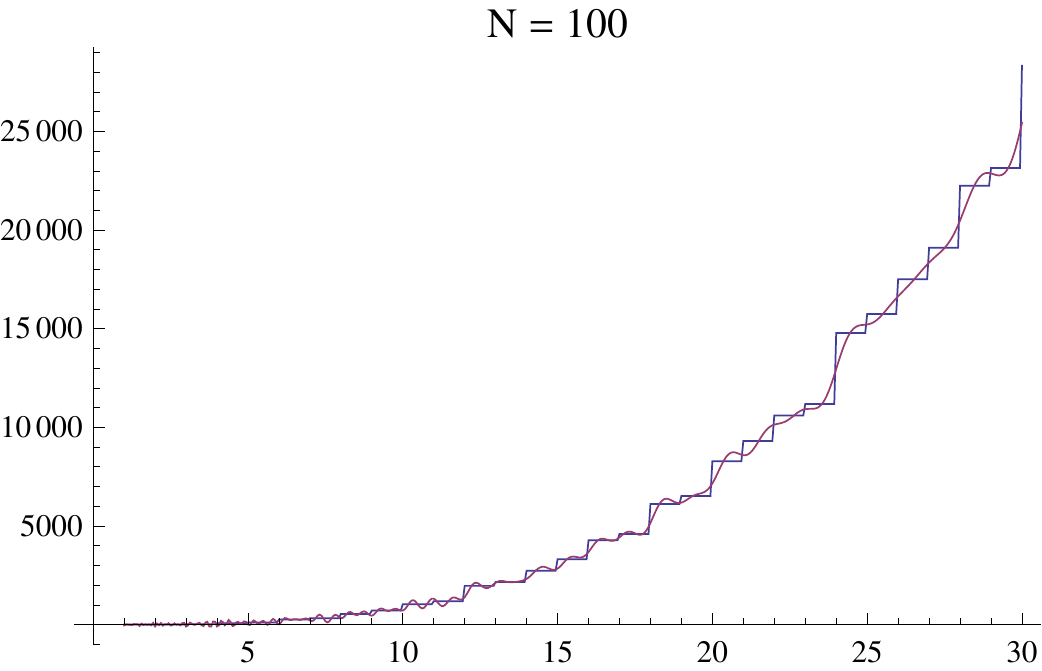}}
    \mbox{\includegraphics[width=\picDblWidth]{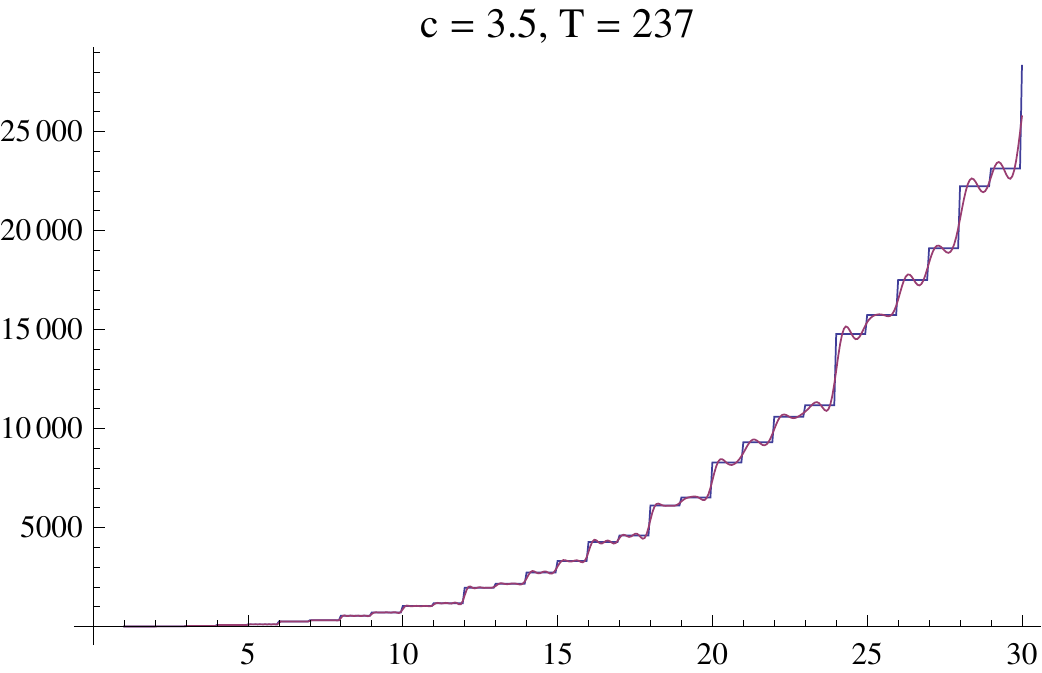}}
  }
  \caption{Approximating the sum of $\sigma(n)^2$ using a sum and an integral}
  \label{fig:SigmaSum2Approximationsx30}
\end{figure*}




\begin{figure*}[ht]
  \mbox{\includegraphics[width=\picDblWidth]{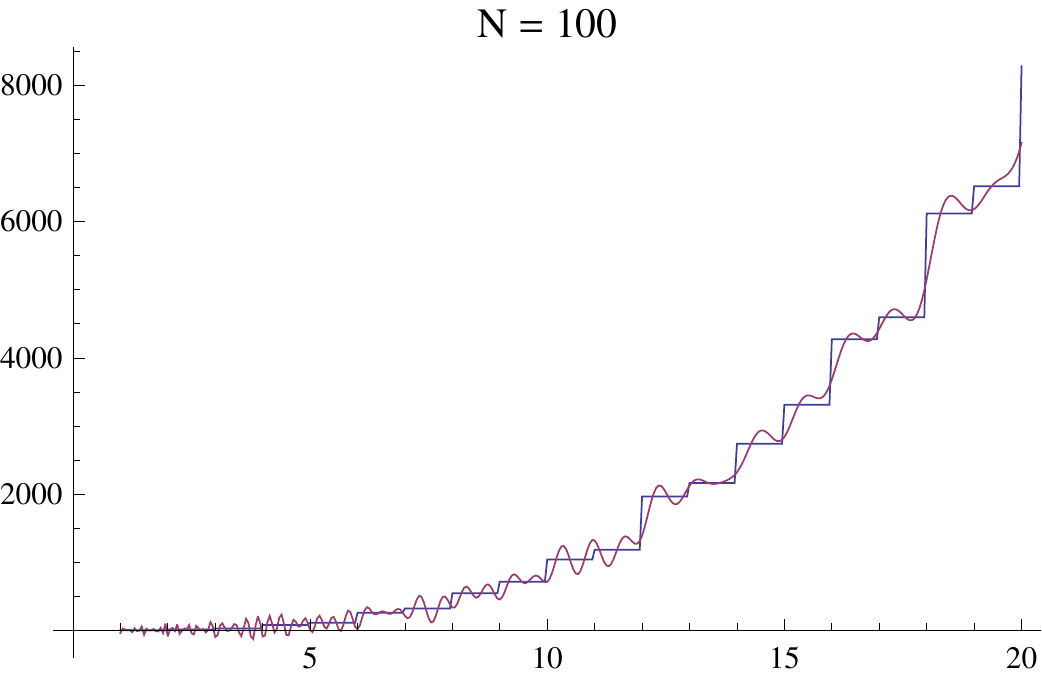}}
  \caption{Approximating the sum of $\sigma(n)^2$ up to $x = 20$ using $N = 100$ pairs of zeros}
  \label{fig:SigmaSum2ApproximationN100x20}
\end{figure*}

\begin{figure*}[ht]
  \mbox{\includegraphics[width=\picDblWidth]{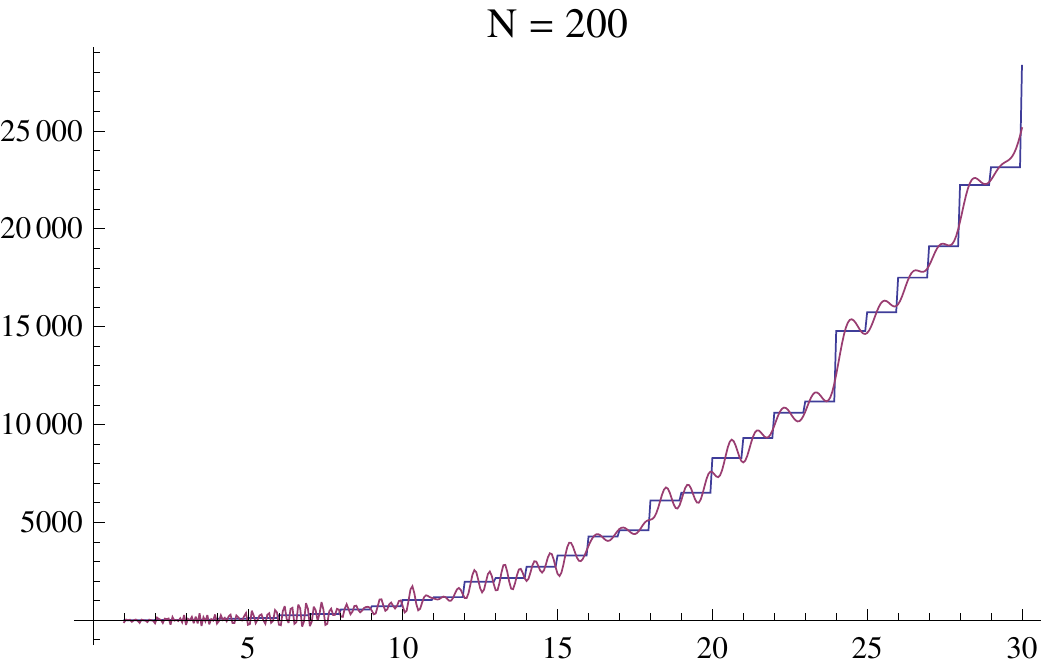}}
  \caption{Approximating the sum of $\sigma(n)^2$ using $N = 200$ pairs of zeros}
  \label{fig:SigmaSum2ApproximationN200}
\end{figure*}

Figure ~\ref{fig:SigmaSum2Approximationsx30} shows the approximation to the step function up to $x = 30$ based on \eqref{E:SigmaSum2Approx} with $N = 100$ pairs of zeta zeros, and the integral \eqref{E:SigmaSum2IntegralApprox} with $c = 3.5$ and $T = 237$.

Figure ~\ref{fig:SigmaSum2ApproximationN100x20} shows the sum in \eqref{E:SigmaSum2Approx}, but going only up to $x = 20$.  Because the scale is different, we can see that the sum is rather ``wavy'' for small $x$.  This raises the question: what would happen if we used more zeta zeros in the sum?

Figure ~\ref{fig:SigmaSum2ApproximationN200}, shows the same sum approximation as in the previous two Figures, but here we use $N = 200$ pairs of zeta zeros.  With more zeros included in the sum, the ``wavyness'' is even more pronounced.  One wonders whether the sum in \eqref{E:SigmaSum2Approx} converges to to the step function as $N$ approaches infinity.



\begin{figure*}[ht]
  \mbox{\includegraphics[width=\picDblWidth]{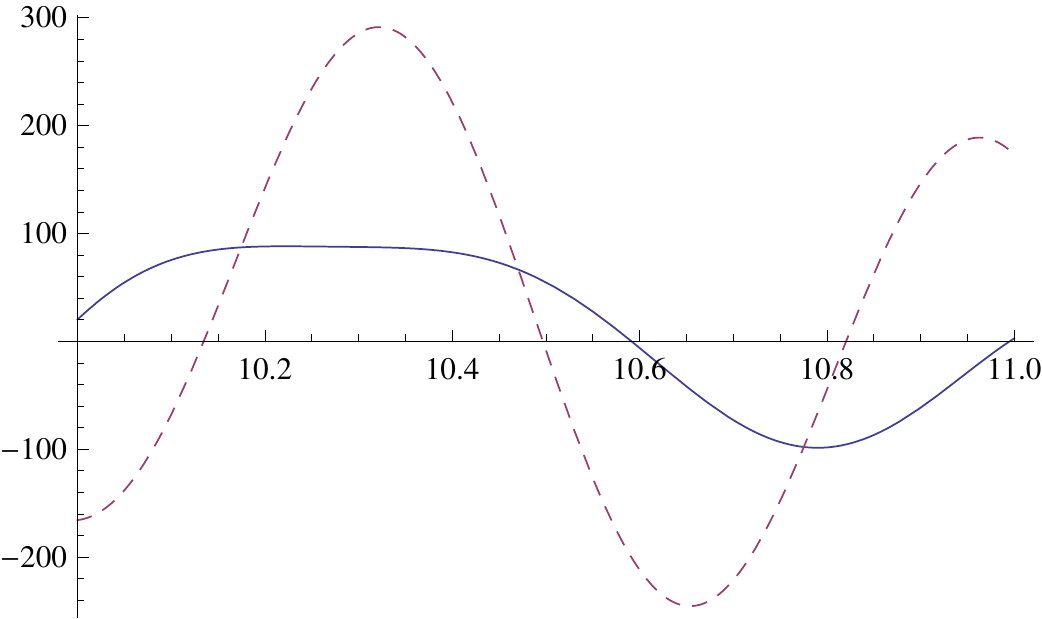}}
  \caption{The sum of the first 89 (solid) and 90 (dashed) zeta terms of \eqref{E:SigmaSum2Approx}}
  \label{fig:SigmaSum2Waves}
\end{figure*}

Figure ~\ref{fig:SigmaSum2Waves} shows a closeup of the sum of the first 89 and 90 zeta terms of \eqref{E:SigmaSum2Approx} (that is, ignoring the initial terms involving $a_1$ through $a_5$).  Between $x = 10$ and $x = 11$, the sum becomes more wavy when the $90^{th}$ term is added to the sum.  At $x = 10.33$, the sum of the first 89 terms (before taking the real part) is about $43.497 + 11.589 i$.  The $90^{th}$ term is about $101.756 + 9.082 i$.  The real part of the $89^{th}$ term is quite a bit larger than the real part of any of the previous terms, or of the sum of the first 89 terms.  This unusually large term is what causes the waves to appear around $x = 11$.  As to why the distortion appears in the form of waves, we'll discuss that in Section \ref{S:Waves}.


\ifthenelse {\boolean{BKMRK}}
  { \section{Computing tau Sums \texorpdfstring{$\tau(n^2)$}{tau(n squared)} and \texorpdfstring{$\tau(n)^2$}{tau(n) squared} } }
  { \section{Computing Tau Sums $\tau(n^2)$ and $\tau(n)^2$} }

$\tau(n)$ is the number of (positive) divisors of the positive integer $n$.  Here we consider two sums involving the $\tau$ function that can be approximated using Perron's formula.

\subsection{The Sum of $\tau(n^2)$}

It can be shown that $\tau(n^2)$ is the number of ordered pairs of positive integers $(i, j)$ such that $lcm(i, j) = n$.

This Dirichlet series holds for $s > 1$ \cite[set $k = 0$ in eq. D-51]{Gould}, \cite[eq. 5.39, p. 237]{McCarthy}, \cite[Equation 1.105, p. 35]{Ivic}:
\[
\sum _{n=1}^{\infty} \frac{\tau (n^2)}{n^s} = \frac{\zeta(s)^3}{\zeta (2s)}.
\]
The integral approximation to the summatory function $T_0(x)$ is
\begin{equation}\label{E:TauSum1IntegralApprox}
T_0(x) \simeq  \frac{1}{2 \pi i}
   \int_{c-i T}^{c+i T} \frac{\zeta(s)^3}{\zeta(2s)} \frac{x^s}{s} \, ds.
\end{equation}
The integrand has a pole of order 3 at $s = 1$, along with poles of order 1 at $s = 0$ and at each $s = \rho/2$, where $\rho$ ranges over all zeta zeros.  The residue at $s = 0$ is $1/4$.  The residue at $s = 1$ is
\begin{align*}
&x
\frac{6 \left(\pi ^4 \left(-3 \gamma _1+1-3 \gamma +3 \gamma ^2\right)+
144 \zeta '(2)^2-12 \pi ^2 \left(\zeta ''(2)-\zeta'(2)+3 \gamma  \zeta '(2)\right)\right)}{\pi ^6} \\
   & \quad  
+ x \log(x)
\frac{6 (3 \gamma -1) \pi ^2-72 \zeta '(2)}{\pi ^4}  \\
   & \quad 
+ x \log(x)^2 \frac{3}{\pi ^2},
\end{align*}
where $\gamma$ is Euler's constant and $\gamma_1 \simeq -.072816$ is the first Stieltjes constant.  The $\gamma_1$ comes from equation \eqref{ZetaAndStieltjesConstants}, the series for the zeta function whose coefficients involve the Stieltjes constants.

The sum over all the residues gives an expression that may approximate the summatory function:
\begin{equation}\label{E:TauSum1Approx}
T_0(x) \simeq
a_1 + a_2 x + a_3 x \log (x) + a_4 x \log ^2(x)
 + 2 \Re\left(\sum _{k=1}^N x^{\frac{\rho _k}{2}} \frac{ \zeta \left(\frac{\rho _k}{2}\right){}^3}{\rho _k \zeta
   '\left(\rho _k\right)}\right)
  + \sum _{k=1}^{M} x^{-\frac{2 k}{2}} \frac{ \zeta \left(-\frac{2 k}{2}\right)^3}{(-2 k) \zeta '(-2 k)},
\end{equation}
where
\[
a_1 = 1/4 = .25,
\]
\[
a_2=\frac{6 \left(\pi ^4 \left(-3 \gamma _1+3 \gamma ^2-3 \gamma
   +1\right)+144 \zeta '(2)^2-12 \pi ^2 \left(\zeta ''(2)+3 \gamma 
   \zeta '(2)-\zeta '(2)\right)\right)}{\pi ^6}\simeq 0.12226,
\]
\[
a_3=\frac{6 (3 \gamma -1) \pi ^2-72 \zeta '(2)}{\pi ^4}\simeq 1.13778,
\]
and
\[
a_4=\frac{3}{\pi ^2}\simeq 0.30396.
\]

The second sum is too small to have any visible effect on the graphs, so when we graph equation \eqref{E:TauSum1Approx}, we will set $M = 0$.

\begin{figure*}[ht]
  \mbox{\includegraphics[width=\picDblWidth]{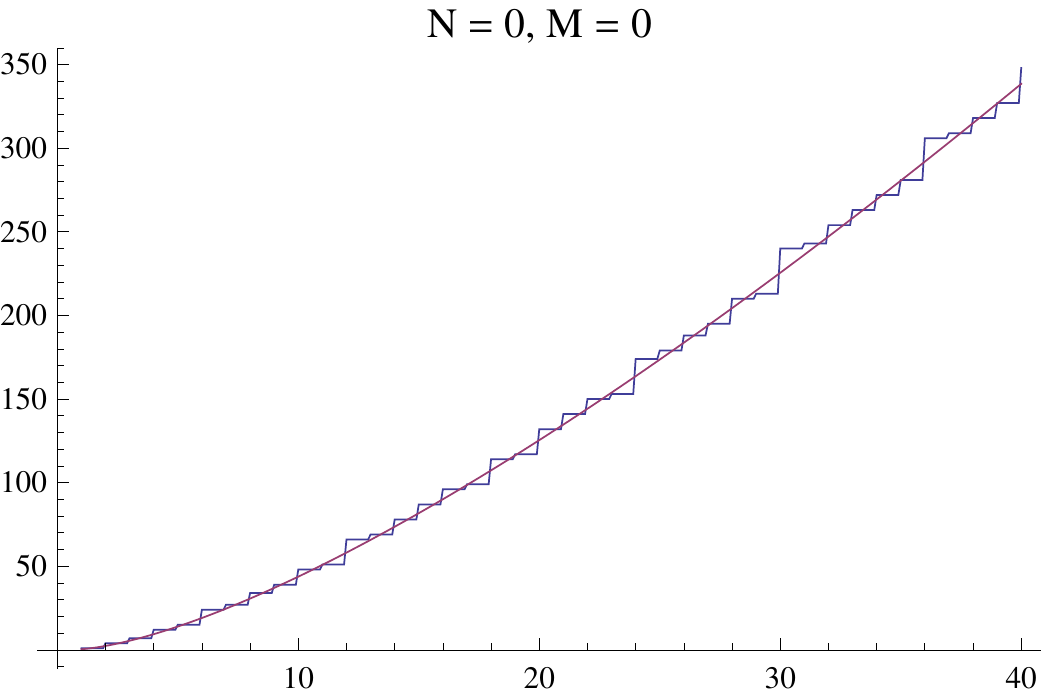}}
  \caption{Approximating the sum of $\tau(n^2)$ with the first four terms of \eqref{E:TauSum1Approx}}
  \label{fig:TauSum1Approx0}
\end{figure*}

\begin{figure*}[ht]
  \centerline{
    \mbox{\includegraphics[width=\picDblWidth]{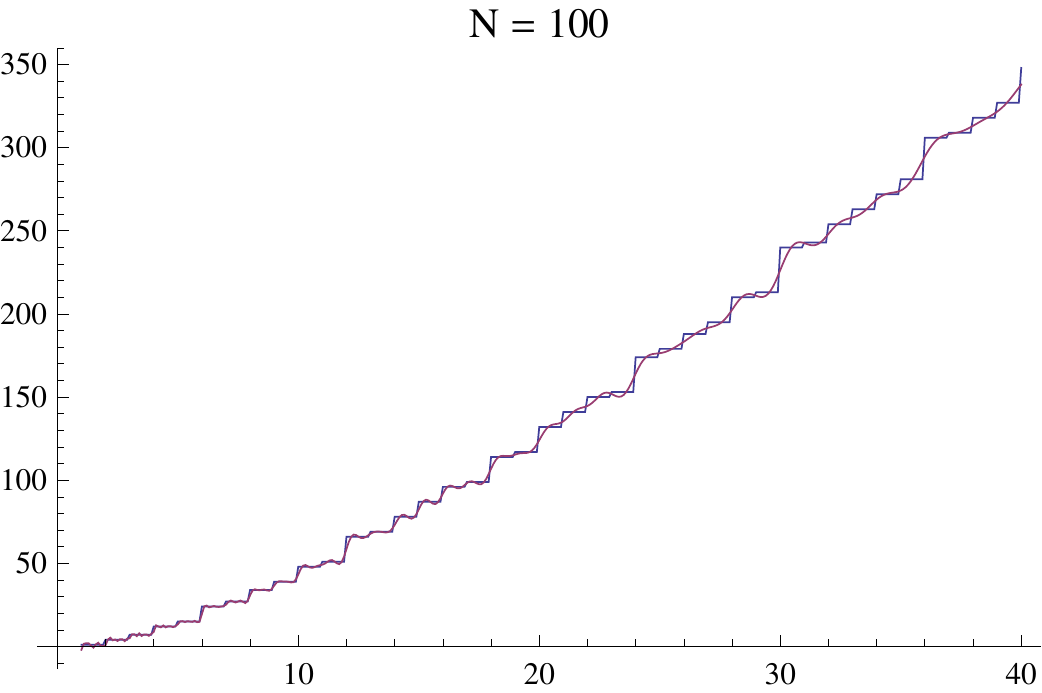}}
    \mbox{\includegraphics[width=\picDblWidth]{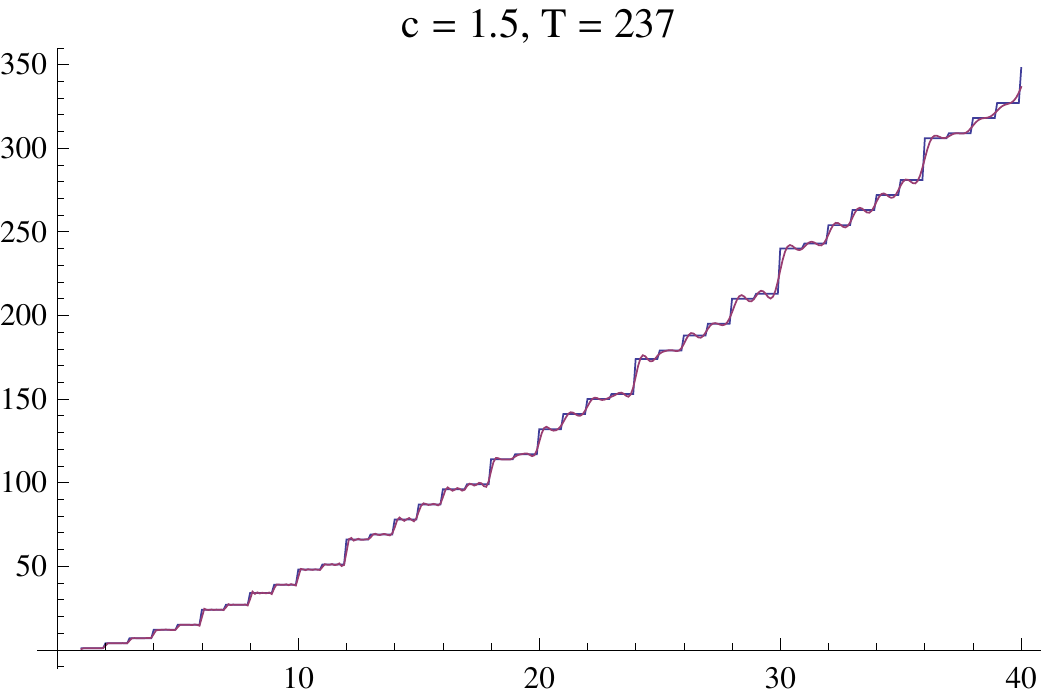}}
  }
  \caption{Approximating the sum of $\tau(n^2)$ using a sum and an integral}
  \label{fig:TauSum1Approximations}
\end{figure*}

The derivation of the expression for the residue at a third-order pole can be quite tedious.  Fortunately, \textit{Mathematica} can easily compute this residue for us:
\begin{verbatim}
    res = Residue[Zeta[s]^3/Zeta[2 s] * x^s/s, {s, 1}]
    resTerms = Collect[res, {Log[x], Log[x]^2}, Simplify].
\end{verbatim}
Then, \verb+resTerms[[3]]+, \verb+resTerms[[2]]+, and \verb+resTerms[[1]]+ are the terms with $x$, $x \log(x)$, and $x \log(x)^2$, respectively.

Figure ~\ref{fig:TauSum1Approx0} shows the approximation we get from setting $N = M = 0$ in equation \eqref{E:TauSum1Approx}.  Figure ~\ref{fig:TauSum1Approximations} shows the approximations using 100 pairs of zeta zeros in \eqref{E:TauSum1Approx}, and the integral in \eqref{E:TauSum1IntegralApprox}.

\subsection{The Sum of $\tau(n)^2$}\label{S:TauSquared}

This Dirichlet series holds for $s > 1$ \cite[set $a = b = 0$ in Theorem 305, p. 256]{HardyAndWright}
\[
\sum _{n=1}^{\infty} \frac{\tau (n)^2}{n^s} = \frac{\zeta(s)^4}{\zeta (2s)}.
\]
This Dirichlet series is similar to the previous one, except that here, the pole at $s = 1$ is of order 4.  The integral approximation to the summatory function $T_0(x)$ is
\begin{equation}\label{E:TauSum2IntegralApprox}
T_0(x) \simeq  \frac{1}{2 \pi i}
   \int_{c-i T}^{c+i T} \frac{\zeta(s)^4}{\zeta(2s)} \frac{x^s}{s} \, ds.
\end{equation}

The residue at $s = 0$ is $-1/8$.  The sum of the residues at $s = 0$ and $s = 1$, and at each $s = \rho/2$, suggests that this approximation might be valid:
\begin{align}
T_0(x) &\simeq
       a_1 + a_2 x + a_3 x\log (x) + a_4 x \log ^2(x) + a_5 x \log ^3(x) \label{E:TauSum2Approx}\\
     & + 2 \Re\left(\sum _{k=1}^N x^{\frac{\rho _k}{2}} \frac{
   \zeta \left(\frac{\rho _k}{2}\right){}^4}{\rho _k \zeta
   '\left(\rho _k\right)}\right)
   + \sum _{k=1}^{M } x^{-\frac{2 k}{2}} \frac{
    \zeta \left(-\frac{2 k}{2}\right)^4}{(-2 k) \zeta '(-2 k)}\notag
\end{align}
where
\[
a_1 = -1/8 = -.125,
\]
\begin{align*}
a_2 = \frac{6}{\pi^8}
& \Big(-4 \pi ^4 \Big((3-12 \gamma _1) \zeta
   '(2)+2 \zeta ^{(3)}(2)-3 \zeta ''(2)+18 \gamma ^2 \zeta '(2)
    -12 \gamma  (\zeta '(2)-\zeta ''(2))\Big) \\
   & \quad
    +\pi ^6 \Big(4 \gamma _1+\gamma  (4-12 \gamma _1)+2 \gamma _2+4
   \gamma ^3-6 \gamma ^2-1\Big) \\
   & \quad
   -1728 \zeta '(2)^3+144 \pi ^2 \zeta
   '(2) (2 \zeta ''(2)+4 \gamma  \zeta '(2)-\zeta
   '(2))\Big) \simeq 0.46032
\end{align*}
where $\gamma_2 \simeq -0.00969$ is the second Stieltjes constant (see equation \eqref{ZetaAndStieltjesConstants}),
\[
a_3=\frac{6 \left(\pi ^4 \left(-4 \gamma _1+6 \gamma ^2-4 \gamma +1\right)+144 \zeta '(2)^2-12 \pi ^2 \left(\zeta ''(2)+4
   \gamma  \zeta '(2)-\zeta '(2)\right)\right)}{\pi ^6}\simeq 0.82327,
\]
\[
a_4=\frac{3 (4 \gamma -1) \pi ^2-36 \zeta '(2)}{\pi ^4}\simeq 0.74434,
\]
and
\[
a_5=\frac{1}{\pi ^2}\simeq 0.10132.
\]

The residue at $s = 1$ accounts for the expressions involving $a_2$, $a_3$, $a_4$, and $a_5$.  The very complicated residue at $s = 1$ can be computed with the \textit{Mathematica} commands
\begin{verbatim}
    res = Residue[Zeta[s]^4/Zeta[2 s] * x^s/s, {s, 1}]
    resTerms = Collect[res, {Log[x], Log[x]^2, Log[x]^3}, Simplify].
\end{verbatim}

\begin{figure*}[ht]
  \mbox{\includegraphics[width=\picDblWidth]{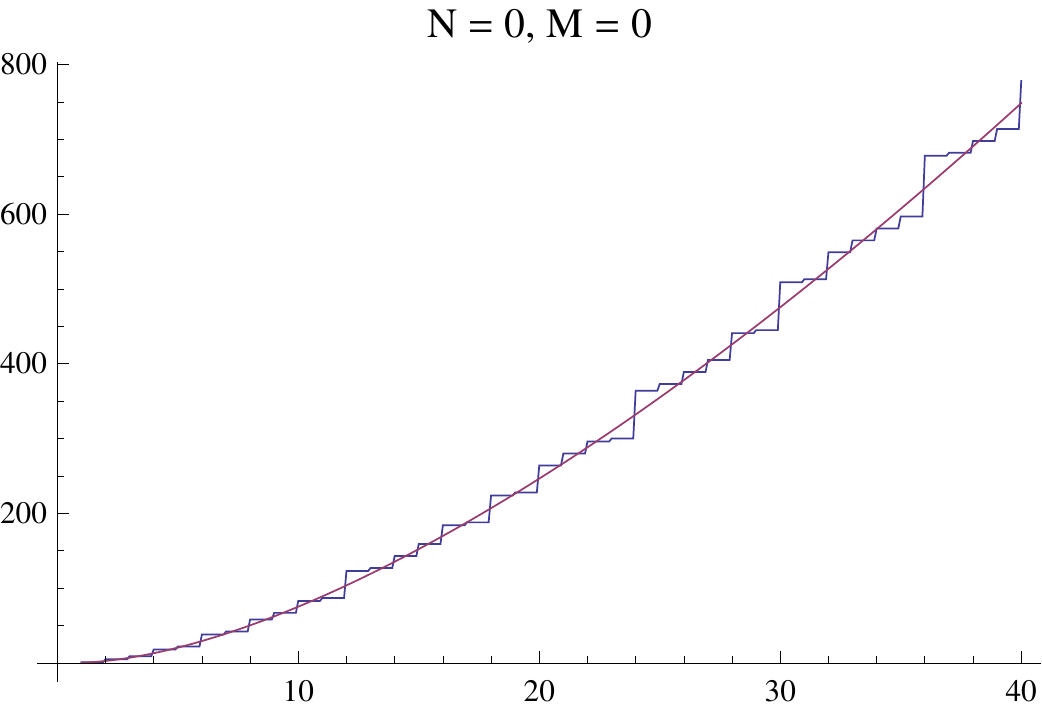}}
  \caption{Approximating the sum of $\tau(n)^2$ with the first five terms of \eqref{E:TauSum2Approx}}
  \label{fig:TauSum2Approx0}
\end{figure*}

\begin{figure*}[ht]
  \centerline{
    \mbox{\includegraphics[width=\picDblWidth]{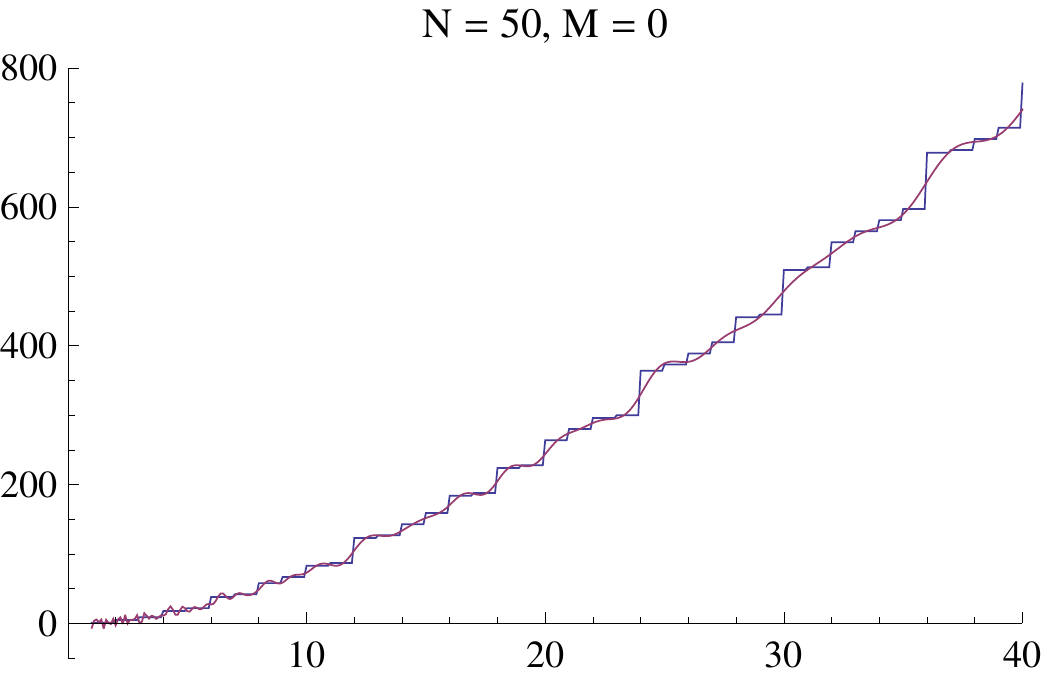}}
    \mbox{\includegraphics[width=\picDblWidth]{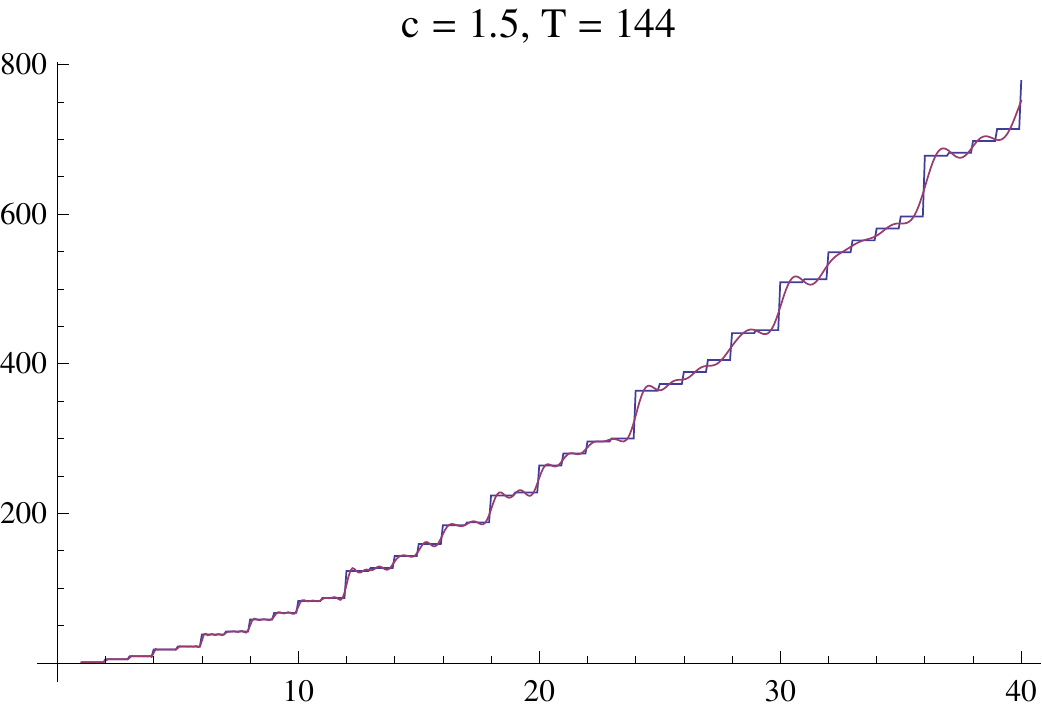}}
  }
  \caption{Approximating the sum of $\tau(n)^2$ using a sum and an integral}
  \label{fig:TauSum2Approximations}
\end{figure*}


\begin{figure*}[ht]
  \centerline{
    \mbox{\includegraphics[width=\picDblWidth]{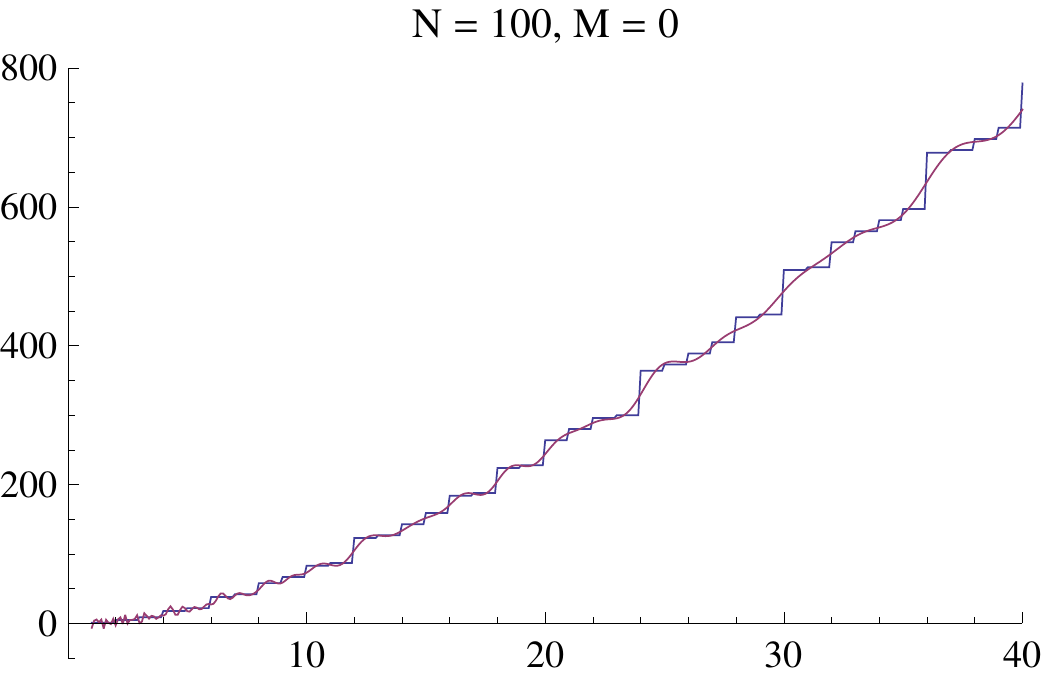}}
    \mbox{\includegraphics[width=\picDblWidth]{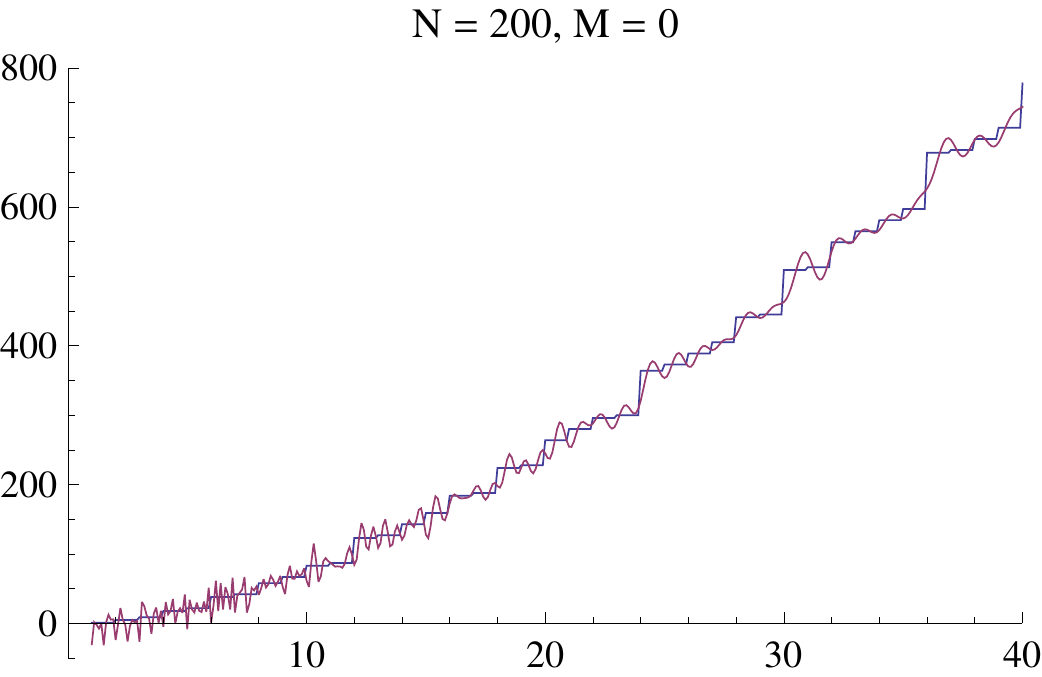}}
  }
  \caption{Approximating the sum of $\tau(n)^2$ using $N = 100$ and $N = 200$ pairs of zeros}
  \label{fig:TauSum2ApproximationsBigN}
\end{figure*}

Figure ~\ref{fig:TauSum2Approximations} shows approximations based on the sum \eqref{E:TauSum2Approx} and the integral \eqref{E:TauSum2IntegralApprox} approximations.  Figure ~\ref{fig:TauSum2ApproximationsBigN} shows the sum approximation with $N = 100$ and $N = 200$.  Notice that this latter approximation is quite wavy, much like the one for $\sigma(n)^2$ in figure ~\ref{fig:SigmaSum2ApproximationN200}.  Figure ~\ref{fig:TauSum2ApproximationsBigN} makes one wonder whether, for any given $x$, the sum even converges as $N$ approaches infinity.

In 1916, Ramanujan \cite{Ramanujan} stated that
\[
\sum _{n \leq x} \tau (n)^2 = a_2 x + a_3 x\log (x) + a_4 x \log ^2(x) + a_5 x \log ^3(x) + O(x^{\frac{3}{5}+\epsilon})
\]
and gave the same values of $a_4$ and $a_5$ stated above; presumably, he also knew the values of $a_2$ and $a_3$.
A proof of this formula appeared in \cite{Wilson}; see also \cite{Suryanarayana}.  The error term can be reduced to $O(x^{\frac{1}{2}+\epsilon})$; see Ivi\'{c}, \cite[Equation 14.30, p. 394]{Ivic}.  Therefore, if the right side of \eqref{E:TauSum2Approx} does, indeed, represent the summatory function, then this error bound is, in effect, a bound on a sum over zeta zeros.  That is, we would have
\[
 \Re\left(\sum _{k=1}^{\infty} x^{\frac{\rho _k}{2}} \frac{
   \zeta \left(\frac{\rho _k}{2}\right){}^4}{\rho _k \zeta
   '\left(\rho _k\right)}\right) = O(x^{\frac{1}{2}+\epsilon}).
\]

\ifthenelse {\boolean{BKMRK}}
  { \section{Tallying \texorpdfstring{$\sigma(n) \tau(n)$}{sigma Times tau}}\label{S:SigmaTau} }
  { \section{Tallying $\sigma(n) \tau(n)$}\label{S:SigmaTau} }

This Dirichlet series for $\sigma(n) \tau(n)$ holds for $s > 2$ \cite[set $a = 1$, $b = 0$ in eq. D-58]{Gould}:

\[
\sum _{n=1}^{\infty} \frac{\sigma(n) \tau(n)}{n^s} = \frac{\zeta(s)^2 \zeta(s - 1)^2}{\zeta(2s - 1)}.
\]

The integral approximation to the summatory function $T_0(x)$ is
\begin{equation}\label{E:SigmaTauIntegralApprox}
T_0(x) \simeq \frac{1}{2 \pi i}
   \int_{c-i T}^{c+i T} \frac{\zeta(s)^2 \zeta(s - 1)^2}{\zeta(2s - 1)} \frac{x^s}{s} \, ds.
\end{equation}

The integrand has a pole of order 1 at $s = 0$, and poles of order 2 at $s = 1$ and at $s = 2$.  There is also a pole at each value of $s$ where $\rho = 2s - 1$, where $\rho$ runs over all zeta zeros.  Therefore, these poles occur at each $s = (\rho + 1)/2$.

The sum of the residues, which we hope will be close to $T_0(x)$, is
\begin{align}
T_0(x) &\simeq a_1 + a_2 x + a_3 x^2 + a_4 x^2 \log (x)  \label{E:SigmaTauSumApprox}\\
     & + 2 \Re\left(\sum _{k=1}^N x^{\frac{\rho _k+1}{2}} \frac{
   \zeta \left(\frac{\rho _k+1}{2}\right)^2
   \zeta \left(\frac{\rho_k-1}{2}\right)^2}
       {2 (\frac{\rho_k+1}{2})
           \zeta '\left(\rho_k\right)}\right)\notag \\
   &+ \sum _{k=1}^{M } x^{\frac{-2 k + 1}{2}} \frac{ \zeta \left(\frac{-2 k + 1}{2}\right)^2
     \zeta \left(\frac{-2 k-1}{2}\right)^2}
     {2 (\frac{-2 k + 1}{2}) \zeta '(-2 k)}.\notag
\end{align}
The apparently spurious 2 in the denominators appears for the same reason as the 2 in \eqref{E:SigmaSum1Approx}.  The second sum is too small to affect our graphs, so we will set $M = 0$ when we graph \eqref{E:SigmaTauSumApprox}.  Here,
\[
a_1=-\frac{1}{48}\simeq 0.020833,
\]
\[
a_2=\frac{1}{2}=0.5,
\]
\[
a_3=\frac{\pi ^2 \left(24 \zeta (3) \zeta '(2)+\pi ^2 \left(4 \gamma
    \zeta (3)-\zeta (3)-4 \zeta '(3)\right)\right)}{144 \zeta (3)^2} \simeq -0.17540,
\]
and
\[
a_4=\frac{\pi ^4}{72 \zeta (3)}\simeq 1.12549.
\]

$a_1$ and $a_2 x$ are the residues at $s = 0$ and $s = 1$, respectively.  $ a_3 x^2 + a_4 x^2 \log (x)$ is the residue at $s = 2$.

\begin{figure*}[ht]
  \mbox{\includegraphics[width=\picDblWidth]{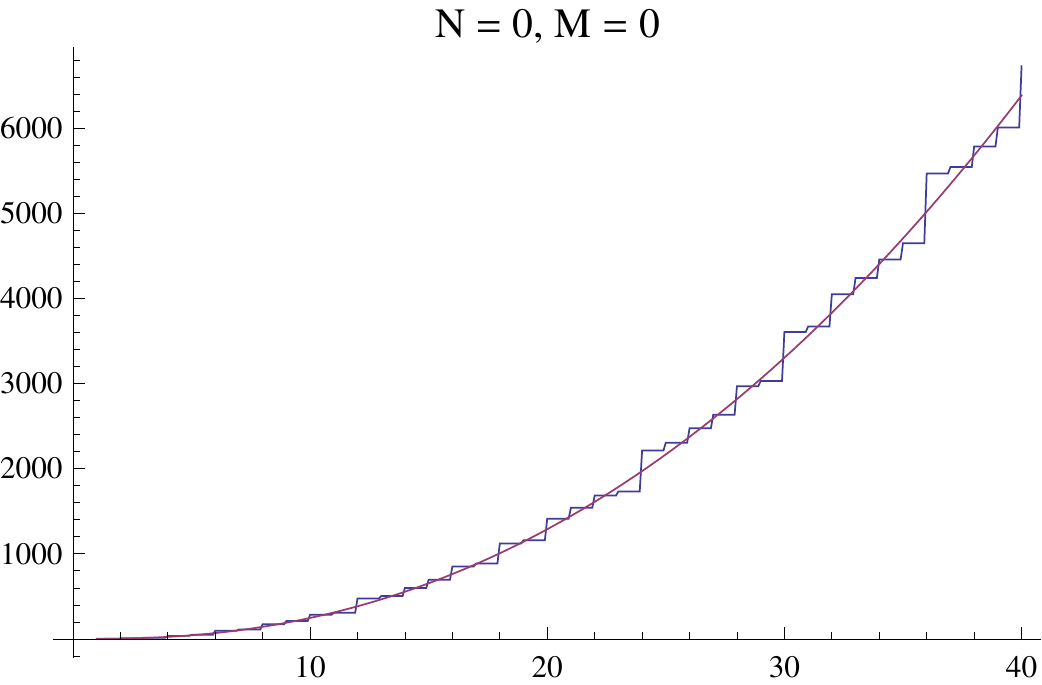}}
  \caption{Approximating the sum of $ \sigma(n) \tau(n)$ with the first four terms of \eqref{E:SigmaTauSumApprox}}
  \label{fig:SigmaTauApprox0}
\end{figure*}

\begin{figure*}[ht]
  \centerline{
    \mbox{\includegraphics[width=\picDblWidth]{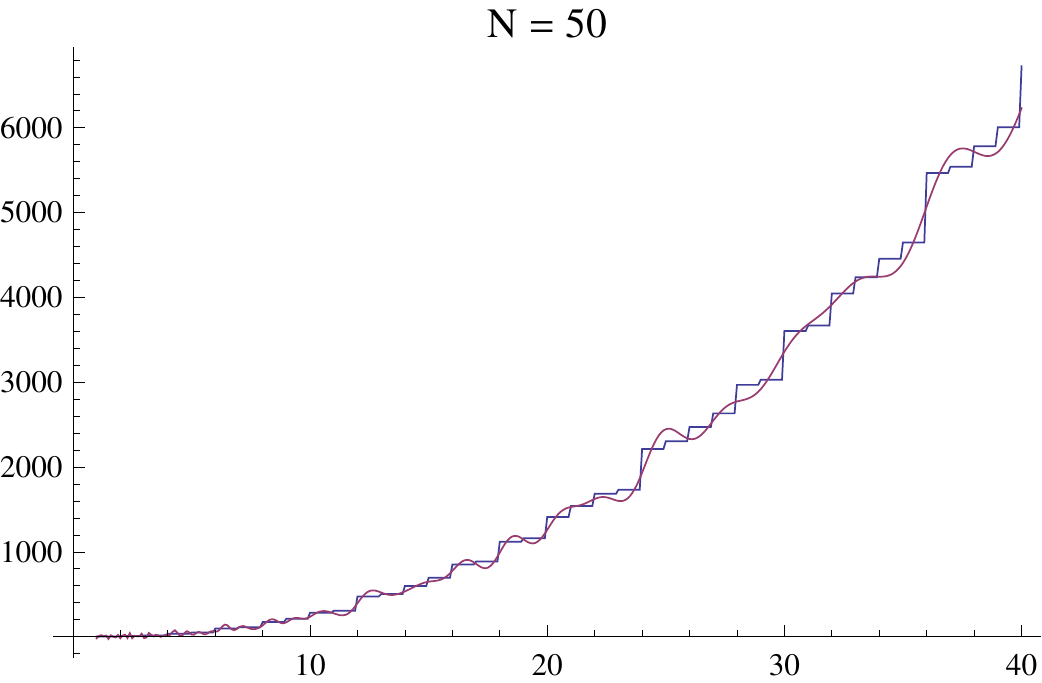}}
    \mbox{\includegraphics[width=\picDblWidth]{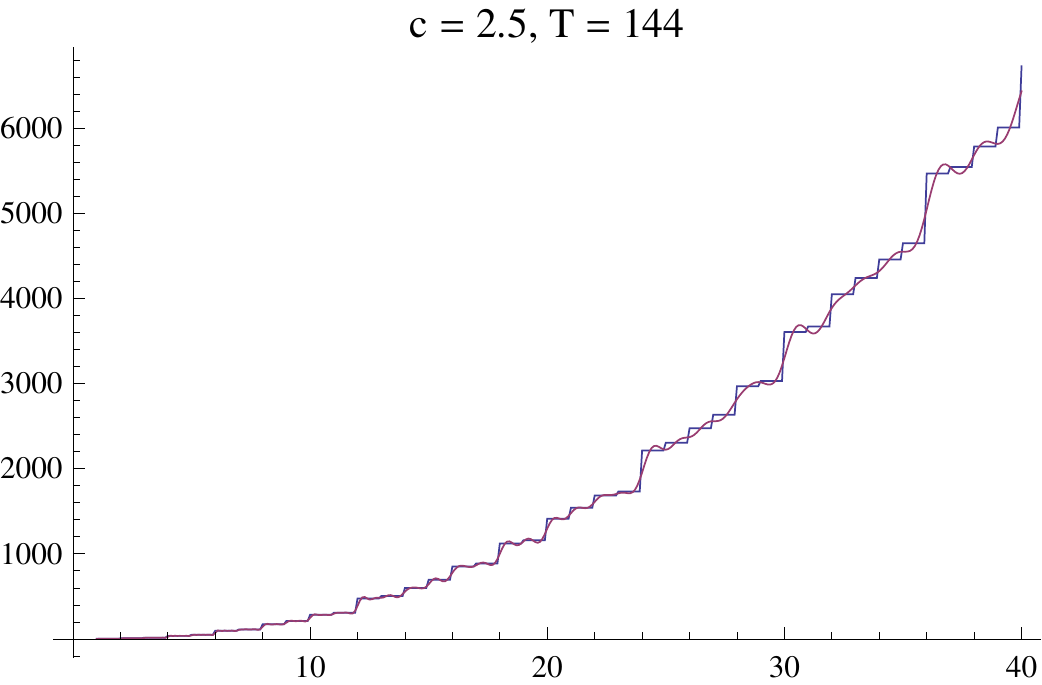}}
  }
  \caption{Approximating the sum of $\sigma(n) \tau(n)$ using a sum and an integral}
  \label{fig:SigmaTauSumGraphs}
\end{figure*}

Figure ~\ref{fig:SigmaTauApprox0} shows the graph of the sum approximation \eqref{E:SigmaTauSumApprox} without using any zeta zeros (i. e., with $N = M = 0$).  Figure ~\ref{fig:SigmaTauSumGraphs} shows the graphs of the sum approximation with 50 pairs of complex zeta zeros, and the integral approximation \eqref{E:SigmaTauIntegralApprox}.

\ifthenelse {\boolean{BKMRK}}
  { \section{Tallying \texorpdfstring{$\lambda(n) \tau(n)$}{lambda Times tau}}\label{S:lambdaTimesTau} }
  { \section{Tallying $\lambda(n) \tau(n)$}\label{S:lambdaTimesTau} }

This Dirichlet series for $\lambda(n) \tau(n)$ holds for $s > 1$ \cite[set $k = 0$ in eq. D-47]{Gould}, \cite[p. 225]{McCarthy}:

\[
\sum _{n=1}^{\infty} \frac{\lambda (n) \sigma (n)}{n^s} = \frac{\zeta (2 s)^2}{\zeta (s)^2}.
\]

The integral approximation to the summatory function $T_0(x)$ is
\begin{equation}\label{E:LambdaTauIntegralApprox}
T_0(x) \simeq \frac{1}{2 \pi i}
   \int_{c-i T}^{c+i T} \frac{\zeta (2 s)^2}{\zeta (s)^2} \frac{x^s}{s} \, ds.
\end{equation}

This integrand has a pole of order 1 at $s = 0$, where the residue is 1.  There is a pole of order 2 where $2s = 1$, that is, at $s = 1/2$.  There is also a pole of order 2 at each $s = \rho$, where $\rho$ runs over all zeta zeros.  Residues at poles of order 2 in the numerator were discussed in Section \ref{S:TallyingTheSquarefreeDivisors}.  There, we computed the residue by taking a limit, using equation \eqref{E:stdResidueFormulaN}.  Here, the residue at $s = 1/2$ can be derived in the same way.  It is
\[
\frac{\sqrt{x} \left(-\zeta \left(\frac{1}{2}\right)+2 \gamma  \zeta
   \left(\frac{1}{2}\right)-\zeta
   '\left(\frac{1}{2}\right)\right)}{\zeta
   \left(\frac{1}{2}\right)^3}+\frac{\sqrt{x} \log (x)}{2 \zeta
   \left(\frac{1}{2}\right)^2}.
\]

Here, for the first time, we encounter poles of order 2 in the denominator of the integrand.  To compute the residues at $s = \rho$, we will use equation \eqref{E:stdResidueFormula2}.  Set $B(s) = \zeta(s)^2$, and set $A(s)$ equal to everything else in the integrand,
\[
A(s) = \zeta (2 s)^2 \frac{x^s}{s}
\]
so that the integrand is equal to $A(s)/B(s)$.  Then we substitute this $A(s)$ and $B(s)$ into equation \eqref{E:stdResidueFormula2}.

The result is a very complicated expression which we shall not write in full here.  The denominator of this expression turns out to be
\begin{equation}\label{E:LambdaTauDenom0}
3 s^2 \left(\zeta (s) \zeta ''(s)+\zeta '(s)^2\right)^2.
\end{equation}
Remember that we want an expression for the residue at $s = \rho$, a zero of zeta.  But notice that when we substitute $s = \rho$, every occurrence of $\zeta(s)$ in \eqref{E:LambdaTauDenom0} becomes 0.  The denominator then simplifies to
\begin{equation}\label{E:LambdaTauDenom}
3 \rho^2 \zeta '(\rho)^4.
\end{equation}

We process the numerator in the same way.  The numerator that results from substituting $A(s)$ and $B(s)$ into \eqref{E:stdResidueFormula2} is very complicated, but can be separated into two terms.  The first term has a factor of $x^s \log(x)$:
\[
3 s x^s \log (x) \zeta (2 s)^2 \left(\zeta (s) \zeta ''(s)+\zeta '(s)^2\right).
\]
Again, because $\zeta(s) = 0$ when we substitute $s = \rho$, this simplifies to
\begin{equation}\label{E:LambdaTauNumer1}
3 \rho x^{\rho} \log (x) \zeta (2 \rho)^2 \zeta '(\rho)^2.
\end{equation}

The second term in the numerator has a factor of $x^\rho$ but without the $\log(x)$:
\begin{align}
-x^s \zeta (2 s) &(3 \zeta '(s)^2 (\zeta (2 s)-4 s \zeta '(2 s))+3 s
   \zeta (2 s) \zeta '(s) \zeta ''(s)\notag \\
&   +\zeta (s) (s \zeta (2 s) \zeta ^{(3)}(s)+3
   \zeta ''(s) (\zeta (2 s)-4 s \zeta '(2 s)))).\notag
\end{align}

Because $\zeta(s)$ will be 0 when we substitute $s = \rho$, this simplifies to
\begin{equation}\label{E:LambdaTauNumer2}
-x^{\rho} \zeta (2 \rho) \left(3 \zeta '(\rho)^2 \left(\zeta (2 \rho)-4 s \zeta '(2 \rho)\right)+3 \rho
   \zeta (2 \rho) \zeta '(\rho) \zeta ''(\rho)\right).
\end{equation}

For the residue at any given zeta zero $\rho$, here's what we end up with.  We get a quotient, the denominator of which is \eqref{E:LambdaTauDenom}.  The numerator has two terms, one of which, \eqref{E:LambdaTauNumer1}, has a factor of $\log(x)$, the other of which, \eqref{E:LambdaTauNumer2}, does not.  We can write the quotient for this zeta zero as
\begin{equation}\label{E:LambdaTauQuot}
\frac{F_1(\rho ) x^{\rho }+F_2(\rho ) x^{\rho } \log (x)}{G(\rho )}.
\end{equation}

There is one of these expressions for each zeta zero.  We sum \eqref{E:LambdaTauQuot} over all zeta zeros.  When we combine this sum with the residues at $s = 0$ and $s = 1/2$, we get the following expression, which may approximate the summatory function, $T_0(x)$:
\begin{equation}\label{E:LambdaTauSumApprox}
T_0(x)\simeq a_1+a_2 \sqrt{x}+a_3 \sqrt{x} \log (x)+2 \Re\left(\sum _{k=1}^N x^{\rho _k}
   \frac{F_1\left(\rho _k\right)+F_2\left(\rho _k\right) \log (x)}{G\left(\rho
   _k\right)}\right),
\end{equation}

where
\[
a_1 = 1,
\]
\[
a_2=\frac{2 \gamma  \zeta \left(\frac{1}{2}\right)-\zeta \left(\frac{1}{2}\right)-\zeta
   '\left(\frac{1}{2}\right)}{\zeta \left(\frac{1}{2}\right)^3}\simeq -1.187104,
\]
and
\[
a_3=\frac{1}{2 \zeta \left(\frac{1}{2}\right)^2}\simeq 0.234452.
\]

The functions $F_1$, $F_2$, and $G$ are
\[
F_1(s)=-\zeta (2 s) \left(s \zeta (2 s) \zeta ''(s)+\zeta '(s) \left(\zeta (2 s)-4 s
   \zeta '(2 s)\right)\right),
\]

\[
F_2(s)=s \zeta (2 s)^2 \zeta '(s),
\]

and

\[
G(s)=s^2 \zeta '(s)^3.
\]

The following \textit{Mathematica} code will perform the above calculations.  This uses the \verb+residueFormula+ function from Section \ref{S:RulesForResidues}:
\begin{verbatim}
A[s] = Zeta[2 s]^2 * x^s/s;
B[s] = Zeta[s]^2;
expr = residueFormula[2, s];
expr2 = Together[expr /. Zeta[s] -> 0];
numer = Collect[Numerator[expr2], Log[x], Simplify];
numer[[2]]/x^s           (* this is F1[s] *)
numer[[1]]/(x^s Log[x])  (* this is F2[s] *)
Denominator[expr2]       (* this is G[s] *)
\end{verbatim}

Note that \eqref{E:LambdaTauSumApprox} has no sum over the real zeta zeros.  This is because, if $s = \rho$ is one of the real zeros -2, -4, ..., then the presence of $2s$ in both $F_1$ and $F_2$ guarantees that every term will be 0, because if $s$ is one of these zeros, then so is $2s$.

\begin{figure*}[ht]
  \mbox{\includegraphics[width=\picDblWidth]{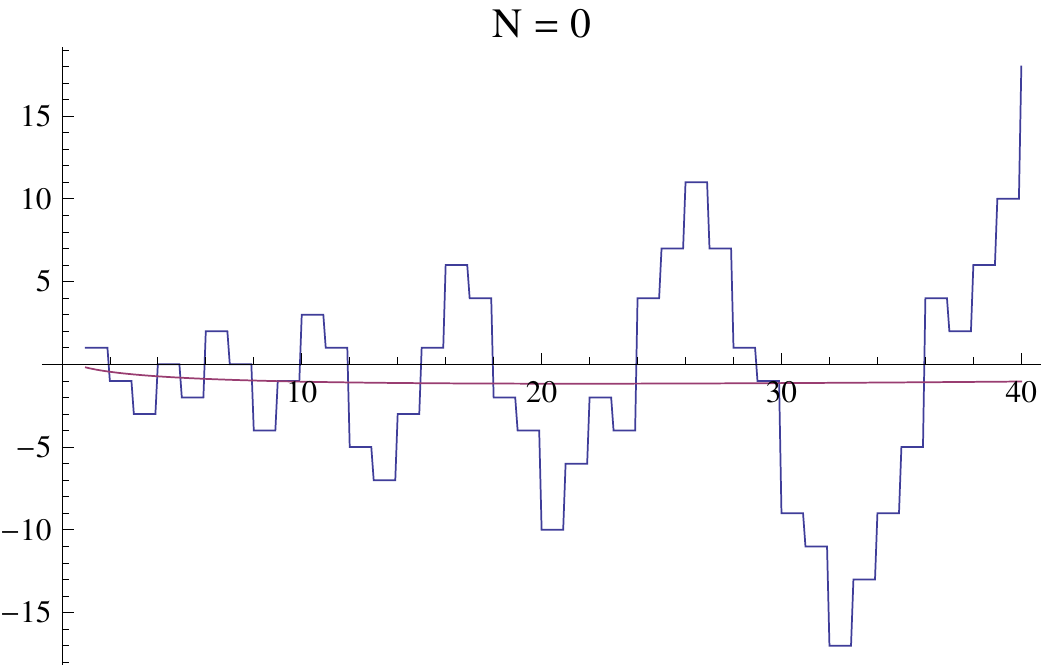}}
  \caption{Approximating the sum of $ \lambda(n) \tau(n)$ with the first three terms of \eqref{E:LambdaTauSumApprox}}
  \label{fig:LambdaTauApprox0}
\end{figure*}

\begin{figure*}[ht]
  \centerline{
    \mbox{\includegraphics[width=\picDblWidth]{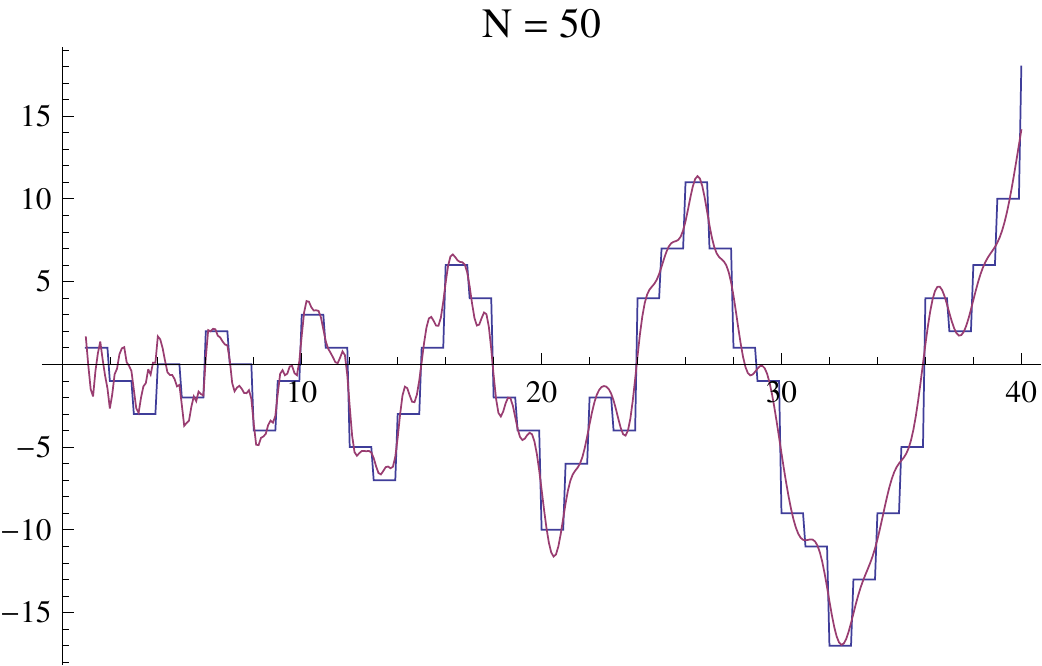}}
    \mbox{\includegraphics[width=\picDblWidth]{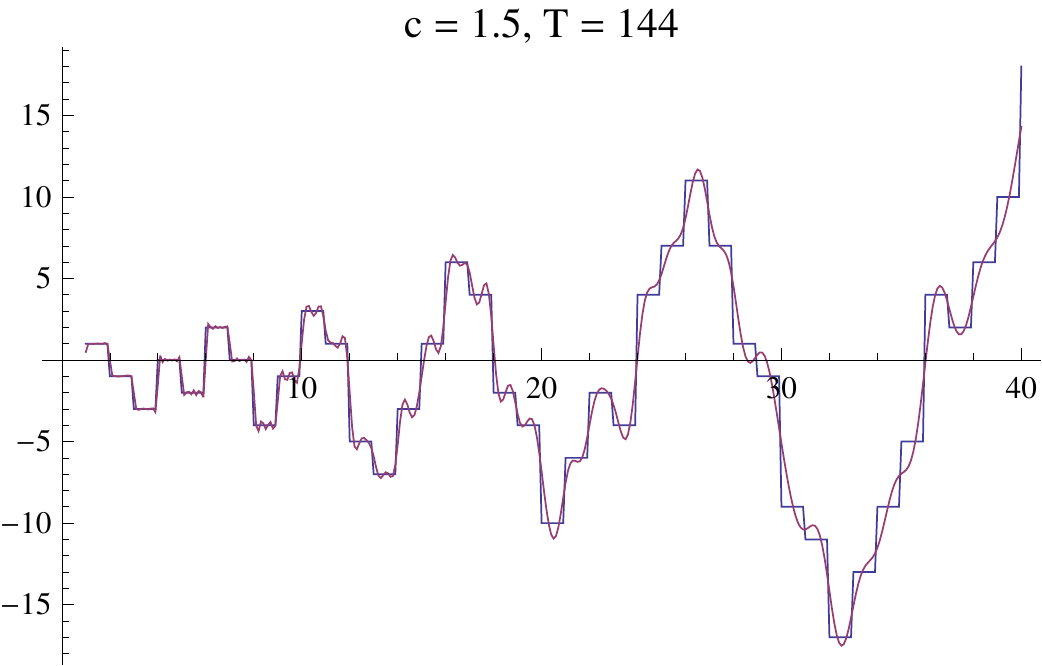}}
  }
  \caption{Approximating the sum of $\lambda(n) \tau(n)$ using a sum and an integral}
  \label{fig:LambdaTauSumGraphs}
\end{figure*}

Figure ~\ref{fig:LambdaTauApprox0} shows the graphs of the sum approximation \eqref{E:LambdaTauSumApprox} without using any zeta zeros.  Because $a_3$ is positive, the expression consisting of the first three terms of \eqref{E:LambdaTauSumApprox} will eventually be positive.  In fact, $a_1+a_2 \sqrt{x}+a_3 \sqrt{x} \log (x) > 0$ if $x \geq \simeq 104.1$.  This suggest an overall tendency for the summatory function $T(x)$ to be positive, but otherwise, the first three terms bear no resemblance to the actual graph of $T(x)$.  Figure ~\ref{fig:LambdaTauSumGraphs} shows the graphs of the sum approximation using 50 pairs of complex zeta zeros, along with the integral approximation \eqref{E:LambdaTauIntegralApprox}.

\ifthenelse {\boolean{BKMRK}}
  { \section{Tallying \texorpdfstring{$\lambda(n) 2^{\nu(n)}$}{lambda Times 2 to the nu} }\label{S:LambdaTwoNu} }
  { \section{Tallying $\lambda(n) 2^{\nu(n)}$}\label{S:LambdaTwoNu} }


Note: see Section \ref{S:TallyingTheSquarefreeDivisors} (``Tallying the Squarefree Divisors'') for a discussion of the sum of $2^{\nu(n)}$.

This Dirichlet series for $\lambda(n) 2^{\nu(n)}$ holds for $s > 1$ \cite[eq. D-21]{Gould}, \cite[p. 227]{McCarthy}:

\[
\sum _{n=1}^{\infty} \frac{\lambda(n) 2^{\nu(n)}}{n^s} = \frac{\zeta (2 s)}{\zeta (s)^2}.
\]

The integral approximation to the summatory function $T_0(x)$ is
\begin{equation}\label{E:LambdaTwoNuIntegralApprox}
T_0(x) \simeq \frac{1}{2 \pi i}
   \int_{c-i T}^{c+i T} \frac{\zeta (2 s)}{\zeta (s)^2} \frac{x^s}{s} \, ds.
\end{equation}

The sum over residues that may approximate the summatory function is
\begin{equation}\label{E:LambdaTwoNuSumApprox}
T_0(x) \simeq a_1 + a_2 \sqrt{x}
 + 2 \Re\left(
       \sum _{k=1}^N x^{\rho _k} \frac{F_1(\rho_k) + F_2(\rho_k) \log(x)}{G(\rho_k)}
     \right),
\end{equation}
where
\[
a_1 = -2
\]
and
\[
a_2 = \frac{1}{\zeta \left(\frac{1}{2}\right)^2} \simeq 0.46890.
\]

The functions $F_1$, $F_2$, and $G$ are:
\[
F_1(s) = -s \zeta (2 s) \zeta ''(s)-\zeta '(s) (\zeta (2 s)-2 s \zeta'(2 s)),
\]
\[
F_2(s) = s \zeta (2 s) \zeta '(s),
\]
and
\[
G(s) = s^2 \zeta '(s)^3.
\]

The following \textit{Mathematica} code will perform these calculations:
\begin{verbatim}
A[s] = Zeta[2 s] * x^s/s;
B[s] = Zeta[s]^2;
expr = residueFormula[2, s];
expr2 = Together[expr /. Zeta[s] -> 0];
numer = Collect[Numerator[expr2], Log[x], Simplify];
numer[[2]]/x^s           (* this is F1[s] *)
numer[[1]]/(x^s Log[x])  (* this is F2[s] *)
Denominator[expr2]       (* this is G[s] *)
\end{verbatim}

\begin{figure*}[ht]
  \centerline{
    \mbox{\includegraphics[width=\picDblWidth]{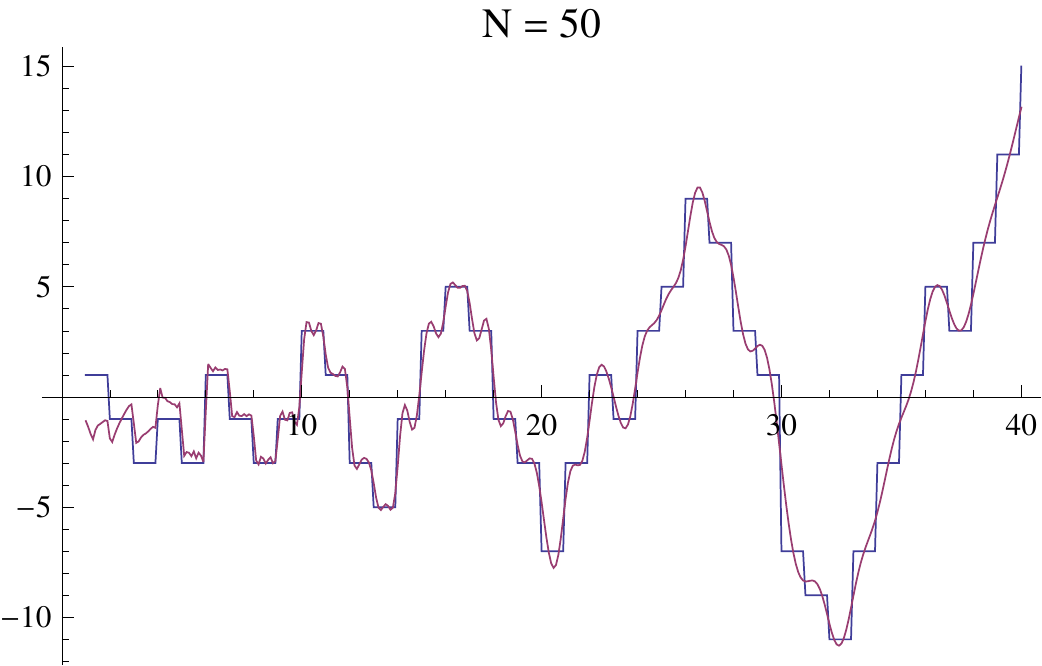}}
    \mbox{\includegraphics[width=\picDblWidth]{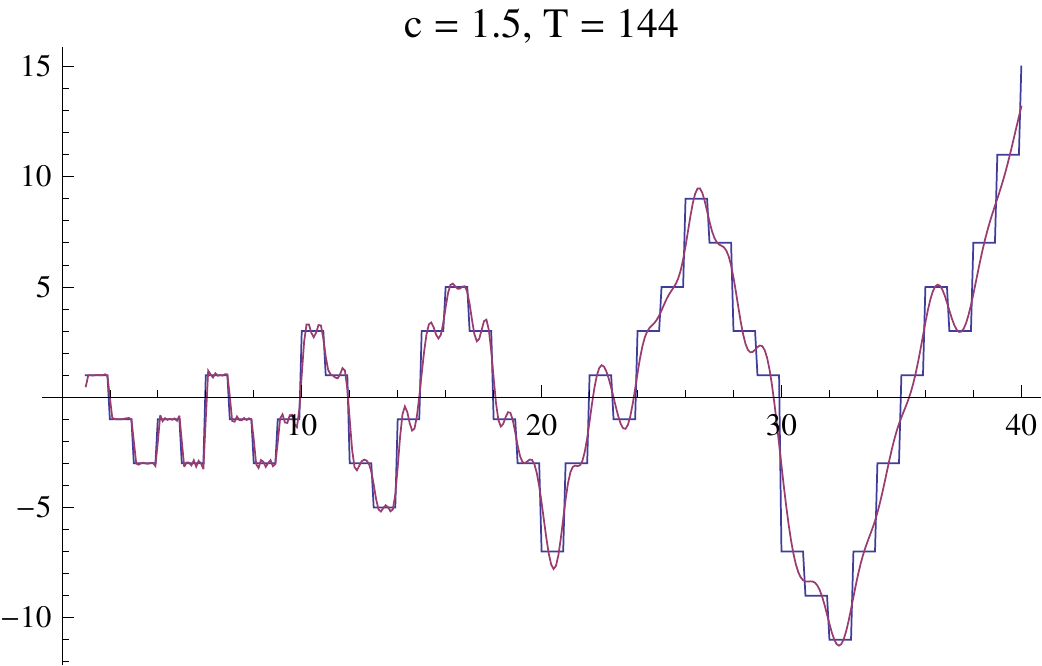}}
  }
  \caption{Approximating the sum of $\lambda(n) 2^{\nu(n)}$ using a sum and an integral}
  \label{fig:LambdaTwoNuGraphs}
\end{figure*}

In Figure ~\ref{fig:LambdaTwoNuGraphs}, notice that, for small $x$, the sum is not very close to $T(x)$.

Note: For this series, the corresponding sum over real zeros appears to diverge.  Here is the sum of the first five terms:
\[
 - \frac{\zeta '(-4)}{x^2 \zeta '(-2)^2}
 - \frac{\zeta '(-8)}{2 x^4 \zeta '(-4)^2}
 - \frac{\zeta '(-12)}{3 x^6 \zeta '(-6)^2}
 - \frac{\zeta '(-16)}{4 x^8 \zeta '(-8)^2}
 - \frac{\zeta '(-20)}{5 x^{10} \zeta '(-10)^2}
\]
\[
= - \frac{8.61152}{x^2} - \frac{65.2338}{x^4} - \frac{605.915}{x^6} - \frac{6409.28}{x^8} - \frac{73829.4}{x^{10}} .
\]

In the previous examples, we have extended the rectangle around which we integrate far to the left, to enclose the real zeros of zeta.  However, it is not necessary to do this.  Aside from the discrepancy at small $x$, it seems sufficient to make the rectangle enclose only the complex zeros and the poles at $s = 0$ and $s = 1/2$.

\ifthenelse {\boolean{BKMRK}}
  { \section{Tallying \texorpdfstring{$\lambda(n) \tau(n^2)$}{lambda Times tau of Squares} }\label{S:LambdaTauOfSquares} }
  { \section{Tallying $\lambda(n) \tau(n^2)$}\label{S:LambdaTauOfSquares} }


This Dirichlet series for $\lambda(n) \tau(n^2)$ holds for $s > 1$ \cite[eq. D-49]{Gould}, \cite[p. 234]{McCarthy}:

\[
\sum _{n=1}^{\infty} \frac{\lambda (n) \tau(n^2)}{n^s} = \frac{\zeta (2 s)^2}{\zeta (s)^3}.
\]

The integral approximation to the summatory function $T_0(x)$ is
\begin{equation}\label{E:LambdaTauOfSquaresIntegralApprox}
T_0(x) \simeq \frac{1}{2 \pi i}
   \int_{c-i T}^{c+i T} \frac{\zeta (2 s)^2}{\zeta (s)^3} \frac{x^s}{s} \, ds.
\end{equation}

The sum over residues that may approximate the summatory function is
\begin{equation}\label{E:LambdaTauOfSquaresSumApprox}
T_0(x) \simeq a_1 + a_2 \sqrt{x} + a_3 \sqrt{x} \log (x) 
 + 2 \Re\left(\sum _{k=1}^N x^{\rho _k}
   \frac{F_1(\rho_k) + F_2(\rho_k) \log(x) + F_3(\rho_k) \log(x)^2}{G(\rho_k)}\right),
\end{equation}
where
\[
a_1 = -2,
\]
\[
a_2=\frac{4 \gamma  \zeta \left(\frac{1}{2}\right)-2 \zeta
   \left(\frac{1}{2}\right)-3 \zeta '\left(\frac{1}{2}\right)}{2
   \zeta \left(\frac{1}{2}\right)^4}\simeq 1.24413,
\]
and
\[
a_3=\frac{1}{2 \zeta \left(\frac{1}{2}\right)^3 }\simeq -0.160544.
\]

The functions $F_1$, $F_2$, $F_3$, and $G$ are:
\begin{align*}
F_1(s)&=
3 s^2 \zeta (2 s)^2 \zeta ''(s)^2+2 \zeta '(s)^2 \left(\zeta (2 s)
   \left(4 s^2 \zeta ''(2 s)+\zeta (2 s)\right)+4 s^2 \zeta '(2
   s)^2-4 s \zeta (2 s) \zeta '(2 s)\right) \\
      & \quad
   -s \zeta (2 s) \zeta '(s)
   \left(s \zeta (2 s) \zeta ^{(3)}(s)-3 \zeta ''(s) \left(\zeta (2
   s)-4 s \zeta '(2 s)\right)\right),
\end{align*}

\[
F_2(s) =
-s \zeta (2 s) \zeta '(s) \left(3 s \zeta (2 s) \zeta ''(s)+2 \zeta
   '(s) \left(\zeta (2 s)-4 s \zeta '(2 s)\right)\right),
\]

\[
F_3(s) = s^2 \zeta (2 s)^2 \zeta '(s)^2,
\]
and
\[
G(s) = 2 s^3 \zeta '(s)^5.
\]

The following \textit{Mathematica} code will perform these calculations:

\begin{verbatim}
A[s] = Zeta[2 s]^2 * x^s/s;
B[s] = Zeta[s]^3;
expr = residueFormula[3, s];
expr2 = Together[expr /. Zeta[s] -> 0];
numer = Collect[Numerator[expr2], {Log[x], Log[x]^2}, Simplify];
numer[[3]]/x^s             (* this is F1[s] *)
numer[[2]]/(x^s Log[x])    (* this is F2[s] *)
numer[[1]]/(x^s Log[x]^2)  (* this is F3[s] *)
Denominator[expr2]         (* this is G[s] *)
\end{verbatim}

\begin{figure*}[ht]
  \centerline{
    \mbox{\includegraphics[width=\picDblWidth]{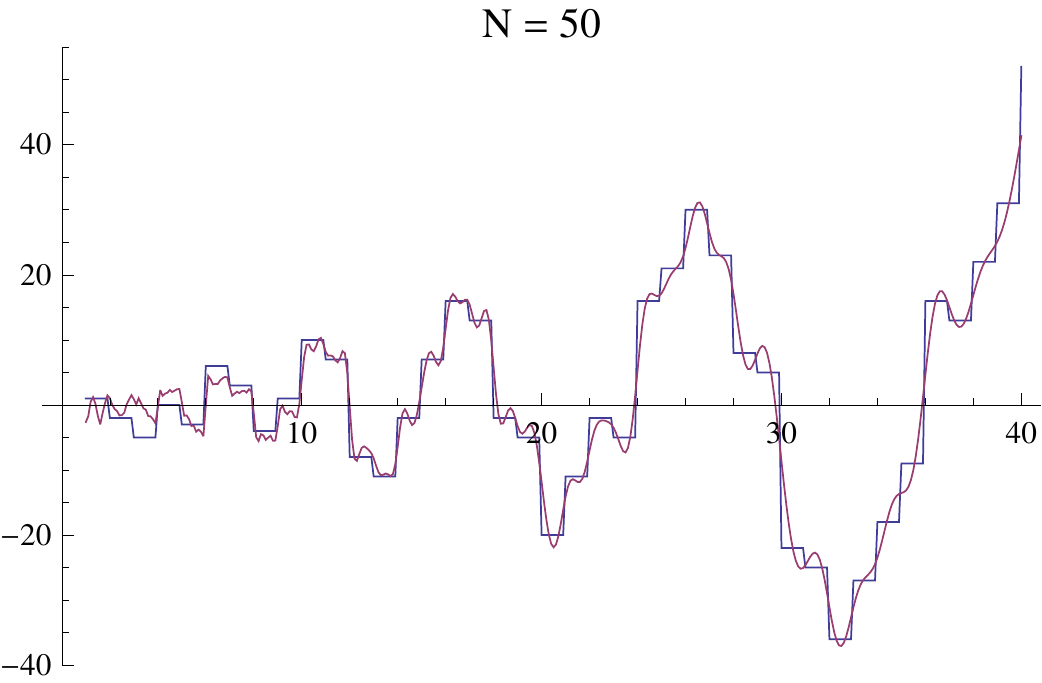}}
    \mbox{\includegraphics[width=\picDblWidth]{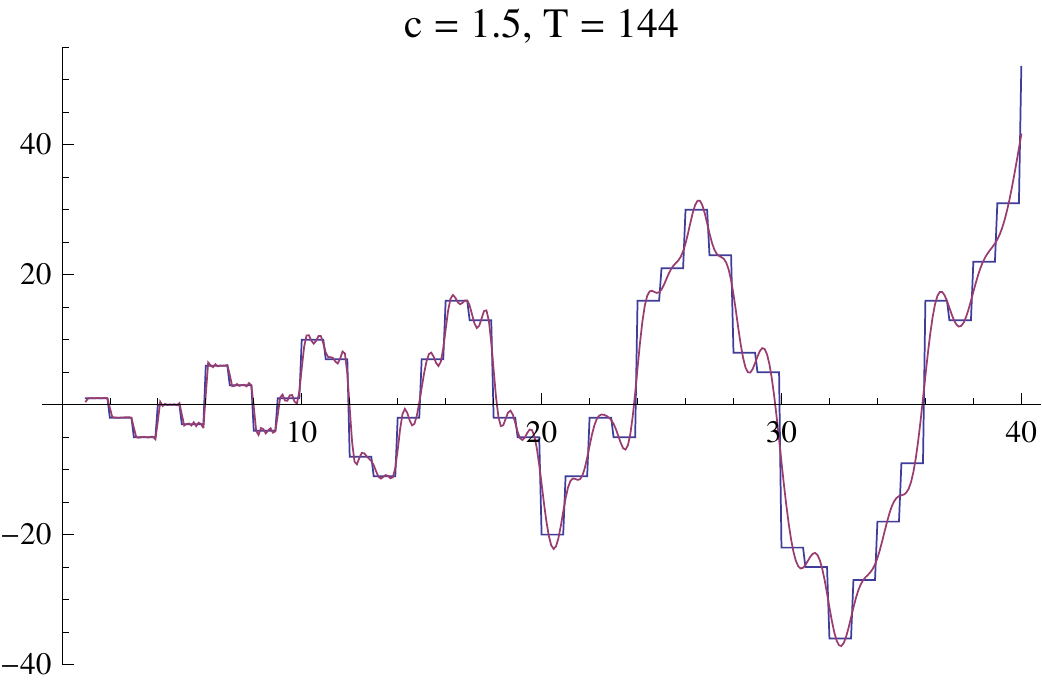}}
  }
  \caption{Approximating the sum of $\lambda(n) \tau(n^2)$ using a sum and an integral}
  \label{fig:LambdaTauOfSquaresGraphs}
\end{figure*}

Because $a_3$ in \eqref{E:LambdaTauOfSquaresSumApprox} is negative, the expression consisting of the first three terms of \eqref{E:LambdaTauOfSquaresSumApprox} will eventually be negative.  In fact, $a_1+a_2 \sqrt{x}+a_3 \sqrt{x} \log (x) < 0$ if $x \geq \simeq 1717.9$.  This suggest an overall tendency for the summatory function $T(x)$ to be negative.  Figure ~\ref{fig:LambdaTauOfSquaresGraphs} shows the graphs of the sum approximation using 50 pairs of complex zeta zeros, along with the integral approximation \eqref{E:LambdaTauOfSquaresIntegralApprox}.

Note: For this series, the corresponding sum over real zeros appears to diverge.  Here is the sum of the first five terms:
\[
 - \frac{2 \zeta '(-4)^2}{x^2 \zeta '(-2)^3}
 - \frac{\zeta'(-8)^2}{x^4 \zeta '(-4)^3}
 - \frac{2 \zeta '(-12)^2}{3 x^6 \zeta'(-6)^3}
 - \frac{\zeta '(-16)^2}{2 x^8 \zeta '(-8)^3}
 - \frac{2 \zeta '(-20)^2}{5 x^{10} \zeta '(-10)^3}
\]
\[
= \frac{4.51601}{x^2} - \frac{135.899}{x^4} + \frac{12996.}{x^6} - \frac{2.73295\times 10^6}{x^8} + \frac{1.03183\times 10^9}{x^{10}} .
\]

\ifthenelse {\boolean{BKMRK}}
  { \section{Tallying \texorpdfstring{$\lambda(n) \tau(n)^2$}{lambda Times tau Squared} }\label{S:LambdaTauSquared} }
  { \section{Tallying $\lambda(n) \tau(n)^2$}\label{S:LambdaTauSquared} }


This Dirichlet series for $\lambda(n) \tau(n)^2$ holds for $s > 1$ \cite[eq. D-50]{Gould}, \cite[p. 234]{McCarthy}:

\[
\sum _{n=1}^{\infty} \frac{\lambda (n) \tau(n)^2}{n^s} = \frac{\zeta (2 s)^3}{\zeta (s)^4}.
\]

The integral approximation to the summatory function $T_0(x)$ is
\begin{equation}\label{E:LambdaTauSquaredIntegralApprox}
T_0(x) \simeq \frac{1}{2 \pi i}
   \int_{c-i T}^{c+i T} \frac{\zeta (2 s)^3}{\zeta (s)^4} \frac{x^s}{s} \, ds.
\end{equation}

The sum over residues that may approximate the summatory function is
\begin{align}\label{E:LambdaTauSquaredSumApprox}
T_0(x)&\simeq a_1 + a_2 \sqrt{x} + a_3 \sqrt{x} \log (x) + a_4 \sqrt{x} \log(x)^2 \notag \\
 & \quad
 + 2 \Re\left(
       \sum _{k=1}^N x^{\rho _k}
       \frac{F_1(\rho_k) + F_2(\rho_k) \log(x) + F_3(\rho_k) \log(x)^2 + F_4(\rho_k) \log(x)^3}{G(\rho_k)}
     \right),
\end{align}
where
\[
a_1 = -2,
\]
\begin{align*}
a_2 = & \frac{1}{2\zeta \left(\frac{1}{2}\right)^6}
\big(
 \left(2-6 \gamma _1\right) \zeta
   \left(\frac{1}{2}\right)^2-6 \gamma  \zeta
   \left(\frac{1}{2}\right) \left(\zeta \left(\frac{1}{2}\right)+2
   \zeta '\left(\frac{1}{2}\right)\right) \\
   & \quad
   +\zeta \left(\frac{1}{2}\right) \left(4 \zeta
   '\left(\frac{1}{2}\right)-\zeta
   ''\left(\frac{1}{2}\right)\right)+6 \gamma ^2 \zeta
   \left(\frac{1}{2}\right)^2+5 \zeta '\left(\frac{1}{2}\right)^2
\big)
   \simeq 2.00358,
\end{align*}
where $\gamma$ is Euler's constant, and $\gamma_1 \simeq -0.072816$ is the first Stieltjes constant,
\[
a_3=-\frac{-3 \gamma  \zeta \left(\frac{1}{2}\right)+\zeta
   \left(\frac{1}{2}\right)+2 \zeta '\left(\frac{1}{2}\right)}{2
   \zeta \left(\frac{1}{2}\right)^5}\simeq -0.510158,
\]
and
\[a_4=\frac{1}{8 \zeta \left(\frac{1}{2}\right)^4}\simeq 0.027484.
\]

The functions $F_1$, $F_2$, $F_3$, $F_4$, and $G$ are:

\begin{align*}
F_1(s)& =
-15 s^3 \zeta (2 s)^2 \zeta ''(s)^3-s \zeta '(s)^2 (12 \zeta
   ''(s) (\zeta (2 s) (4 s^2 \zeta ''(2 s)+\zeta (2
   s)) \\
   & \quad
   +4 s^2 \zeta '(2 s)^2-4 s \zeta (2 s) \zeta '(2s)) \\
     & \quad
   +s \zeta (2 s) (s \zeta (2 s) \zeta ^{(4)}(s)-4
   \zeta ^{(3)}(s) (\zeta (2 s)-4 s \zeta '(2
   s))))+5 s^2 \zeta (2 s) \zeta '(s) \zeta ''(s) \\
   & \quad
   (2 s \zeta (2 s) \zeta ^{(3)}(s)-3 \zeta ''(s) (\zeta
   (2 s)-4 s \zeta '(2 s)))-2 \zeta '(s)^3 (-12 s
   \zeta '(2 s) (2 s^2 \zeta ''(2 s)+\zeta (2 s)) \\
   & \quad
   +12 s^2 \zeta '(2 s)^2+\zeta (2 s) (-8 s^3 \zeta ^{(3)}(2 s)+12 s^2
   \zeta ''(2 s)+3 \zeta (2 s))),
\end{align*}

\begin{align*}
F_2(s)& =
s \zeta '(s) (15 s^2 \zeta (2 s)^2 \zeta ''(s)^2+6 \zeta '(s)^2
   (\zeta (2 s) (4 s^2 \zeta ''(2 s)+\zeta (2 s))+4
   s^2 \zeta '(2 s)^2 \\
   & \quad
   - 4 s \zeta (2 s) \zeta '(2 s)) - 4 s \zeta
   (2 s) \zeta '(s) (s \zeta (2 s) \zeta ^{(3)}(s)-3 \zeta
   ''(s) (\zeta (2 s)-4 s \zeta '(2 s)))),
\end{align*}

\[
F_3(s) =
-3 s^2 \zeta (2 s) \zeta '(s)^2 \left(2 s \zeta (2 s) \zeta
   ''(s)+\zeta '(s) \left(\zeta (2 s)-4 s \zeta '(2 s)\right)\right),
\]

\[
F_4(s) = s^3 \zeta (2 s)^2 \zeta '(s)^3,
\]
and
\[
G(s) = 6 s^4 \zeta '(s)^7.
\]

The following \textit{Mathematica} code will perform these calculations:

\begin{verbatim}
A[s] = Zeta[2 s]^2 * x^s/s;
B[s] = Zeta[s]^4;
expr = residueFormula[4, s];
expr2 = Together[expr /. Zeta[s] -> 0];
numer = Collect[Numerator[expr2], {Log[x], Log[x]^2, Log[x]^3}, Simplify];
numer[[4]]/x^s             (* this is F1[s] *)
numer[[3]]/(x^s Log[x])    (* this is F2[s] *)
numer[[2]]/(x^s Log[x]^2)  (* this is F3[s] *)
numer[[1]]/(x^s Log[x]^3)  (* this is F4[s] *)
Denominator[expr2]         (* this is G[s] *)
\end{verbatim}

\begin{figure*}[ht]
  \centerline{
    \mbox{\includegraphics[width=\picDblWidth]{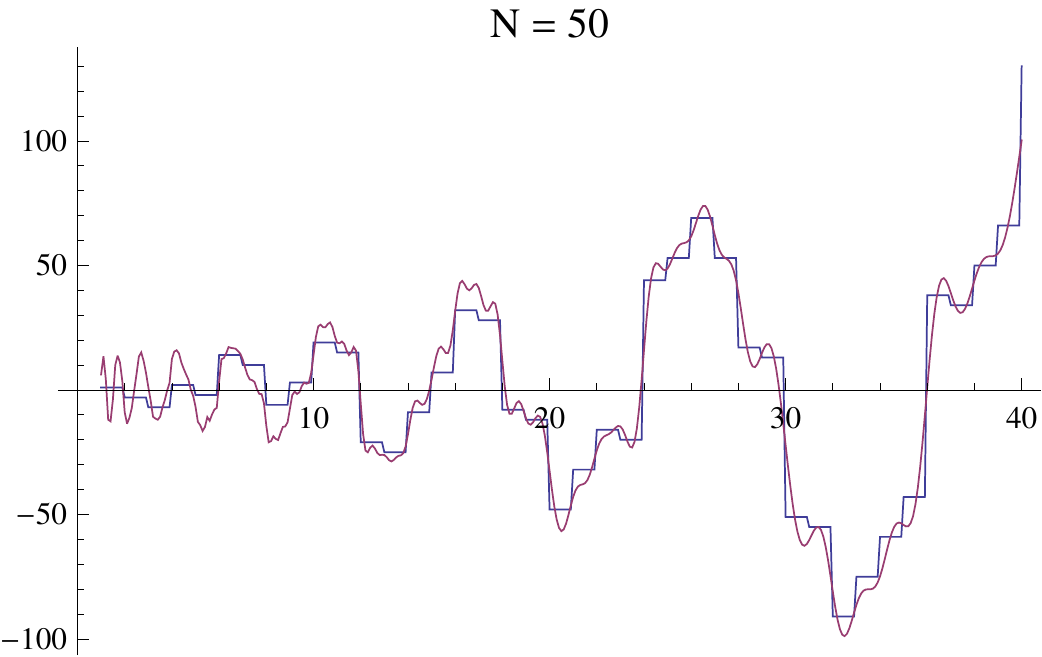}}
    \mbox{\includegraphics[width=\picDblWidth]{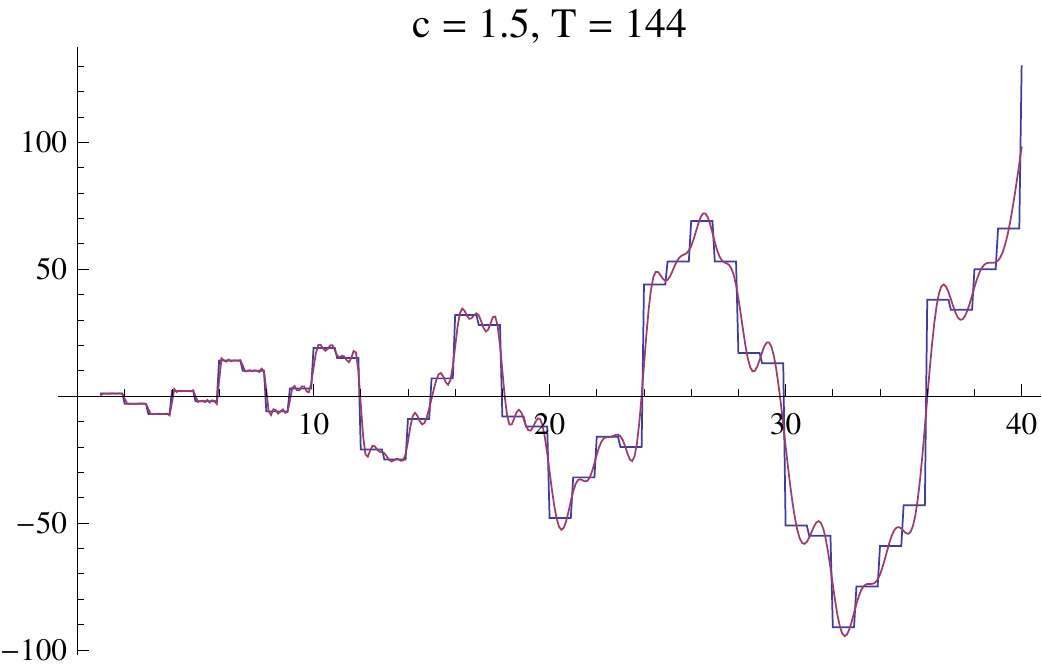}}
  }
  \caption{Approximating the sum of $\lambda(n) \tau(n)^2$ using a sum and an integral}
  \label{fig:LambdaTauSquaredGraphs}
\end{figure*}

Note: For this series, the corresponding sum over real zeros appears to diverge.  Here is the sum of the first five terms:
\[
- \frac{4 \zeta '(-4)^3}{x^2 \zeta '(-2)^4}
- \frac{2 \zeta '(-8)^3}{x^4 \zeta '(-4)^4}
- \frac{4 \zeta '(-12)^3}{3 x^6 \zeta '(-6)^4}
- \frac{\zeta '(-16)^3}{x^8 \zeta '(-8)^4}
- \frac{4 \zeta '(-20)^3}{5 x^{10} \zeta '(-10)^4}
\]
\[
= - \frac{2.36826}{x^2} - \frac{283.112}{x^4} - \frac{278745.}{x^6}
  - \frac{1.16534\times 10^9}{x^8} - \frac{1.44207\times 10^{13}}{x^{10}} .
\]

\ifthenelse {\boolean{BKMRK}}
  { \section{Tallying \texorpdfstring{$\lambda(n) \sigma(n)$}{lambda Times sigma} }\label{S:LambdaSigma} }
  { \section{Tallying $\lambda(n) \sigma(n)$}\label{S:LambdaSigma} }

This Dirichlet series for $\lambda(n) \sigma(n)$ holds for $s > 2$ \cite[set $k = 1$ in eq. D-47]{Gould}, \cite[p. 232]{McCarthy}:
\[
\sum _{n=1}^{\infty } \frac{\lambda (n) \sigma (n)}{n^s}
 = \frac{\zeta (2 s) \zeta (2s-2)}{\zeta (s) \zeta (s-1)}
\]

The integral approximation to the summatory function $T_0(x)$ is
\begin{equation}\label{E:LambdaSigmaIntegralApprox}
T_0(x) \simeq \frac{1}{2 \pi i}
   \int_{c-i T}^{c+i T} \frac{\zeta (2 s) \zeta (2s-2)}{\zeta (s) \zeta (s-1)} \frac{x^s}{s} \, ds.
\end{equation}

The integrand has poles at $s = 1/2$ and at $s = 3/2$.  The pole at $s = 1/2$ gives rise to the term involving $a_1$ in equation \eqref{E:LambdaSigmaSumApprox} below. The pole at $s = 3/2$ gives rise to the term involving $a_2$.

The integrand also has poles at each $s = \rho_k$ and each $s = \rho_k + 1$.  To compute the residues at these poles, we will proceed in two steps.  First, for the residue at each $s = \rho_k$, we apply equation \eqref{E:stdResidueFormula1} with
\[
A(s) = \frac{\zeta (2 s) \zeta (2s-2)}{\zeta (s-1)} \frac{x^s}{s}
\]
and
\[
B(s) = \zeta(s).
\]
Equation equation \eqref{E:stdResidueFormula1} tells us that the residue is
\[
\frac{A(s)}{B'(s)} = x^s \frac{\zeta (2 s) \zeta (2s-2)}{s \zeta (s-1) \zeta'(s)}.
\]
When we evaluate this at $s = \rho_k$, we get
\[
x^\rho_k \frac{\zeta (2 \rho_k) \zeta (2 \rho_k - 2)}{\rho_k \zeta (\rho_k-1) \zeta'(\rho_k)}.
\]

For the residue at each $s = \rho_k + 1$, we take
\[
A(s) = \frac{\zeta (2 s) \zeta (2s-2)}{\zeta (s)} \frac{x^s}{s}
\]
and
\[
B(s) = \zeta(s-1).
\]
Equation equation \eqref{E:stdResidueFormula1} tells us that the residue is
\[
\frac{A(s)}{B'(s)} = x^s \frac{\zeta (2 s) \zeta (2s-2)}{s \zeta(s) \zeta'(s-1)}.
\]
When we evaluate this at $s = \rho_k+1$, we get
\[
x^{\rho_k+1} \frac{\zeta (2 (\rho_k+1)) \zeta (2 \rho_k)}{(\rho_k+1) \zeta(\rho_k+1) \zeta'(\rho_k)}.
\]

Therefore, the formula for the sum of the residues, which may approximate $T_0(x)$, is
\begin{equation}\label{E:LambdaSigmaSumApprox}
T_0(x) \simeq a_1 \sqrt{x} + a_2 x^{3/2}
 + 2 \Re\left(\sum _{k=1}^N \left(x^{\rho _k}
   F_1\left(\rho _k\right)+x^{1+\rho _k} F_2\left(\rho _k\right)\right)\right),
\end{equation}

where

\[
a_1=-\frac{1}{12 \zeta \left(-\frac{1}{2}\right) \zeta \left(\frac{1}{2}\right)}\simeq
   -0.274495,
\]

\[
a_2=\frac{\zeta (3)}{3 \zeta \left(\frac{1}{2}\right) \zeta
   \left(\frac{3}{2}\right)}\simeq -0.105029,
\]

\[
F_1(s)=\frac{\zeta (2 s) \zeta (2 s-2)}{s \zeta (s-1) \zeta '(s)},
\]

and

\[
F_2(s)=\frac{\zeta (2 s) \zeta (2 s+2)}{(s+1) \zeta (s+1) \zeta '(s)}.
\]

Note the presence of $\zeta(2s)$ in both $F_1$ and $F_2$.  If $s = -2 k$ is a real zeta zero, then so is $2s = -4 k$.  Therefore, the sum over the real zeta zeros vanishes.

The following \textit{Mathematica} code will perform these calculations.  Parts A and B calculate the residues at $s = \rho_k$ and $s = \rho_k + 1$, respectively.
\begin{verbatim}
(* part A. integrand has a pole at each zero of Zeta[s] *)
(* set B[s] = Zeta[s]; set A[s] = everything else in the integrand *)
A[s] = Zeta[2 s] Zeta[2s - 2] / Zeta[s - 1] * x^s/s ;
B[s] = Zeta[s] ;
residueFormula[1, s] / x^s (* = F1[s] *)

(* part B. integrand has a pole at s = (1 + each zeta zero) *)
(* set B[s] = Zeta[s - 1]; set A[s] = everything else *)
A[s] = Zeta[2 s] Zeta[2s - 2] / Zeta[s] * x^s/s ;
B[s] = Zeta[s - 1] ;
residueFormula[1, s]/x^s /. s -> s+1 (* = F2[s] *)
\end{verbatim}

\begin{figure*}[ht]
  \centerline{
    \mbox{\includegraphics[width=\picDblWidth]{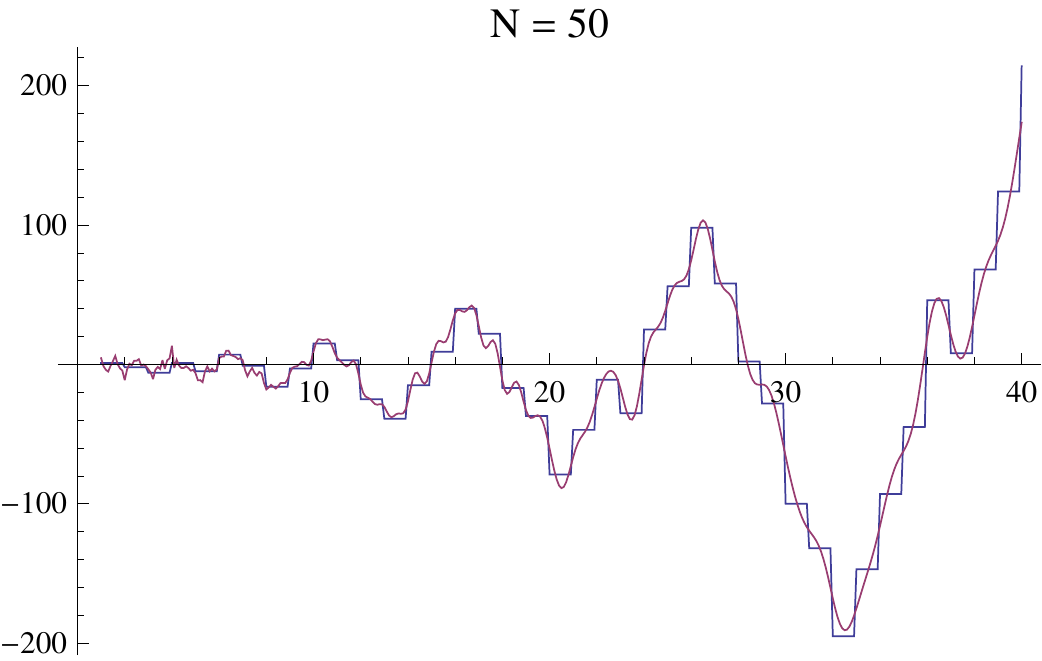}}
    \mbox{\includegraphics[width=\picDblWidth]{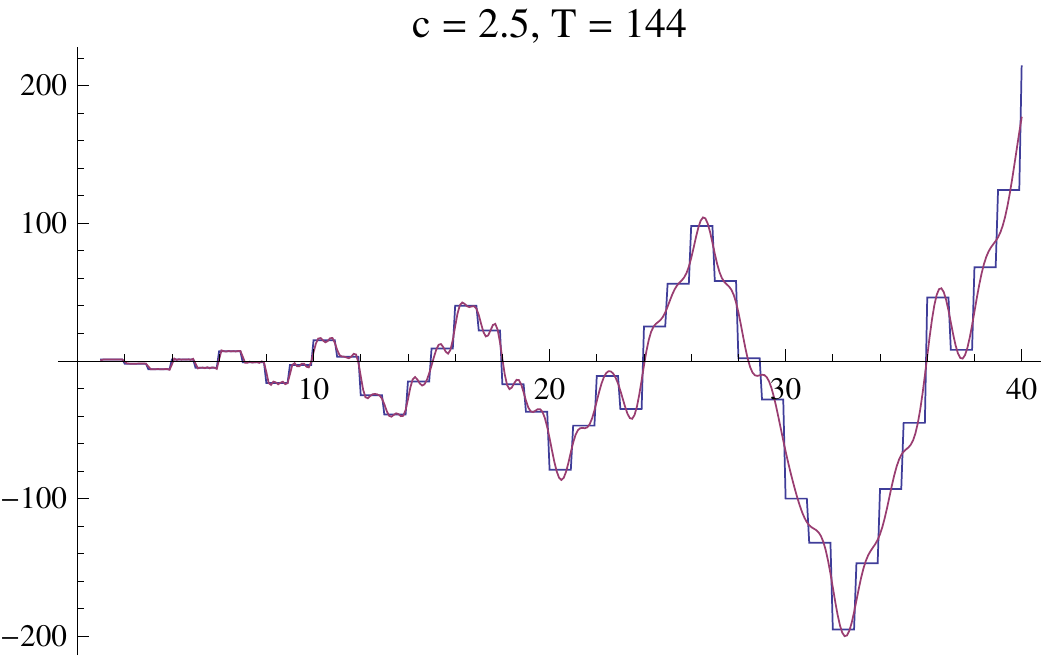}}
  }
  \caption{Approximating the sum of $\lambda(n) \sigma(n)$ using a sum and an integral}
  \label{fig:LambdaSigmaGraphs}
\end{figure*}

Figure ~\ref{fig:LambdaSigmaGraphs} shows the summatory function of $\lambda(n) \sigma(n)$ and an approximation using the first 50 pairs of zeta zeros.

\begin{figure*}[ht]
    \mbox{\includegraphics[width=\picDblWidth]{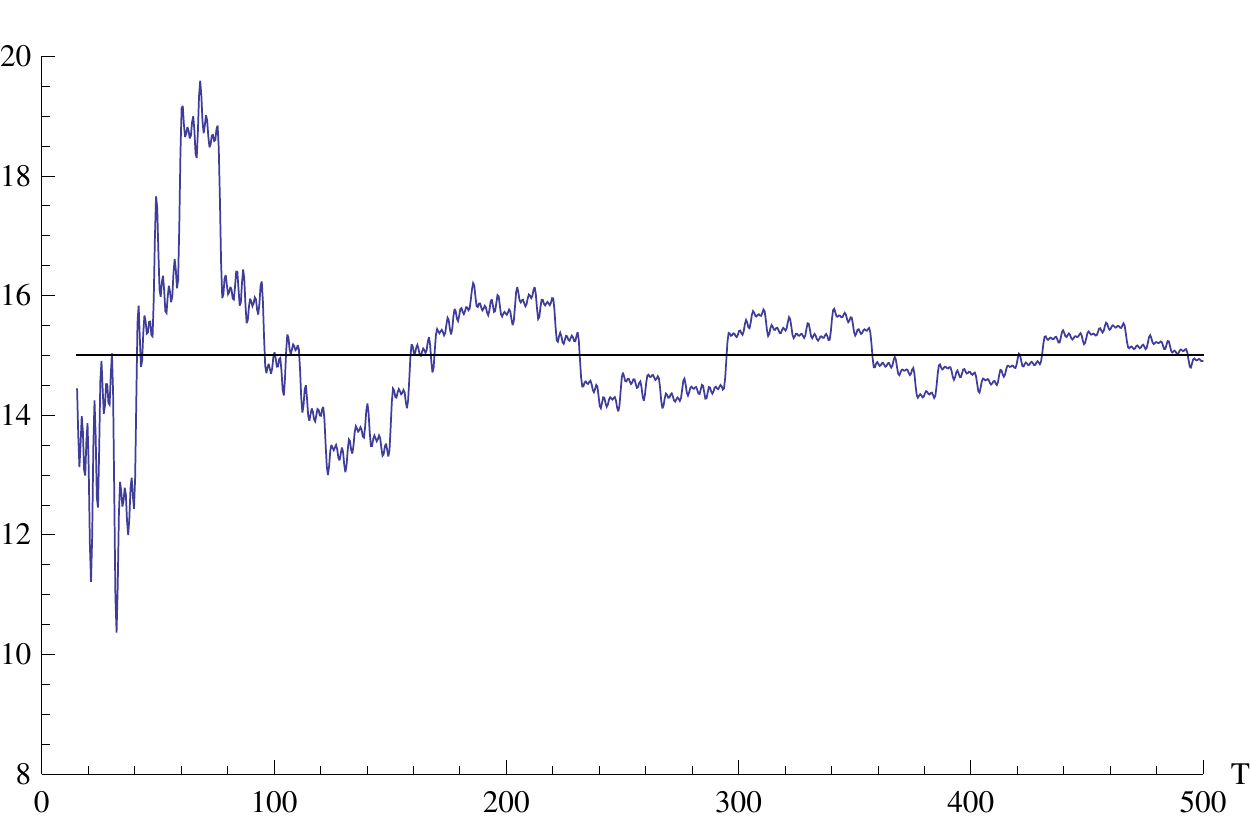}}
  \caption{$x = 10.5$, $c = 2.1$}
  \label{fig:LambdaSigmaIntVsT}
\end{figure*}

\begin{figure*}[ht]
  \mbox{\includegraphics[width=\picDblWidth]{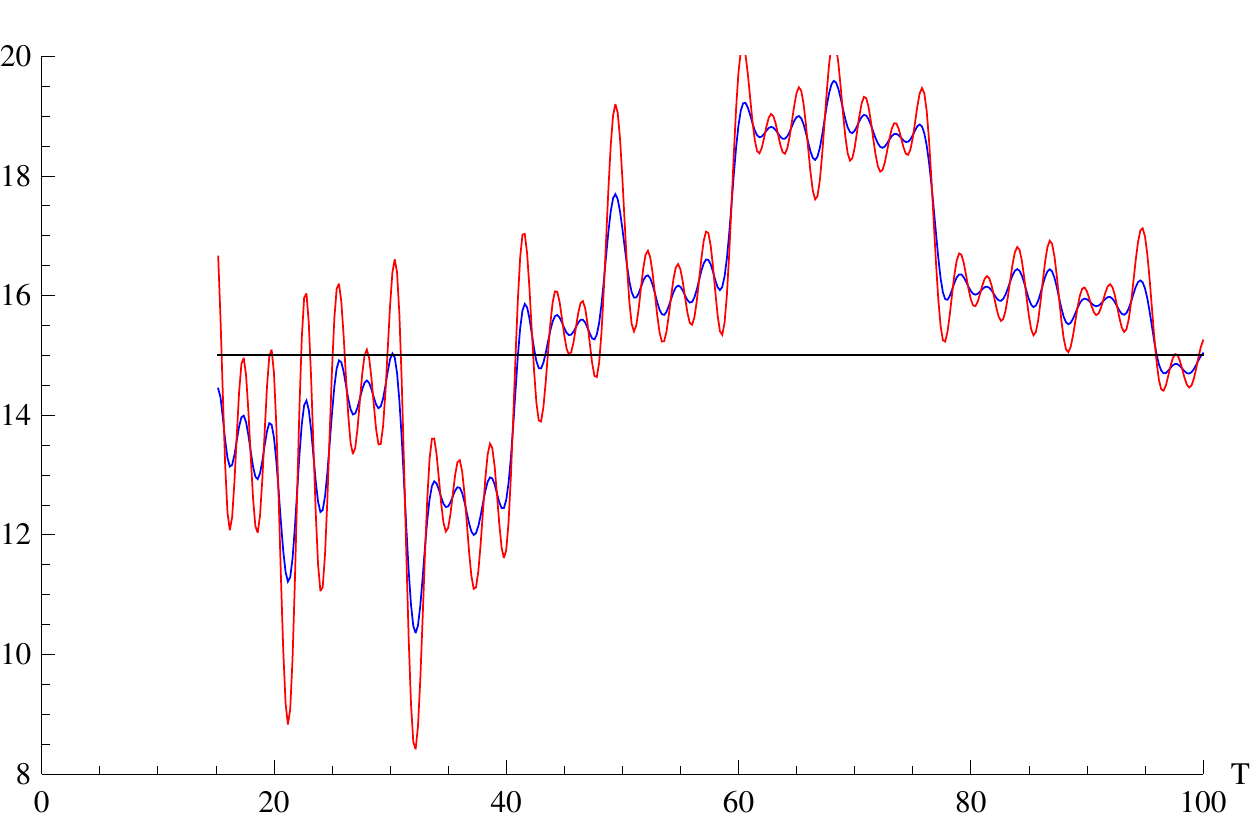}}
    \caption{$x = 10.5:  c = 2.1$ (blue) and $c = 2.5$ (red)}
  \label{fig:LambdaSigmac1Vc2}
\end{figure*}

Here is another way to visualize how close the Perron integral \eqref{E:LambdaSigmaIntegralApprox} is to the summatory function.  Take $x = 10.5$.  The summatory function at $x = 10.5$ is 15.  Set $c = 2.1$ in equation \eqref{E:LambdaSigmaIntegralApprox}.  Figure ~\ref{fig:LambdaSigmaIntVsT} shows how close the integral is to 15 as a function of $T$.

Again taking $x = 10.5$, figure ~\ref{fig:LambdaSigmac1Vc2} shows the integrals with $c = 2.1$ and $c = 2.5$.

\ifthenelse {\boolean{BKMRK}}
  { \section{Tallying \texorpdfstring{$\lambda(n) \sigma(n)^2$}{lambda Times sigma Squared} } }
  { \section{Tallying $\lambda(n) \sigma(n)^2$} }

This Dirichlet series for $\lambda(n) \sigma(n)^2$ holds for $s > 3$ \cite[set $h = k = 1$ in eq. D-46]{Gould}, \cite[p. 232]{McCarthy}:

\[
\sum _{n=1}^{\infty } \frac{\lambda (n) \sigma(n)^2}{n^s}
 = \frac{\zeta (2 s) \zeta (2 s - 2) \zeta(2s - 4)}{\zeta (s) \zeta (s-2) \zeta (s - 1)^2}.
\]

The integral approximation to the summatory function $T_0(x)$ is
\begin{equation}\label{E:LambdaSigmaSquaredIntegralApprox}
T_0(x) \simeq \frac{1}{2 \pi i}
   \int_{c-i T}^{c+i T}
    \frac{\zeta (2 s) \zeta (2 s - 2) \zeta(2s - 4)}{\zeta (s) \zeta (s-2) \zeta (s - 1)^2}
    \frac{x^s}{s} \, ds.
\end{equation}

The integrand has poles at $s = 1/2$, at $s = 3/2$, and at $s = 5/2$.  These residues give rise to the terms involving $a_1$, $a_2$, and $a_3$ in equation \eqref{E:LambdaSigmaSquaredSumApprox} below.

The integrand also has first-order poles at each $s = \rho_k$ and at $s = \rho_k + 2$.  These residues at these poles are computed the same way we computed the residues in Section \ref{S:LambdaSigma}, above.  These residues account for the terms
\[
x^{\rho_k} F_1(\rho_k) + x^{\rho_k + 2} F_2(\rho_k)
\]
in equation \eqref{E:LambdaSigmaSquaredSumApprox}.

The integrand also has a second-order pole at each $s = \rho_k + 1$.  To compute the residue at these poles, we use the standard formula \eqref{E:stdResidueFormula2} with
\[
A(s) = \frac{\zeta (2 s) \zeta (2s-2) \zeta(2s-4)}{\zeta (s) \zeta(s-2)} \frac{x^s}{s}
\]
and
\[
B(s) = \zeta(s-1)^2.
\]
This accounts for
\[
x^{\rho_k + 1} \frac{F_3(\rho_k) + F_4(\rho_k) \log(x)}{G(\rho_k)}
\]
in the next equation.

Therefore, the sum over the residues, which we hope will approximate the summatory function, is
\begin{align}\label{E:LambdaSigmaSquaredSumApprox}
T_0(x)&\simeq a_1 x^{1/2} + a_2 x^{3/2} + a_3 x^{5/2} \notag \\
 & \quad
 + 2 \Re\left(
       \sum _{k=1}^N
       x^{\rho_k} F_1(\rho_k) + x^{\rho_k + 2} F_2(\rho_k) + x^{\rho_k + 1} \frac{F_3(\rho_k) + F_4(\rho_k) \log(x)}{G(\rho_k)}
     \right),
\end{align}
where
\[
a_1=-\frac{1}{1440 \zeta \left(-\frac{3}{2}\right) \zeta
   \left(-\frac{1}{2}\right)^2 \zeta \left(\frac{1}{2}\right)}\simeq
   -0.431757,
\]
\[
a_2=-\frac{\zeta (3)}{36 \zeta \left(-\frac{1}{2}\right) \zeta
   \left(\frac{1}{2}\right)^2 \zeta \left(\frac{3}{2}\right)}\simeq
   0.02883,
\]
and
\[
a_3=\frac{\zeta (3) \zeta (5)}{5 \zeta \left(\frac{1}{2}\right)
   \zeta \left(\frac{3}{2}\right)^2 \zeta
   \left(\frac{5}{2}\right)}\simeq -0.018646.
\]
$F_1$, $F_2$, $F_3$, $F_4$, and $G$ are given by

\[
F_1(s)=\frac{\zeta (2 s) \zeta (2 s-4) \zeta (2 s-2)}{s \zeta (s-2)
   \zeta (s-1)^2 \zeta '(s)},
\]

\[
F_2(s)=\frac{\zeta (2 s) \zeta (2 s+2) \zeta (2 s+4)}{(s+2) \zeta
   (s+1)^2 \zeta (s+2) \zeta '(s)},
\]

\begin{align*}
F_3(s) &=
  -(s+1) \zeta (s-1) \zeta (2 s) \zeta (s+1) \zeta (2 s-2)
   \zeta (2 s+2) \zeta ''(s) \\
  & \quad
   -\zeta '(s) ((s+1) \zeta (2 s)
   \zeta (s+1) \zeta (2 s-2) \zeta (2 s+2) \zeta '(s-1) \\
  & \quad
   +\zeta (s-1)
   ((s+1) \zeta (2 s) \zeta (2 s-2) \zeta (2 s+2) \zeta'(s+1) \\
  & \quad
   +\zeta (s+1) (\zeta (2 s+2) (\zeta (2 s-2)
   (\zeta (2 s)-2 (s+1) \zeta '(2 s)) \\
  & \quad
   -2 (s+1) \zeta (2 s)
   \zeta '(2 s-2))-2 (s+1) \zeta (2 s) \zeta (2 s-2) \zeta '(2
   (s+1))))),
\end{align*}

\[
F_4(s)=(s+1) \zeta (s-1) \zeta (2 s) \zeta (s+1) \zeta (2 s-2) \zeta
   (2 s+2) \zeta '(s),
\]
and
\[
G(s)=(s+1)^2 \zeta (s-1)^2 \zeta (s+1)^2 \zeta '(s)^3.
\]

The following \textit{Mathematica} code will perform these calculations.  Parts A, B, and C calculate the residues at $s = \rho_k$, $s = \rho_k + 2$, and $s = \rho_k + 1$, respectively.
\smaller
\begin{verbatim}
(* part A. integrand has a pole at each zero of Zeta[s] *)
(* set B[s] = Zeta[s]; set A[s] = everything else in the integrand *)
A[s] = Zeta[2 s] Zeta[2 s - 2] Zeta[2 s - 4] / (Zeta[s - 2] Zeta[s - 1]^2) * x^s/s ;
B[s] = Zeta[s] ;
residueFormula[1, s] / x^s  (* = F1[s] *)

(* part B. integrand has a pole at s = (2 + each zeta zero) *)
(* set B[s] = Zeta[s - 2]; set A[s] = everything else *)
A[s] = Zeta[2 s] Zeta[2 s - 2] Zeta[2 s - 4] / (Zeta[s] Zeta[s - 1]^2) * x^s/s ;
B[s] = Zeta[s - 2] ;
residueFormula[1, s]/x^s  /. s -> s+2  (* = F2[s] *)

(* part C. integrand has a pole of order 2 at s = (1 + each zeta zero) *)
(* set B[s] = Zeta[s - 1]^2; set A[s] = everything else *)
A[s] = Zeta[2 s] Zeta[2 s - 2] Zeta[2 s - 4] / (Zeta[s] Zeta[s - 2]) * x^s/s ;
B[s] = Zeta[s - 1]^2 ;
expr = residueFormula[2, s];
(* we will evaluate expr only at s - 1 = (zeta zero), so remove all Zeta[s-1] *)
expr2 = Together[expr /. Zeta[s-1] -> 0];
(* compute the residue at s = (1 + each zeta zero), so make this substitution *)
expr3 = expr2 /. s -> s+1;
numer = Collect[Numerator[expr3], Log[x], Simplify];
numer[[2]]/x^(s+1)          (* = F3[s] *)
numer[[1]]/(x^(s+1) Log[x]) (* = F4[s] *)
Denominator[expr3]          (* = G[s] *)
\end{verbatim}
\larger

Note the presence of $\zeta(2s)$ in $F_1$, $F_2$, $F_3$, and $F_4$.  If $s = -2 k$ is a real zeta zero, then so is $2s = -4 k$.  Therefore, the sum over the real zeta zeros vanishes.

\begin{figure*}[ht]
  \centerline{
    \mbox{\includegraphics[width=\picDblWidth]{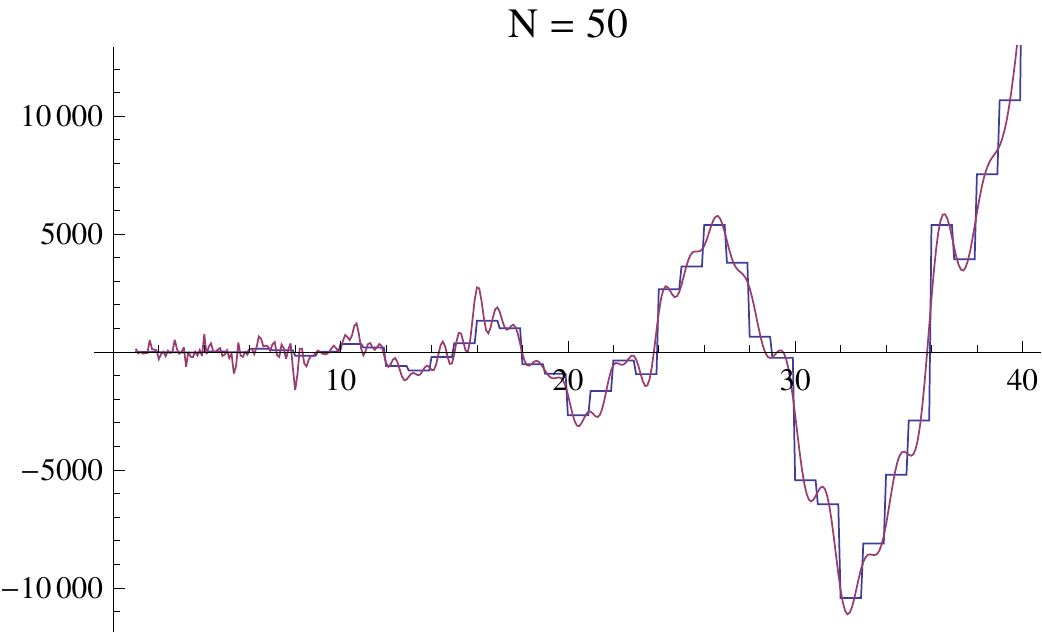}}
    \mbox{\includegraphics[width=\picDblWidth]{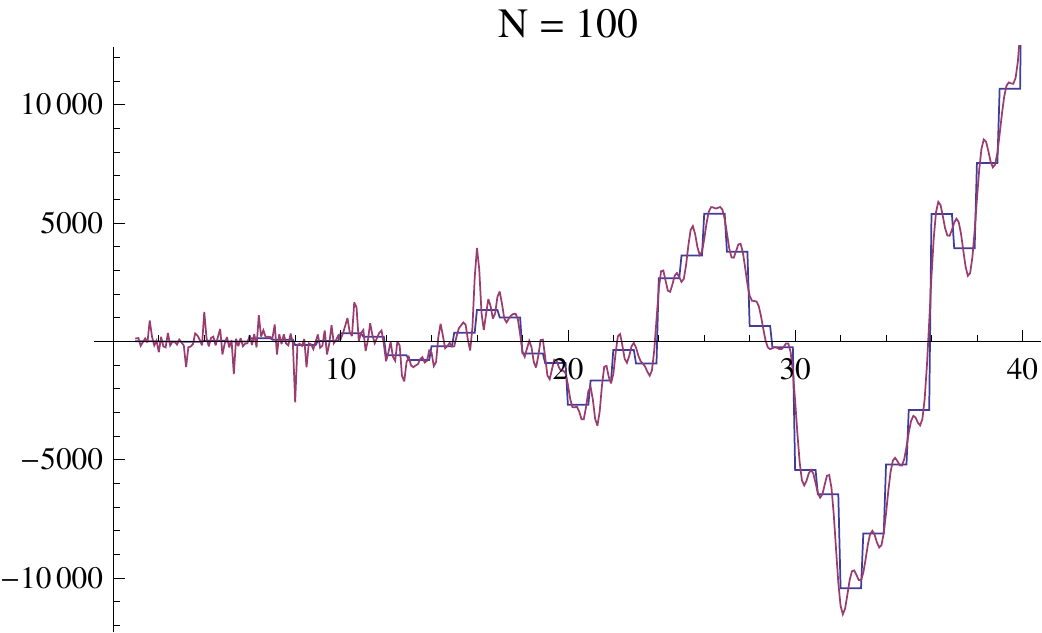}}
  }
  \caption{Approximating the sum of $\lambda(n) \sigma(n)^2$ using equation \eqref{E:LambdaSigmaSquaredSumApprox} with $N = 50$ and $N = 100$}
  \label{fig:LambdaSigmaSquaredSumGraphs}
\end{figure*}

\begin{figure*}[ht]
  \centerline{
    \mbox{\includegraphics[width=\picDblWidth]{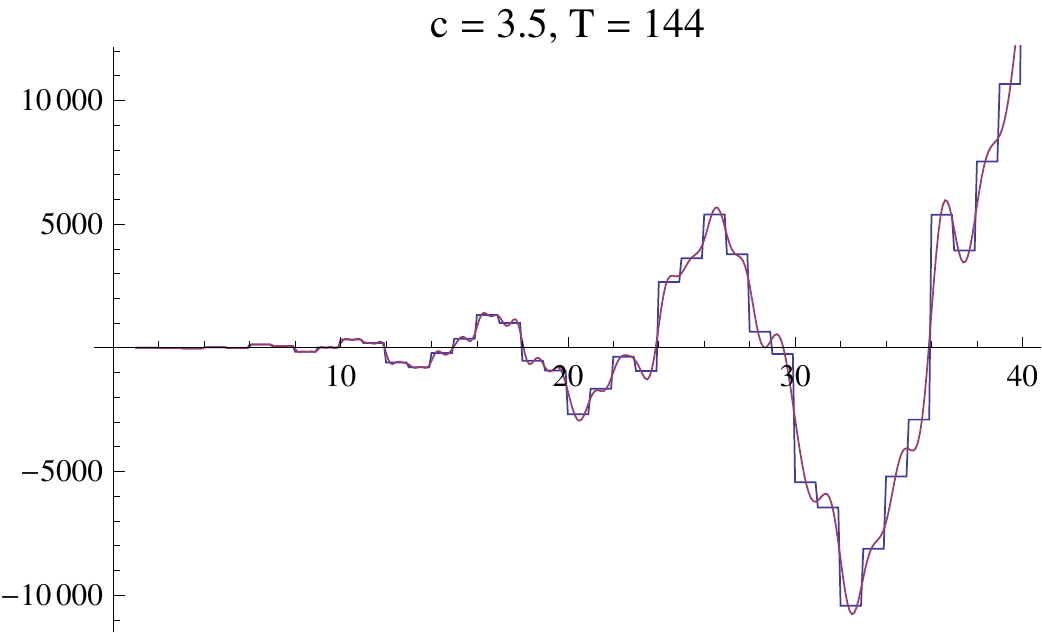}}
    \mbox{\includegraphics[width=\picDblWidth]{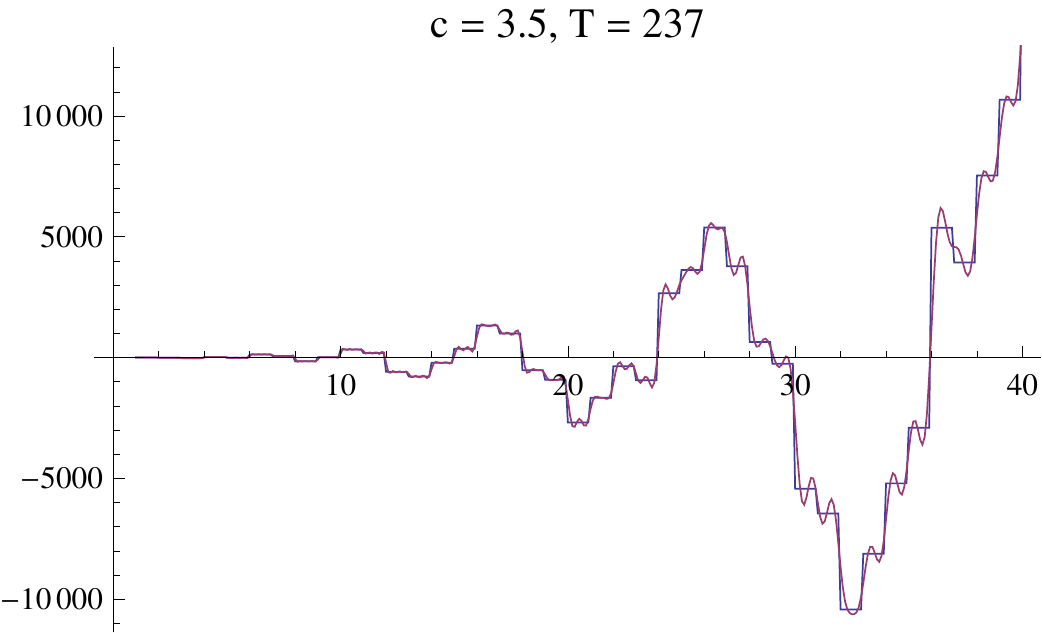}}
  }
  \caption{Approximating the sum of $\lambda(n) \sigma(n)^2$ using the integral in \eqref{E:LambdaSigmaSquaredIntegralApprox}}
  \label{fig:LambdaSigmaSquaredIntegralGraphs}
\end{figure*}

Figure ~\ref{fig:LambdaSigmaSquaredSumGraphs} shows the summatory function of $\lambda(n) \sigma(n)^2$ and the approximations that come from equation \eqref{E:LambdaSigmaSquaredSumApprox} using the first 50 pairs and 100 of zeta zeros.  Although both approximations can reproduce the large swings in the step function, neither approximation is particularly good at reproducing the details of the step function.  Further, the approximation does not appear to be improved by using 100 pairs of zeros instead of 50.  (As usual, using more zeros does enable the sum to reproduce the general shape of the step function for larger $x$).  On the other hand, Figure ~\ref{fig:LambdaSigmaSquaredIntegralGraphs} shows that the integral approximation \eqref{E:LambdaSigmaSquaredIntegralApprox} works quite well.

\ifthenelse {\boolean{BKMRK}}
  { \section{Tallying \texorpdfstring{$\lambda(n) \sigma(n) \tau(n)$}{lambda Times sigma Times tau} }\label{S:LambdaSigmaTau} }
  { \section{Tallying $\lambda(n) \sigma(n) \tau(n)$}\label{S:LambdaSigmaTau} }

This Dirichlet series for $\lambda(n) \sigma(n) \tau(n)$ holds for $s > 2$ \cite[eq. D-46]{Gould}, \cite[p. 232]{McCarthy}
\[
\sum _{n=1}^{\infty } \frac{\lambda (n) \sigma (n) \tau (n)}{n^s}
 = \frac{\zeta (2 s)^2 \zeta (2 s-2)^2}{\zeta (2 s-1) \zeta (s)^2 \zeta (s-1)^2}.
\]

The summatory function is
\[
T(x) = \sum _{k=1}^x \lambda (n) \sigma (n) \tau (n).
\]

The integral approximation to the modified summatory function $T_0(x)$ is
\begin{equation}\label{E:LambdaSigmaTauIntegralApprox}
T_0(x) \simeq \frac{1}{2 \pi i}
   \int_{c-i T}^{c+i T}
    \frac{\zeta (2 s)^2 \zeta (2 s-2)^2}{\zeta (2 s-1) \zeta (s)^2 \zeta (s-1)^2}
    \frac{x^s}{s} \, ds.
\end{equation}

The sum over the residues, which we hope will approximate the summatory function, is

\begin{align}\label{E:LambdaSigmaTauSumApprox}
& T_0(x) \simeq
 a_1 \sqrt{x} + a_2 \sqrt{x} \log (x) + a_3 x^{3/2} + a_4 x^{3/2} \log (x) \\
 & \quad
 + 2\Re\left(\sum _{k=1}^N \left(x^{\frac{1+\rho _k}{2}} F_1\left(\rho _k\right)+x^{\rho
   _k} \frac{F_2\left(\rho _k\right)+F_3\left(\rho _k\right) \log (x)}{G_1\left(\rho
   _k\right)}+x^{1+\rho _k} \frac{F_4\left(\rho _k\right)+F_5\left(\rho _k\right) \log
   (x)}{G_2\left(\rho _k\right)}\right)\right) \notag
\end{align}
where
\[
a_1=\frac{\zeta \left(-\frac{1}{2}\right) \left(\zeta \left(\frac{1}{2}\right) (-24 \log
   (A)+3+\log (2 \pi ))+\zeta '\left(\frac{1}{2}\right)\right)+\zeta
   \left(\frac{1}{2}\right) \zeta '\left(-\frac{1}{2}\right)-2 \gamma  \zeta
   \left(\frac{1}{2}\right) \zeta \left(-\frac{1}{2}\right)}{72 \zeta
   \left(-\frac{1}{2}\right)^3 \zeta \left(\frac{1}{2}\right)^3}\simeq 0.321773,
\]
where $\gamma$ is Euler's constant and $A \simeq 1.282427$ is Glaisher's constant defined by
\[
\log (A)=\frac{1}{12}-\zeta '(-1),
\]
\[
a_2=-\frac{1}{144 \zeta \left(-\frac{1}{2}\right)^2 \zeta
   \left(\frac{1}{2}\right)^2}\simeq -0.075348,
\]

\begin{align*}
a_3&=\frac{2 \zeta (3)}{3 \pi ^4 \zeta \Big(\frac{1}{2}\Big)^3 \zeta
    \Big(\frac{3}{2}\Big)^3} \\
   & \quad  
   \Big(-18 \zeta \Big(\frac{1}{2}\Big) \zeta
    \Big(\frac{3}{2}\Big) \zeta (3) \zeta '(2)+\pi ^2 \Big(6 \gamma \zeta
    \Big(\frac{1}{2}\Big) \zeta \Big(\frac{3}{2}\Big) \zeta (3)  \\
   & \quad 
   -3 \zeta \Big(\frac{3}{2}\Big) \zeta (3) \zeta
    '\Big(\frac{1}{2}\Big)-\zeta \Big(\frac{1}{2}\Big) \Big(3 \zeta
    (3) \zeta '\Big(\frac{3}{2}\Big)+\zeta \Big(\frac{3}{2}\Big)
    \Big(\zeta (3)-6 \zeta'(3)\Big)\Big)\Big)\Big) \\
    & \shoveleft{\simeq -0.002403},
\end{align*}

\[
a_4=\frac{\zeta (3)^2}{\pi ^2 \zeta \left(\frac{1}{2}\right)^2 \zeta
   \left(\frac{3}{2}\right)^2}\simeq 0.010059,
\]

\[
F_1(s)=\frac{\zeta (s-1)^2 \zeta (s+1)^2}{(s+1) \zeta \left(\frac{s-1}{2}\right)^2 \zeta
   \left(\frac{s+1}{2}\right)^2 \zeta '(s)},
\]

\begin{align*}
F_2(s)&=-\zeta (2 s) \zeta (2 s-2) \\
& \quad
 \Big(s \zeta (s-1) \zeta (2 s) \zeta (2 s-2) \zeta (2 s-1) \zeta ''(s) +\zeta '(s) \\
& \quad
 \Big(2 s \zeta (2 s) \zeta (2 s-2) \zeta (2 s-1) \zeta '(s-1) +\zeta (s-1) \\
& \quad
 \Big(\zeta (2 s) \Big(\zeta (2 s-2) \Big(\zeta (2 s-1)+2 s \zeta
   '(2 s-1)\Big) \\
& \quad
   -4 s \zeta (2 s-1) \zeta '(2 s-2)\Big)
   -4 s \zeta (2 s-2) \zeta (2
   s-1) \zeta '(2 s)\Big)\Big)\Big),
\end{align*}

\[
F_3(s)=s \zeta (s-1) \zeta (2 s)^2 \zeta (2 s-2)^2 \zeta (2 s-1) \zeta '(s),
\]

\[
G_1(s)=s^2 \zeta (s-1)^3 \zeta (2 s-1)^2 \zeta '(s)^3,
\]

\begin{align*}
F_4(s)&=\zeta (2 s) \zeta (2 (s+1)) \\
& \quad
 \Big(-(s+1) \zeta (2 s) \zeta (s+1) \zeta (2 (s+1)) \zeta (2
   s+1) \zeta ''(s)-\zeta '(s) \\
 & \quad
    \Big(\zeta (2 s) \Big(2 (s+1) \zeta (2 (s+1)) \zeta (2 s+1) \zeta '(s+1)+\zeta (s+1) \\
 & \quad
    \Big(\zeta (2 (s+1)) \Big(\zeta (2 s+1)+2 (s+1) \zeta '(2 s+1)\Big)-4
   (s+1) \zeta (2 s+1) \zeta '(2 (s+1))\Big)\Big) \\
 & \quad
   -4 (s+1) \zeta (s+1) \zeta (2 (s+1)) \zeta (2s+1) \zeta '(2 s)\Big)\Big),
\end{align*}

\[
F_5(s)=(s+1) \zeta (2 s)^2 \zeta (s+1) \zeta (2 (s+1))^2 \zeta (2 s+1) \zeta '(s),
\]
and, finally,
\[
G_2(s)=(s+1)^2 \zeta (s+1)^3 \zeta (2 s+1)^2 \zeta '(s)^3.
\]

The following \textit{Mathematica} code will perform these calculations:

\smaller
\begin{verbatim}
(* part A. integrand has a pole at each zero of Zeta[2s - 1] *)
(* set B[s] = Zeta[2s - 1]; set A[s] = everything else in the integrand *)
A[s] = Zeta[2s]^2 Zeta[2s - 2]^2 / (Zeta[s]^2 Zeta[s - 1]^2) * x^s/s ;
B[s] = Zeta[2s - 1] ;
(* compute the residue at s = (1 + (zeta zero))/2, so make this substitution *)
residueFormula[1, s] / x^s  /. s -> (1 + s)/2  (* = F1[s] *)

(* part B. integrand has a pole of order 2 at each zero of Zeta[s] *)
(* set B[s] = Zeta[s]^2; set A[s] = everything else in the integrand *)
A[s] = Zeta[2s]^2 Zeta[2s - 2]^2 / (Zeta[2s - 1] Zeta[s - 1]^2) * x^s/s ;
B[s] = Zeta[s]^2 ;
exprB = residueFormula[2, s];
(* we will evaluate exprB only at s = zero of zeta function, so remove Zeta[s] *)
exprB2 = Together[exprB /. Zeta[s] -> 0];
numerB = Collect[Numerator[exprB2], Log[x], Simplify];
numerB[[2]]/x^s           (* = F2[s] *)
numerB[[1]]/(x^s Log[x])  (* = F3[s] *)
Denominator[exprB2]       (* = G1[s] *)

(* part C. integrand has a pole of order 2 at each zero of Zeta[s - 1] *)
(* set B[s] = Zeta[s-1]^2; set A[s] = everything else in the integrand *)
A[s] = Zeta[2s]^2 Zeta[2s - 2]^2 / (Zeta[2s - 1] Zeta[s]^2) * x^s/s ;
B[s] = Zeta[s - 1]^2 ;
exprC = residueFormula[2, s];
(* we evaluate exprC only at s-1 = zero of zeta function, so remove Zeta[s-1] *)
exprC2 = Together[exprC /. Zeta[s-1] -> 0];
(* we want the residue at s = (1 + each zeta zero); so make this substitution *)
exprC3 = exprC2 /. s -> s+1;
numerC = Collect[Numerator[exprC3], Log[x], Simplify];
numerC[[2]]/x^(s+1)           (* = F4[s] *)
numerC[[1]]/(x^(s+1) Log[x])  (* = F5[s] *)
Denominator[exprC3]           (* = G2[s] *)
\end{verbatim}
\larger

What about the sum over the real zeta zeros?  This sum would have the same form as the terms in \eqref{E:LambdaSigmaTauSumApprox}, but with $\rho_k$ replaced with $-2 k$.  But because of the $\zeta(2s)$ that is present in $F_2$, $F_3$, $F_4$, and $F_5$, we see that those parts of this sum would be zero.  However, the part of the sum involving $F_1(s)$ evaluated at $s = -2 k$, that is,
\[
\sum _{k=1}^{M} x^{ \frac{1 + (-2k)}{2} } F_1(-2k)
\]
is quite badly behaved: the coefficients $F_1(-2k)$ rapidly grow large.  The first ten are approximately $0.56426$, $-0.96909$, $6.478$, $-90.794$, $2197.1$, $-82446$, $4.4633 \cdot 10^6$, $-3.3062 \cdot 10^8$,  $3.2185 \cdot 10^{10}$, and $-3.9889 \cdot 10^{12}$.  As we saw in section ~\ref{S:LambdaTwoNu}, if the sum over real zeros does not improve the accuracy of the approximation to $T(x)$, we may be able to omit those terms.  All we have to do is integrate around a rectangle in the complex plane that does not enclose the real zeta zeros.

\begin{figure*}[ht]
  \centerline{
    \mbox{\includegraphics[width=\picDblWidth]{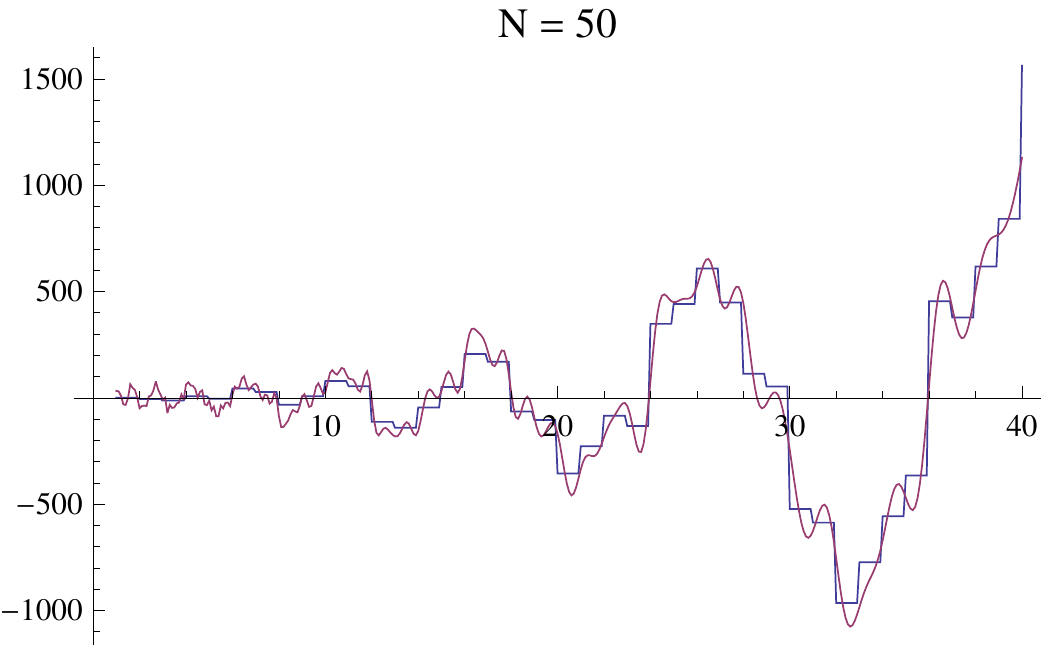}}
    \mbox{\includegraphics[width=\picDblWidth]{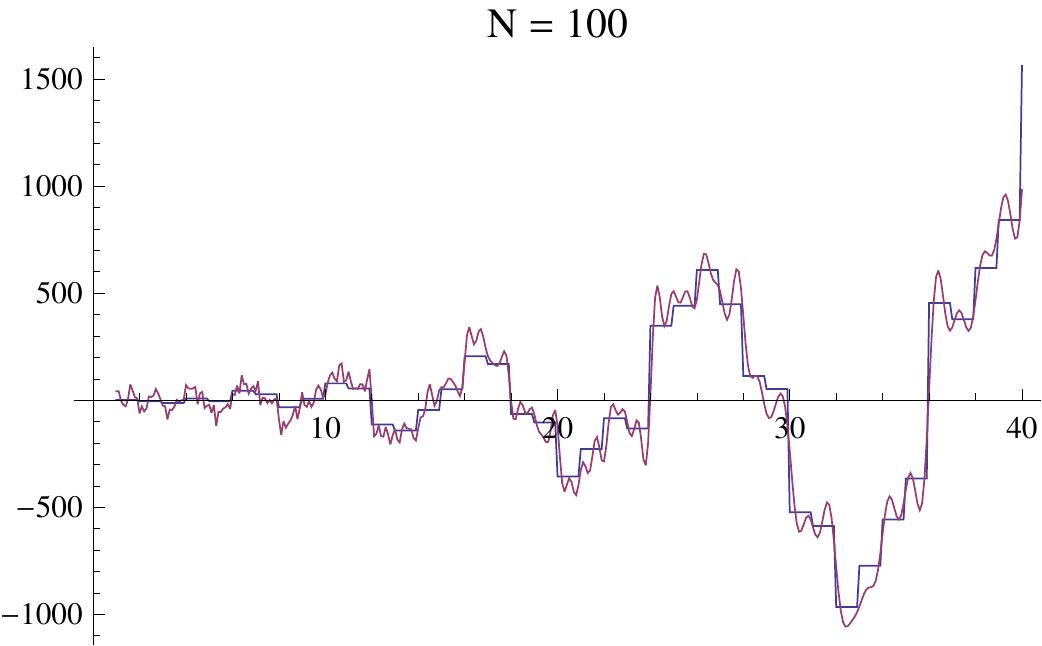}}
  }
  \caption{Approximating the sum of $\lambda(n) \sigma(n) \tau(n)$ using \eqref{E:LambdaSigmaTauSumApprox} with $N = 50$ and $N = 100$}
  \label{fig:LambdaSigmaTauSumGraphs}
\end{figure*}

\begin{figure*}[ht]
  \centerline{
    \mbox{\includegraphics[width=\picDblWidth]{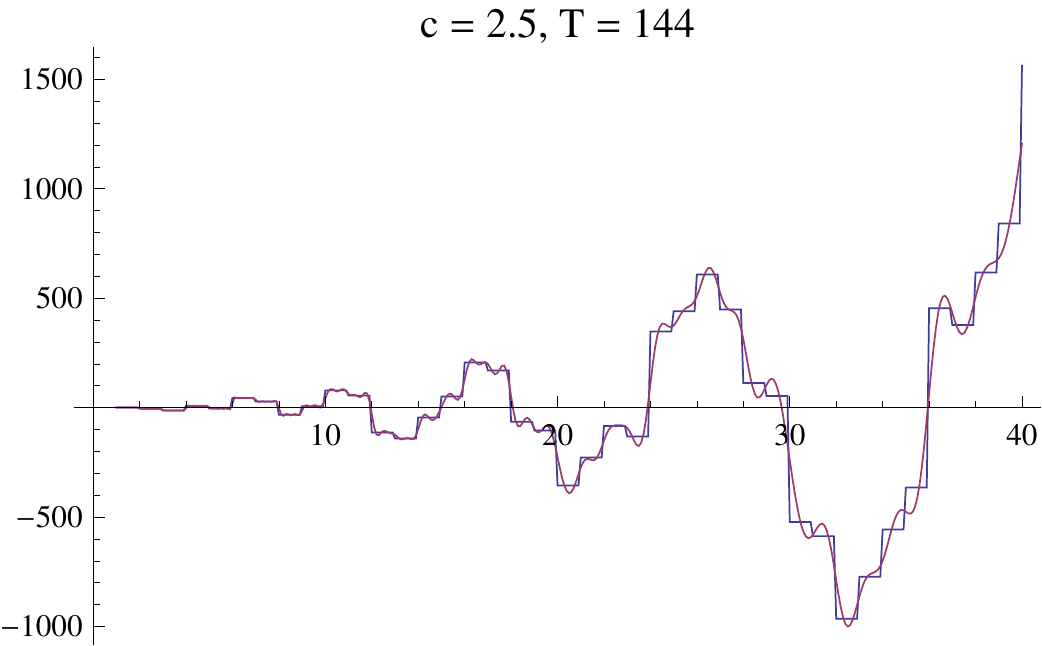}}
    \mbox{\includegraphics[width=\picDblWidth]{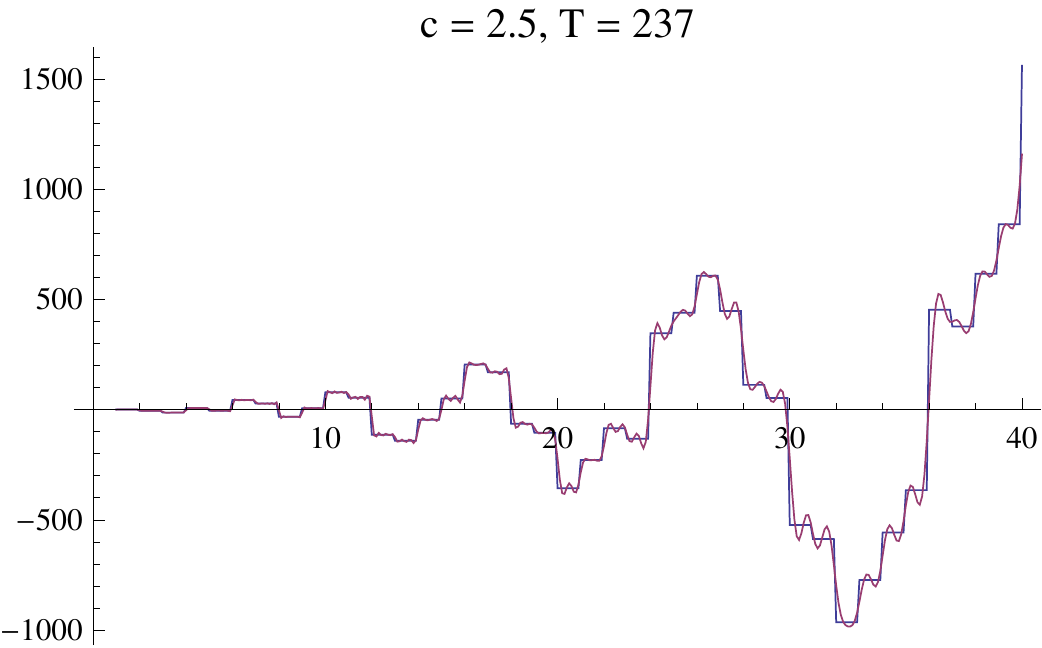}}
  }
  \caption{Approximating the sum of $\lambda(n) \sigma(n) \tau(n)$ using the integral in \eqref{E:LambdaSigmaTauIntegralApprox}}
  \label{fig:LambdaSigmaTauIntegralGraphs}
\end{figure*}

Figure ~\ref{fig:LambdaSigmaTauSumGraphs} shows the summatory function of $\lambda(n) \sigma(n) \tau(n)$ and the approximations that come from equation \eqref{E:LambdaSigmaTauSumApprox} using the first 50 pairs and 100 of zeta zeros.  Neither sum is particularly good at reproducing the details of the step function.  Figure ~\ref{fig:LambdaSigmaTauIntegralGraphs} shows that the integral approximation \eqref{E:LambdaSigmaTauIntegralApprox} works quite well.

\section{Tallying Greatest Common Divisors}\label{S:TallyGCD}

Define the function $P(n)$ to be
\[
P(n) = \sum_{k=1}^{n} \gcd(k, n).
\]

The following Dirichlet series for $P(n)$ holds for $s > 2$ \cite[eq. 15]{Toth}:
\[
\sum _{n=1}^{\infty } \frac{P(n)}{n^s} = \frac{\zeta(s - 1)^2}{\zeta(s)}.
\]
This Dirichlet series comes from the fact that $P(n)$ is the convolution of $(\phi * Id$)(n) where $Id(n) = n$ is the identity function, so the dirichlet series for $P$ is the product of the Dirichlet series of $\phi$ and $Id$:
\[
\sum _{n=1}^{\infty } \frac{P(n)}{n^s} = \frac{\zeta(s - 1)}{\zeta(s)}  \sum _{n=1}^{\infty } \frac{n}{n^s}
= \frac{\zeta(s - 1)}{\zeta(s)}  \zeta(s - 1).
\]
Here, we used the Dirichlet series for $\phi$, given by \eqref{E:PhiDirichletSeries}.

The summatory function is
\[
T(x) = \sum _{k=1}^x P(n).
\]

The integral approximation to the summatory function $T(x)$ is
\begin{equation}\label{E:SumGCDIntegralApprox}
T_0(x) \simeq \frac{1}{2 \pi i}
   \int_{c-i T}^{c+i T} \frac{\zeta(s - 1)^2}{\zeta(s)} \frac{x^s}{s} \, ds.
\end{equation}

The integrand has a pole of order 1 at $s = 0$, and a pole of order 2 at $s = 2$.  There is also a pole at each zeta zero.  The residue at $s = 0$ is $-1/72$.  The residue at $s = 2$ is
\[
\frac{3 x^2 \left((4 \gamma -1) \pi ^2-12 \zeta '(2)\right)}{2 \pi
   ^4}+\frac{3 x^2 \log (x)}{\pi ^2}.
\]

The sum of the residues, which we hope will be close to $T_0(x)$, is
\begin{equation}\label{E:SumGCDSumApprox}
T_0(x) \simeq
 a_1 + a_2 x^2+a_3 x^2 \log (x)+
 2 \Re\left(
     \sum _{k=1}^N x^{\rho _k} \frac{ \zeta \left(\rho_k-1\right){}^2}{\rho _k \zeta '\left(\rho _k\right)}
   \right)
   + \sum _{k=1}^M x^{-2 k} \frac{ \zeta (-2 k-1)^2}{(-2 k) \zeta '(-2 k)}
\end{equation}
where
\[
a_1 = -1/72 \simeq -0.013889,
\]
\[
a_2 = \frac{3 \left((4 \gamma -1) \pi ^2-12 \zeta '(2)\right)}{2 \pi ^4} \simeq 0.37217,
\]
and
\[
a_3 = 3/\pi^2 \simeq 0.30396.
\]

\begin{figure*}[ht]
  \mbox{\includegraphics[width=\picDblWidth]{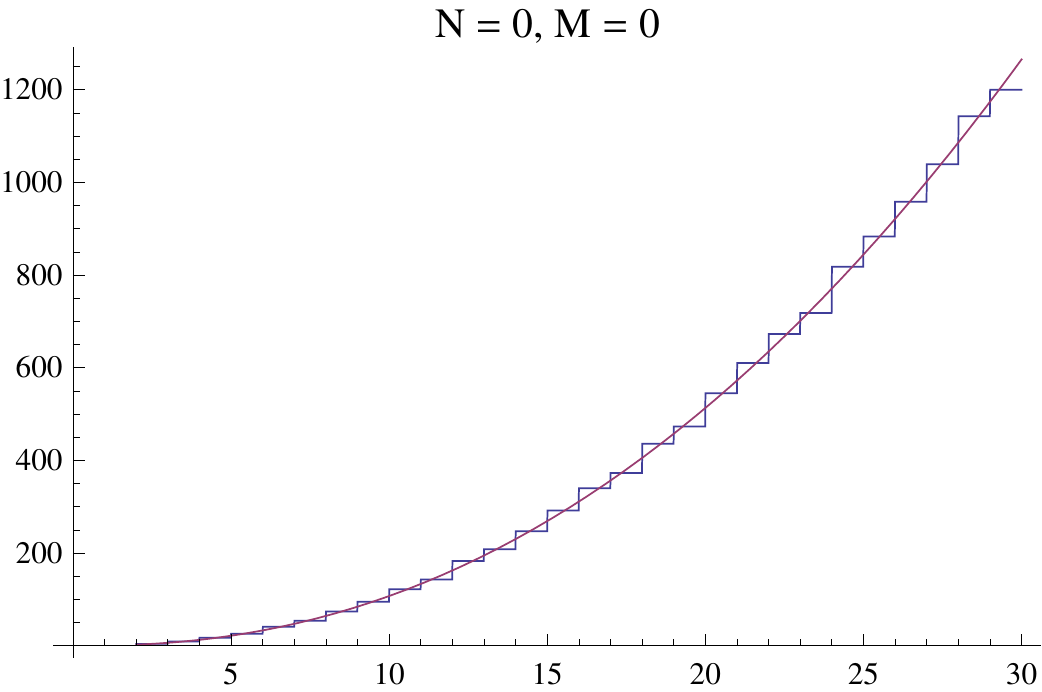}}
  \caption{Approximating the sum of $ P(n)$ with the first three terms of \eqref{E:SumGCDSumApprox}}
  \label{fig:GCDApprox0}
\end{figure*}

\begin{figure*}[ht]
  \centerline{
    \mbox{\includegraphics[width=\picDblWidth]{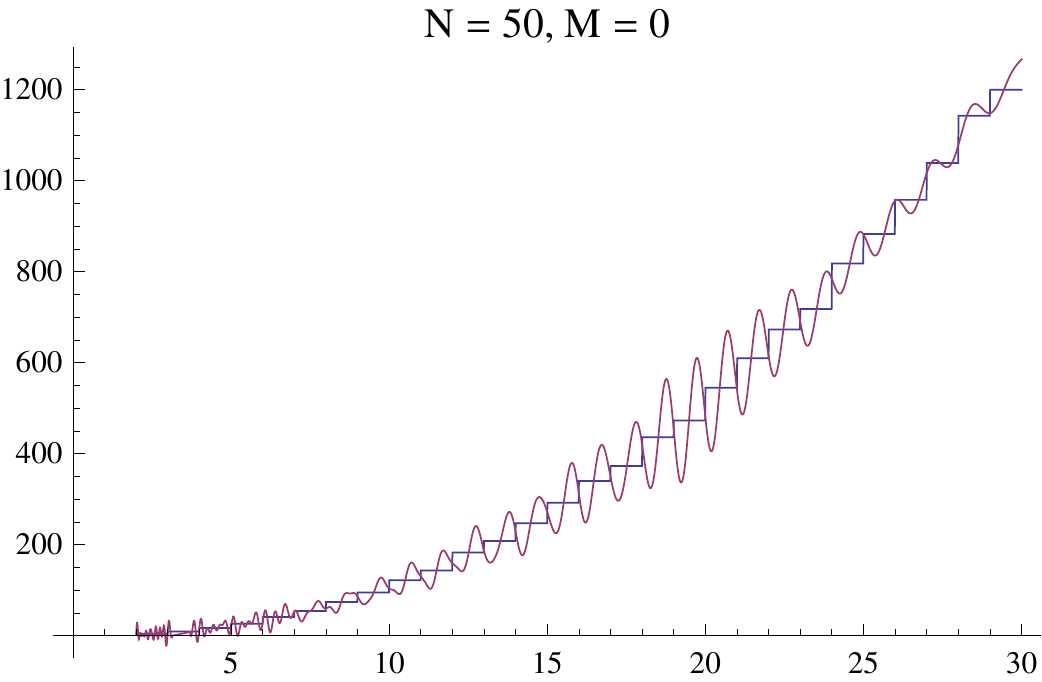}}
    \mbox{\includegraphics[width=\picDblWidth]{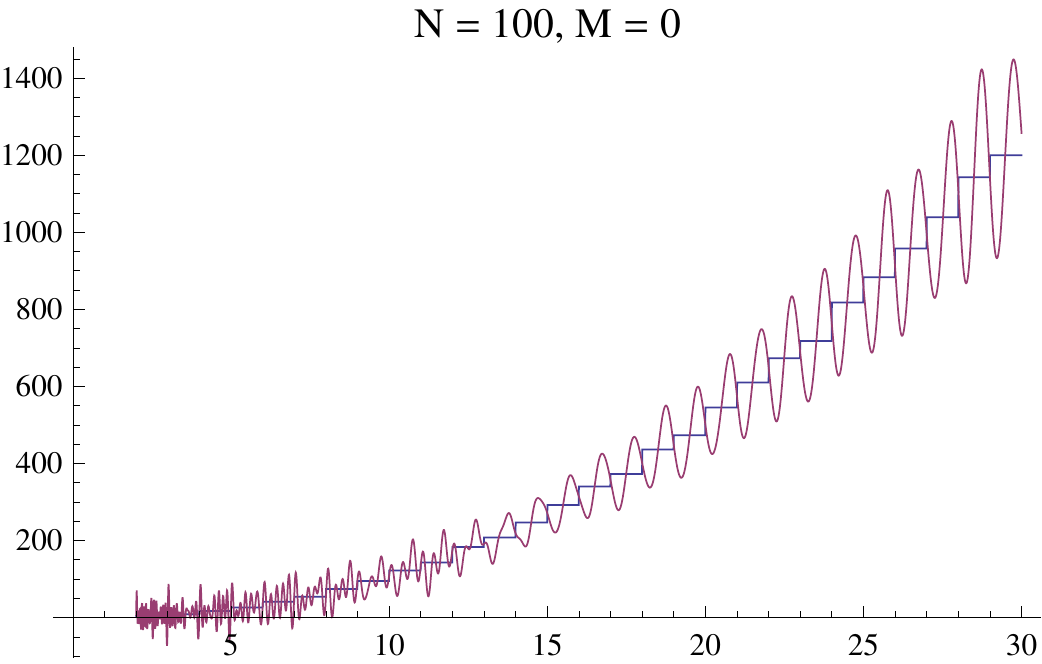}}
  }
  \caption{Approximating the sum of $P(n)$ using 50 and 100 pairs of zeta zeros}
  \label{fig:GCDSumGraphs}
\end{figure*}

\begin{figure*}[ht]
  \mbox{\includegraphics[width=\picDblWidth]{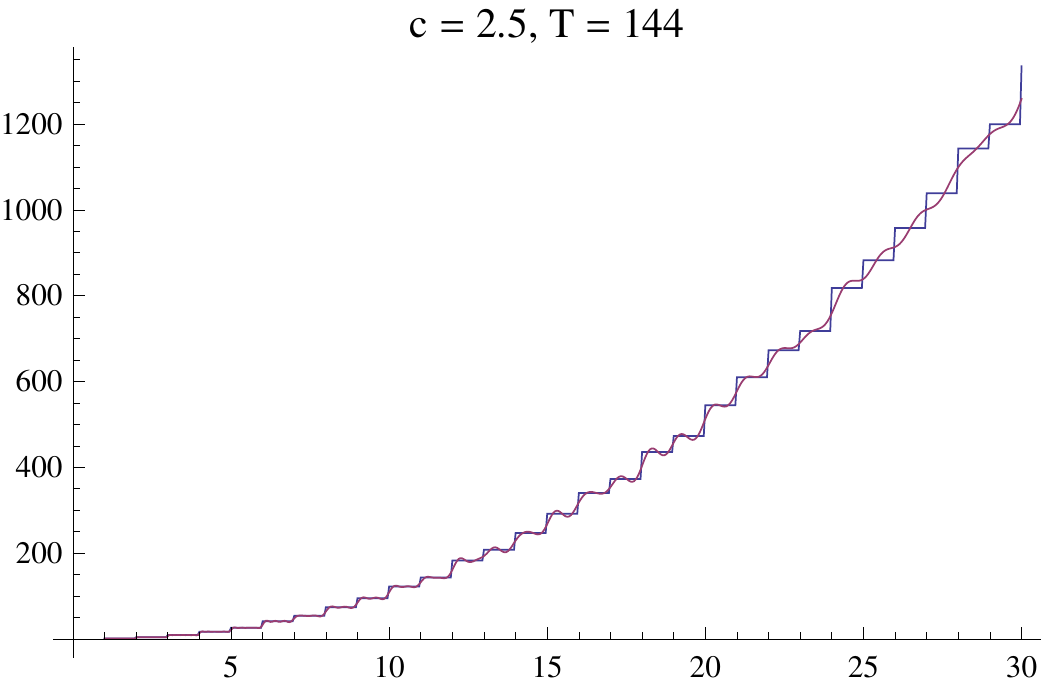}}
  \caption{Approximating the sum of $P(n)$ using an integral}
  \label{fig:GCDIntegralGraphs}
\end{figure*}

Figure ~\ref{fig:GCDApprox0} shows the graph of the sum approximation \eqref{E:SumGCDSumApprox} without using any zeta zeros.  Figure ~\ref{fig:GCDSumGraphs} shows the graphs of the sum approximation \eqref{E:SumGCDSumApprox} using 50 and 100 pairs of complex zeta zeros.  Figure ~\ref{fig:GCDIntegralGraphs} shows the integral approximation \eqref{E:SumGCDIntegralApprox}.

The sums using 50 zeta zeros is pretty poor.  The obvious thing to try is to use more zeros.  But that only makes the graph worse.
If one makes a table of values of, say, the first 50 coefficients in the first sum in equation \eqref{E:SumGCDSumApprox}, that is,
\[
\frac{ \zeta \left(\rho_k-1\right){}^2}{\rho _k \zeta '\left(\rho _k\right)}
\]
for $1 \leq k \leq 50$, one sees that neither the real nor the imaginary parts of the coefficients are getting small.  This is in contrast to previous cases where the sum yielded a good approximation.  It may be that the first sum in \eqref{E:SumGCDSumApprox} does not even converge.

\begin{figure*}[ht]
  \centerline{
    \mbox{\includegraphics[width=\picDblWidth]{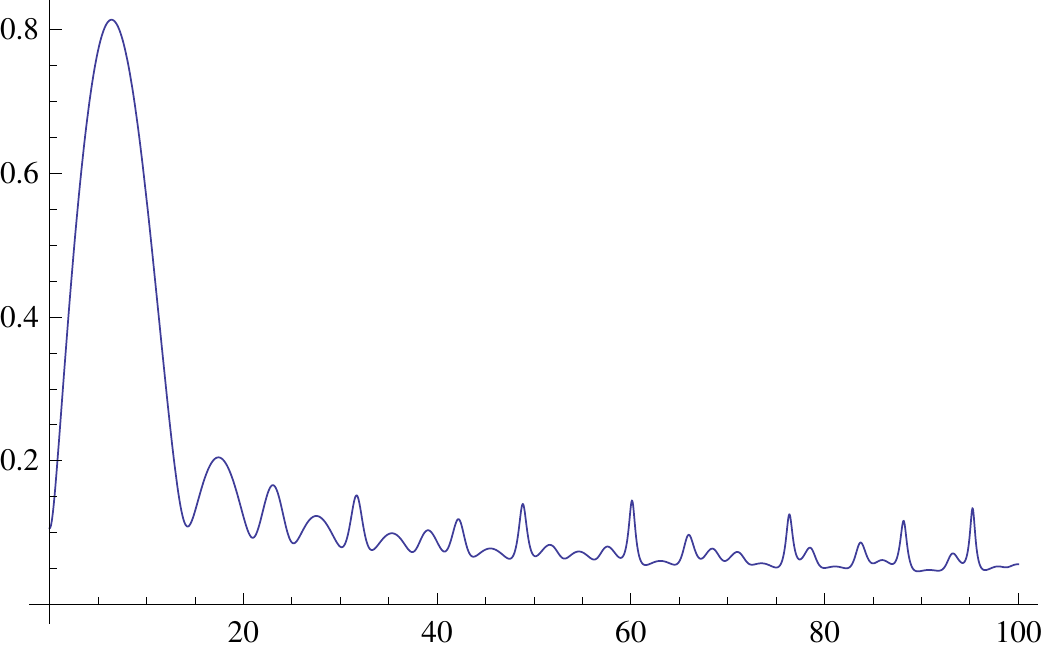}}
    \mbox{\includegraphics[width=\picDblWidth]{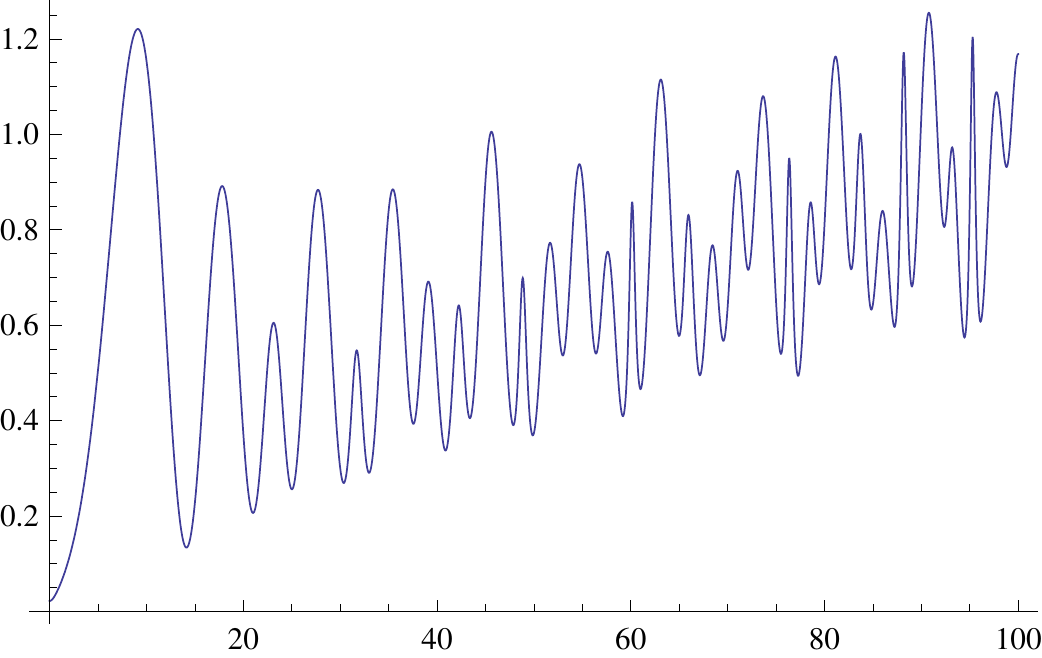}}
  }
  \caption{Left: function \eqref{E:PhiZeta}; right: function \eqref{E:GCDZeta}}
  \label{fig:PhiAndGCDZetaGraphs}
\end{figure*}

To see what's going on, let's compare the two functions
\[
\frac{ \zeta(\rho_k-1)}{\rho _k \zeta'(\rho _k)}
\]
from \eqref{E:PhiSum}, the sum for the summatory function of $\phi(n)$, and
\[
\frac{ \zeta(\rho_k-1){}^2}{\rho _k \zeta'(\rho _k)},
\]
from the summatory function of this section.  In particular, let's graph the function
\[
\frac{ \zeta(s-1)}{s \zeta'(s)}
\]
where $s$ lies on the line $1/2 + i t$, that is,
\begin{equation}\label{E:PhiZeta}
\frac{ \zeta(1/2 + i t - 1)}{(1/2 + i t) \zeta '(1/2 + i t)}
\end{equation}
for, say, $0 \leq t \leq 100$.  Let's also graph
\begin{equation}\label{E:GCDZeta}
\frac{ \zeta(1/2 + i t - 1)^2}{(1/2 + i t) \zeta'(1/2 + i t)}
\end{equation}
over the same interval.

\bigskip


It is also stated \cite[Theorem 1.1]{Bordelles} that, for $x > 1$ and every $\epsilon > 0$,
\[
T(x) = a_2 x^2+a_3 x^2 \log (x) + O_{\epsilon}(x^{1 + \theta + \epsilon})
\]
where $a_2$ and $a_3$ have the values given earlier in this section, and where $1/4 \leq \theta \leq 131/416$ is the exponent that appears in the Dirichlet divisor problem.

\section{Powerful (``power-full'') Numbers}\label{S:Powerful}

We now consider integers that, in some sense, are the opposite of squarefree integers.  Let $k$ be an integer $\ge 2$.  If the prime factorization of $n$ is
\[
n=p_1^{a_1} \cdot p_2^{a_2} \cdots p_r^{a_r},
\]
then we say that $n$ is $k$-full if every exponent $a_j \geq k$.  (We also define 1 to be $k$-full for all $k \ge 2$).  The characteristic function of the $k$-full numbers is
\[
 f_k(n)=
  \begin{cases}
    1  &\text{if $n$ is $k$-full,}\\
    0 &\text{otherwise.}
  \end{cases}
\]

If $k = 2$, then we call the 2-full numbers ``square-full'', and we have this Dirichlet series \cite[p. 33]{Ivic}:
\[
\sum _{n=1}^{\infty } \frac{f_2(n)}{n^s} = \frac{\zeta(2 s) \zeta(3 s)}{\zeta(6 s)},
\]
which converges for $s > 1/2$.

The summatory function, which counts the square-full numbers up to $x$, is
\[
T(x) = \sum _{n=1}^x f_2(n).
\]

The integral approximation to the summatory function $T_0(x)$ is
\begin{equation}\label{E:PowerfullIntegralApprox}
T_0(x) \simeq \frac{1}{2 \pi i}
   \int_{c-i T}^{c+i T} \frac{\zeta(2 s) \zeta(3 s)}{\zeta(6 s)} \frac{x^s}{s} \, ds.
\end{equation}

The integrand has poles of order 1 at $s = 0$, at $s = 1/3$, and at $s = 1/2$, where the residues are $-1/2$,
\[
x^{1/3} \frac{ \zeta \left(\frac{2}{3}\right)}{\zeta(2)},
\]
and
\[
x^{1/2} \frac{ \zeta \left(\frac{3}{2}\right)}{\zeta(3)},
\]
respectively.

From these, from the residues at each zero of $\zeta(6s)$, that is, at $s = \rho/6$, we get the following sum, which we hope approximates the summatory function $T_0(x)$:
\begin{equation}\label{E:PowerfullSumApprox}
T_0(x) \simeq
 a_1 + a_2 x^{1/3} + a_3 x^{1/2}
 + 2 \Re\left(\sum _{k=1}^N x^{\frac{\rho_k}{6}}
       \frac{ \zeta(\frac{\rho_k}{2}) \zeta(\frac{\rho_k}{3}) } {\rho_k \zeta'(\rho_k)} \right)
 + \sum_{k=1}^M x^{\frac{-2k}{6}}
       \frac{ \zeta(\frac{-2k}{2}) \zeta(\frac{-2k}{3}) } {(-2k) \zeta'(-2k)}
\end{equation}
where $a_1 = -1/2$,
\[
 a_2 = \frac{ \zeta \left(\frac{2}{3}\right)}{\zeta(2)} \simeq -1.48795,
\]
and
\[
a_3 = \frac{ \zeta \left(\frac{3}{2}\right)}{\zeta (3)} \simeq 2.17325.
\]

\begin{figure*}[ht]
  \mbox{\includegraphics[width=\picDblWidth]{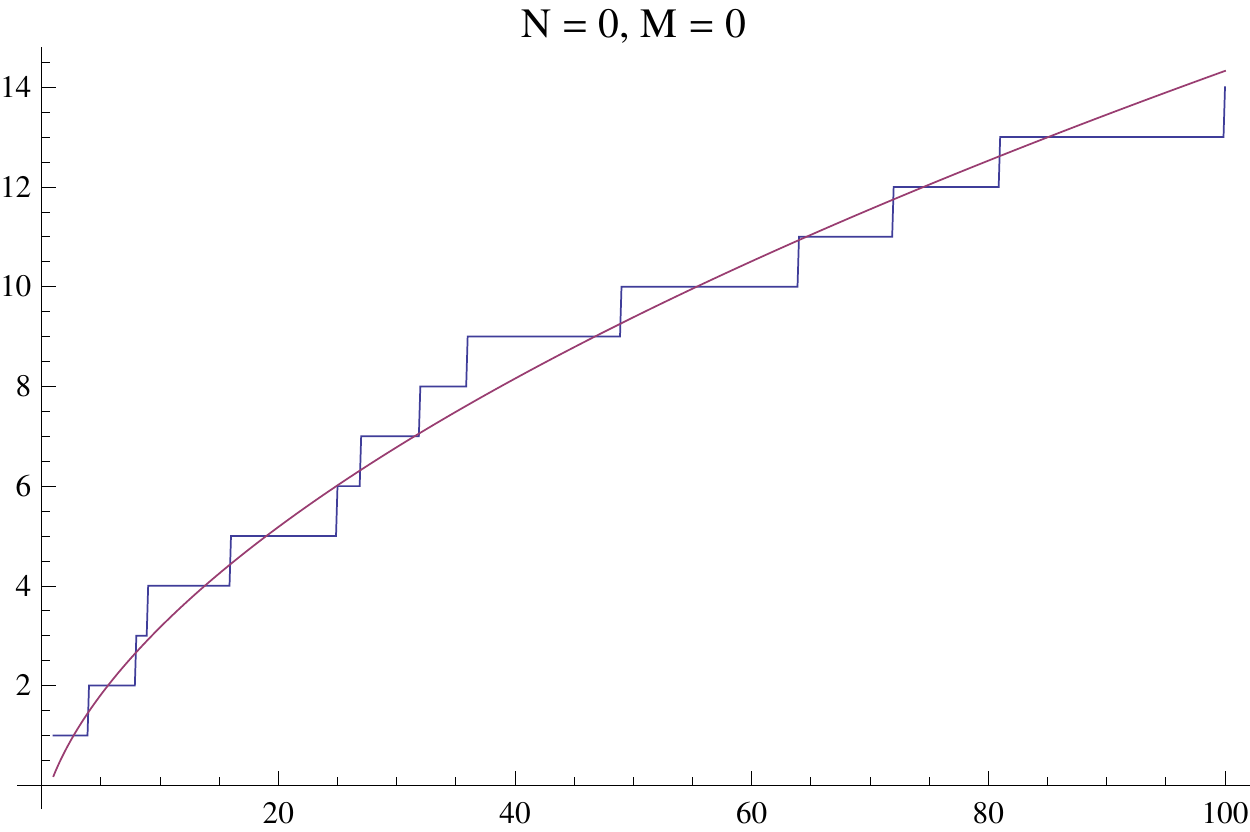}}
  \caption{Counting squarefull numbers using \eqref{E:PowerfullSumApprox} with $N = M = 0$}
  \label{fig:SquareFullSumGraph}
\end{figure*}

\begin{figure*}[ht]
  \centerline{
    \mbox{\includegraphics[width=\picDblWidth]{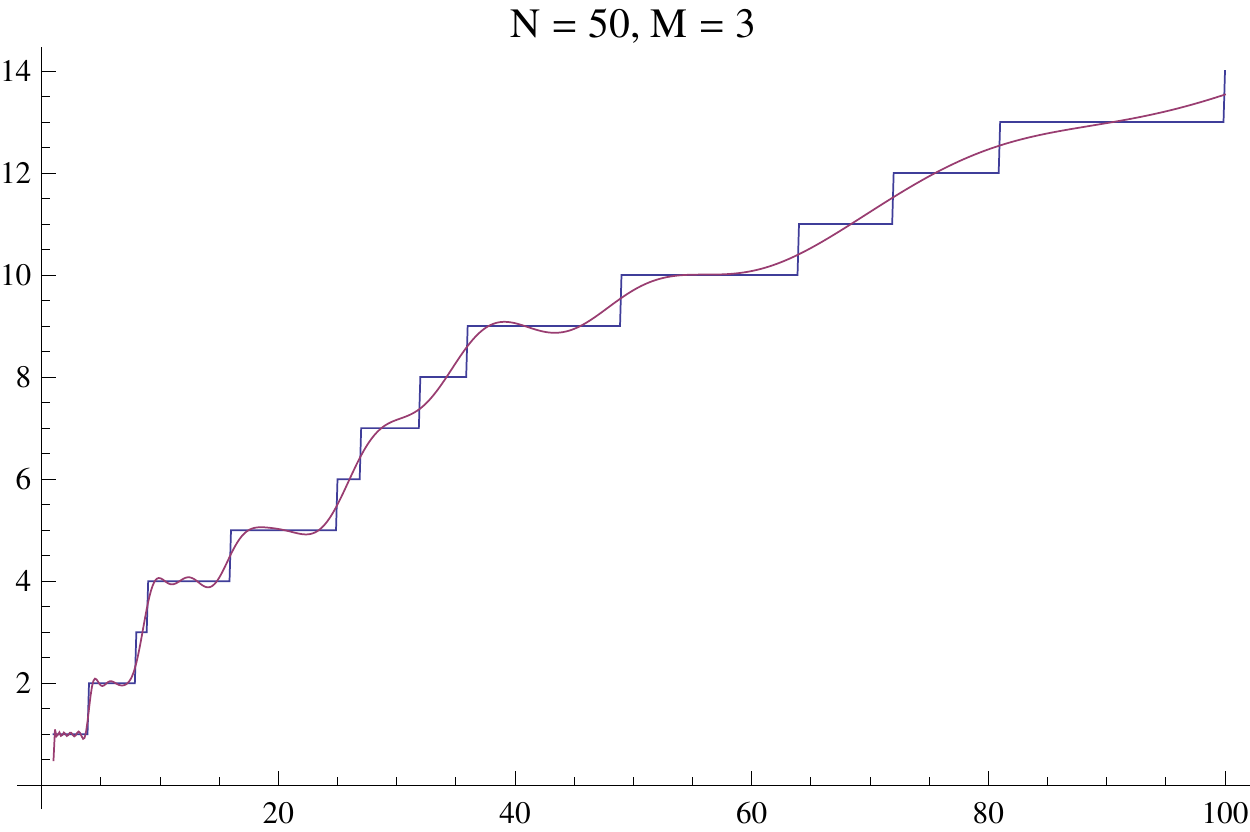}}
    \mbox{\includegraphics[width=\picDblWidth]{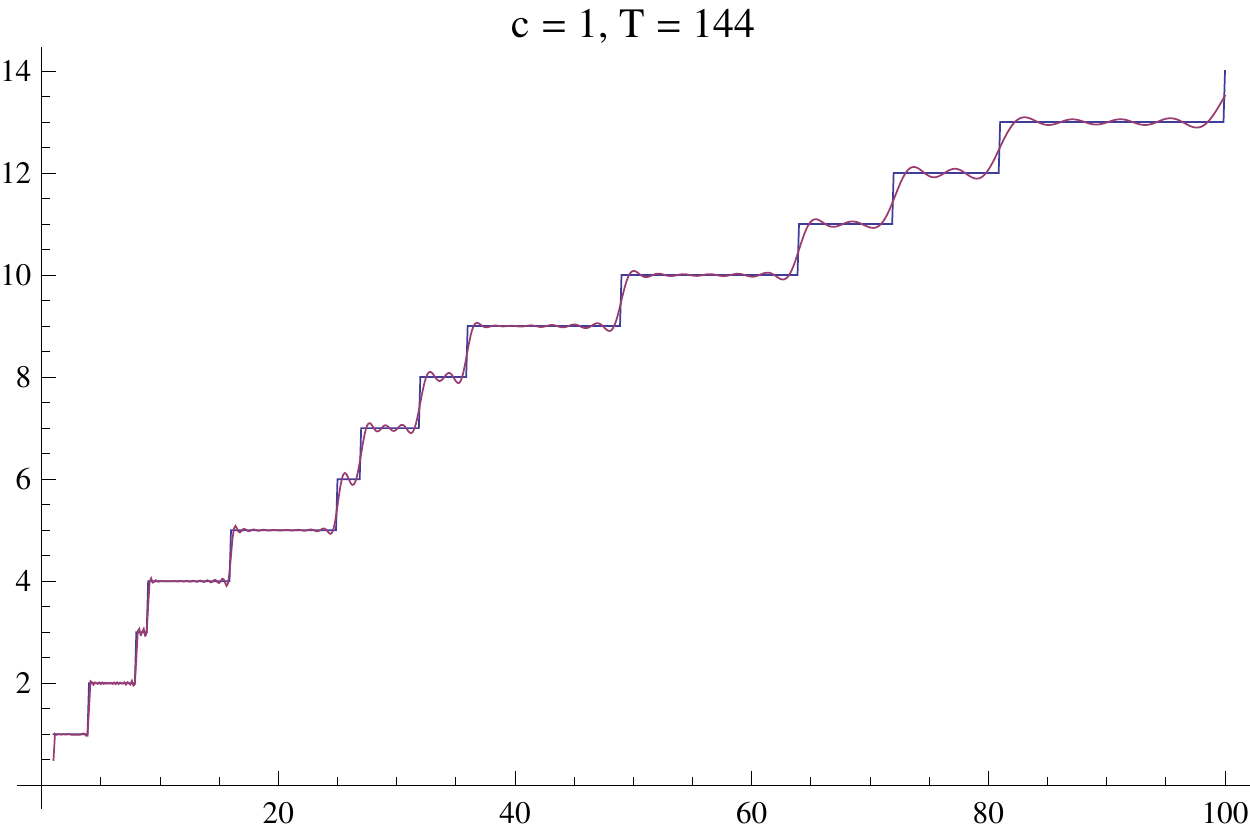}}
  }
  \caption{Counting squarefull numbers using \eqref{E:PowerfullSumApprox} (left) and \eqref{E:PowerfullIntegralApprox} (right)}
  \label{fig:SquareFullSumIntegralGraphs}
\end{figure*}

\begin{figure*}[ht]
  \centerline{
    \mbox{\includegraphics[width=\picDblWidth]{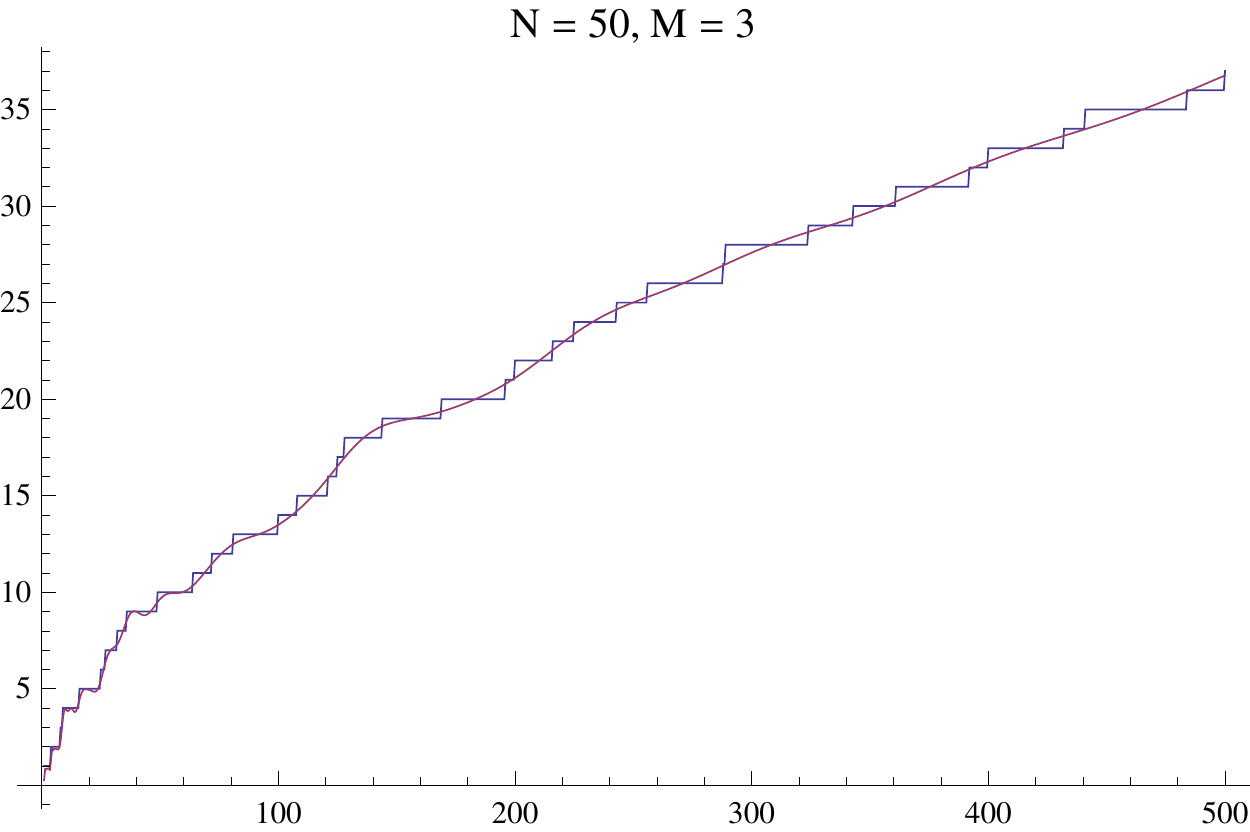}}
    \mbox{\includegraphics[width=\picDblWidth]{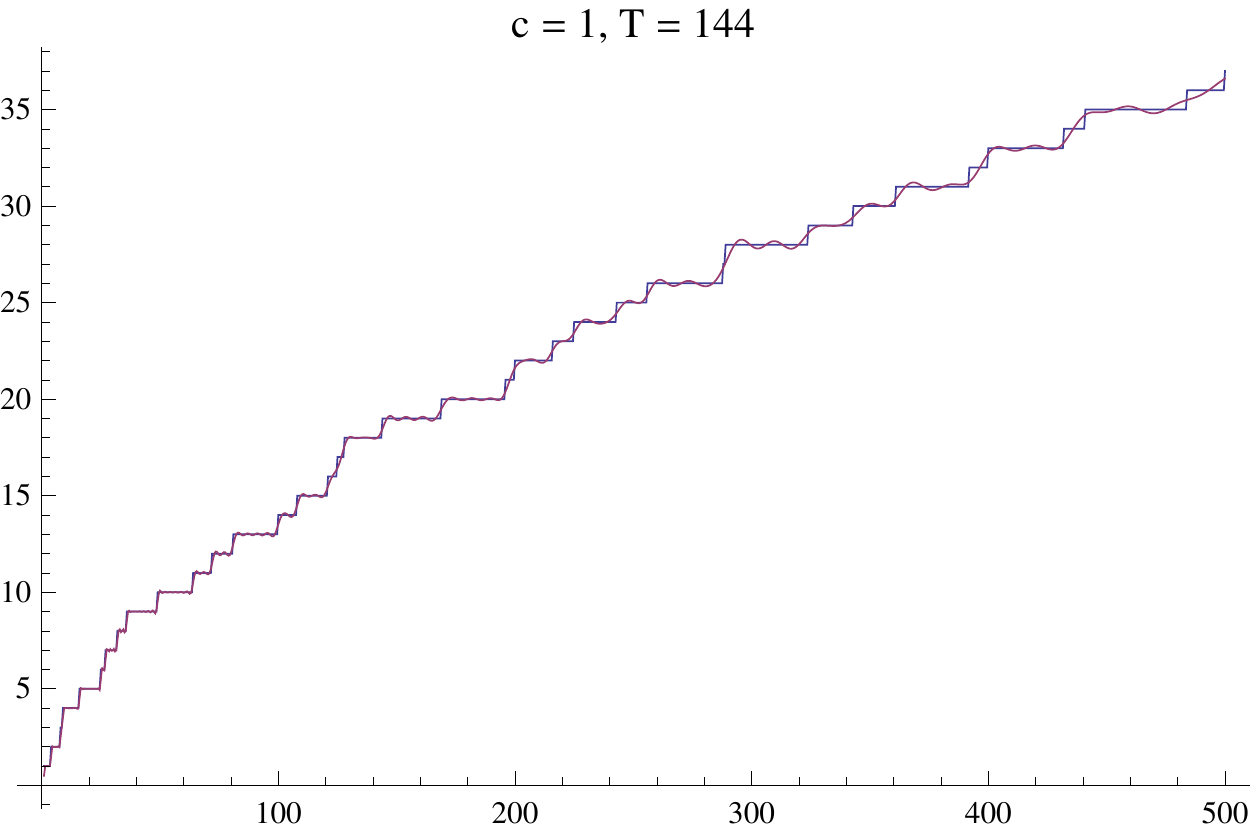}}
  }
  \caption{Counting squarefull numbers for $x \leq 500$}
  \label{fig:SquareFullSumIntegralGraphs2}
\end{figure*}

Figure ~\ref{fig:SquareFullSumGraph} shows the summatory function and the approximation that comes from setting $N = M = 0$ in equation \eqref{E:PowerfullSumApprox}, that is, using only the terms $a_1 + a_2 x^{1/3} + a_3 x^{1/2}$.
The left side of Figure ~\ref{fig:SquareFullSumIntegralGraphs} shows the approximation using $N = 50$ pairs of complex zeta zeros and $M = 3$ real zeros.  Without including the 3 real zeros, the approximation would be visibly, but slightly, too low.
The right side of Figure ~\ref{fig:SquareFullSumIntegralGraphs} shows the approximation based on the integral approximation \eqref{E:PowerfullIntegralApprox}.
Figure ~\ref{fig:SquareFullSumIntegralGraphs2} shows the same approximations, but for a larger range of $x$.

Bateman and Grosswald \cite{Bateman} proved that the number of squarefull numbers up to x is
\[
x^{1/3} \frac{\zeta(2/3)}{\zeta(2)} + x^{1/2} \frac{\zeta(3/2)}{\zeta(3)} + O( x^{1/6}e^{-a \omega(x)} )
\]
where $a$ is a positive constant, and
\[
\omega(x) = \log(x)^{4/7} \log(\log(x))^{-3/7}.
\]
The approximation in equation \eqref{E:PowerfullSumApprox} with $N = M = 0$ consists of the first two terms of this estimate, plus the constant $-1/2$, but we have not proved that the difference between our estimate and the summatory function is small.

\section{What's Really Going on Here?  An Analogy With Fourier Series}\label{S:Waves}

When you stop to think about it, Fourier series are pretty amazing.  By taking a linear combination of a countable number of ``wavy'' functions (sines and/or cosines), we can get a function that is a straight line over an entire interval that contains an uncountable number of points!

Fourier series do have one limitation, however: they are periodic.  That is, the sum of sines and cosines will represent a desired function over only a finite interval, say, between 0 and $\pi$, at which point the sum will repeat.

We will see that, like Fourier series, the sums we computed above are also linear combinations of ``wavy'' functions.  But in one sense, our sums are even more amazing than Fourier series. The step functions we are trying to approximate are quite irregular and are not periodic. Nevertheless, our sums often appear to approximate the step functions all the way out to infinity!

Consider the sums we've obtained over pairs of complex zeta zeros, for example, \eqref{E:psiZetaZeroSumProved}, \eqref{E:MertensApproximation}, \eqref{E:QxSum}, \eqref{E:PhiSum}, \eqref{E:LxSum}, or \eqref{E:TwoNuSumApprox}.

A typical term in one of these sums has the form $c_k x^{\rho_k}$ where $c_k$ is a complex coefficient, say $c_k = a_k + i b_k$.  Let's write the $k^{th}$ complex zeta zero as $\rho_k = 1/2 + i t_k$.  Then

\begin{equation*}
x^{\rho_k} = x^{1/2 + i t_k} = x^{1/2} \cdot x^{i t_k} = x^{1/2} \cdot \exp(i t_k \log(x))
\end{equation*}
\begin{equation*}
= x^{1/2} \cdot [ \cos(t_k \log(x)) + i \sin(t_k \log(x)) ].
\end{equation*}

Then
\begin{equation*}
c_k x^{\rho_k} = x^{1/2} \cdot (a_k + i b_k) \cdot [ \cos(t_k \log(x)) + i \sin(t_k \log(x)) ].
\end{equation*}

We will want the real part of this product, which is:
\begin{equation*}
\Re( c_k x^{\rho_k} ) = x^{1/2} \cdot [a_k \cos(t_k \log(x)) - b_k \sin(t_k \log(x))]
\end{equation*}
\begin{equation*}
=  a_k \cdot x^{1/2} \cos(t_k \log(x)) - b_k \cdot x^{1/2} \sin(t_k \log(x)).
\end{equation*}

To make this more concrete, let's consider the first such term, with $\rho_1 \simeq 1/2 + 14.135 i$, and, for comparison, the tenth such term, with $\rho_{10} \simeq 1/2 + 49.774 i$, and let's take the coefficients $c_1$ and $c_{10}$ to be 1.  Then the real parts of the products are
\[
x^{1/2} \cos(14.135 \log(x))
\]
and
\[
x^{1/2} \cos(49.774 \log(x)).
\]

Here are the graphs of these two functions.  Both curves are bounded by the envelopes $\pm x^{1/2}$.
\begin{figure*}[ht]
  \centerline{
    \mbox{\includegraphics[width=\picDblWidth]{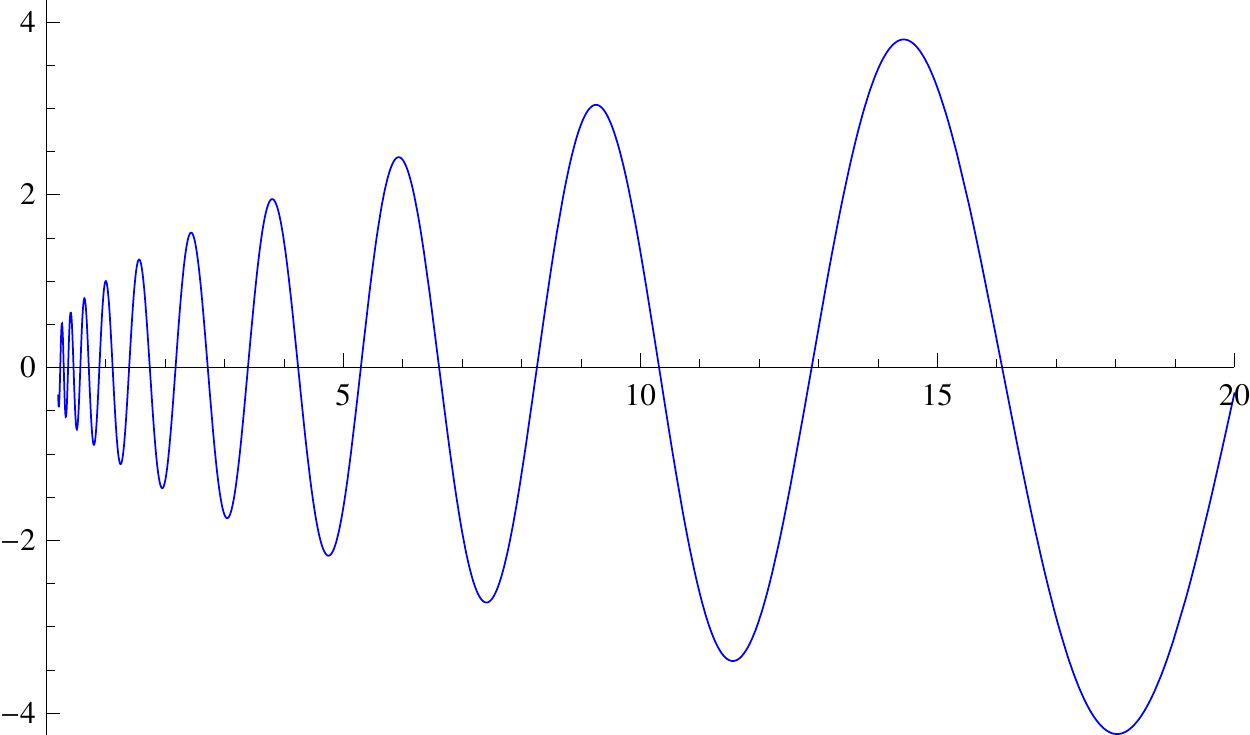}}
    \mbox{\includegraphics[width=\picDblWidth]{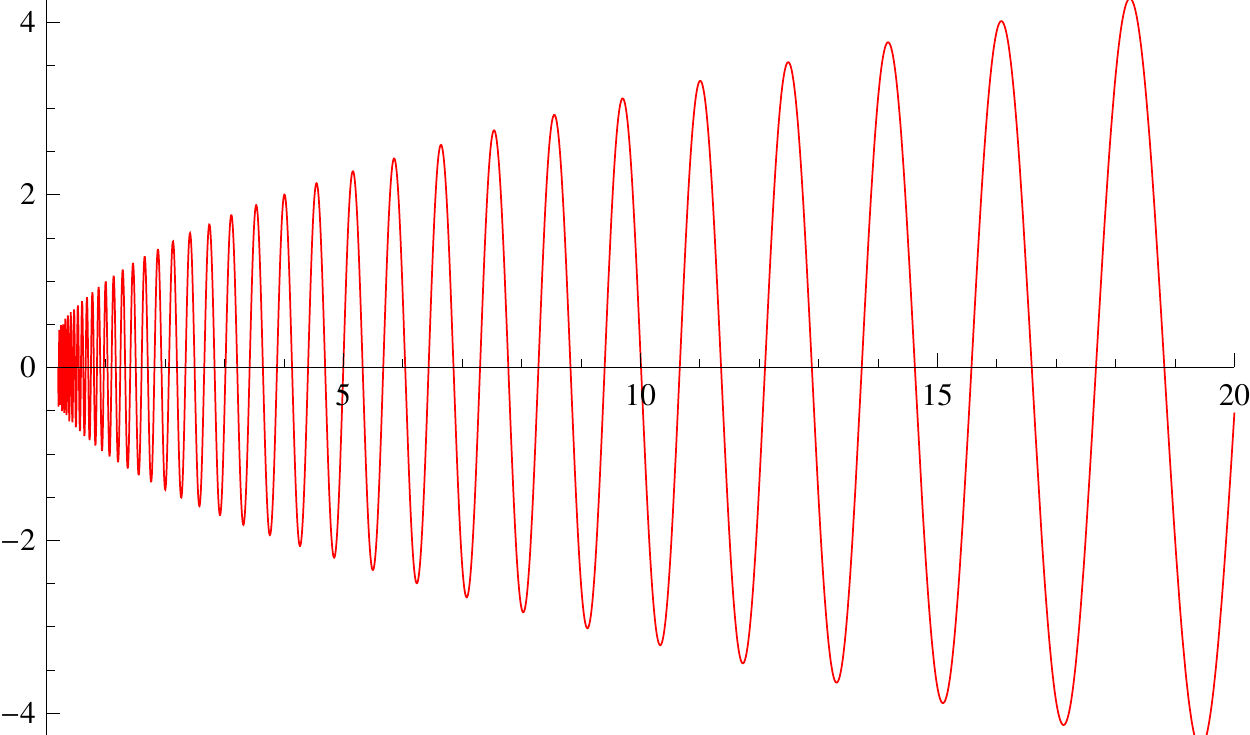}}
  }
  \caption{$x^{1/2} \cos(14.135 \log(x))$ and $x^{1/2} \cos(49.774 \log(x))$}
  \label{fig:WavyGraphs1and10}
\end{figure*}

Again, the approximating sums are linear combinations of functions that look like these.  We can think of the functions
\[
x^{1/2} \cos(t_k \log(x))
\]
and
\[
x^{1/2} \sin(t_k \log(x))
\]
as \emph{basis} functions for our approximating sums, in the same sense that $\sin(n x)$ and $\cos(n x)$ are the basis functions for Fourier series.

This analogy also helps explain why there appears to be a Gibbs-type phenomenon in some of the graphs above.

\section{Counting Primes in Arithmetic Progressions}\label{S:CountingPrimesInAP}



In section \ref{S:rVonM}, we saw how to use the first $N$ pairs of zeros of the Riemann zeta function to count the primes up to $x$.  Here, we will make a first pass at using zeros of Dirichlet $L$-functions to count primes in arithmetic progressions.

Let $q$ and $a$ be integers with $q  > 0$. If $q$ and $a$ are relatively prime, then the arithmetic progression $qn + a$ where $n = 0, 1, 2, 3, \dots$ contains an infinite number of primes.  This is Dirichlet's Theorem \cite[p. 76]{LeVeque}.  The number of such primes that are less than or equal to $x$ is a function of $q$, $a$, and $x$ that is usually denoted by $\pi_{q,a}(x)$.  For an introduction to characters and Dirichlet $L$-functions, see \cite[Chapter 6]{LeVeque}.

Equations \eqref{E:arithProg1} through \eqref{E:arithProg4} below come from equations (6) through (8) in \cite{Martin}.  See \cite{RubinsteinAndSarnak} for an estimate of the error term for these equations.

\begin{equation}\label{E:arithProg1}
\pi_{q,a}(x) =
 \frac{\pi(x)}{\phi(q)} + E(x, q, a) \frac{x^{1/2}}{\phi(q) \log(x)}
\end{equation}
where
\begin{equation}\label{E:arithProg2}
E(x, q, a) = -c(q, a) - \sum_{\chi \neq \chi_0} \chi(a)^{*} E(x, \chi),
\end{equation}

\begin{equation}\label{E:arithProg3}
c(q, a) = -1 + \#(\text{square roots of }a, \text{ mod }q),
\end{equation}
and
\begin{equation}\label{E:arithProg4}
E(x, \chi) =  \sum_{k=1}^{2 N} \frac{x^{i \rho_k}}{1/2 + i \rho_k}.
\end{equation}

In equation \eqref{E:arithProg1}, the main term is $\frac{\pi(x)}{\phi(q)}$.  Except for $c(q, a)$, the rest of equation \eqref{E:arithProg1} depends on zeros of $L$-functions.

In equation \eqref{E:arithProg2}, the asterisk denotes the complex conjugate, and we sum over all characters $\chi \pmod{q}$, except for the principal character, $\chi_0$.

In equation \eqref{E:arithProg3}: if $a$ has no square roots $\pmod{q}$, that is, if $a$ is a quadratic nonresidue $\pmod{q}$, then $c(q,a) = -1$.  Otherwise, $a$ must have at least two square roots $\pmod{q}$, which makes $c(q, a) \geq 1$.  This causes a bias in ``prime number races'': $\pi_{q,a}(x)$ often exceeds $\pi_{q,b}(x)$ if $b$ is a square $\pmod{q}$ and $a$ is not.  For example, $\pi_{3,2}(x) > \pi_{3,1}(x)$ for all $x < 608981813029$.  See \cite{GranvilleAndMartin} for more details on this fascinating topic.

In equation \eqref{E:arithProg4}, for the character $\chi \pmod{q}$, we sum over complex zeros of the corresponding Dirichlet $L$-function.  We include $2N$ such zeros in the sum: namely, the $N$ zeros having the smallest positive imaginary part, and the $N$ having the least negative imaginary part.  Note that these zeros do not necessarily occur in conjugate pairs.  $\rho_k$ is the imaginary part of the $k^{th}$ zero of the $L$-function that corresponds to the character $\chi$.


\begin{figure*}[ht]
  \mbox{\includegraphics[width=4.20in]{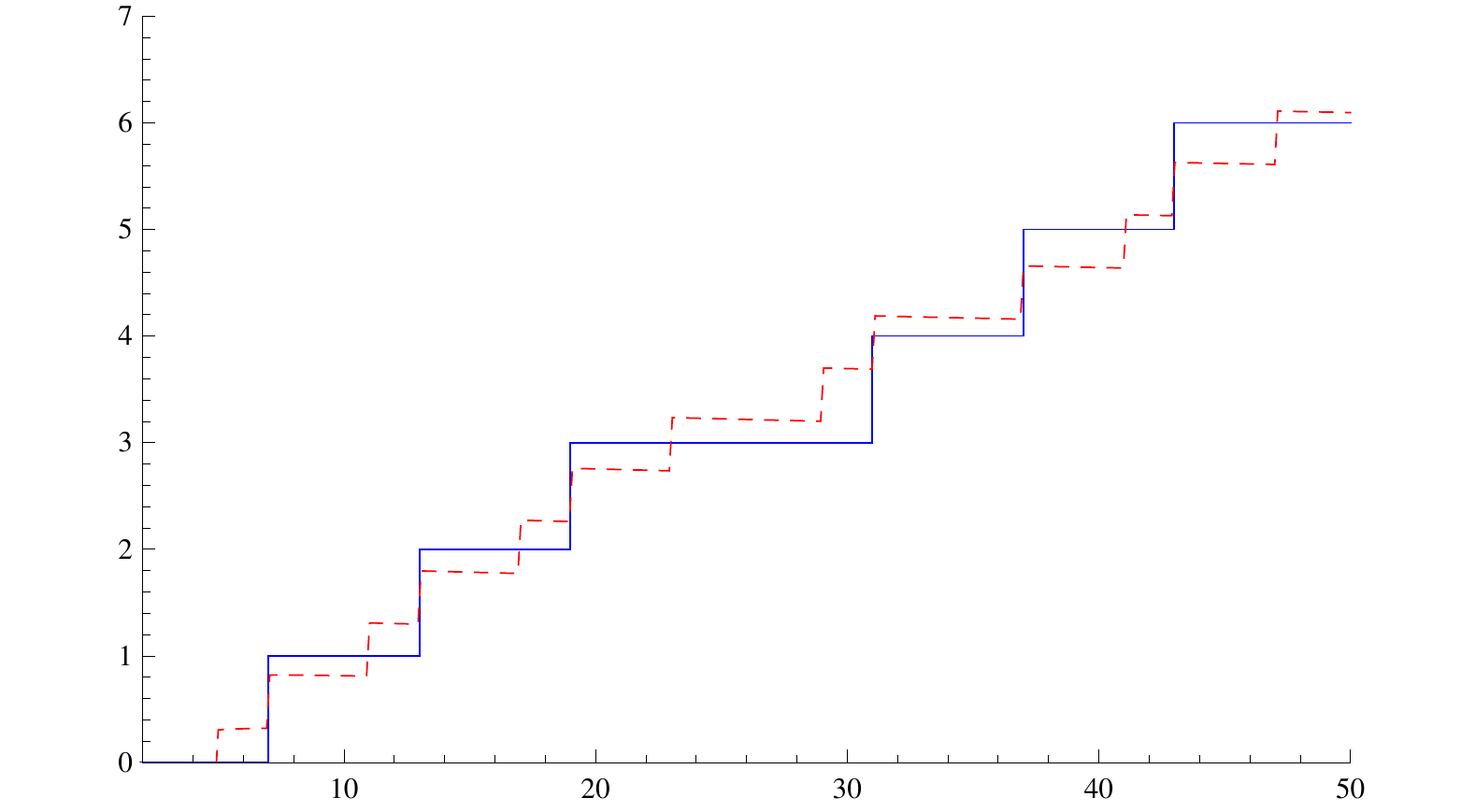}}
  \caption{$\pi_{3,1}(x)$ (solid), and an approximation \eqref{E:arithProg5} using no zeros of $L$-functions (dashed)}
  \label{fig:Pix31N0}
\end{figure*}


The solid curve in figure ~\ref{fig:Pix31N0} is $\pi_{3,1}(x)$.  The dashed curve is an approximation to $\pi_{3,1}(x)$ using equations \eqref{E:arithProg1} -- \eqref{E:arithProg4} with no Dirichlet $L$ zeros taken into account.  That is, the dashed curve is the approximation

\begin{equation}\label{E:arithProg5}
\pi_{q,a}(x) \simeq
 \frac{\pi(x)}{\phi(q)} - \frac{c(q, a)}{\phi(q)} \frac{x^{1/2}}{\log(x)}.
\end{equation}

(This approximation looks like a step function in spite of the second term on the right, because the second term is relatively small and doesn't vary much for $x$ between 2 and 50).  Note that this approximation uses $\pi(x)$, so the dashed curve jumps up at \emph{every} prime.

\begin{figure*}[ht]
  \mbox{\includegraphics[width=4.20in]{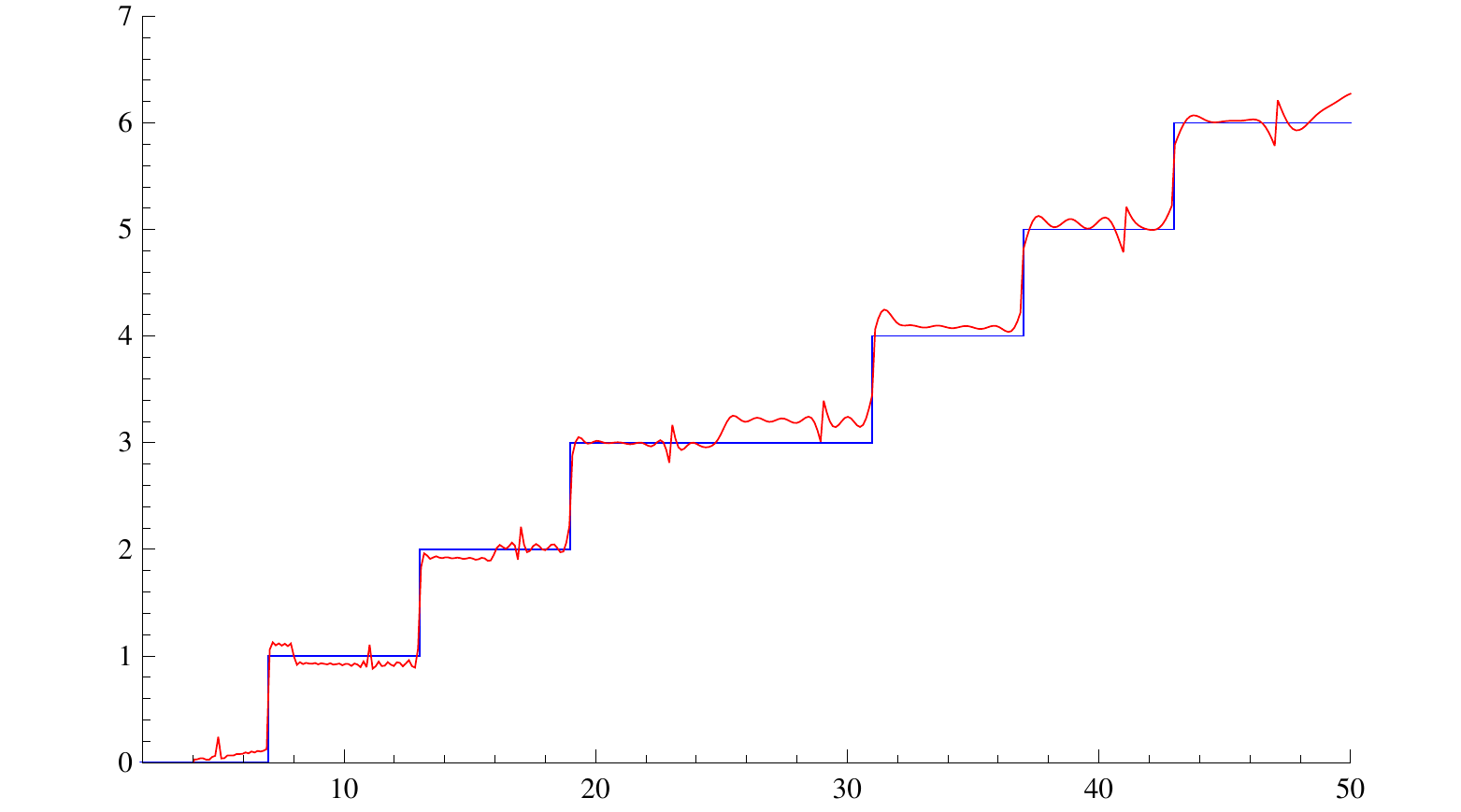}}
  \caption{$\pi_{3,1}(x)$ and an approximation using 200 zeros of $L$-functions}
  \label{fig:Pix31N100}
\end{figure*}

Figure ~\ref{fig:Pix31N100} shows $\pi_{3,1}(x)$ and an approximation using the first $2N = 200$ zeros of $L$-functions in equations \eqref{E:arithProg1} -- \eqref{E:arithProg4}.  Note the prominent glitches at $x$ = 17, 23, 29, 41, and 47.  These are places where $\pi(x)$ jumps, but $\pi_{3,1}(x)$ does not.

Warning: there is a bit of empirical fakery in this graph!  The pair of arithmetic progressions $\{3n + 1, 3n + 2\}$ omits one prime, namely, 3.  The value for $\pi_{3,1}(x)$ that comes out of equation \eqref{E:arithProg1} is too large by about $1/2$.  The same is true for $\pi_{3,2}(x)$.  So, to account for the missing prime among the two progressions $\{3n + 1, 3n + 2\}$, the values shown in figure ~\ref{fig:Pix31N100} are $1/2$ less than the values we get from equation \eqref{E:arithProg1}.  In spite of this legerdemain, the sums over zeros of $L$-functions provide a decent approximation to the counts of primes in these arithmetic progressions.

\begin{figure*}[ht]
  \mbox{\includegraphics[width=4.20in]{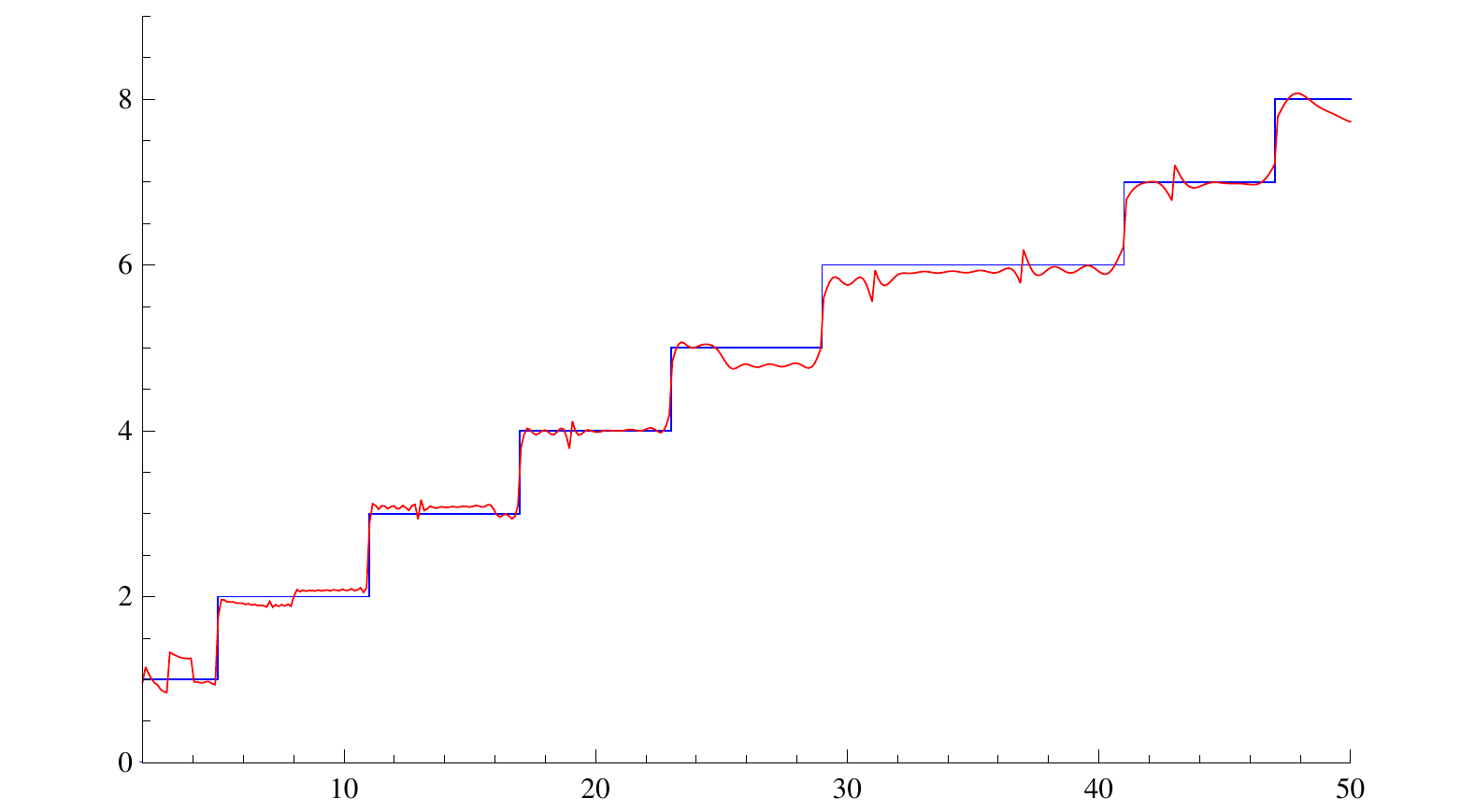}}
  \caption{$\pi_{3,2}(x)$ and an approximation using 200 zeros of $L$-functions}
  \label{fig:Pix32N100}
\end{figure*}

\begin{figure*}[ht]
  \centerline{
    \mbox{\includegraphics[width=\picDblWidth]{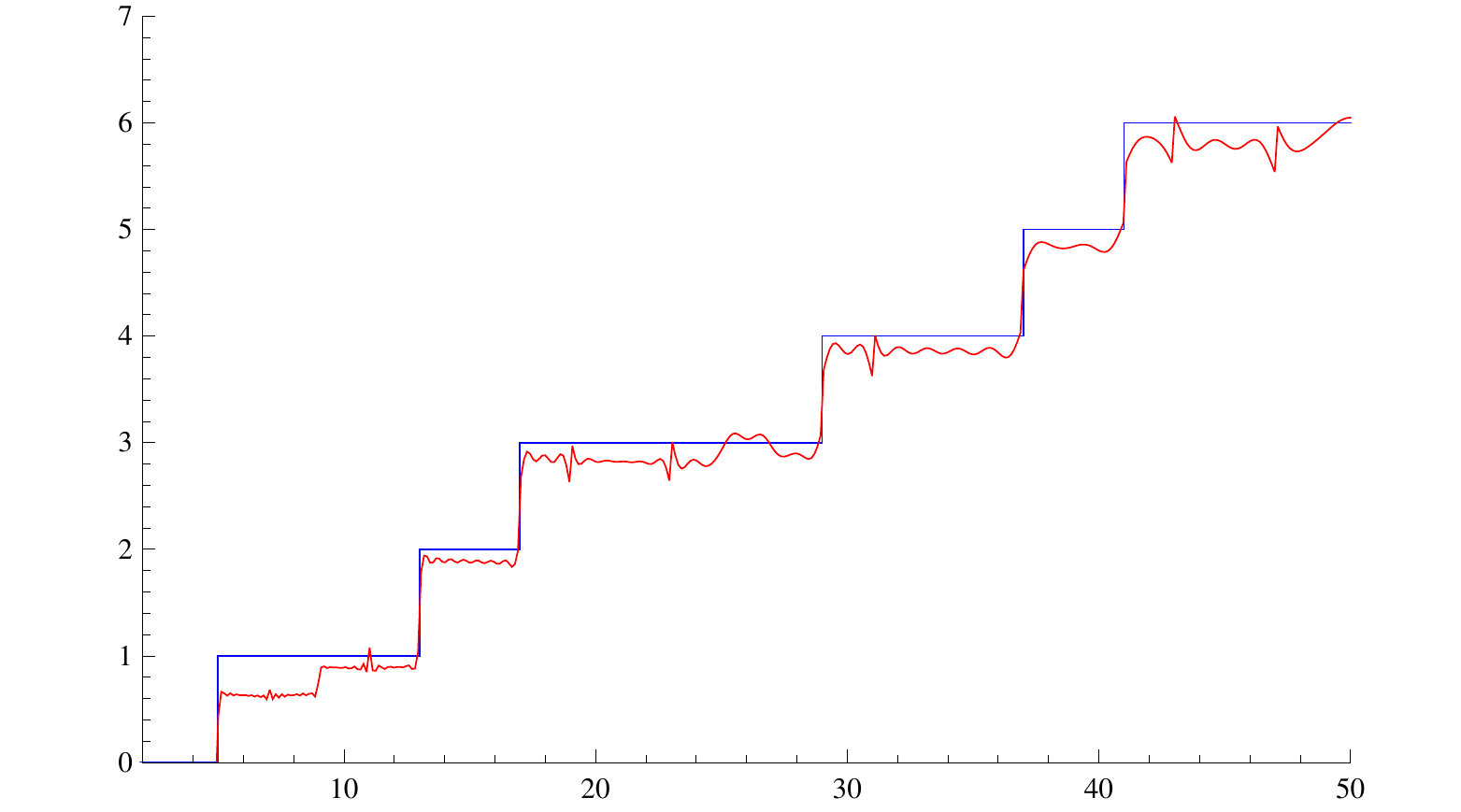}}
    \mbox{\includegraphics[width=\picDblWidth]{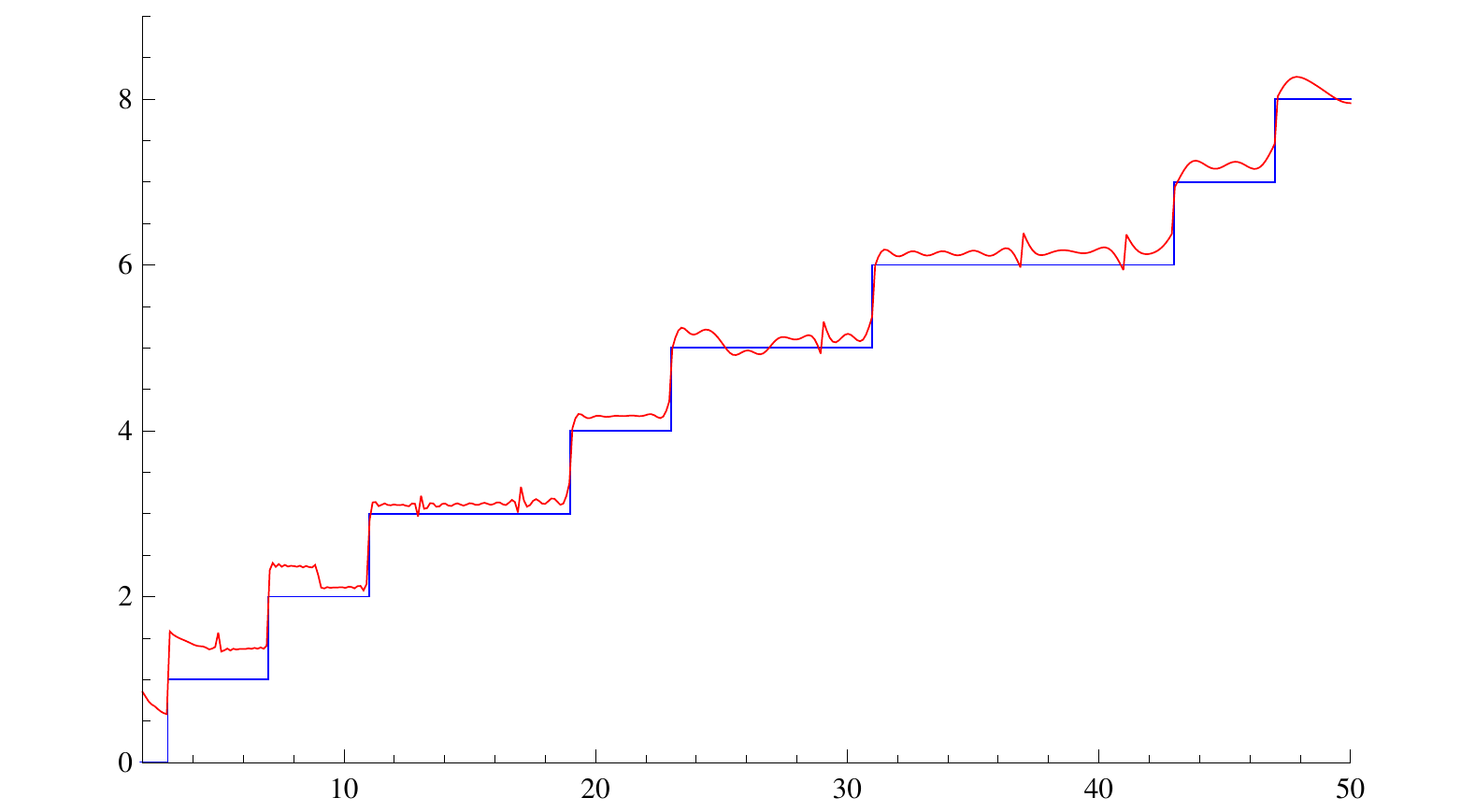}}
  }
  \caption{$\pi_{4,1}(x)$ (left) and $\pi_{4,3}(x)$ (right)}
  \label{fig:Pix4aN100}
\end{figure*}

Figure ~\ref{fig:Pix4aN100} shows $\pi_{4,1}(x)$ and $\pi_{4,3}(x)$ and their approximations using $2N = 200$ zeros of $L$-functions.  In order to account for the fact that this pair of progressions omits the prime 2, the values shown in the graphs are $1/2$ less than the values we get from equation \eqref{E:arithProg1}.

\begin{figure*}[ht]
  \centerline{
    \mbox{\includegraphics[width=\picDblWidth]{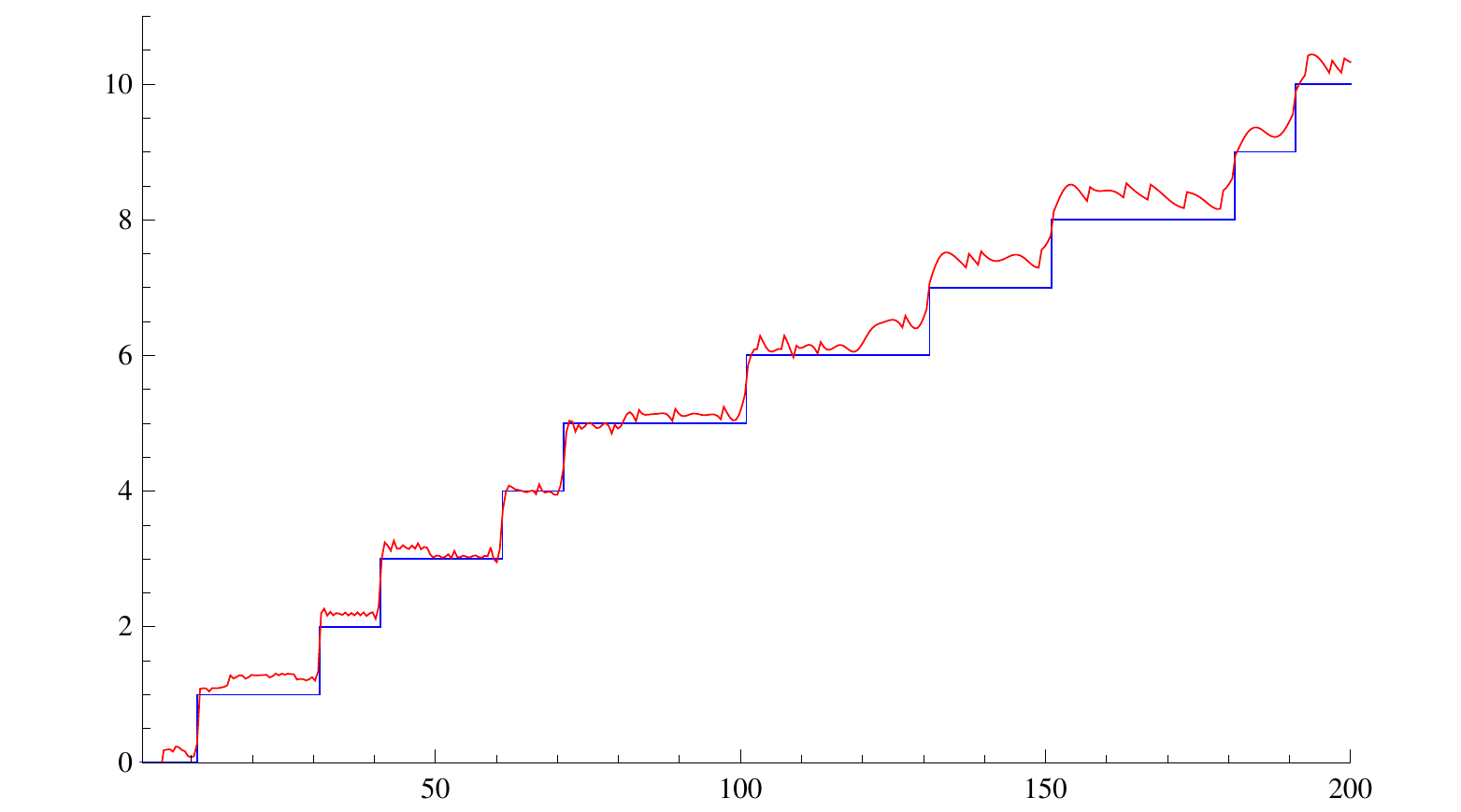}}
    \mbox{\includegraphics[width=\picDblWidth]{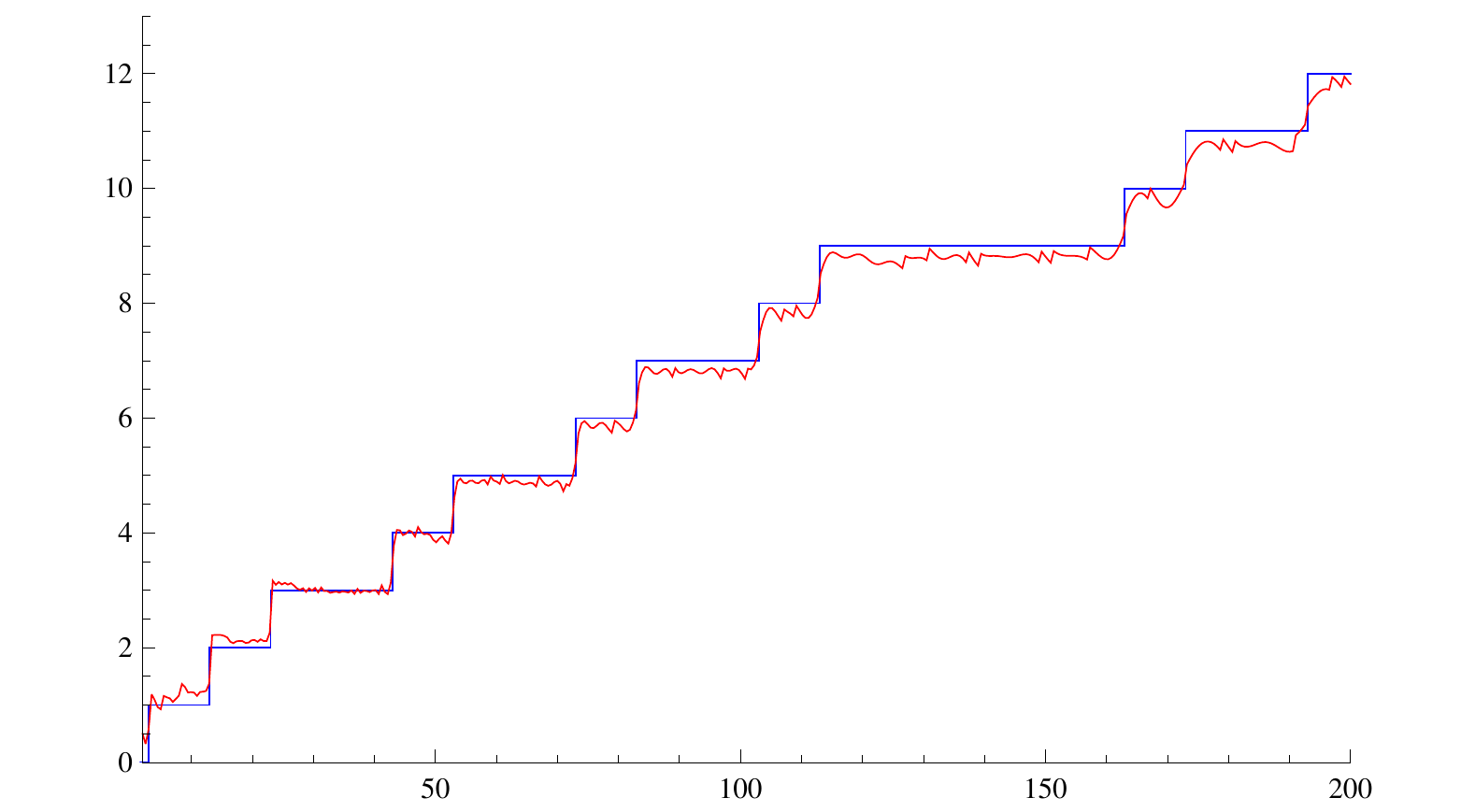}}
  }
  \caption{$\pi_{10,1}(x)$ (left) and $\pi_{10,3}(x)$ (right)}
  \label{fig:Pix10a13N100}
\end{figure*}

\begin{figure*}[ht]
  \centerline{
    \mbox{\includegraphics[width=\picDblWidth]{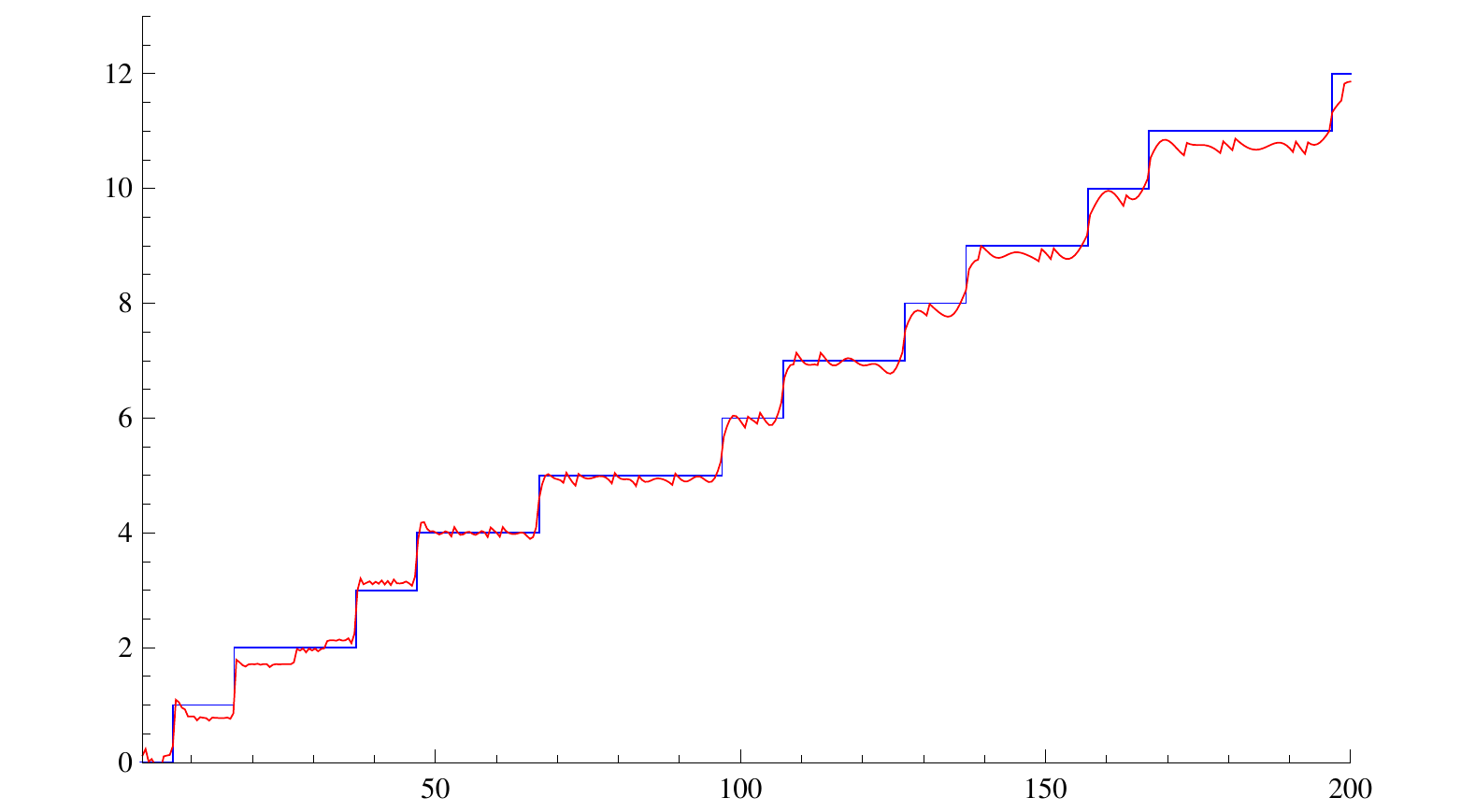}}
    \mbox{\includegraphics[width=\picDblWidth]{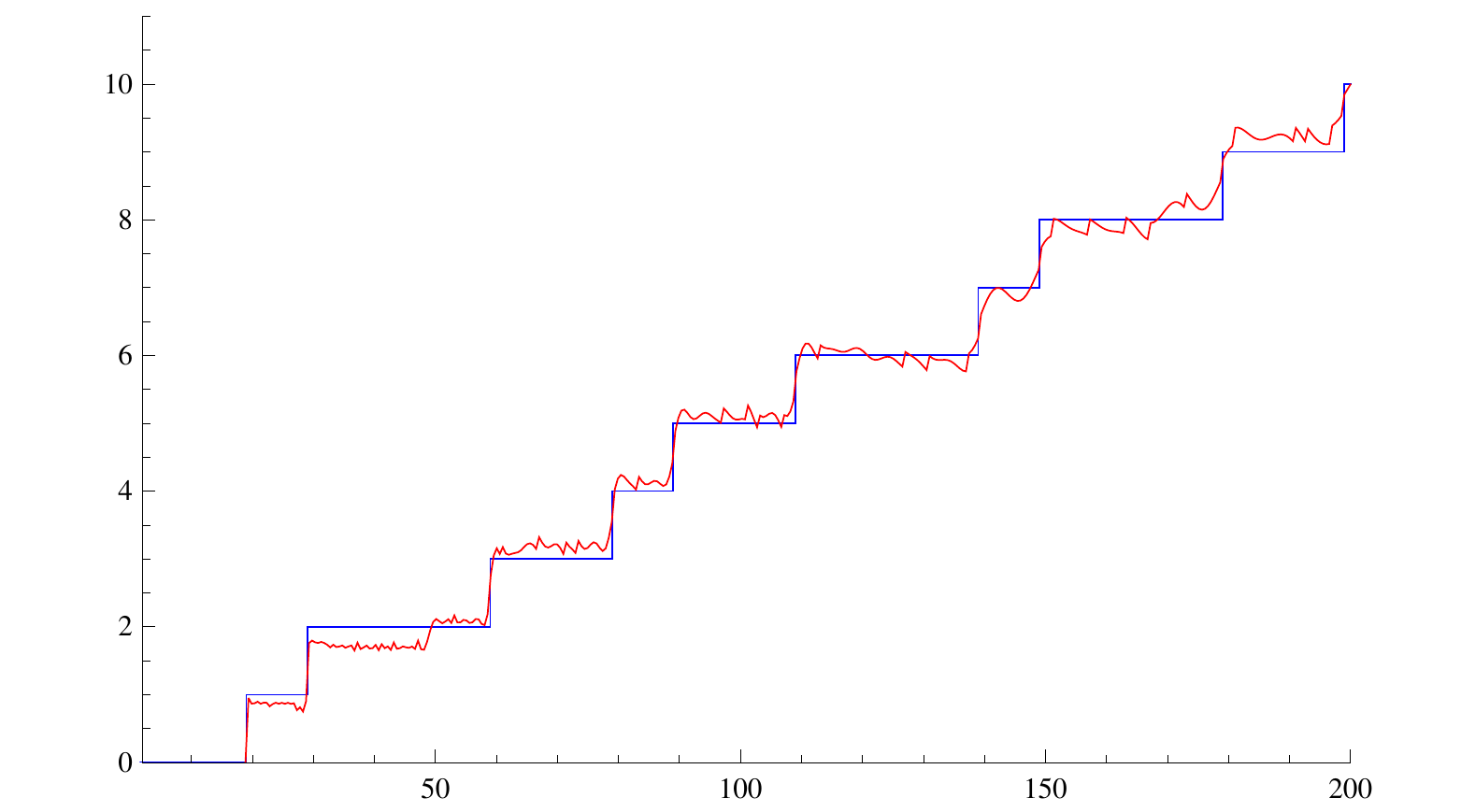}}
  }
  \caption{$\pi_{10,7}(x)$ (left) and $\pi_{10,9}(x)$ (right)}
  \label{fig:Pix10a79N100}
\end{figure*}

In figures ~\ref{fig:Pix10a13N100} and ~\ref{fig:Pix10a79N100}, we use $2N = 200$ zeros of $L$-functions to approximate $\pi_{10,1}(x)$, $\pi_{10,3}(x)$, $\pi_{10,7}(x)$, and $\pi_{10,9}(x)$, the four progressions $\mod 10$, up to $x = 200$.  These four progressions, taken together, omit two primes, namely, 2 and 5.  For the progressions with $a = 1, 3, 7, 9$, we adjust the values by, respectively, $-1/4$, $0$, $-3/4$, and $-1$.  These are empirical guesses.  It is not clear why these values seem to work better than, say, $-1/2$ for all of the progressions.

\section{Questions For Further Research}\label{S:ResearchQuestions}

The answers to the following questions are not known to this author.

Can the Riemann-von Mangoldt formula (equations \eqref{E:r-vonm} - \eqref{E:r-vonm2}), which counts primes using a sum over zeta zeros, be obtained by applying Perron's formula?

In section ~\ref{S:TauSquared} that considers the sum of $\tau(n^2)$, when we use 200 zeros in the sum, the graph seems worse in some ways than when we use only 100 zeros.
The same thing happens in section ~\ref{S:SigmaTau} with $\sigma(n) \tau(n)$, in section ~\ref{S:LambdaSigmaTau} with $\lambda(n) \sigma(n) \tau(n)$, and in section ~\ref{S:TallyGCD} where we tally GCD's.

Why does this happen?  Can we predict in which cases will this occur?  Would the approximation get worse and worse if we used even more zeros?

In section ~\ref{S:LambdaTwoNu} on $\lambda(n) 2^{\nu(n)}$, section ~\ref{S:LambdaTauOfSquares} on $\lambda(n) \tau(n^2)$ and section ~\ref{S:LambdaTauSquared} on $\lambda(n) \tau(n)^2$, and  ~\ref{S:LambdaSigmaTau} on $\lambda(n) \sigma(n) \tau(n)$, if we extend the rectangle to the left to include some real zeros of zeta, the sum of the residues appears to diverge.  Can we prove this?  For which Dirichlet series will this happen?

In section ~\ref{S:Waves} on the analogy with Fourier series: is there anything more we can say that's useful or interesting?

In section ~\ref{S:CountingPrimesInAP} on counting primes in arithmetic progressions, why, for $q = 10$, do we need different offsets for $a$ = 1, 3, 7, and 9?

In several cases, there are provable theorems which say (more or less) that a sum over zeta zeros does, in fact, approach the respective arithmetical function.  These cases are:
\begin{itemize}
 \item equations \eqref{E:r-vonm} - \eqref{E:r-vonm2} (counting primes),
 \item equation \eqref{E:psiZetaZeroSumProved} (computing $\psi(x)$),
 \item equation \eqref{E:provedSumM0} (computing the Mertens function),
 \item equation \eqref{E:CountingSquareFree} (counting squarefree integers) (a similar result holds for cubefree integers),
 \item equation \eqref{E:WiertelakEq} (counting squarefree divisors).
\end{itemize}
Can similar theorems be proven for the other arithmetical functions that are considered here?

\bigskip

\noindent\emph{Email: rjbaillie@frii.com}


\bigskip


\begin{thebibliography}{99}

\bibitem{BartzMx2}
  K. M. Bartz,
  \emph{On some complex explicit formulae connected with the M\"{o}bius function, II},
  Acta Arithmetica, vol. 57 (1991), pp. 295--305.

\bibitem{BartzSquareFree}
  K. M. Bartz,
  \emph{On some connections between zeta-zeros and square-free integers},
  Monatshefte f\"{u}r Mathematik, vol 114 no. 1, March 1992, 15--34.
  Available at\\
\url{http://www.digizeitschriften.de/index.php?id=resolveppn&PPN=GDZPPN002488116}\\
  (note that pages 18 and 19 are missing from the online version).

\bibitem{BartzCubeFree}
  K. M. Bartz,
  \emph{On some connections between zeta-zeros and 3-free integers},
  Functiones et Approximatio, 28, (2000), 233--236.  Available at
  \url{http://www.staff.amu.edu.pl/~fa/XXVIII/fa-28-1-233.pdf}

\bibitem{Bateman}
  P. T. Bateman and E. Grosswald,
  \emph{On a Theorem of Erd\"{o}s and Szekeres},
  Illinois Journal of Mathematics, vol. 2, issue 1, pp. 88--98, (1958).

\bibitem{Bordelles}  
  O. Bordell\`{e}s,
  \emph{A Note on the Average Order of the gcd-sum Function},
  Journal of Integer Sequences, vol. 10, (2007), Article 07.3.3.  Available at\\
\url{http://www.cs.uwaterloo.ca/journals/JIS/VOL10/Bordelles/bordelles90.pdf}

\bibitem{Conrey}
  J. B. Conrey,
  \emph{More than two fifths of the zeros of the Riemann zeta function are on the critical line},
  J. Reine angew. Math. (1989), 399: 1--16.

\bibitem{Davenport}
  H. Davenport,
  \emph{Multiplicative Number Theory},
  Springer: New York, 3rd edition, 2000.

\bibitem{Finch}
  S. Finch,
  \emph{Mathematical Constants},
  Cambridge, 2003.
  See also the \emph{Errata and Addenda to Mathematical Constants},
  at \url{http://algo.inria.fr/csolve/erradd.pdf}


\bibitem{Gould}
  H. W. Gould and T. Shonhiwa, 
  \emph{A Catalog of Interesting Dirichlet Series},
  available at\\ 
\url{http://www.math-cs.ucmo.edu/~mjms/2008.1/hwgould.pdf}

\bibitem{Gourdon}
  X. Gourdon,
  \emph{The $10^{13}$ First Zeros of the Riemann Zeta Function, and Zeros Computation at Very Large Height},
  available at\\
\url{http://numbers.computation.free.fr/Constants/Miscellaneous/zetazeros1e13-1e24.pdf}

\bibitem{GranvilleAndMartin}
  A. Granville and G. Martin,
  \emph{Prime Number Races}
  American Mathematical Monthly, vol. 113 no. 1, 2006, pp. 1–-33,
  available at
 \url{http://mathdl.maa.org/images/upload_library/22/Ford/granville1.pdf}

\bibitem {HardyAndLittlewood}
  G. H. Hardy and J. E. Littlewood, (1921),
  \emph{The zeros of Riemann's zeta-function on the critical line},
  Math. Z. (1921) 10: 283--317

\bibitem{HardyAndWright}
  G. H. Hardy and E. M. Wright,
  \emph{An Introduction to the Theory of Numbers},
  Oxford, 4th edition, 1960.

\bibitem{Ingham}
  A. E. Ingham,
  \emph{The Distribution of Prime Numbers}.
  Cambridge University Press, 1990.

\bibitem{Ivic}
  A. Ivi\'{c},
  \emph{The Riemann Zeta-function},
  New York: Wiley, 1985.

\bibitem{Kaplan}
  W. Kaplan,
  \emph{Advanced Calculus},
  Addison-Wesley, 1959.


\bibitem{KotnikandteRiele}
  T. Kotnik and H. J. J. te Riele,
  \emph{The Mertens Conjecture Revisited},
  Lecture Notes in Computer Science 4076, Proceedings of the 7th Algorithmic Number Theory Symposium,
  Berlin, 2006, pp. 156–-167.

\bibitem{LeVeque}
  W. J. LeVeque,
  \emph{Topics in Number Theory},
  vol. 2, Reading, MA: Addison–Wesley, 1961.

\bibitem{Martin}
  G. Martin,
  \emph{Asymmetries in the Shanks-Renyi Prime Number Race},
  available at\\
\url{http://front.math.ucdavis.edu/0010.5086}

\bibitem{McCarthy}
  P. J. McCarthy,
  \emph{Introduction to Arithmetical Functions},
  New York: Springer-Verlag, 1986.

\bibitem {MontgomeryAndVaughan}
  H. L. Montgomery and R. C. Vaughan,
  \emph{Multiplicative Number Theory: I. Classical Theory},
  Cambridge: Cambridge University Press, 2007

\bibitem {OdlyzkoAndteRiele}
  A. M. Odlyzko and H. J. J. te Riele,
  \emph{Disproof of the Mertens Conjecture},
  J. reine angew. Math. 357, 1985, pp. 138–-160.

\bibitem{Perron}
  O. Perron,
  \emph{Zur Theorie der Dirichletschen Reihen},
  J. reine angew. Math. 134, 1908, pp. 95–-143.

\bibitem{Suryanarayana}
  D. Suryanarayana and R. S. R. C. Rao,
  \emph{On an Asymptotic Formula of Ramanujan},
  Mathematica Scandinavica, vol. 32 (1973) pp. 258--264.
  Available at \url{http://www.mscand.dk/article.php?id=2131}

\bibitem{Ramanujan}
  S. Ramanujan,
  \emph{Some Formulae in the Analytic Theory of Numbers},
  Messenger of Mathematics, vol. XLV, 1916, pp. 81--84.

\bibitem{Rekos}
  M. R\c{e}ko\'{s},
  \emph{On some complex explicit formulae connected with Euler's $\phi$ function},
  Functiones et Approximatio, 29, (2001), 113--124.
  Available at \url{http://www.staff.amu.edu.pl/~fa/XXIX/fa-29-1-113.pdf}

\bibitem {RieselBook}
  H. Riesel, \textit{Prime Numbers and Computer Methods for Factorization},
  2nd ed., Boston: Birkhauser, 1994.

\bibitem {RieselPaper}
  H. Riesel and G. Gohl,
  \emph{Calculations Related to Riemann's Prime Number Formula},
  Mathematics of Computation, 24(112), 1970, pp. 969--983.

\bibitem {RubinsteinAndSarnak}
  M. Rubinstein and P. Sarnak,
  \emph{Chebyshev's Bias},
  Experimental Mathematics vol. 3 no. 3 (1994), pp. 173-–197.

\bibitem{Sloane}
  N. J. A. Sloane,
  \emph{The On-Line Encyclopedia of Integer Sequences}. Available at\\
\url{http://www.research.att.com/~njas/sequences}.

\bibitem {Tanaka1}
  M. Tanaka,
  \emph{A Numerical Investigation on Cumulative Sum of the Liouville Function},
  Tokyo J. Math. 3, 187--189, 1980.

\bibitem {Tanaka2}
  M. Tanaka,
  \emph{On the M\"{o}bius and Allied Functions},
  Tokyo J. Math. 3, 215--218, 1980.

\bibitem {TitchmarshFirstEd}
  E. C. Titchmarsh,
  \emph{The Theory of the Riemann Zeta-Function},
  First Edition, Oxford University Press, 1951.

\bibitem {Titchmarsh}
  E. C. Titchmarsh,
  \emph{The Theory of the Riemann Zeta-Function},
  Second Edition, Oxford University Press, 1986.

\bibitem {Toth}
  L. T\'{o}th,
  \emph{A Survey of Gcd-Sum Functions},
  Journal of Integer Sequences, Vol. 13, Issue 8, 2010, Available at
  \url{http://www.cs.uwaterloo.ca/journals/JIS/VOL13/Toth/toth10.pdf}

\bibitem {Wagon}
  S. Wagon,
  \emph{Mathematica in Action},
  3rd ed., New York: Springer, 2010, pp. 512--519.

\bibitem{Wiertelak}
  K. Wiertelak,
  \emph{On some connections between zeta-zeros and square-free divisors of an integer},
  Functiones et Approximatio, 31, (2003), 133--145.  Available at\\
\url{http://www.staff.amu.edu.pl/~fa/XXXI/fa-31-1-133.pdf}

\bibitem{Wilson}
  B. Wilson,
  \emph{Proofs of Some Formulae Enunciated by Ramanujan},
  Proc. London Math. Soc., vol. 21, 1923, pp. 235--255.

\end{thebibliography}
\end{document}